\documentclass[amsbook,12pt]{article}
\usepackage{epsfig, amsfonts,amssymb, amsmath}
\usepackage[usenames]{color}
\usepackage{picinpar}
\newcommand{\nl}{\mbox{}\\}

\begin{document}
\thispagestyle{empty}
%
%
\mbox{} \vspace{-2.250cm} \\
\begin{center}
{\Large \bf
Dois Problemas em Equa\c c\~oes Diferenciais Parciais} \\
\nl
\mbox{} \vspace{-0.100cm} \\
{\large \sc Paulo R. Zingano} \\
\mbox{} \vspace{-0.400cm} \\
{\small Departamento de Matem\'atica Pura e Aplicada} \\
\mbox{} \vspace{-0.650cm} \\
{\small Universidade Federal do Rio Grande do Sul} \\
\mbox{} \vspace{-0.650cm} \\
{\small Porto Alegre, RS 91509-900, Brazil} \\
\nl
\mbox{} \vspace{-0.250cm} \\
\textit{Trabalho in\'edito submetido para avalia\c c\~ao} \\
\textit{para fins de progress\~ao funcional \`a Classe E} \\
\nl
\mbox{} \vspace{-0.250cm} \\
{\bf Abstract} \\
\mbox{} \vspace{-0.500cm} \\
\begin{minipage}[t]{12.250cm}
{\small
\mbox{} \hspace{+0.350cm}
In this work,
we examine two important problems
in the theory of nonlinear partial
differential equations.
In {\small \sc Part I},
we propose and solve a more general
and complete version of the
celebrated Leray's problem
for the incompressible Navier-Stokes equations
in $ \mathbb{R}^{3} \!$,
which in its simplest form
was suggested by J.$\;$Leray
in 1934
(and solved only in the 1980s by
T.$\;$ Kato, K.$\;$Masuda
and other authors).
A number of related new results
of clear interest
to the theory of Leray's solutions
are also given here. \\
%
%
\mbox{} \hspace{+0.350cm}
In {\small \sc Part II},
which is independent of
{\small \sc Part I}
and can be read separately,
we introduce an important new
collection of problems
concerning global existence results
and blow-up phenomena for solutions
of conservative advection-diffusion
equations in $ \mathbb{R}^{n} $
where heterogeneities in the lower order
terms tend to destabilize the solution
(everywhere or in certain regions),
strongly competing with the
viscous dissipation effects
to determine the overall solution behavior.
This work extends some
recent contributions led by the author
\cite{Zingano2010, Zingano2011}
(see also
\cite{BarrionuevoOliveiraZingano2014, BrazMeloZingano2015})
where the analysis was confined
to one-dimensional problems,
primarily in the case of
linear advection,
for which solutions
are automatically global
and can only misbehave
by increasing unboundedly
as $ t \rightarrow \infty $.
Here,
we consider the much more challenging
case of superlinear advection
(and arbitrary dimension),
which may cause finite-time blow-up
in several important spaces.
We then point out a new kind of
phenomena --- one that may be
properly named ``anti-Fujita"
for its vivid contrast to
the type of blow-up behavior
discovered by Fujita in the 1960s,
and which has been investigated ever since ---
that has apparently
been completely overlooked
in the literature.
For better clarity,
we restrict ourselves here
to the case
of linear nondegenerate diffusion,
but similar properties
and behavior are also to be found
in much more general diffusion phenomena,
as will be reported shortly.
}
\end{minipage}
\end{center}
\newpage
%
%
\mbox{} \vspace{-0.550cm} \\
%
%
\thispagestyle{empty}
{\Large \bf \'Indice} \\
\nl
{\sc Parte I: Problema de Leray} \\
\mbox{} \vspace{-0.200cm} \\
1. Introdu\c c\~ao
\dotfill $\;$\mbox{ }1 \\
2. Preliminares, I
\dotfill $\;$\mbox{ }7 \\
3. Preliminares, II
\dotfill $\;$12 \\
4. Prova de (1.12$a$)
\dotfill $\;$16 \\
5. Prova de (1.12$b$)
\dotfill $\;$23 \\
Ap\^endice A
\dotfill $\;$28 \\
\nl
{\sc Parte II: Problema de Exist\^encia Global
para Equa\c c\~oes de Advec\c c\~ao-Difus\~ao Conservativas} \\
\mbox{} \vspace{-0.200cm} \\
1. Introdu\c c\~ao
\dotfill $\;$31 \\
2. Preliminares
\dotfill $\;$35 \\
3. Prova de (1.10)
\dotfill $\;$40 \\
4. Condi\c c\~oes de exist\^encia global
\dotfill $\;$48 \\
\nl
{\sc Refer\^encias}
\dotfill $\;$51 \\
\nl
\nl
{\large \bf Contribui\c c\~oes Principais} \\
\nl
{\sc Parte I} \\
\mbox{} \vspace{-0.400cm} \\
Teorema 3.1
\dotfill $\;$12 \\
Teorema 4.1
\dotfill $\;$16 \\
Teorema 5.1
\dotfill $\;$23 \\
Teorema A.1
\dotfill $\;$30 \\
\nl
{\sc Parte II} \\
\mbox{} \vspace{-0.400cm} \\
Teorema A
\dotfill $\;$33 \\
Teorema B
\dotfill $\;$34 \\
Teorema 3.1
\dotfill $\;$40 \\
Teorema 4.1
\dotfill $\;$48 \\
\newpage
%
%
\mbox{} \vspace{-2.225cm} \\
%
%
\begin{center}
{\Large \sc Parte I} \\
\nl
\mbox{} \vspace{-0.350cm} \\
{\Large \bf Problema de Leray} \\
\end{center}
\setcounter{page}{1}
\mbox{} \vspace{-0.250cm} \\
%
%
%
%

{\bf 1. Introdu\c c\~ao} \\

Em seu trabalho seminal \cite{Leray1934},
Leray
construiu solu\c c\~oes (fracas) globais
de e\-ner\-gia finita
$ \mbox{\boldmath $u$}(\cdot,t) \in $
$ {\displaystyle
L^{\infty}([\;\!0, \infty), \:\!L^{2}_{\sigma}(\mathbb{R}^{3}))
\cap
C_{\!\;\!w}([\;\!0, \infty), \:\!L^{2}(\mathbb{R}^{3}))
\cap
L^{2}([\;\!0, \infty), \;\!\dot{H}^{1}(\mathbb{R}^{3}))
} $
para as equa\c c\~oes de Navier-Stokes
em $ \mathbb{R}^{3} \!\:\!$, \\
\mbox{} \vspace{-0.600cm} \\
\begin{equation}
\tag{1.1$a$}
\mbox{\boldmath $u$}_{t} +\,
\mbox{\boldmath $u$} \!\;\!\cdot\!\;\! \nabla\;\!
\mbox{\boldmath $u$}
\,+\,
\nabla p
\,=\;
\nu \,
\Delta \mbox{\boldmath $u$},
\qquad
\nabla \!\cdot\!\;\! \mbox{\boldmath $u$}(\cdot,t)
\,=\,0,
\end{equation}
\mbox{} \vspace{-0.900cm} \\
\begin{equation}
\tag{1.1$b$}
\mbox{\boldmath $u$}(\cdot,0) \,=\,
\mbox{\boldmath $u$}_0 \in L^{2}_{\sigma}(\mathbb{R}^{3}),
\end{equation}
\mbox{} \vspace{-0.200cm} \\
onde
$ \nu > 0 $ \'e constante,
e
$ {\displaystyle
L^{2}_{\sigma}(\mathbb{R}^{3})
} $
denota o espa\c co
de fun\c c\~oes
$ {\displaystyle
\;\!
\mbox{\bf u} =
(\:\! \mbox{u}_{\mbox{}_{1}} \!\:\!,
 \:\! \mbox{u}_{\mbox{}_{2}} \!\:\!,
 \:\! \mbox{u}_{\mbox{}_{3}} \!\;\!)
\in L^{2}(\mathbb{R}^{3})
\;\!
} $
com
$ {\displaystyle
\;\!
\nabla \!\cdot \mbox{\bf u}
= 0
} $
em $ {\cal D}^{\prime}(\mathbb{R}^{3}) $.
Ademais,
estas solu\c c\~oes
reproduzem o estado inicial
$ \mbox{\boldmath $u$}_0 \!\;\!$
em $ L^{2} \!\;\!$
(i.e.,
$ {\displaystyle
\|\, \mbox{\boldmath $u$}(\cdot,t) - \mbox{\boldmath $u$}_{0} \;\!
\|_{\mbox{}_{\scriptstyle L^{2}(\mathbb{R}^{3})}}
\!\!\,\!\rightarrow 0
\;\!
} $
ao $\;\! t \mbox{\footnotesize $\,\searrow\,$} 0 $)
e satisfazem
a desigualdade de energia \\
\mbox{} \vspace{-0.600cm} \\
\begin{equation}
\tag{1.2}
\|\,\mbox{\boldmath $u$}(\cdot,t) \,
\|_{\mbox{}_{\scriptstyle L^{2}(\mathbb{R}^{3})}}^{\:\!2}
\:\!+\;
2 \;\!\nu \!\!\;\!
\int_{0}^{\;\!\mbox{\mbox{\footnotesize $t$}}} \!
\|\, D \:\!\mbox{\boldmath $u$}(\cdot,s) \,
\|_{\mbox{}_{\scriptstyle L^{2}(\mathbb{R}^{3})}}^{\:\!2}
\;\!ds
\;\leq\;
\|\,\mbox{\boldmath $u$}_{0} \;\!
\|_{\mbox{}_{\scriptstyle L^{2}(\mathbb{R}^{3})}}^{\:\!2}
\end{equation}
\mbox{} \vspace{-0.100cm} \\
para todo $ t > 0 $
\cite{Galdi2000, KreissLorenz1989, Leray1934, Serrin1963}.
Enquanto a {\em unicidade\/}
das solu\c c\~oes de Leray
corres\-pon\-dentes a um estado inicial
$ \mbox{\boldmath $u$}_0 \in L^{2}_{\sigma}(\mathbb{R}^{3}) $
qualquer
permanece fundamentalmente em aberto
at\'e hoje,
Leray mostrou que
existe um instante
de tempo
$ {\displaystyle
\;\!
0 < \mbox{\small $T$}_{\!\:\!\ast\ast}
\!\;\!< \infty
\:\!
} $
(de\-pen\-dendo
dos dados $ \nu, \;\!\mbox{\boldmath $u$}_0 $ fornecidos)
tal que
$ {\displaystyle
\;\!
\mbox{\boldmath $u$} \in
C^{\infty}(\mathbb{R}^{3} \!\times\!\:\!
[\,\mbox{\small $T$}_{\!\:\!\ast\ast}\!\;\!, \infty))
} $\footnote{%
%
%
%
\'E conhecido tamb\'em
que
se tem
$ {\displaystyle
\;\!
T_{\!\:\!\ast\ast} \!\;\! <\;\!
K \;\! \nu^{-\,5} \,
\|\, \mbox{\boldmath $u$}_0 \;\!
\|_{{\scriptstyle L^{2}(\mathbb{R}^{3})}}^{\:\!4}
\!
} $
para certa constante absoluta
$ K > 0 $, \linebreak
\mbox{} \vspace{-0.500cm} \\
ver e.g.$\;$\cite{KreissHagstromLorenzZingano2003, %
LorenzZingano2012, Leray1934}.
($\:\!$No {\sc Ap\^endice A} a seguir,
melhoraremos o valor dado para $K\!\:\!$,
mostrando que
$ \:\!K \!\;\!< 0.000\,753\,026 \:\!$.)
Condi\c c\~oes adicionais
sobre o estado inicial
$ \mbox{\boldmath $u$}_0 \in L^{2}_{\sigma}(\mathbb{R}^{3}) $
--- por exemplo,
$ \mbox{\boldmath $u$}_0 \in H^{s}_{\sigma}(\mathbb{R}^{3}) $,
$ s > 3/2 $ ---
garantem al\'em disso
$ {\displaystyle
\mbox{\boldmath $u$} \in
C^{\infty}(\mathbb{R}^{3} \!\times\!\;\!(\:\!0,\mbox{\small $T$}_{\!\:\!\ast}])
} $
para certo
$ {\displaystyle
0 < \mbox{\small $T$}_{\!\:\!\ast} \leq
\mbox{\small $T$}_{\!\:\!\ast\ast}
} $
\cite{Yuan2008}.
%
%
%
}
%
%
e \\
\mbox{} \vspace{-0.550cm} \\
\begin{equation}
\tag{1.3}
%
%
\mbox{\boldmath $u$}(\cdot,t) \in
L^{\infty}_{\tt loc}
([\;\!\mbox{\small $T$}_{\!\:\!\ast\ast},\infty),
H^{m}(\mathbb{R}^{3}))
\end{equation}
\mbox{} \vspace{-0.150cm} \\
para cada $ m \geq 0 $,
onde $ H^{m}(\mathbb{R}^{3}) $
denota o espa\c co de Sobolev
das fun\c c\~oes
(neste caso, com valores em $ \mathbb{R}^{3} $)
em $ L^{2}(\mathbb{R}^{3}) $
com derivadas (espaciais) fracas de ordem at\'e $m$ \linebreak
em $ L^{2}(\mathbb{R}^{3}) $.
Um problema b\'asico importante
deixado aberto por Leray
em 1934
foi
(denotando
$ {\displaystyle
\;\!
W(t) := \|\, \mbox{\boldmath $u$}(\cdot,t) \,
\|_{\mbox{}_{\scriptstyle L^{2}(\mathbb{R}^{3})}}^{\:\!2}
\!\:\!
} $,
como em \cite{Leray1934}):  \\
\mbox{} \vspace{+0.075cm} \\
\mbox{} \hspace{+0.050cm}
\fbox{%
\begin{minipage}{14.00cm}
\nl
\mbox{} \vspace{-0.400cm} \\
\textit{%
\mbox{}$\;\;$J'ignore si W(t) tend n\'ecessairement vers $\:\!0$
quand $\,t$ augmente ind\'efiniment,
} \\
\mbox{} \hspace{+9.750cm}
J.$\;$Leray (\cite{Leray1934}, p.$\;$248) \\
\mbox{} \vspace{-0.900cm} \\
\end{minipage}
}

%
%
\mbox{} \vspace{-0.550cm} \\
ou seja,
se vale (ou n\~ao)
que \\
\mbox{} \vspace{-0.750cm} \\
\begin{equation}
\tag{1.4}
\lim_{t\,\rightarrow\,\infty}
\:
\|\, \mbox{\boldmath $u$}(\cdot,t) \,
\|_{\mbox{}_{\scriptstyle L^{2}(\mathbb{R}^{3})}}
=\; 0.
\end{equation}
\mbox{} \vspace{-0.175cm} \\
Esta quest\~ao somente foi resolvida
(positivamente) 50 anos mais tarde
por Kato
\cite{Kato1984} \linebreak
e subsequentemente
tamb\'em por outros autores
\cite{KajikiyaMiyakawa1986, Masuda1984, Wiegner1987}.
V\'arios desenvolvimentos e extens\~oes
importantes de (1.4)
vem sendo estabelecidos
(ver e.g.$\;$\cite{BenameurSelmi2012, BrazLorenzMeloZingano2014, %
KreissHagstromLorenzZingano2003, OliverTiti2000, %
SchonbekWiegner1996, SchutzZinganoZingano2014}
e a discuss\~ao abaixo).
Em particular,
uma prova extremamente simples
para (1.4)
foi obtida em
\cite{SchutzZinganoZingano2014},
com base no m\'etodo em
\cite{KreissHagstromLorenzZingano2003},
utilizando somente t\'ecnicas
j\'a conhecidas em 1934!$\:\!$\footnote{%
%
%
Como em \cite{KreissHagstromLorenzZingano2003},
esta prova faz uso tradicional de t\'ecnicas cl\'assicas
como estimativas de energia e transformadas de Fourier,
mas \'e v\'alida apenas em dimens\~ao $ n = 2, 3$.
Em contraste, a prova de (1.4) em \cite{Kato1984}
pode ser estendida a $ n = 4 $,
e o argumento (muito envolvente) desenvolvido em
\cite{Wiegner1987}
consegue estabelecer (1.4) para $ n \geq 2 $
qualquer.
}
%
%
$\!$($\,\!$Para uma descri\c c\~ao detalhada
do m\'etodo em \cite{SchutzZinganoZingano2014},
ver \cite{Perusato2014}.) \linebreak
\mbox{} \vspace{-0.900cm} \\

Dado $ t_0 \!\geq 0 $,
\'e natural que se tente aproximar
as solu\c c\~oes $ \mbox{\boldmath $u$}(\cdot,t) $
de (1.1)
para $ t > t_0 $
pelas solu\c c\~oes
$ {\displaystyle
\mbox{\boldmath $v$}(\cdot,t) =
e^{\nu \:\! \Delta \:\!(t - t_0)} \;\!
\mbox{\boldmath $u$}(\cdot,t_0)
} $
dos problemas lineares
associados \\
\mbox{} \vspace{-0.625cm} \\
\begin{equation}
\tag{1.5}
\mbox{\boldmath $v$}_{t}
\,=\;
\nu \,
\Delta \mbox{\boldmath $v$},
\qquad
\mbox{\boldmath $v$}(\cdot,t_0) \,=\,
\mbox{\boldmath $u$}(\cdot,t_0).
\end{equation}
\mbox{} \vspace{-0.200cm} \\
Para estas solu\c c\~oes,
\'e f\'acil obter
v\'arias estimativas de decaimento,
como e.g. \\
\mbox{} \vspace{-0.550cm} \\
\begin{equation}
\tag{1.6$a$}
\lim_{t\,\rightarrow\,\infty}
\:
\|\, \mbox{\boldmath $v$}(\cdot,t) \,
\|_{\mbox{}_{\scriptstyle L^{2}(\mathbb{R}^{n})}}
\,=\;
0,
\end{equation}
\mbox{} \vspace{-0.750cm} \\
\begin{equation}
\tag{1.6$b$}
\lim_{t\,\rightarrow\,\infty}
\:
t^{\mbox{}^{\;\! \frac{\scriptstyle n}{\scriptstyle 4}}}
\,\!
\|\, \mbox{\boldmath $v$}(\cdot,t) \,
\|_{\mbox{}_{\scriptstyle L^{\infty}(\mathbb{R}^{n})}}
\,=\;
0,
\end{equation}
\mbox{} \vspace{-0.750cm} \\
\begin{equation}
\tag{1.6$c$}
\lim_{t\,\rightarrow\,\infty}
\:
t^{\mbox{}^{\;\! \frac{\scriptstyle s}{\scriptstyle 2}}}
\,\!
\|\, \mbox{\boldmath $v$}(\cdot,t) \,
\|_{\mbox{}_{\scriptstyle \dot{H}^{\!\;\!s}(\mathbb{R}^{n})}}
\,=\;
0,
\qquad
s \geq 0,
\end{equation}
\mbox{} \vspace{-0.125cm} \\
v\'alidas para todo $n$.
(Em (1.6$c$) acima,
$ \!\:\!\dot{H}^{\!\;\!s}\!\;\!(\mathbb{R}^{n}) $
denota o espa\c co de Sobolev homo\-g\^eneo
formado pelas fun\c c\~oes
$ {\displaystyle
\mbox{\boldmath $v$} = (v_{\scriptscriptstyle 1}\!\;\!,\!...,v_{n})
\in L^{2}(\mathbb{R}^{n})
} $
tais que
$ {\displaystyle
\;\!
| \cdot |^{\,\!s}_{\mbox{}_{2}} \:\!
|\,\hat{\mbox{\boldmath $v$}}(\cdot)\,|_{\mbox{}_{2}}
\in\!\;\! L^{2}(\mathbb{R}^{n})
} $,
onde
$ \hat{\mbox{\boldmath $v$}}(\cdot) $
denota a transformada de Fourier de
$ \mbox{\boldmath $v$}(\cdot) $,
com norma
$ {\displaystyle
\,\!
\| \cdot
\|_{\mbox{}_{\scriptstyle \dot{H}^{\!\:\!s}\!\;\!(\mathbb{R}^{n})}}
\!
} $
dada por \\
\mbox{} \vspace{-0.500cm} \\
\begin{equation}
\tag{1.7}
\|\, \mbox{\boldmath $v$} \,
\|_{\mbox{}_{\scriptstyle \dot{H}^{\!\:\!s}\!\;\!(\mathbb{R}^{n})}}
=\,
\Bigl\{%
\int_{\mathbb{R}^{n}} \!\!
|\,\xi\,|_{\mbox{}_{2}}^{\:\!2\:\!s}
\;\!
|\, \hat{\mbox{\boldmath $v$}}(\xi)\,|_{\mbox{}_{2}}^{2}
\,d\xi
\:\Bigr\}^{\!\!\;\!1/2}
\!\!.
\end{equation}
\mbox{} \vspace{-0.050cm} \\
Podemos assim esperar
que estimativas similares a (1.6)
sejam tamb\'em v\'alidas
para as solu\c c\~oes de Leray
$ \mbox{\boldmath $u$}(\cdot,t) $
de (1.1),
ao menos para $ n = 3 $,
dado
$ \mbox{\boldmath $ u $}_0 \in L^{2}_{\sigma}(\mathbb{R}^{n}) $
qualquer.
Este \'e essencialmente o
{\small \sc Problema de Leray}
para as equa\c c\~oes de Navier-Stokes,
com sua vers\~ao mais b\'asica
dada em (1.4) acima.
Outras quest\~oes
se p\~oem aqui naturalmente;
por exemplo,
a respeito do comportamento
da diferen\c ca
(ou erro)
$ \mbox{\boldmath $u$}(\cdot,t) - \mbox{\boldmath $v$}(\cdot,t) $
para $ t \gg 1 $.
Na norma $L^{2}\!\;\!$,
esta quest\~ao sobre o erro
foi respondida por Wiegner
\cite{Wiegner1987},
tendo-se,
para $ t_0 \geq 0 $ qualquer,  \\
\mbox{} \vspace{-0.550cm} \\
\begin{equation}
\tag{1.8}
\lim_{t\,\rightarrow\,\infty}
\,
t^{\;\! \frac{\scriptstyle n - 2}{\scriptstyle 4}}
\;\!
\|\, \mbox{\boldmath $u$}(\cdot,t) -\,
e^{\:\!\nu \:\! \Delta \:\!(t - t_0)}
\;\!\mbox{\boldmath $u$}(\cdot,t_0) \,
\|_{\mbox{}_{\scriptstyle L^{2}(\mathbb{R}^{n})}}
\:\!=\;
0.
\end{equation}
\mbox{} \vspace{-0.150cm} \\
Como com (1.4), uma prova mais simples
para (1.8) foi recentemente dada
em \cite{SchutzZinganoZingano2014},
supondo $ n \leq 3 $,
tendo-se tamb\'em mostrado
os resultados correspondentes a (1.6$b$): \linebreak
\mbox{} \vspace{-0.550cm} \\
\begin{equation}
\tag{1.9$a$}
\lim_{t\,\rightarrow\,\infty} \,
t^{\;\! \frac{\scriptstyle n}{\scriptstyle 4}}
\;\!
\|\, \mbox{\boldmath $u$}(\cdot,t) \,
\|_{\mbox{}_{\scriptstyle L^{\infty}(\mathbb{R}^{n})}}
\:\!=\;
0,
\end{equation}
\mbox{} \vspace{-0.750cm} \\
\begin{equation}
\tag{1.9$b$}
\lim_{t\,\rightarrow\,\infty}
\,
t^{\;\! \frac{\scriptstyle n - 1}{\scriptstyle 2}}
\;\!
\|\, \mbox{\boldmath $u$}(\cdot,t) -\,
e^{\:\!\nu \:\! \Delta \:\!(t - t_0)}
\;\!\mbox{\boldmath $u$}(\cdot,t_0) \,
\|_{\mbox{}_{\scriptstyle L^{\infty}(\mathbb{R}^{n})}}
\:\!=\;
0,
\end{equation}
\mbox{} \vspace{-0.100cm} \\
em dimens\~ao $ n = 2, 3 $
(cf.$\;$\cite{SchutzZinganoZingano2014}, Section 4).
No presente trabalho,
vamos estabelecer as estimativas (mais dif\'\i ceis)
correspondentes a (1.6$c$)
no caso das solu\c c\~oes de Leray \linebreak
do problema (1.1),
al\'em de estimativas similares sobre
o erro (ver (1.10$b$) abaixo).
%
%
Os resultados s\~ao simples
de descrever:
dados
$ \mbox{\boldmath $ u $}_0 \!\;\!\in L^{2}_{\sigma}(\mathbb{R}^{n}) $,
$ t_0 \!\;\!\geq 0 $
quaisquer,
\mbox{tem-se} \\
\mbox{} \vspace{-0.575cm} \\
\begin{equation}
\tag{1.10$a$}
\lim_{t\,\rightarrow\,\infty}
\:
t^{\mbox{}^{\scriptstyle \frac{\scriptstyle s}{2} }}
\|\, \mbox{\boldmath $u$}(\cdot,t) \,
\|_{\mbox{}_{\scriptstyle \dot{H}^{\!\:\!s}\!\;\!(\mathbb{R}^{n})}}
\!\;\!=\:0,
\end{equation}
\mbox{} \vspace{-0.720cm} \\
\begin{equation}
\tag{1.10$b$}
\lim_{t\,\rightarrow\,\infty}
\:
t^{\mbox{}^{\scriptstyle
\frac{{\scriptstyle n} - 2}{4} \;\!+\;\!\frac{\scriptstyle s}{2}
}}
\;\!
\|\, \mbox{\boldmath $u$}(\cdot,t) \:\!-\:\!
e^{\:\!\nu \:\!\Delta (t - t_0)} \mbox{\boldmath $u$}(\cdot,t_0) \,
\|_{\mbox{}_{\scriptstyle \dot{H}^{\!\:\!s}\!\;\!(\mathbb{R}^{n})}}
\!\;\!=\:0,
\end{equation}
\mbox{} \vspace{-0.120cm} \\
para todo $ s \geq 0 $ ($s$ real),
e $ \;\!n = 2, 3 $.
Note-se que (1.4), (1.8) e (1.9)
tornam-se agora
consequ\^encias
simples
de (1.10),
que pode assim
ser considerada
como a forma geral completa
do pro\-blema de Leray,
resolvida neste trabalho
(exceto que,
no caso especial
$ n = 2 $,
tem-se
(1.10$a$) derivada
previamente em
\cite{BenameurSelmi2012},
usando um procedimento diferente).
A obten\c c\~ao de
(1.10$b$)
\'e particularmente delicada,
e utiliza
a estimativa \\
\mbox{} \vspace{-0.450cm} \\
\begin{equation}
\tag{1.11}
\|\, e^{\:\!\nu \:\!\Delta (t - t_0)}
\mbox{\boldmath $u$}(\cdot,t_0) \;\!-\,
e^{\:\!\nu \:\!\Delta (\:\!t - t_1\!\;\!)}
\mbox{\boldmath $u$}(\cdot,t_1) \,
\|_{\mbox{}_{\scriptstyle \dot{H}^{\!\:\!s}\!\;\!(\mathbb{R}^{n})}}
\!\;\!\leq\;\!
K\;\! \nu^{\scriptstyle \,-\,
\frac{1}{2} \,-\, \gamma} \,
( t_1 \!\;\!-\;\! t_0 )^{\scriptstyle \!\;\!\frac{1}{2}}
\;\!
\|\, \mbox{\boldmath $u$}_0 \;\!
\|_{\mbox{}_{\scriptstyle L^{2}(\mathbb{R}^{n})}}^{\:\!2}
(t - t_{1} \!\;\!)^{\mbox{}^{\scriptstyle \!\;\!-\, \gamma}}
%
%
\end{equation}
\mbox{} \vspace{-0.250cm} \\
para $ t > t_1 > t_0 \geq 0 $ arbitr\'arios,
derivada na Se\c c\~ao 3 abaixo,
onde
$ \gamma = n/4 \;\!+\;\! s/2 $,
e $ \:\!K \!= K\!\;\!(n,s) > 0 \;\!$
\'e uma constante
(cujo valor
depende apenas dos par\^ametros $\;\!n, s $
considerados).
%
%
Leitores interessados prioritariamente
nas novas contribui\c c\~oes
do presente trabalho
podem neste ponto
consultar diretamente
os seguintes resultados:
Teorema 3.1,
Teorema 4.1 (e Lema 4.1),
Teorema 5.1
e
Teorema A.1 (Ap\^endice A).
Tamb\'em pode ser conveniente
rever rapidamente os
Teoremas 4.2 e 4.3
e
os resultados
revisados
na Se\c c\~ao 2
a seguir. \\
\nl
%
%
%
%
{\bf Observa\c c\~ao 1.1.}
Tomando-se $ s = m $ (inteiro) em (1.10),
resulta
(por (1.7) e do fato \linebreak
de se ter$\;\!$\footnote{%
%
%
Para a defini\c c\~ao de
$ {\displaystyle
\;\!
\| \, D^{m} \:\!\mbox{\boldmath $u$}(\cdot,t) \,
\|_{\scriptstyle L^{2}(\mathbb{R}^{n})}
\!\:\!
} $
e outras normas aqui usadas,
ver (1.19), (1.20) a seguir.
De modo geral,
o s\'\i mbolo $ D^{m} $
refere-se coletivamente a todas
as derivadas espaciais de ordem $m$,
enquanto
$ \:\!D^{\alpha} \!\:\!$
denota uma derivada particular,
correspondente ao multi-\'\i ndice $\:\!\alpha\:\!$
indicado.
%
%
%
}
$ {\displaystyle
\|\, \mbox{\boldmath $u$}(\cdot,t) \,
\|_{\mbox{}_{\scriptstyle \dot{H}^{\!\:\!m}\!\;\!(\mathbb{R}^{n})}}
\!\;\!=\;\!
\|\, D^{m} \:\!\mbox{\boldmath $u$}(\cdot,t) \,
\|_{\mbox{}_{\scriptstyle L^{2}(\mathbb{R}^{n})}}
\!\:\!
} $,
por Parseval): \\
\mbox{} \vspace{-0.475cm} \\
\begin{equation}
\tag{1.12$a$}
\lim_{t\,\rightarrow\,\infty}
\:
t^{\mbox{}^{\scriptstyle \frac{\scriptstyle m}{2} }}
\|\, D^{m} \:\!\mbox{\boldmath $u$}(\cdot,t) \,
\|_{\mbox{}_{\scriptstyle L^{2}(\mathbb{R}^{n})}}
\!\;\!=\:0,
\end{equation}
\mbox{} \vspace{-0.750cm} \\
\begin{equation}
\tag{1.12$b$}
\lim_{t\,\rightarrow\,\infty}
\:
t^{\mbox{}^{\scriptstyle
\frac{{\scriptstyle n} - 2}{4} \;\!+\;\!\frac{|\;\!{\scriptstyle \alpha}\;\!|}{2}
}}
\;\!
\|\, D^{\alpha} \:\! \mbox{\boldmath $u$}(\cdot,t) \:\!-\:\!
D^{\alpha} \:\![\:
e^{\:\!\nu \:\!\Delta (t - t_0)} \mbox{\boldmath $u$}(\cdot,t_0)
\;\!] \,
\|_{\mbox{}_{\scriptstyle L^{2}(\mathbb{R}^{n})}}
\!\;\!=\:0,
\end{equation}
\mbox{} \vspace{-0.150cm} \\
para todo $ m \geq 0 $,
$ {\displaystyle
\;\!
\alpha = (\:\! \alpha_{\mbox{}_{1}}\!\,\!, \:\!
\alpha_{\mbox{}_{2}}\!\;\!,\!\;\!...,\alpha_{\mbox{}_{\scriptstyle n}})
} $,
sendo
$ {\displaystyle
|\,\alpha\,| =
\alpha_{\mbox{}_{1}}\!\:\!+
\alpha_{\mbox{}_{2}}\!\:\!+ ... +
\alpha_{\mbox{}_{\scriptstyle n}}
\!\;\!
} $
a ordem de $ \alpha $. \linebreak
Assim, (1.10) descreve o comportamento
de $ \mbox{\boldmath $u$}(\cdot,t) $,
$ \mbox{\boldmath $u$}(\cdot,t) - \mbox{\boldmath $v$}(\cdot,t) $
e suas derivadas (espaciais) de qualquer ordem.
Na verdade, (1.10) e (1.12)
s\~ao equivalentes:
tendo-se \linebreak
(1.12),
obt\'em-se (1.10)
aplicando-se a
propriedade de interpola\c c\~ao \\
\mbox{} \vspace{-0.475cm} \\
\begin{equation}
\tag{1.13}
\|\, \mbox{\boldmath $u$}(\cdot,t) \,
\|_{\mbox{}_{\scriptstyle \dot{H}^{\!\:\!s}\!\;\!(\mathbb{R}^{n})}}
\leq\:
\|\, \mbox{\boldmath $u$}(\cdot,t) \,
\|_{\mbox{}_{\scriptstyle \dot{H}^{\!\:\!s_{\mbox{}_{1}}}\!\;\!(\mathbb{R}^{n})}}
  ^{\mbox{}^{\scriptstyle \:\!\alpha_{\mbox{}_{1}}}}
\:\!
\|\, \mbox{\boldmath $u$}(\cdot,t) \,
\|_{\mbox{}_{\scriptstyle \dot{H}^{\!\:\!s_{\mbox{}_{2}}}\!\;\!(\mathbb{R}^{n})}}
  ^{\mbox{}_{\scriptstyle \:\!\alpha_{\mbox{}_{2}}}}
\!,
\qquad
s_{\mbox{}_{1}} <\:\! s <\:\! s_{\mbox{}_{2}}
\!\;\!,
\end{equation}
\mbox{} \vspace{-0.110cm} \\
onde
$ {\displaystyle
\;\!
\alpha_{\mbox{}_{1}} \!=\:\!
\theta_{\mbox{}_{1}} \:\!s \:\!/
s_{\mbox{}_{1}}
\!\;\!
} $,
$ {\displaystyle
\:\!
\alpha_{\mbox{}_{2}} \!=\:\!
\theta_{\mbox{}_{2}} \:\!s \:\!/
s_{\mbox{}_{2}}
\!\;\!
} $,
sendo
$ {\displaystyle
\;\!
\theta_{\mbox{}_{1}} \!\:\!,
\:\!\theta_{\mbox{}_{2}}
\in (\:\!0, 1)
\:\!
} $
dados por
$ {\displaystyle
\;\!
s^{-\;\!1} \!\:\!=
\theta_{\mbox{}_{1}}
s_{\scriptscriptstyle 1}^{-\;\!1}
\!\;\!+\;\!
\theta_{\mbox{}_{2}}
s_{\scriptscriptstyle 2}^{-\;\!1}
\!\:\!
} $. \linebreak
A obten\c c\~ao de (1.10)
nas Se\c c\~oes 2, 3 abaixo
ser\'a feita
considerando-se a forma (1.12)
destas propriedades.
Se desejado,
seria tamb\'em suficiente
derivar o resultado
no caso particular $ \nu = 1 $;
uma vez obtido,
os resultados
(1.10), (1.11), (1.12)
para $ \nu > 0 $ geral
decorreriam ent\~ao
de argumentos simples de escala
(dado que,
sendo $ \mbox{\boldmath $u$}(x,t) $,
$ p(x,t) $ \linebreak
uma solu\c c\~ao de Leray
do sistema (1.1)
para dado $ \nu > 0 $,
ent\~ao
$ {\displaystyle
\mbox{\boldmath $U$}(x,t) :=
\mbox{\boldmath $u$}(\nu \;\!x, \nu \;\!t)
} $,
$ {\displaystyle
P(x,t) := p(\nu \;\!x, \nu \;\!t)
} $
definem uma solu\c c\~ao de Leray
para (1.1) com $ \nu = 1 $). \\
%
%
\nl
%
%
%
%
{\bf Observa\c c\~ao 1.2.}
A an\'alise e resultados a seguir
podem tamb\'em ser adaptados/es\-ten\-didos para
o pro\-blema de Navier-Stokes
com for\c cas externas,
ou seja, \\
\mbox{} \vspace{-0.600cm} \\
\begin{equation}
\tag{1.14$a$}
\mbox{\boldmath $u$}_{t} +\,
\mbox{\boldmath $u$} \!\;\!\cdot\!\;\! \nabla\;\!
\mbox{\boldmath $u$}
\,+\,
\nabla p
\,=\;
\nu \,
\Delta \mbox{\boldmath $u$}
\,+\,
\mbox{\boldmath $f$}(\cdot,t),
\qquad
\nabla \!\cdot\!\;\! \mbox{\boldmath $u$}(\cdot,t)
\,=\,0,
\end{equation}
\mbox{} \vspace{-0.900cm} \\
\begin{equation}
\tag{1.14$b$}
\mbox{\boldmath $u$}(\cdot,0) \,=\,
\mbox{\boldmath $u$}_0 \in L^{2}_{\sigma}(\mathbb{R}^{3}),
\end{equation}
\mbox{} \vspace{-0.200cm} \\
onde se sup\~oe,
no caso mais simples,
$ {\displaystyle
\mbox{\boldmath $f$} \in
C^{\infty}(\mathbb{R}^{n} \!\times \!\:\![\,0, \infty))
\cap\,
L^{2}(\mathbb{R}^{n} \!\times \!\:\![\,0, \infty))
} $
satisfa\-zendo (1.15) abaixo:
considerando-se a proje\c c\~ao
de Helmholtz
$ {\displaystyle
\;\!
\mbox{\boldmath $g$}(\cdot,t)
\:\!=\,
\mathbb{P}_{\mbox{}_{\!H}}[\;\!\mbox{\boldmath $f$}(\cdot,t)\,]
\;\!
} $
de $ \mbox{\boldmath $f$}(\cdot,t) $
em $ L^{2}_{\sigma}(\mathbb{R}^{n}) $,
a suposi\c c\~ao \\
\mbox{} \vspace{-0.500cm} \\
\begin{equation}
\tag{1.15}
\int_{\:\!0}^{\infty}
\!\!\!\;\!
(1 + t)^{\:\!m/2} \:
\|\,D^{m} \,\!\mbox{\boldmath $g$}(\cdot,t) \,
\|_{\mbox{}_{\scriptstyle L^{2}(\mathbb{R}^{n})}}
\;\!dt
\;<\, \infty,
\qquad
\forall \;\, m \geq 0
\end{equation}
\mbox{} \vspace{-0.120cm} \\
permite a validade
de (1.6), (1.10), (1.11), (1.12) acima,
onde agora
$ {\displaystyle
\;\!
\mbox{\boldmath $v$}(\cdot,t)
\equiv
\mbox{\boldmath $v$}(\cdot,t\:\!; \:\! t_0)
} $
deve ser definida
como a solu\c c\~ao
(\'unica)
em
$ L^{\infty}(\:\![\,t_0,\infty), L^{2}_{\sigma}(\mathbb{R}^{n})) $
do problema \\
\mbox{} \vspace{-0.600cm} \\
\begin{equation}
\tag{1.16}
\mbox{\boldmath $v$}_{t}
\;=\;
\nu \;\! \Delta \mbox{\boldmath $v$}
\,+\,
\mbox{\boldmath $g$}(\cdot,t),
\qquad
\mbox{\boldmath $v$}(\cdot,t_0) \,=\,
\mbox{\boldmath $u$}(\cdot,t_0).
\end{equation}
\mbox{} \vspace{-0.250cm} \\
Hip\'oteses mais fracas sobre
$ \mbox{\boldmath $g$}(\cdot,t) $
podem tamb\'em ser adotadas no lugar de (1.15),
com resultados correspondentes mais fracos;
por exemplo,
tendo-se apenas \\
\mbox{} \vspace{-0.500cm} \\
\begin{equation}
\tag{1.17}
\int_{\:\!0}^{\infty}
\!\!\!\;\!
(1 + t)^{\:\!m/2} \:
\|\,D^{m} \,\!\mbox{\boldmath $g$}(\cdot,t) \,
\|_{\mbox{}_{\scriptstyle L^{2}(\mathbb{R}^{n})}}
\;\!dt
\;<\, \infty,
\qquad
m = 0, 1,
\end{equation}
\mbox{} \vspace{-0.150cm} \\
obt\'em-se
(adaptando-se a prova em \cite{SchutzZinganoZingano2014},
como feito em \cite{BrazLorenzMeloZingano2014})
que as estimativas
(1.4), (1.6$a$), (1.6$b$), (1.8) e (1.9) acima
permanecem v\'alidas.
Para mais detalhes,
ver \cite{BrazLorenzMeloZingano2014}.
Em \cite{Wiegner1987},
obt\'em-se (1.4) e (1.8)
para as solu\c c\~oes de (1.14)
supondo-se (1.17) para
$ m = 0 $ apenas,
e adicionalmente \\
\mbox{} \vspace{-0.650cm} \\
\begin{equation}
\tag{1.18}
\limsup_{t\,\rightarrow\,\infty}
\;
t^{\;\!n/4 \,+\, 1/2} \:
\|\, \mbox{\boldmath $g$}(\cdot,t) \,
\|_{\mbox{}_{\scriptstyle L^{n}(\mathbb{R}^{n})}}
\;\!<\;
\infty.
\end{equation}
\mbox{} \vspace{-1.050cm} \\
\nl

Descrevendo sucintamente o conte\'udo
das se\c c\~oes seguintes,
apresentamos na Se\-\c{c}\~ao~2
v\'arios resultados conhecidos
que ser\~ao relevantes na
deriva\c c\~ao das estimativas
(1.10), (1.11) e (1.12) a seguir.
Esta discuss\~ao \'e baseada
em \cite{KreissHagstromLorenzZingano2002, %
KreissHagstromLorenzZingano2003, Leray1934}.
Na Se\c c\~ao~3,
estabelecemos (1.11),
que ser\'a importante
para simplificar a
prova de (1.10$b$), (1.12$b$)
mais adiante.
A Se\c c\~ao 4
\'e inteiramente voltada
\`a obten\c c\~ao
das estimativas (1.10$a$) e (1.12$a$),
enquanto a Se\c c\~ao~5
\'e dedicada a
(1.10$b$), (1.12$b$).
Com exce\c c\~oes ocasionais, \linebreak
apresentaremos os detalhes
no caso espec\'\i fico $ n = 3 $
apenas, j\'a que
as provas
corres\-pondentes
em dimens\~ao $ n = 2 $
podem ser feitas
seguindo um procedimento
inteiramente an\'alogo
(e, em alguns casos,
bem mais simples). \\
\nl
\mbox{} \vspace{-0.300cm} \\
%
%
%
%
%
{\bf Nota\c c\~ao.}
Como j\'a visto acima,
usaremos (em geral) letras em negrito
para gran\-de\-zas vetoriais,
e.g.
$ {\displaystyle
\mbox{\boldmath $u$}(x,t)
=
} $
$ {\displaystyle
(\:\! u_{\mbox{}_{\!\:\!1}}\!\;\!(x,t),...,
 \:\! u_{\mbox{}_{\!\;\!\scriptstyle n}}\!\;\!(x,t) \:\!)
} $,
denotando
por $ | \cdot |_{\mbox{}_{2}} $
(ou simplesmente $ | \cdot | $)
a norma Euclideana em $ \mathbb{R}^{n}\!\;\!$,
ver e.g.$\;$(1.7).
$\!$Como \'e usual,
$ \nabla p \equiv \nabla p(\cdot,t) $
denota o gradiente (espacial) de $ \;\!p(\cdot,t) $,
$ D_{\!\;\!j} \!\;\!=\:\! \partial / \partial x_{\!\;\!j} \!\;\! $,
e
$ {\displaystyle
\:\!
\nabla \!\cdot \mbox{\boldmath $u$}
\:\!=
  D_{\mbox{}_{\!\:\!1}} u_{\mbox{}_{\!\:\!1}} \!\;\!+
  ... +
  D_{\mbox{}_{\scriptstyle \!\;\!n}} \:\! u_{\mbox{}_{\scriptstyle \!\;\!n}} )
} $
\'e o divergente (espacial) de
$ \:\!\mbox{\boldmath $u$}(\cdot,t) $;
analogamente,
$ {\displaystyle
\mbox{\boldmath $u$} \cdot \nabla \mbox{\boldmath $u$}
\;\!=\;\!
u_{\mbox{}_{\!\:\!1}} D_{\mbox{}_{\!\:\!1}} \mbox{\boldmath $u$}
+ ... +
u_{\mbox{}_{\scriptstyle \!\:\!n}}
D_{\mbox{}_{\scriptstyle \!\:\!n}} \mbox{\boldmath $u$}
} $,
etc. \linebreak
$ {\displaystyle
\| \;\!\cdot\;\!
\|_{\mbox{}_{\scriptstyle L^{q}(\mathbb{R}^{n})}}
\!\;\!
} $,
$ 1 \leq q \leq \infty $,
denota a norma tradicional
do espa\c co de Lebesgue
$ L^{q}(\mathbb{R}^{n}) $,
pondo-se,
para $ 1 \leq q < \infty $:
\newpage
\mbox{} \vspace{-1.250cm} \\
\begin{equation}
\tag{1.19$a$}
\|\, \mbox{\boldmath $u$}(\cdot,t) \,
\|_{\mbox{}_{\scriptstyle L^{q}(\mathbb{R}^{n})}}
\;\!=\;
\Bigl\{\,
\sum_{i\,=\,1}^{n} \int_{\mathbb{R}^{n}} \!
|\:u_{i}(x,t)\,|^{q} \;\!dx
\,\Bigr\}^{\!\!\:\!1/q}
\end{equation}
\mbox{} \vspace{-0.750cm} \\
\begin{equation}
\tag{1.19$b$}
\|\, D \mbox{\boldmath $u$}(\cdot,t) \,
\|_{\mbox{}_{\scriptstyle L^{q}(\mathbb{R}^{n})}}
\;\!=\;
\Bigl\{\,
\sum_{i, \,j \,=\,1}^{n} \int_{\mathbb{R}^{n}} \!
|\, D_{\!\;\!j} \;\!u_{i}(x,t)\,|^{q} \;\!dx
\,\Bigr\}^{\!\!\:\!1/q}
\end{equation}
\mbox{} \vspace{-0.750cm} \\
\begin{equation}
\tag{1.19$c$}
\|\, D^{2} \mbox{\boldmath $u$}(\cdot,t) \,
\|_{\mbox{}_{\scriptstyle L^{q}(\mathbb{R}^{n})}}
\;\!=\;
\Bigl\{\!\!
\sum_{\mbox{} \;\;i, \,j, \,\ell \,=\,1}^{n}
\!\;\! \int_{\mathbb{R}^{n}} \!
|\, D_{\!\;\!j} \:\!D_{\ell} \, u_{i}(x,t)\,|^{q} \;\!dx
\,\Bigr\}^{\!\!\:\!1/q}
\end{equation}
\mbox{} \vspace{-0.000cm} \\
e,
mais geralmente,
para $ m \geq 1 $ qualquer: \\
\mbox{} \vspace{-0.650cm} \\
\begin{equation}
\tag{1.19$d$}
\|\, D^m \:\! \mbox{\boldmath $u$}(\cdot, t)\,
\|_{\mbox{}_{\scriptstyle L^{q}(\mathbb{R}^{n})}}
=
\biggl(\;\! \sum_{i\,=\,1}^{n}
\sum_{\;\!j_{\mbox{}_{1}}=\,1}^{n}
\!\cdots\!
\sum_{j_{\mbox{}_{m}}=\,1}^{n}
\int_{\mbox{}_{\mbox{\scriptsize $\!\;\!\mathbb{R}^{n}$}}}
\!\!\!\!\;\!
|\;\!D_{\scriptstyle \!\;\!j_{\mbox{}_{1}}}
\!\!\:\!\cdots
D_{\scriptstyle \!\;\!j_{\mbox{}_{m}}}
u_{i}(x,t)\:|^q\, dx\:\!\biggr)^{\!\!\:\!1/q}\!\!\!\!,
\end{equation}
\mbox{} \vspace{-0.150cm} \\
denotando-se por
$ {\displaystyle
\;\!
\|\, \mbox{\boldmath $u$}(\cdot,t) \,
\|_{\mbox{}_{\scriptstyle L^{\infty}(\mathbb{R}^{n})}}
\!=\;\!
\max \, \bigl\{\,
\|\,u_{i}(\cdot,t)\,
\|_{\mbox{}_{\scriptstyle L^{\infty}(\mathbb{R}^{n})}}
\!\!: \, 1 \leq i \leq n
\,\bigr\}
} $
o supremo (essencial)
de
$ \mbox{\boldmath $u$}(\cdot,t) $,
e similarmente
para
$ {\displaystyle
\;\!
\|\, D \;\!\mbox{\boldmath $u$}(\cdot,t) \,
\|_{\mbox{}_{\scriptstyle L^{\infty}(\mathbb{R}^{n})}}
\!\,\!
} $,
$ {\displaystyle
\;\!
\|\, D^{2} \:\!\mbox{\boldmath $u$}(\cdot,t) \,
\|_{\mbox{}_{\scriptstyle L^{\infty}(\mathbb{R}^{n})}}
\!\,\!
} $,
etc.
(Com estas defini\c c\~oes,
tem-se
$ {\displaystyle
\;\!
\|\, \mbox{\boldmath $u$}(\cdot,t) \,
\|_{\mbox{}_{\scriptstyle L^{q}(\mathbb{R}^{n})}}
\!\;\!\rightarrow\;\!
\|\, \mbox{\boldmath $u$}(\cdot,t) \,
\|_{\mbox{}_{\scriptstyle L^{\infty}(\mathbb{R}^{n})}}
\!\;\!
} $
ao
$ q \rightarrow \infty $,
assim como,
mais geralmente,
$ {\displaystyle
\|\, D^m \:\! \mbox{\boldmath $u$}(\cdot,t) \,
\|_{\mbox{}_{\scriptstyle L^{q}(\mathbb{R}^{n})}}
\!\;\!\rightarrow \;\!
\|\, D^m \:\! \mbox{\boldmath $u$}(\cdot,t) \,
\|_{\mbox{}_{\scriptstyle L^{\infty}(\mathbb{R}^{n})}}
\!\;\!
} $,
para todo $m$.)$\:\!$\footnote{%
%
%
Mais seriamente,
conv\'em observar que,
com as defini\c c\~oes (1.19), (1.20),
se uma desigualdade de tipo Nirenberg-Gagliardo
$ {\displaystyle
\|\,{\sf u}\,\|_{\mbox{}_{\scriptstyle L^{\!q}}}
\!\leq\!\;\!K\;\!
\|\,{\sf u}\,\|_{\mbox{}_{\scriptstyle L^{\!r_{\mbox{}_{\!\:\!1}}}}}^{1 - \theta}
\|\;\!\nabla {\sf u}\,
\|_{\mbox{}_{\scriptstyle L^{\!r_{\mbox{}_{\!\:\!2}}}}}^{\:\!\theta}
\!
} $,
$ 0 \leq \theta \leq 1 $,
valer para fun\c c\~oes {\em escalares\/} $ {\sf u} $ \linebreak
($K \!> 0$ constante),
ent\~ao ela ser\'a
automaticamente v\'alida
para fun\c c\~oes
{\em vetoriais\/} $ \mbox{\boldmath $u$} $
com a {\em mesma\/} \\
constante $K\!\;\!$
do caso escalar.
Ademais,
tem-se
$ {\displaystyle
\;\!
\|\, D^m \:\! \mbox{\boldmath $u$}(\cdot,t)\,
\|_{\mbox{}_{\scriptstyle L^{\!q}}}
\!\:\!\leq\:\!
\|\, D^m \:\! \mbox{\boldmath $u$}(\cdot,t)\,
\|_{\mbox{}_{\scriptstyle L^{\!q_{\mbox{}_{\!\:\!1}}}}}^{1 - \theta}
\|\, D^m \:\! \mbox{\boldmath $u$}(\cdot,t)\,
\|_{\mbox{}_{\scriptstyle L^{\!q_{\mbox{}_{\!\:\!2}}}}}^{\:\!\theta}
\!\!\;\!
} $
se \linebreak
$ 1/q = (1 - \theta)/q_{\mbox{}_{1}} \!+ \theta/q_{\mbox{}_{2}} $,
$ 0 \leq \theta \leq 1 $,
e assim por diante.
}
%
%
%
Ocasionalmente,
resulta tamb\'em conveniente
usar a seguinte defini\c c\~ao alterna\-ti\-va
para a norma do sup de $ \mbox{\boldmath $u$}(\cdot,t) $, \\
\mbox{} \vspace{-0.675cm} \\
\begin{equation}
\tag{1.20}
\|\, \mbox{\boldmath $u$}(\cdot,t) \,
\|_{\mbox{}_{\scriptstyle \infty}}
=\;
\mbox{ess}\,\sup\; \bigl\{\:
|\,\mbox{\boldmath $u$}(x,t) \,|_{\mbox{}_{2}}
\!\!\:\!: \: x \in \mathbb{R}^{n}
\,\bigr\}.
\end{equation}
\mbox{} \vspace{-0.250cm} \\
Podemos tamb\'em utilizar
$ \|\,\mbox{\boldmath $u$}(\cdot,t)\,\|_{\mbox{}_{\scriptstyle L^q}}\!\;\!$
no lugar de
$\;\! \|\,\mbox{\boldmath $u$}(\cdot,t)\,
\|_{\mbox{}_{\scriptstyle L^q(\mathbb{R}^{n})}}\!\;\!$,
etc,
por simpli\-ci\-dade.
Constantes ser\~ao usualmente representadas
pelas letras
$ C \!\,\!$, $\!\:\!c$, $\!\:\!K$;
escrevemos
$ \:\!C(\lambda_{\mbox{}_{1}}\!\;\!, ...,
\lambda_{\mbox{}_{\scriptstyle k}})\!\;\! $
para observar que o valor da constante
$C$ em quest\~ao
depende apenas \linebreak
dos par\^ametros
$ \;\!\{\;\!\lambda_{\mbox{}_{1}}\!\;\!, ...,
\lambda_{\mbox{}_{\scriptstyle k}} \}\;\!$
indicados
(a menos que explicitamente mencionado
em contr\'ario).
Por conveni\^encia e economia,
usamos tipicamente o mesmo s\'\i mbolo
para denotar constantes com diferentes
valores num\'ericos
(por exemplo,
escrevemos $ C^2 \!\;\!$
ou $ 10 \,C + 1 $,
$ \mbox{cosh}\;C $, etc,
nova\-mente como $C\!\;\!$,
e assim por diante),
como usualmente feito na literatura. \\
%
%
\mbox{} \vspace{-0.450cm} \\
\nl
%
%
%
%
%
{\bf Agradecimentos.}
Parte das contribui\c c\~oes
feitas na {\sc Parte I}
n\~ao teria provavelmente sido
pos\-s\'\i vel
sem v\'arias ideias e m\'etodos
introduzidos em
\cite{KreissHagstromLorenzZingano2002, %
KreissHagstromLorenzZingano2003}.
O autor \'e especialmente grato a
Thomas Hagstrom e Jens Lorenz
pelas in\'umeras discuss\~oes
ocorridas
durante sua visita \`a
Universidade do Novo M\'exico
em 2001$\;\!$-$\;\!$2002.

%
%
%
%
%
\newpage
%
%

%
\mbox{} \vspace{-2.000cm} \\
%

%
%

{\bf 2. Preliminares, I} \\
\setcounter{section}{2}
\mbox{} \vspace{-0.450cm} \\

Nesta se\c c\~ao,
reunimos por conveni\^encia
v\'arios resultados b\'asicos
dados em
\cite{KreissHagstromLorenzZingano2002, %
KreissHagstromLorenzZingano2003, Leray1934}
que ter\~ao papel relevante na deriva\c c\~ao
das estimativas (1.10), (1.11), (1.12) \linebreak
nas se\c c\~oes seguintes.
No texto,
restringiremos nossa aten\c c\~ao
ao caso (fundamental) \linebreak
de dimens\~ao $ n = 3 $,
mas todos os argumentos usados
podem ser facilmente estendidos/adaptados
de modo
a se aplicarem a $ n = 2 $
igualmente,
com apenas pequenas mudan\c cas
\'obvias.
Em todo o texto,
$ \mbox{\boldmath $u$}(\cdot,t) $
sempre denotar\'a
uma solu\c c\~ao
de Leray (dada, qualquer)
para
as equa\c c\~oes (1.1),
mesmo que nada seja dito explicitamente. \linebreak
\mbox{} \vspace{-1.000cm} \\

Para a constru\c c\~ao
das solu\c c\~oes de Leray
$ \mbox{\boldmath $u$}(\cdot,t) $
do problema (1.1),
ver e.g.$\;$\cite{Galdi2000, Leray1934}
e a discuss\~ao abaixo.
Estas solu\c c\~oes
foram obtidas em \cite{Leray1934}
usando um procedimento
de re\-gu\-lariza\c c\~ao engenhoso
revisado a seguir.
Tomando
$ {\displaystyle
\;\!
G \in C^{\infty}_{0}(\mathbb{R}^{n})
\;\!
} $
n\~ao negativa (qualquer)
com
$ \!\;\!\int_{\mathbb{^R}^{3}} \!\:\!G(x) \,dx \:\!=\:\! 1 $,
e
definindo
$ {\displaystyle
\;\!
\bar{\mbox{\boldmath $u$}}_{\mbox{}_{\scriptstyle \!\;\!0, \,\delta}}
\!\;\!(\cdot)
\in C^{\infty}(\mathbb{R}^{3})
} $
pela
convolu\c c\~ao de
$ {\displaystyle
\;\!
\mbox{\boldmath $u$}_{0}(\cdot)
\;\!
} $
com
$ {\displaystyle
\;\!
G_{\mbox{}_{\scriptstyle \!\delta}}(x)
\:\!=\:\!
\delta^{\;\!-\,n} \;\!G(x/\delta)
} $,
$ \;\!\delta > 0 $,
introduz-se
$ {\displaystyle
\;\!
\mbox{\boldmath $u$}_{\mbox{}_{\scriptstyle \!\;\!\delta}},
\,
p_{\mbox{}_{\scriptstyle \!\;\!\delta}}
\in
C^{\infty}(\:\!\mathbb{R}^{3} \!\times\!\;\! [\,0, \infty\:\!)\:\!)
} $
como a solu\c c\~ao
(\'unica, cl\'assica, globalmente definida e em $L^{2}$)
do problema regularizado \\
\mbox{} \vspace{-0.575cm} \\
\begin{equation}
\tag{2.1$a$}
\mbox{\small $ {\displaystyle \frac{\partial}{\partial \;\!t} }$}
\,\mbox{\boldmath $u$}_{\mbox{}_{\scriptstyle \!\;\!\delta}}
\:\!+\:
\bar{\mbox{\boldmath $u$}}_{\mbox{}_{\scriptstyle \!\;\!\delta}}
(\cdot,t) \!\;\!
\cdot\!\;\! \nabla\;\!
\mbox{\boldmath $u$}_{\mbox{}_{\scriptstyle \!\;\!\delta}}
\:\!+\,
\nabla \:\!p_{\mbox{}_{\scriptstyle \!\;\!\delta}}
\;=\;
\nu \,
\Delta \:\!\mbox{\boldmath $u$}_{\mbox{}_{\scriptstyle \!\;\!\delta}},
\qquad
\nabla \!\cdot\!\;\!
\mbox{\boldmath $u$}_{\mbox{}_{\scriptstyle \!\;\!\delta}}(\cdot,t)
\,=\,0,
\end{equation}
\mbox{} \vspace{-0.950cm} \\
\begin{equation}
\tag{2.1$b$}
\mbox{\boldmath $u$}_{\mbox{}_{\scriptstyle \!\;\!\delta}}(\cdot,0)
\,=\,
\bar{\mbox{\boldmath $u$}}_{\mbox{}_{\scriptstyle \!\;\!0, \,\delta}}
\!\;\!:=\,
G_{\mbox{}_{\scriptstyle \!\delta}}
\!\:\!\ast
\mbox{\boldmath $u$}_{0}
\;\!
\in\!
\bigcap_{m\,=\,1}^{\infty}
\! H^{m}(\mathbb{R}^{3}),
\end{equation}
\mbox{} \vspace{-0.100cm} \\
onde
$ {\displaystyle
\;\!
\bar{\mbox{\boldmath $u$}}_{\mbox{}_{\scriptstyle \!\;\!\delta}}
\!\;\!(\cdot,t)
\!\;\!:=\:\!
G_{\mbox{}_{\scriptstyle \!\delta}} \!\,\!\ast
{\mbox{\boldmath $u$}}_{\mbox{}_{\scriptstyle \!\;\!\delta}}
\!\;\!(\cdot,t)
} $.
Em \cite{Leray1934},
Leray logrou mostrar que,
para uma sequ\^encia
$ {\displaystyle
\delta^{\;\!\prime}
\!\rightarrow
\:\!0
} $
apropriada,
tem-se a
converg\^encia fraca (em $L^{2}(\mathbb{R}^{3})$) \\
\mbox{} \vspace{-0.600cm} \\
\begin{equation}
\tag{2.2}
\mbox{\boldmath $u$}_{\mbox{}_{\scriptstyle \!\;\!\delta^{\;\!\prime}}}(\cdot,t)
\,\rightharpoonup \,
\mbox{\boldmath $u$}(\cdot,t)
\quad \;\,
\mbox{as } \;\,\delta^{\;\!\prime}
\!\rightarrow
\:\!0,
\qquad \;\;\;
\forall \;\, t \geq 0,
\end{equation}
\mbox{} \vspace{-0.200cm} \\
i.e.,
$ {\displaystyle
\mbox{\boldmath $u$}_{\mbox{}_{\scriptstyle \!\;\!\delta^{\;\!\prime}}}(\cdot,t)
\,\rightarrow \,
\mbox{\boldmath $u$}(\cdot,t)
\;\!
} $
fracamente em $L^{2}(\mathbb{R}^{3}) $,
para cada $ \;\! t \geq 0 $
(ver \cite{Leray1934}, p.$\;$237),
com \linebreak
$ {\displaystyle
\mbox{\boldmath $u$}(\cdot,t)
\in
L^{\infty}([\;\!0, \infty \:\!), L^{2}_{\sigma}(\mathbb{R}^{3}))
\cap
L^{2}([\;\!0, \infty \:\!), \mbox{$\stackrel{.}{H}$}\mbox{}^{1}(\mathbb{R}^{3}))
\cap
C^{0}_{w}([\;\!0, \infty \:\!), L^{2}(\mathbb{R}^{3}))
\:\!
} $
cont\'\i nua em $L^{2}(\mathbb{R}^{3})\!\;\! $
no instante $ t = 0 $
e resolvendo
a equa\c c\~ao (1.1$a$)
no sentido de distribui\c c\~oes.
Ademais, (1.2)
\'e satisfeita
para todo $ t \geq 0 $,
de modo que,
em particular,
tem-se \\
\mbox{} \vspace{-0.475cm} \\
\begin{equation}
\tag{2.3}
\int_{0}^{\:\!\infty} \!\!\!
\|\, D \mbox{\boldmath $u$}(\cdot,t) \,
\|_{\mbox{}_{\scriptstyle L^{2}(\mathbb{R}^{3})}}^{\:\!2}
\:\!dt
\;\leq\;
\mbox{\small $ {\displaystyle
\frac{1}{\:\!2\;\!\mbox{\normalsize $\nu$}\:\!} }$} \:
\|\: \mbox{\boldmath $u$}_0 \,
\|_{\mbox{}_{\scriptstyle L^{2}(\mathbb{R}^{3})}}^{\:\!2}
\!.
\end{equation}
\mbox{} \vspace{-0.150cm} \\
Outra propriedade importante
obtida em \cite{Leray1934}
\'e que
$ {\displaystyle
\;\!
\mbox{\boldmath $u$} \in
C^{\infty}(\:\!\mathbb{R}^{3} \!\times \!\;\!
[\,\mbox{\small $T$}_{\!\;\!\ast\ast}\!\;\!, \infty\:\!)\:\!)
\;\!
} $
para certo
$ \mbox{\small $T$}_{\!\:\!\ast\ast} \!\gg 1 $,
com
$ {\displaystyle
D^{m} \mbox{\boldmath $u$}(\cdot,t)
\in L^{\infty}_{\tt loc}
(\:\![\, \mbox{\small $T$}_{\!\:\!\ast\ast}\!\;\!,
 \infty \:\!), L^{2}(\mathbb{R}^{3})\:\!)
\;\!
} $
para cada $ \:\! m \geq 1 $.
Este fato \linebreak
permite simplificar
significativamente
o argumento desenvolvido
para (1.10)$\,-\,$(1.12)
mais adiante.
Outros resultados
importantes
referem-se \`a proje\c c\~ao
de Helmholtz
de
$ {\displaystyle
-\,
\mbox{\boldmath $u$}(\cdot,t)
\!\;\!\cdot\!\;\! \nabla \:\!
\mbox{\boldmath $u$}(\cdot,t)
} $
em
$ L^{2}_{\sigma}(\mathbb{R}^{3}) $,
ou seja,
o campo
$ \:\!\mbox{\boldmath $Q$}(\cdot,t) \in L^{2}_{\sigma}(\mathbb{R}^{3}) \:\!$
dado por \\
\mbox{} \vspace{-0.600cm} \\
\begin{equation}
\tag{2.4}
\mbox{\boldmath $Q$}(\cdot,t)
\;\! := \;
-\:
\mbox{\boldmath $u$}(\cdot,t)
\!\;\!\cdot\!\;\! \nabla \:\!
\mbox{\boldmath $u$}(\cdot,t)
\,-\;\! \nabla \:\!p\:\!(\cdot,t),
\qquad
\mbox{a.e. }\;\! t > 0.
\end{equation}
%
%
As estimativas que precisaremos
de
$ \:\!\mbox{\boldmath $Q$}(\cdot,t) \in L^{2}_{\sigma}(\mathbb{R}^{3}) \:\!$
acima
s\~ao revisadas a seguir. \linebreak
\nl
%
%
%
%
{\bf Proposi\c c\~ao 2.1.}
\textit{%
Para quase todo $ \;\! s > 0 $
$($e todo $\;\!s \geq \mbox{\small $T$}_{\!\;\!\ast\ast})$,
tem-se
} \\
\mbox{} \vspace{-0.600cm} \\
\begin{equation}
\tag{2.5}
\|\: e^{\;\!\nu \,\!\Delta
(\:\!\mbox{\footnotesize $t$} \;\!-\, \mbox{\footnotesize $s$})}
\,\! \mbox{\boldmath $Q$}(\cdot,s) \,
\|_{\mbox{}_{\scriptstyle L^{2}(\mathbb{R}^{3})}}
\leq\:
K \,
\nu^{\mbox{}^{\scriptstyle \!-\,\frac{3}{4} }}
\:\!
( t - s )^{\mbox{}^{\scriptstyle \! -\,\frac{3}{4} }}
\:\!
\|\, \mbox{\boldmath $u$}(\cdot,s) \,
\|_{\mbox{}_{\scriptstyle L^{2}(\mathbb{R}^{3})}}
\:\!
\|\, D \mbox{\boldmath $u$}(\cdot,s) \,
\|_{\mbox{}_{\scriptstyle L^{2}(\mathbb{R}^{3})}}
\end{equation}
\mbox{} \vspace{-0.200cm} \\
\textit{%
para todo $ \,t > s $,
onde
$ \:\!K \!\:\!=\:\! (\:\! 8 \:\!\pi )^{-\,3/4} \!\:\!$.
} \\
%
%
\nl
%
%
{\small
{\bf Prova:}
O argumento a seguir
\'e adaptado de
\cite{KreissHagstromLorenzZingano2002, %
KreissHagstromLorenzZingano2003}.
Seja
$ {\displaystyle
\;\!
\mathbb{F}\:\![\;\!f\;\!]
\equiv
\hat{f}
\;\!
} $
a transformada de Fourier
de uma dada fun\c c\~ao
$ \!\;\!f \!\;\!\in L^{1}(\mathbb{R}^{3}) $,
viz., \\
\mbox{} \vspace{-0.525cm} \\
\begin{equation}
\tag{2.6}
\mathbb{F}\:\![\;\!f\;\!]\:\!(k)
\:\equiv\;
\hat{f}(k)
\,:=\;
(\:\!2\:\!\pi)^{-\,3/2} \!\!
\int_{\mbox{}_{\scriptstyle \!\:\!\mathbb{R}^{3}}}
\!\!\:\!
e^{\mbox{\scriptsize $- \!\;\!
\stackrel{\mbox{\tiny $\circ$}}{\mbox{\sf \i$\!\;\!$\i}}
\!\:\! k \!\;\!\cdot\!\;\! x$}}
f(x) \: dx,
\qquad
k \in \mathbb{R}^{3}
%
%
\end{equation}
\mbox{} \vspace{-0.150cm} \\
(onde
$ {\displaystyle
\:\!
\stackrel{\mbox{\tiny $\circ$}}{\mbox{\sf \i$\!\;\!$\i}}
\mbox{}\!\!\;\!\mbox{}^{2}
\!\;\!=\;\! - \;\!1
} $).
Dada
$ {\displaystyle
\;\!
\mbox{\bf v}(\cdot,s) =
(\:\! \mbox{v}_{\mbox{}_{\!\;\!1}}\!\;\!(\cdot,s),
      \mbox{v}_{\mbox{}_{\!\;\!2}}\!\;\!(\cdot,s),
      \mbox{v}_{\mbox{}_{\!\;\!3}}\!\;\!(\cdot,s)\:\! )
\in
L^{1}(\mathbb{R}^{3}) \cap L^{2}(\mathbb{R}^{3})
\:\!
} $
arbitr\'aria,
ob\-t\'em-se,
pela identidade de Parseval, \\
\mbox{} \vspace{-0.520cm} \\
\begin{equation}
\notag
\begin{split}
\|\:  e^{\;\!\nu \:\!\Delta \mbox{\scriptsize $(t-s)$}} \;\!
\mbox{\bf v}(\cdot,s) \,
\|_{\mbox{}_{\scriptstyle L^{2}(\mathbb{R}^{3})}}^{\:\!2}
\;\!&=\;\;\!
\|\; \mathbb{F}\;\![\, e^{\;\! \nu \:\!\Delta \mbox{\scriptsize $(t-s)$}}
\;\!\mbox{\bf v}(\cdot,s) \: ] \,
\|_{\mbox{}_{\scriptstyle L^{2}(\mathbb{R}^{3})}}^{\:\!2} \\
&=\,
\int_{\mathbb{R}^{3}} \!\!
e^{-\,2\,\nu \;\!|\,\mbox{\scriptsize $k$}\,
|_{\mbox{}_{\!\;\!2}}^{\:\!\scriptstyle 2} \:\!(t \,-\, s)}
\;\!
|\, \hat{\mbox{\bf v}}(k,s) \,
|_{\mbox{}_{\!\;\!2}}^{\:\!2}
\, dk \\
&\leq\;
\|\, \hat{\mbox{\bf v}}(\cdot,s) \,
\|_{\mbox{}_{\scriptstyle \infty}}^{\:\!2} \!\!\;\!
\int_{\mathbb{R}^{3}} \!\!
e^{-\,2\,\nu \;\!|\,\mbox{\scriptsize $k$}\,
|_{\mbox{}_{\!\;\!2}}^{\:\!\scriptstyle 2} \:\!(t \,-\, s)}
\, dk \\
&=\;
\Bigl(\;\! \frac{\;\!\pi\;\!}{2} \,\Bigr)^{\!\!\;\!3/2}
\!\;\!
\nu^{-\,3/2}
\,( t - s )^{- \,3/2}
\,
\|\, \hat{\mbox{\bf v}}(\cdot,s) \,
\|_{\mbox{}_{\scriptstyle \infty}}^{\:\!2}
\!\;\!,
\end{split}
\end{equation}
\mbox{} \vspace{-0.300cm} \\
ou seja, \\
\mbox{} \vspace{-0.800cm} \\
\begin{equation}
\tag{2.7}
\|\:  e^{\;\!\nu\:\!\Delta \mbox{\scriptsize $(t-s)$}} \;\!
\mbox{\bf v}(\cdot,s) \,
\|_{\mbox{}_{\scriptstyle L^{2}(\mathbb{R}^{3})}}
\;\!\leq\;\;\!
\Bigl(\;\! \frac{\;\!\pi\;\!}{2} \,\Bigr)^{\!\!\;\!3/4}
\!\;\!
\nu^{-\,3/4} \,( t - s )^{- \,3/4}
\,
\|\, \hat{\mbox{\bf v}}(\cdot,s) \,
\|_{\mbox{}_{\scriptstyle \infty}}
\!\;\!,
\end{equation}
\mbox{} \vspace{-0.175cm} \\
onde
$ {\displaystyle
\;\!
|\;\!\cdot\;\!|_{\mbox{}_{2}}
} $
denota a norma Euclideana
em $ \mathbb{R}^{3}\!\;\!$,
e
$ {\displaystyle
\;\!
\|\, \hat{\mbox{\bf v}}(\cdot,s) \,
\|_{\mbox{}_{\scriptstyle \infty}}
\!=\,
\sup \;\{\:
|\, \hat{\mbox{\bf v}}(k,s) \,|_{\mbox{}_{2}}
\!: \, k\,\in\,\mathbb{R}^{3} \:\!\}
} $.
Como ser\'a mostrado abaixo,
(2.5) segue de
uma aplica\c c\~ao direta de (2.7)
a
$ {\displaystyle
\;\!
\mbox{\bf v}(\cdot,s) =
\mbox{\boldmath $Q$}(\cdot,s)
} $.
Para isso,
\'e preciso
que se estime
$ {\displaystyle
\;\!
\|\, \hat{\mbox{\boldmath $Q$}}(\cdot,s) \,
\|_{\mbox{}_{\scriptstyle \infty}}
\!\;\!
} $:
como tem-se
$ {\displaystyle
\;\!
\mathbb{F}\:\![\,\nabla \!\:\!
P(\cdot,s)\,]\:\!(k) \;\!=\;\;\!
\stackrel{\mbox{\tiny o}}{\mbox{\sf \i$\!\;\!$\i}} \!\:\!
\hat{p}(k,s) \;\!k
\,
} $
e
$ {\displaystyle
\:\!
\mbox{\small $\sum$}_{j\,=\,1}^{3} k_{j}
\hat{Q}_{j}(k,s)
=\;\!0
} $
(pois
$ \nabla \!\:\!\cdot \mbox{\boldmath $Q$}(\cdot,s) =\;\! 0 $),
segue que
$ {\displaystyle
\;\!
\mathbb{F}\:\![\,\nabla \!\:\!
P(\cdot,s)\,]\:\!(k)
} $
e
$ {\displaystyle
\:\!
\hat{\mbox{\boldmath $Q$}}(k,s)
\;\!
} $
s\~ao vetores
ortogonais em $\mathbb{C}^{3}\!$,
para todo $k \in \mathbb{R}^{3}\!$.
Lembrando,
por (2.4),
que
$ {\displaystyle
\;\!
\hat{\mbox{\boldmath $Q$}}(k,s)
+\:\!
\mathbb{F}\:\![\,\nabla \!\:\!
P(\cdot,s)\;\!]\:\!(k)
=
} $ \linebreak
$ {\displaystyle
-\:
\mathbb{F}\;\![\, \mbox{\boldmath $u$}(\cdot,s) \!\;\!\cdot \!\;\!\nabla
\mbox{\boldmath $u$}(\cdot,s)\;\!]\:\!(k)
} $,
obt\'em-se \\
\mbox{} \vspace{-0.600cm} \\
\begin{equation}
\tag{2.8}
|\,\hat{\mbox{\boldmath $Q$}}(k,s) \,|_{\mbox{}_{2}}
\,\leq\;\,\!
|\: \mathbb{F}\;\![\, \mbox{\boldmath $u$}(\cdot,s)\!\;\!\cdot \!\;\!\nabla
\mbox{\boldmath $u$}\;\!(\cdot,s)\;\!]\:\!(k)
\:|_{\mbox{}_{2}}
\end{equation}
\mbox{} \vspace{-0.275cm} \\
para todo $\;\! k \in \mathbb{R}^{3}\!$,
de modo que \\
\mbox{} \vspace{-0.600cm} \\
\begin{equation}
\tag{2.9}
\|\:
\hat{\mbox{\boldmath $Q$}}(\cdot,s)
\,\|_{\mbox{}_{\scriptstyle \infty}}
\;\!\leq\;
\|\;
\mathbb{F}\;\![\, \mbox{\boldmath $u$}
\!\;\!\cdot \!\;\!\nabla
\mbox{\boldmath $u$}\,]\:\!(\cdot,s)
\:\|_{\mbox{}_{\scriptstyle \infty}}.
\end{equation}
\mbox{} \vspace{-0.200cm} \\
Por outro lado,
tem-se,
para $ 1 \leq i \leq 3 $, \\
\mbox{} \vspace{-0.750cm} \\
\begin{equation}
\notag
\begin{split}
|\; \mathbb{F}\;\![\, \mbox{\boldmath $u$}(\cdot,s)
\!\;\!\cdot \!\;\!\nabla
\;\!u_{{\scriptstyle i}}(\cdot,s)\,]\:\!(k) \:|
\;\,&\leq\;
\sum_{j\,=\,1}^{3} \;\!
|\; \mathbb{F}\;\![\, u_{{\scriptstyle j}}(\cdot,s) \;\!
D_{{\scriptstyle \!\;\!j}} \:\! u_{{\scriptstyle i}}(\cdot,s)
\,]\:\!
(k) \:| \\
&\leq\;
(\:\!2\:\!\pi)^{-\,3/2} \;\!
\sum_{j\,=\,1}^{3} \;\!
\|\, u_{{\scriptstyle j}}(\cdot,s) \;\!
D_{{\scriptstyle \!\;\!j}} \:\! u_{{\scriptstyle i}}(\cdot,s) \,
\|_{\mbox{}_{\scriptstyle L^{1}(\mathbb{R}^{3})}} \\
&\leq\;
(\:\!2\:\!\pi)^{-\,3/2} \,
\|\, \mbox{\boldmath $u$}(\cdot,s) \,
\|_{\mbox{}_{\scriptstyle L^{2}(\mathbb{R}^{3})}}
\;\!
\|\, \nabla \!\;\! u_{{\scriptstyle i}}(\cdot,s) \,
\|_{\mbox{}_{\scriptstyle L^{2}(\mathbb{R}^{3})}}
\!\;\!,
\end{split}
\end{equation}
\mbox{} \vspace{+0.125cm} \\
usando Cauchy-Schwarz.
$\!$($\;\!$Aqui,
como sempre,
$ D_{\!\;\!j} \!\;\!=\,\! \partial/\partial x_{\!\;\!j} \!\;\!$.)
Isso fornece \\
\mbox{} \vspace{-0.550cm} \\
\begin{equation}
\tag{2.10}
\|\; \mathbb{F}\;\![\, \mbox{\boldmath $u$} \!\;\!\cdot \!\;\!\nabla
\;\!\mbox{\boldmath $u$}\,]\:\!(\cdot,s) \:\|_{\mbox{}_{\scriptstyle \infty}}
\:\leq\;
(\:\!2\:\!\pi)^{-\,3/2} \,
\|\, \mbox{\boldmath $u$}(\cdot,s) \,
\|_{\mbox{}_{\scriptstyle L^{2}(\mathbb{R}^{3})}}
\;\!
\|\, D \mbox{\boldmath $u$}(\cdot,s) \,
\|_{\mbox{}_{\scriptstyle L^{2}(\mathbb{R}^{3})}}
\!\:\!.
\end{equation}
\mbox{} \vspace{-0.200cm} \\
De (2.7), (2.9) and (2.10),
obt\'em-se (2.5),
concluindo a prova da Proposi\c c\~ao 2.1.
}
\mbox{} \hfill $\Box$ \\
%
%
\mbox{} \vspace{-0.750cm} \\

Note-se que,
repetindo o argumento acima
para as solu\c c\~oes
das equa\c c\~oes regularizadas (2.1),
obt\'em-se
de modo an\'alogo que \\
\mbox{} \vspace{-0.600cm} \\
\begin{equation}
\tag{2.11}
\|\: e^{\;\!\nu \:\!
\mbox{\scriptsize $\Delta$} \mbox{\scriptsize $(t - s)$}} \;\!
\mbox{\boldmath $Q$}_{\mbox{}_{\scriptstyle \!\delta}}\!\;\!(\cdot,s) \,
\|_{\mbox{}_{\scriptstyle L^{2}(\mathbb{R}^{3})}}
\leq\,
K \,
\nu^{\mbox{}^{\scriptstyle \!-\, \frac{3}{4}}}
(t - s)^{\mbox{}^{\scriptstyle \!\!-\, \frac{3}{4}}}
\|\, \mbox{\boldmath $u$}_{\mbox{}_{\scriptstyle \!\:\!\delta}}
\!\;\!(\cdot,s) \,
\|_{\mbox{}_{\scriptstyle L^{2}(\mathbb{R}^{3})}}
\;\!
\|\, D \mbox{\boldmath $u$}_{\mbox{}_{\scriptstyle \!\:\!\delta}}
\!\;\!(\cdot,s) \,
\|_{\mbox{}_{\scriptstyle L^{2}(\mathbb{R}^{3})}}
\end{equation}
\mbox{} \vspace{-0.250cm} \\
para cada $\;\!t > s > 0 $,
sendo
$ K \!\;\!=\:\! (\:\!8 \:\!\pi)^{-\,3/4} \!$,
como antes,
e \\
\mbox{} \vspace{-0.600cm} \\
\begin{equation}
\tag{2.12}
\mbox{\boldmath $Q$}_{\mbox{}_{\scriptstyle \!\delta}}
\!\:\!(\cdot,s)
\;=\; -\:
\bar{\mbox{\boldmath $u$}}_{\mbox{}_{\scriptstyle \!\:\!\delta}}\!\:\!(\cdot,s)
\!\;\!\cdot \!\;\!\nabla
\mbox{\boldmath $u$}_{\mbox{}_{\scriptstyle \!\:\!\delta}}\!\:\!(\cdot,s)
\,-\;\!
\nabla p_{\mbox{}_{\scriptstyle \!\delta}}\!\;\!(\cdot,s).
\end{equation}
\mbox{} \vspace{-0.250cm} \\
A estimativa (2.11)
\'e muito \'util,
visto que
%
%
as solu\c c\~oes regularizadas
$ {\displaystyle
\:\!
\mbox{\boldmath $u$}_{\mbox{}_{\scriptstyle \!\:\!\delta}}\!\:\!(\cdot,t)
\:\!
} $
definidas em~(2.1)
satisfazem
a desigualdade de energia \\
\mbox{} \vspace{-0.675cm} \\
\begin{equation}
\tag{2.13}
\|\,\mbox{\boldmath $u$}_{\mbox{}_{\scriptstyle \!\:\!\delta}}
\!\:\!(\cdot,t) \,
\|_{\mbox{}_{\scriptstyle L^{2}(\mathbb{R}^{3})}}^{\:\!2}
+\;
2 \, \nu \!\!\;\!
\int_{0}^{\;\!\mbox{\mbox{\footnotesize $t$}}} \!\!\;\!
\|\, D \mbox{\boldmath $u$}_{\mbox{}_{\scriptstyle \!\:\!\delta}}
\!\:\!(\cdot,s) \,
\|_{\mbox{}_{\scriptstyle L^{2}(\mathbb{R}^{3})}}^{\:\!2}
\;\!ds
\;\;\!\leq\:
\|\, \mbox{\boldmath $u$}_{0} \,
\|_{\mbox{}_{\scriptstyle L^{2}(\mathbb{R}^{3})}}^{\:\!2}
\end{equation}
\mbox{} \vspace{-0.100cm} \\
para todo
$\;\! t > 0 $
(e todo
$ \delta > 0 $),
de modo que
$ {\displaystyle
\;\!
\|\,\mbox{\boldmath $u$}_{\mbox{}_{\scriptstyle \!\:\!\delta}}
\!\:\!(\cdot,t) \,
\|_{\mbox{}_{\scriptstyle L^{2}(\mathbb{R}^{3})}}
\!\;\!
} $,
$
\int_{0}^{\;\!\mbox{\mbox{\footnotesize $t$}}} \!
\|\, D \mbox{\boldmath $u$}_{\mbox{}_{\scriptstyle \!\:\!\delta}}
\!\:\!(\cdot,s) \,
\|_{\mbox{}_{\scriptstyle L^{2}(\mathbb{R}^{3})}}^{\:\!2}
\!\;\!ds
$ \linebreak
\mbox{} \vspace{-0.550cm} \\
podem ser estimadas
independentemente de $ \:\!\delta > 0 $.
Isso ser\'a usado no Teorema 3.1
(Se\c c\~ao~3)
para mostrar que
a escolha particular
de
$\;\!t_0 \!\;\!\geq 0 \;\!$
ao se definir
as aproxima\-\c c\~oes
%
%
(1.5)
n\~ao \'e relevante
com respeito \`as
propriedades
(1.10$b$), (1.12$b$).
Conv\'em tamb\'em
generalizar
a estimativa (2.5)
para derivadas
de ordem superior.
No caso da equa\c c\~ao do calor,
ser\'a \'util
lembrarmos
a seguinte estimativa (bem conhecida): \\
\mbox{} \vspace{-0.650cm} \\
\begin{equation}
\tag{2.14}
\|\, D^{\alpha} \:\![\,
e^{\;\!\nu \:\!\Delta \tau} \,\!
\mbox{u} \,]\,
\|_{\mbox{}_{\scriptstyle L^{2}(\mathbb{R}^{n})}}
\;\!\leq\;
K\!\;\!(n,\:\!m)
\;
\|\: \mbox{u} \:
\|_{\mbox{}_{\scriptstyle L^{r}(\mathbb{R}^{n})}}
\:\!
(\;\!\nu \;\!\tau\:\!)^{\mbox{}^{\scriptstyle
\!\! -\, \frac{\scriptstyle n}{2}
\left( \frac{1}{\scriptstyle r} \;\!-\;\! \frac{1}{2} \right)
\,-\, \frac{|\:\!{\scriptstyle \alpha}\:\!|}{2} }}
%
\end{equation}
\mbox{} \vspace{-0.125cm} \\
para todo $ \tau > 0 $,
e quaisquer
$ \alpha $ (multi-\'\i ndice),
$ 1 \leq r \leq 2 $,
$ \mbox{u} \in L^{r}(\mathbb{R}^{n}) $
considerados,
$ n \geq 1 $ arbitr\'ario,
e onde
$ \;\!m = |\;\!\alpha\:\!| $.
($\:\!$Para uma deriva\c c\~ao de (2.14),
ver e.g.$\;$\cite{KreissLorenz1989, %
LorenzZingano2012}.) \linebreak
%
%
%
%
%
{\bf Proposi\c c\~ao 2.2.}
\textit{%
%
%
Para quase todo $ \;\! s > 0 $
$($e todo $\;\!s \geq \mbox{\small $T$}_{\!\;\!\ast\ast})$,
tem-se
} \\
\mbox{} \vspace{-0.600cm} \\
\begin{equation}
\tag{2.15}
\|\, D^{\alpha} \{\, e^{\;\!\nu \,\!\Delta
(\:\!\mbox{\footnotesize $t$} \;\!-\, \mbox{\footnotesize $s$})}
\,\! \mbox{\boldmath $Q$}(\cdot,s) \,\}\,
\|_{\mbox{}_{\scriptstyle L^{2}(\mathbb{R}^{3})}}
\leq\:
K\!\;\!(m) \,\;\!
\nu^{\mbox{}^{\scriptstyle \!-\,\gamma }}
\:\!
( t - s )^{\mbox{}^{\scriptstyle \! -\,\gamma}}
\:\!
\|\, \mbox{\boldmath $u$}(\cdot,s) \,
\|_{\mbox{}_{\scriptstyle L^{2}(\mathbb{R}^{3})}}
\:\!
\|\, D \mbox{\boldmath $u$}(\cdot,s) \,
\|_{\mbox{}_{\scriptstyle L^{2}(\mathbb{R}^{3})}}
\end{equation}
\mbox{} \vspace{-0.500cm} \\
\textit{%
para todo $ \,t > s $,
onde
$ m = |\;\!\alpha\;\!| $,
$ \gamma = m/2 + 3/4 $,
e $ \:\!K\!\;\!(m) $
depende apenas de $\;\!m$.
} \\
%
%
\mbox{} \vspace{-0.050cm} \\
%
%
{\small
{\bf Prova:}
Este resultado \'e uma consequ\^encia direta
de (2.5) e (2.14) acima.
De fato, tem-se: \\
\mbox{} \vspace{-0.150cm} \\
\mbox{} \hspace{+1.175cm}
$ {\displaystyle
\|\, D^{\alpha} \{\, e^{\;\!\nu \,\!\Delta
(\:\!\mbox{\footnotesize $t$} \;\!-\, \mbox{\footnotesize $s$})}
\,\! \mbox{\boldmath $Q$}(\cdot,s) \,\}\,
\|_{\mbox{}_{\scriptstyle L^{2}(\mathbb{R}^{3})}}
\leq
} $ \\
\mbox{} \vspace{-0.100cm} \\
\mbox{} \hspace{+2.750cm}
$ {\displaystyle
\leq \:
K\!\:\!(m)
\,\;\!
\nu^{\mbox{}^{\scriptstyle \!-\, \frac{m}{2} }}
(t - s)^{\mbox{}^{\scriptstyle \!-\, \frac{m}{2} }}
\|\: e^{\;\!\nu \,\!\Delta
(\:\!\mbox{\footnotesize $t$} \;\!-\, \mbox{\footnotesize $s$})/2}
\,\! \mbox{\boldmath $Q$}(\cdot,s) \,\}\,
\|_{\mbox{}_{\scriptstyle L^{2}(\mathbb{R}^{3})}}
} $
\mbox{} \hfill
\mbox{[}$\;\!$por (2.14)$\;\!$\mbox{]} \\
\mbox{} \vspace{-0.120cm} \\
\mbox{} \hspace{+2.750cm}
$ {\displaystyle
\leq \:
K\!\:\!(m)
\,\;\!
\nu^{\mbox{}^{\scriptstyle \!-\, \frac{m}{2} \,-\, \frac{3}{4} }}
(t - s)^{\mbox{}^{\scriptstyle \!-\, \frac{m}{2} \,-\, \frac{3}{4} }}
\|\, \mbox{\boldmath $u$}(\cdot,s) \,
\|_{\mbox{}_{\scriptstyle L^{2}(\mathbb{R}^{3})}}
\;\!
\|\, D \:\!\mbox{\boldmath $u$}(\cdot,s) \,
\|_{\mbox{}_{\scriptstyle L^{2}(\mathbb{R}^{3})}}
\!\;\!
} $.
\mbox{} \hfill
\mbox{[}$\;\!$(2.5)$\;\!$\mbox{]} \\
}
\mbox{} \hfill $ \Box $ \\
%
\mbox{} \vspace{-0.650cm} \\

Para fins do pr\'oximo resultado
a ser revisado nesta se\c c\~ao,
dado na Proposi\c c\~ao~2.3,
precisaremos das seguintes
desigualdades elementares
de Nirenberg-Gagliardo
(\mbox{\small SNG})
para
fun\c c\~oes
$ {\displaystyle
\;\!
\mbox{u} \in H^{2}(\mathbb{R}^{3})
} $
quaisquer: \\
\mbox{} \vspace{-0.600cm} \\
\begin{equation}
\tag{2.16$a$}
\|\: \mbox{u} \:\|_{\mbox{}_{\scriptstyle \infty}}
\;\!\leq\,
K_{\mbox{}_{\!\;\!0}} \,
\|\: \mbox{u} \:
\|_{\mbox{}_{\scriptstyle L^{2}(\mathbb{R}^{3})}}^{\:\!1/4}
\;\!
\|\: D^{2} \mbox{u} \:
\|_{\mbox{}_{\scriptstyle L^{2}(\mathbb{R}^{3})}}^{\:\!3/4}
\!\:\!,
\qquad
K_{\mbox{}_{\!\;\!0}} \!\;\!<\;\! 0.678,
\end{equation}
\mbox{} \vspace{-0.150cm} \\
ver
e.g.$\;$\mbox{\{}$\,$\cite{Taylor2011}, Proposition 2.4, p.$\;$5$\,$\mbox{\}},
ou \mbox{\{}$\,$\cite{Schutz2008}, Teorema 4.5.1, p.$\;$52$\,$\mbox{\}};
$\;\!$e\\
\mbox{} \vspace{-0.600cm} \\
\begin{equation}
\tag{2.16$b$}
\|\, D \:\!\mbox{u} \:
\|_{\mbox{}_{\scriptstyle L^{2}(\mathbb{R}^{3})}}
\leq\,
K_{\mbox{}_{\!\;\!1}} \;\!
\|\: \mbox{u} \:
\|_{\mbox{}_{\scriptstyle L^{2}(\mathbb{R}^{3})}}^{\:\!1/2}
\;\!
\|\: D^{2} \mbox{u} \:
\|_{\mbox{}_{\scriptstyle L^{2}(\mathbb{R}^{3})}}^{\:\!1/2}
\!\!\;\!,
\qquad
K_{\mbox{}_{\!\;\!1}} \!\;\!=\;\! 1,
\end{equation}
\mbox{} \vspace{-0.150cm} \\
facilmente obtida
usando a transformada de Fourier.
De (2.16$a$), (2.16$b$),
obt\'em-se \\
\mbox{} \vspace{-0.600cm} \\
\begin{equation}
\tag{2.17}
\|\: \mbox{u} \:\|_{\mbox{}_{\scriptstyle \infty}}
\;\!
\|\, D \,\!\mbox{u} \:
\|_{\mbox{}_{\scriptstyle L^{2}(\mathbb{R}^{3})}}^{\:\!1/2}
\leq\,
K_{\mbox{}_{\!\;\!2}} \;\!
\|\: \mbox{u} \:
\|_{\mbox{}_{\scriptstyle L^{2}(\mathbb{R}^{3})}}^{\:\!1/2}
\;\!
\|\: D^{2} \mbox{u} \:
\|_{\mbox{}_{\scriptstyle L^{2}(\mathbb{R}^{3})}}
\!\:\!,
\quad \;\,
K_{\mbox{}_{\!\;\!2}} \!\;\!=\;\!
K_{\mbox{}_{\!\;\!0}} \;\!
K_{\mbox{}_{\!\;\!1}}^{\:\!1/2}
\!\:\!< 1.
\end{equation}
\mbox{} \vspace{-0.150cm} \\
Lembramos tamb\'em a defini\c c\~ao
de $ \mbox{\small $T$}_{\!\:\!\ast\ast} $
dada pela propriedade (1.3) na Se\c c\~ao~1. \\
\nl
%
%
%
%
{\bf Proposi\c c\~ao 2.3.}
\textit{%
Seja
$ \,\mbox{\boldmath $u$}(\cdot,t) $
solu\c c\~ao de Leray
para $\;\!(1.1)$.
Ent\~ao,
existe
$ \,t_{\!\;\!\ast\ast} \!\;\!\geq \mbox{\small $T$}_{\!\:\!\ast\ast}$
$($com $\:\!t_{\!\;\!\ast\ast}\!\:\!$
dependendo da solu\c c\~ao
$\;\!\mbox{\boldmath $u$}$\/$)$
suficientemente grande
tal que
$ {\displaystyle
\|\, D \mbox{\boldmath $u$}(\cdot,t) \,
\|_{{\scriptstyle L^{2}(\mathbb{R}^{3})}}
\!\:\!
} $
\'e uma fun\c c\~ao suave e
monotonicamente decrescente
de $ \, t $
no intervalo
$\:\! [\,t_{\!\;\!\ast\ast}\!\;\!, \:\!\mbox{\small $\infty$}\:\!)$.
} \\
%
%
%
\mbox{} \vspace{-0.000cm} \\
{\small
{\bf Prova:}
O argumento abaixo
\'e adaptado da prova de
\mbox{\{}$\,$\cite{KreissHagstromLorenzZingano2002},
Lemma 2.2$\,$\mbox{\}}.
Considere-se
$ {\displaystyle
\;\!
t_{0}
\geq\;\!
\mbox{\small $T$}_{\!\:\!\ast\ast}
} $
(a ser escolhido abaixo),
e seja
$ \;\! t > t_{0} $.
Aplicando
$ {\displaystyle
D_{\mbox{}_{\scriptstyle \!\ell}} \!\;\!=\:\!
\partial/\partial \:\!x_{\mbox{}_{\scriptstyle \!\ell}}
} $
\`a primeira equa\c c\~ao em (1.1$a$),
tomando o produto escalar com
$ {\displaystyle
\:\!D_{\mbox{}_{\scriptstyle \!\ell}}
\mbox{\boldmath $u$}(\cdot,t)
} $
e
integrando em
$ {\displaystyle
\mathbb{R}^{3} \!\times\!\;\!
[\, t_{0}, \;\!t\;\!]
} $,
obt\'em-se,
somando em $ 1 \leq \ell \leq 3 $, \\
\newpage
\mbox{} \vspace{-0.750cm} \\
\mbox{} \hspace{+3.500cm}
$ {\displaystyle
\|\, D \mbox{\boldmath $u$}(\cdot,t) \,
\|_{\mbox{}_{\scriptstyle L^{2}(\mathbb{R}^{3})}}^{\:\!2}
+\:
2 \: \nu
\!\!\:\!
\int_{\mbox{\footnotesize $ t_{0} $}}
    ^{\mbox{\footnotesize $\:\!t$}}
\!
\|\, D^{2} \mbox{\boldmath $u$}(\cdot,s) \,
\|_{\mbox{}_{\scriptstyle L^{2}(\mathbb{R}^{3})}}^{\:\!2}
ds
\;\;\!=
} $ \\
\mbox{} \vspace{+0.020cm} \\
\mbox{} \hspace{+1.000cm}
$ {\displaystyle
=\;\:\!
\|\, D \mbox{\boldmath $u$}(\cdot,t_{0}) \,
\|_{\mbox{}_{\scriptstyle L^{2}(\mathbb{R}^{3})}}^{\:\!2}
\!\;\!+\:
2
\sum_{i, \, j, \, \ell}
\int_{\mbox{\footnotesize $ t_{0} $}}
    ^{\mbox{\footnotesize $\:\!t$}}
\!
\int_{\mathbb{R}^{3}}
\!\!\!\;\!
u_{i}(x,s) \;\!
D_{\mbox{}_{\scriptstyle \!\!\;\!\ell}} u_{j}(x,s)
\;\!
D_{\scriptstyle \!j}
D_{\mbox{}_{\scriptstyle \!\!\;\!\ell}}
u_{i}(x,s) \,
dx \: ds
} $ \\
\mbox{} \vspace{+0.050cm} \\
\mbox{} \hspace{+1.000cm}
$ {\displaystyle
\leq\;
\|\, D \mbox{\boldmath $u$}(\cdot,t_{0}) \,
\|_{\mbox{}_{\scriptstyle L^{2}(\mathbb{R}^{3})}}^{\:\!2}
+\:
2 \!\!\:\!
\int_{\mbox{\footnotesize $ t_{0} $}}
    ^{\mbox{\footnotesize $\:\!t$}}
\!
\|\, \mbox{\boldmath $u$}(\cdot,s) \,
\|_{\mbox{}_{\scriptstyle \infty}}
\:\!
\|\, D \mbox{\boldmath $u$}(\cdot,s) \,
\|_{\mbox{}_{\scriptstyle L^{2}(\mathbb{R}^{3})}}
\|\, D^{2} \mbox{\boldmath $u$}(\cdot,s) \,
\|_{\mbox{}_{\scriptstyle L^{2}(\mathbb{R}^{3})}}
\;\!
ds
} $ \\
\mbox{} \vspace{+0.025cm} \\
\mbox{} \hspace{+1.000cm}
$ {\displaystyle
\leq\;
\|\, D \mbox{\boldmath $u$}(\cdot,t_{0}) \,
\|_{\mbox{}_{\scriptstyle L^{2}(\mathbb{R}^{3})}}^{\:\!2}
+\:
2 \!\!\:\!
\int_{\mbox{\footnotesize $ t_{0} $}}
    ^{\mbox{\footnotesize $\:\!t$}}
\!
\|\, \mbox{\boldmath $u$}(\cdot,s) \,
\|_{\mbox{}_{\scriptstyle L^{2}(\mathbb{R}^{3})}}^{\:\!1/2}
\:\!
\|\, D \mbox{\boldmath $u$}(\cdot,s) \,
\|_{\mbox{}_{\scriptstyle L^{2}(\mathbb{R}^{3})}}^{\:\!1/2}
\|\, D^{2} \mbox{\boldmath $u$}(\cdot,s) \,
\|_{\mbox{}_{\scriptstyle L^{2}(\mathbb{R}^{3})}}^{\:\!2}
\;\!
ds
} $, \\
\mbox{} \vspace{+0.100cm} \\
por (2.17),
lembrando
(1.19) e (1.20).
Em particular,
tem-se \\
\mbox{} \vspace{-0.150cm} \\
\mbox{} \hspace{+3.300cm}
$ {\displaystyle
\|\, D \mbox{\boldmath $u$}(\cdot,t) \,
\|_{\mbox{}_{\scriptstyle L^{2}(\mathbb{R}^{3})}}^{\:\!2}
+\:
2 \, \nu \!\!\;\!
\int_{\mbox{\footnotesize $ t_{0} $}}
    ^{\mbox{\footnotesize $\:\!t$}}
\!
\|\, D^{2} \mbox{\boldmath $u$}(\cdot,s) \,
\|_{\mbox{}_{\scriptstyle L^{2}(\mathbb{R}^{3})}}^{\:\!2}
ds
\;\;\!\leq
} $ \\
\mbox{} \vspace{-0.700cm} \\
\mbox{} \hfill (2.18) \\
\mbox{} \vspace{-0.400cm} \\
\mbox{} \hspace{+0.300cm}
$ {\displaystyle
\leq\;
\|\, D \mbox{\boldmath $u$}
(\cdot,t_{0}) \,
\|_{\mbox{}_{\scriptstyle L^{2}(\mathbb{R}^{3})}}^{\:\!2}
+\:
2 \!
\int_{\mbox{\footnotesize $ t_{0} $}}
    ^{\mbox{\footnotesize $\:\!t$}}
\Bigl[\;
\|\, \mbox{\boldmath $u$}_{0} \,
\|_{\mbox{}_{\scriptstyle L^{2}(\mathbb{R}^{3})}}
\:\!
\|\, D \mbox{\boldmath $u$}(\cdot,s) \,
\|_{\mbox{}_{\scriptstyle L^{2}(\mathbb{R}^{3})}}
\;\!
\Bigr]^{\!1/2}
\|\, D^{2} \mbox{\boldmath $u$}(\cdot,s) \,
\|_{\mbox{}_{\scriptstyle L^{2}(\mathbb{R}^{3})}}^{\:\!2}
ds
} $ \\
\mbox{} \vspace{+0.150cm} \\
para todo
$ {\displaystyle
\;\!
t \geq t_{0}
} $.
$\!$Seja ent\~ao
$ \;\! t_{0} \geq t_{\!\;\!\ast} \;\!$
tal que,
por (1.2):
$ {\displaystyle
\;\!
\|\, \mbox{\boldmath $u$}_{0} \,
\|_{\mbox{}_{\scriptstyle L^{2}(\mathbb{R}^{3})}}
\,\!
\|\, D \mbox{\boldmath $u$}(\cdot,t_{0}) \,
\|_{\mbox{}_{\scriptstyle L^{2}(\mathbb{R}^{3})}}
\!\!\;\!< \nu
} $. \linebreak
\mbox{} \vspace{-0.550cm} \\
De fato,
com esta escolha,
segue de (2.18)
que \\
\mbox{} \vspace{-0.625cm} \\
\begin{equation}
\tag{2.19$a$}
\|\, \mbox{\boldmath $u$}_{0} \,
\|_{\mbox{}_{\scriptstyle L^{2}(\mathbb{R}^{3})}}
\,\!
\|\, D \mbox{\boldmath $u$}(\cdot,s) \,
\|_{\mbox{}_{\scriptstyle L^{2}(\mathbb{R}^{3})}}
<\;\! \nu
\qquad
\forall \;\;\!s \geq t_0.
\end{equation}
\mbox{} \vspace{-0.200cm} \\
\mbox{$[\,$}Prova
de (2.19$a$): sendo falso,
existiria
$\;\! t_{1} \!\;\!> t_{0} \;\!$
tal que
$ {\displaystyle
\;\!
\|\, \mbox{\boldmath $u$}_{0} \,
\|_{\mbox{}_{\scriptstyle L^{2}(\mathbb{R}^{3})}}
\,\!
\|\, D \mbox{\boldmath $u$}(\cdot,s) \,
\|_{\mbox{}_{\scriptstyle L^{2}(\mathbb{R}^{3})}}
\!< \nu
} $ \linebreak
\mbox{} \vspace{-0.530cm} \\
para todo
$ \;\!t_0 \leq s < t_{1}$,
com
$ {\displaystyle
\;\!
\|\, \mbox{\boldmath $u$}_{0} \,
\|_{\mbox{}_{\scriptstyle L^{2}(\mathbb{R}^{3})}}
\,\!
\|\, D \mbox{\boldmath $u$}(\cdot,t_{1}) \,
\|_{\mbox{}_{\scriptstyle L^{2}(\mathbb{R}^{3})}}
\!= \;\!\nu
} $.
Tomando $\;\! t = t_{1} $
em (2.18), \linebreak
\mbox{} \vspace{-0.500cm} \\
resultaria
$ {\displaystyle
\;\!
\|\, D \mbox{\boldmath $u$}(\cdot,t_1) \,
\|_{\mbox{}_{\scriptstyle L^{2}(\mathbb{R}^{3})}}
\!\;\!\leq
\|\, D \mbox{\boldmath $u$}(\cdot,t_0) \,
\|_{\mbox{}_{\scriptstyle L^{2}(\mathbb{R}^{3})}}
\!\:\!
} $,
e, ent\~ao,
$ {\displaystyle
\;\!
\|\, \mbox{\boldmath $u$}_0 \,
\|_{\mbox{}_{\scriptstyle L^{2}(\mathbb{R}^{3})}}
\!\;\!
\|\, D \mbox{\boldmath $u$}(\cdot,t_1) \,
\|_{\mbox{}_{\scriptstyle L^{2}(\mathbb{R}^{3})}}
\!\,\!\leq
} $ \\
\mbox{} \vspace{-0.530cm} \\
$ {\displaystyle
\|\, \mbox{\boldmath $u$}_0 \,
\|_{\mbox{}_{\scriptstyle L^{2}(\mathbb{R}^{3})}}
\|\, D \mbox{\boldmath $u$}(\cdot,t_0) \,
\|_{\mbox{}_{\scriptstyle L^{2}(\mathbb{R}^{3})}}
\!< \;\!\nu
} $.
$\!$Esta contradi\c c\~ao mostra (2.19$a$),
como afirmado.
\mbox{\footnotesize $\Box$}$\,]$ \linebreak
%
\mbox{} \vspace{-0.700cm} \\

De (2.18) e (2.19$a$),
segue que \\
\mbox{} \vspace{-0.650cm} \\
\begin{equation}
\tag{2.19$b$}
\|\, D \mbox{\boldmath $u$}(\cdot,t) \,
\|_{\mbox{}_{\scriptstyle L^{2}(\mathbb{R}^{3})}}^{\:\!2}
+\:
2 \, \gamma
\!\!\;\!
\int_{\mbox{\footnotesize $ t_{2} $}}
    ^{\mbox{\footnotesize $\:\!t$}}
\!
\|\, D^{2} \mbox{\boldmath $u$}(\cdot,s) \,
\|_{\mbox{}_{\scriptstyle L^{2}(\mathbb{R}^{3})}}^{\:\!2}
ds
\;\;\!\leq\;
\|\, D \mbox{\boldmath $u$}
(\cdot,t_{2}) \,
\|_{\mbox{}_{\scriptstyle L^{2}(\mathbb{R}^{3})}}^{\:\!2}
\end{equation}
\mbox{} \vspace{-0.150cm} \\
para todo
$\;\! t \geq t_{2} \geq t_{0} $,
onde
$ {\displaystyle
\,
\gamma \:\! := \:
\nu \,-\: \|\, \mbox{\boldmath $u$}(\cdot,t_0) \,
\|_{\mbox{}_{\scriptstyle L^{2}(\mathbb{R}^{3})}}^{\:\!1/2}
\|\, D \mbox{\boldmath $u$}(\cdot,t_0) \,
\|_{\mbox{}_{\scriptstyle L^{2}(\mathbb{R}^{3})}}^{\:\!1/2}
\!
} $
\'e uma constante \linebreak
positiva.
Isto conclui a prova,
em vista de resultados
cl\'assicos de regularidade
de $ \mbox{\boldmath $u$}(\cdot,t) $ \linebreak
(ver e.g.$\;$\cite{KreissLorenz1989, %
Leray1934, LorenzZingano2012}),
tendo-se
$ \:\!t_{\!\;\!\ast\ast} \!\:\!=\;\! t_{0} \:\!$
com
$\;\!t_{0} \geq \mbox{\small $T$}_{\!\:\!\ast\ast} $
escolhido em (2.19$a$) acima.
}
\mbox{} \hfill $\Box$ \\
%
\mbox{} \vspace{-0.550cm} \\
%
%
%
%
%
\nl
{\bf Observa\c c\~ao 2.1.}
Como mostrado em
\cite{KreissHagstromLorenzZingano2002},
uma consequ\^encia da prova acima
\'e que
tem-se
$ {\displaystyle
\,
t_{\ast\ast} \!\;\!<
0.212 \cdot\:\! \nu^{-\;\!5} \;\!
\|\, \mbox{\boldmath $u$}_{0} \;\!
\|_{\scriptstyle L^{2}(\mathbb{R}^{3})}
  ^{\:\!4}
\!\;\!
} $
sempre.
Um argumento mais elaborado
desenvolvido
no \mbox{\small \sc Ap\^endice A}
do presente texto
produz a estimativa mais fina \\
\mbox{} \vspace{-0.075cm} \\
\mbox{} \hspace{+3.450cm}
\fbox{%
\begin{minipage}{7.175cm}
\nl
\mbox{} \vspace{-0.400cm} \\
$ {\displaystyle
\mbox{} \;\;\,
t_{\ast\ast}
\,<\;
0.000\,753\,026 \;\!\cdot\;\!
\nu^{\!\;\!-\,5} \:
\|\, \mbox{\boldmath $u$}_{0} \;\!
\|_{\mbox{}_{\scriptstyle L^{2}(\mathbb{R}^{3})}}
  ^{\mbox{}^{\scriptstyle \:\!4}}
\!\:\!
} $. \\
\end{minipage}
}
\mbox{} \vspace{-1.050cm} \\
\mbox{} \hfill (2.20) \\
%

%
%

{\bf 3. Preliminares, II} \\
\setcounter{section}{3}

Nesta se\c c\~ao,
vamos utilizar os resultados
revisados acima para estabelecer (1.11).
Al\'em disso,
vamos tamb\'em
obter (1.12$a$) para
$ m = 1 $ (Teorema 3.2 abaixo)
e $ m = 0 $ (Teorema 3.3),
revisando as provas dadas
em \cite{KreissHagstromLorenzZingano2002, %
KreissHagstromLorenzZingano2003}
e \cite{SchutzZinganoZingano2014},
respectivamente. \\
\mbox{} \vspace{-0.200cm} \\
%
%
%
%
%
%
\mbox{} \hspace{-0.800cm}
\fbox{%
\begin{minipage}[t]{16.000cm}
\mbox{} \vspace{-0.450cm} \\
\mbox{} \hspace{+0.300cm}
\begin{minipage}[t]{15.000cm}
\mbox{} \vspace{+0.100cm} \\
{\bf Teorema 3.1.}
\textit{%
Seja
$ \;\!\mbox{\boldmath $u$}(\cdot,t) $
solu\c c\~ao de
Leray para $\;\!(1.1)$.
Dados
$\;\! \tilde{t}_0 \!\;\!> t_0 \!\;\!\geq 0 $,
tem-se
} \linebreak
\mbox{} \vspace{-0.675cm} \\
\begin{equation}
\tag{3.1}
\|\, D^{\alpha} \mbox{\boldmath $v$}(\cdot,t) -
  D^{\alpha} \tilde{\mbox{\boldmath $v$}}(\cdot,t) \,
\|_{\mbox{}_{\scriptstyle L^{2}(\mathbb{R}^{3})}}
\leq\;\!
K\!\;\!(m)
\;\!\;\!
\nu^{\mbox{}^{\scriptstyle \!
-\,\left(\frac{5}{4} \;\!+\;\!\frac{\scriptstyle m}{2} \right)}}
\;\!
\|\, \mbox{\boldmath $u$}_0 \,
\|_{\mbox{}_{\scriptstyle L^{2}(\mathbb{R}^{3})}}^{\:\!2}
(\tilde{t}_{0} \!\;\!- t_0 )^{\mbox{}^{\scriptstyle
\frac{1}{2} }}
(t - \,\!\tilde{t}_0 )^{\mbox{}^{\scriptstyle \!
-\,\left(\frac{3}{4} \;\!+\;\!\frac{\scriptstyle m}{2} \right)}}
\end{equation}
\mbox{} \vspace{-0.350cm} \\
\textit{%
para todo $\,t > \tilde{t}_0 $,
onde
$ {\displaystyle
\,
\mbox{\boldmath $v$}(\cdot,t)
\:\!=\;\!
e^{\:\!\mbox{\scriptsize $\Delta$} (\:\!\mbox{\footnotesize $t$}
\;\!-\,\mbox{\footnotesize $t_0$})}
\:\!\mbox{\boldmath $u$}(\cdot,t_0)
} $,
$ {\displaystyle
\:\!
\tilde{\mbox{\boldmath $v$}}(\cdot,t)
\:\!=\;\!
e^{\:\!\mbox{\scriptsize $\Delta$} (\:\!\mbox{\footnotesize $t$}
\;\!-\,\mbox{\footnotesize $\tilde{t}_0$})}
\:\!\mbox{\boldmath $u$}(\cdot,\tilde{t}_0)
} $,
$ {\displaystyle
\;\!
m = |\;\!\alpha\;\!|
} $. \\
}
\end{minipage}
\end{minipage}
}
%
%
%
\nl
%
%
\mbox{} \vspace{-0.000cm} \\
{\small
{\bf Prova:}
Come\c camos escrevendo
$ \:\!\mbox{\boldmath $v$}(\cdot,t) \:\!$
na forma \\
\mbox{} \vspace{-0.550cm} \\
\begin{equation}
\notag
\mbox{\boldmath $v$}(\cdot,t)
\;=\;\:\!
e^{\;\!\nu \,\!\Delta \:\!
(\;\!\mbox{\footnotesize $t$} \;\!-\,\mbox{\footnotesize $t_0$})} \;\!
[\,\mbox{\boldmath $u$}(\cdot,t_0) -
\mbox{\boldmath $u$}_{\mbox{}_{\scriptstyle \!\:\!\delta}}
\!\;\!(\cdot,t_0) \,]
\:+\:
e^{\;\!\nu \,\!\Delta \:\!
(\;\!\mbox{\footnotesize $t$} \;\!-\,\mbox{\footnotesize $t_0$})} \;\!
\mbox{\boldmath $u$}_{\mbox{}_{\scriptstyle \!\:\!\delta}}
\!\;\!(\cdot,t_0),
\qquad
t > t_0,
\end{equation}
\mbox{} \vspace{-0.200cm} \\
sendo
$ {\displaystyle
\;\!
\mbox{\boldmath $u$}_{\mbox{}_{\scriptstyle \!\:\!\delta}}
\!\;\!(\cdot,t)
\;\!
} $
dada em (2.1),
$ \delta > 0 $.
Como \\
\mbox{} \vspace{-0.650cm} \\
\begin{equation}
\notag
\mbox{\boldmath $u$}_{\mbox{}_{\scriptstyle \!\:\!\delta}}
\!\;\!(\cdot,t_0)
\;=\;\:\!
e^{\;\!\nu \,\!\Delta \:\!\mbox{\footnotesize $t_0$}}
\;\!\bar{\mbox{\boldmath $u$}}_{0, \,\delta}
\:+\!\:\!
\int_{\:\!\mbox{\footnotesize $0$}}^{\;\!\mbox{\footnotesize $t_0$}}
\!\!
e^{\;\!\nu \,\!\Delta \:\!(\;\!\mbox{\footnotesize $t_0$} \;\!-\,
\mbox{\footnotesize $s$})} \;\!
\mbox{\boldmath $Q$}_{\mbox{}_{\scriptstyle \!\delta}}
\!\;\!(\cdot,s)
\,ds,
\end{equation}
\mbox{} \vspace{-0.125cm} \\
onde
$ {\displaystyle
\;\!
\bar{\mbox{\boldmath $u$}}_{\mbox{}_{\scriptstyle \!\;\!0, \,\delta}}
\!\;\!=\;\!
G_{\mbox{}_{\scriptstyle \!\delta}}
\!\:\!\ast
\mbox{\boldmath $u$}_{0}
} $,
$ {\displaystyle
\:\!
\mbox{\boldmath $Q$}_{\mbox{}_{\scriptstyle \!\delta}}
\!\:\!(\cdot,s)
= -\,
\bar{\mbox{\boldmath $u$}}_{\mbox{}_{\scriptstyle \!\:\!\delta}}\!\:\!(\cdot,s)
\!\;\!\cdot \!\;\!\nabla
\mbox{\boldmath $u$}_{\mbox{}_{\scriptstyle \!\:\!\delta}}\!\:\!(\cdot,s)
\;\!-\:\!
\nabla p_{\mbox{}_{\scriptstyle \!\delta}}\!\;\!(\cdot,s)
} $,
dados em (2.1$b$) e (2.12),
obt\'em-se \\
\mbox{} \vspace{-0.725cm} \\
\begin{equation}
\notag
\mbox{\boldmath $v$}(\cdot,t)
\;=\;\:\!
e^{\;\!\nu \,\!\Delta \:\!(\;\!\mbox{\footnotesize $t$} \;\!-\,
\mbox{\footnotesize $t_0$})} \;\!
[\,\mbox{\boldmath $u$}(\cdot,t_0) -
\:\!\mbox{\boldmath $u$}_{\mbox{}_{\scriptstyle \!\:\!\delta}}
\!\;\!(\cdot,t_0) \,]
\;+\;
e^{\;\!\nu \,\!\Delta \mbox{\footnotesize $t$}} \:\!
\bar{\mbox{\boldmath $u$}}_{0,\,\delta}
\;+
\int_{\:\!\mbox{\footnotesize $0$}}^{\;\!\mbox{\footnotesize $t_0$}}
\!\!
e^{\;\!\nu \,\!\Delta \:\!(\;\!\mbox{\footnotesize $t$} \;\!-\,
\mbox{\footnotesize $s$})} \;\!
\mbox{\boldmath $Q$}_{\mbox{}_{\scriptstyle \!\delta}}
\!\:\!(\cdot,s)
\; ds,
\end{equation}
\mbox{} \vspace{-0.200cm} \\
para $ \;\!t > t_0 $.
Analogamente,
tem-se,
para $\;\! t > \tilde{t}_0 $: \\
\mbox{} \vspace{-0.725cm} \\
\begin{equation}
\notag
\tilde{\mbox{\boldmath $v$}}(\cdot,t)
\;=\;\:\!
e^{\;\!\nu \,\!\Delta \:\!(\;\!\mbox{\footnotesize $t$} \;\!-\,
\mbox{\footnotesize $\tilde{t}_0$})} \;\!
[\,\mbox{\boldmath $u$}(\cdot,\tilde{t}_0) -
\:\!\mbox{\boldmath $u$}_{\mbox{}_{\scriptstyle \!\:\!\delta}}
\!\;\!(\cdot,\tilde{t}_0) \,]
\;+\;
e^{\;\!\nu \,\!\Delta \mbox{\footnotesize $t$}} \:\!
\bar{\mbox{\boldmath $u$}}_{0,\,\delta}
\;+
\int_{\:\!\mbox{\footnotesize $0$}}^{\;\!\mbox{\footnotesize $\tilde{t}_0$}}
\!\!
e^{\;\!\nu \,\!\Delta \:\!(\;\!\mbox{\footnotesize $t$} \;\!-\,
\mbox{\footnotesize $s$})} \;\!
\mbox{\boldmath $Q$}_{\mbox{}_{\scriptstyle \!\delta}}
\!\:\!(\cdot,s)
\; ds.
\end{equation}
\mbox{} \vspace{-0.150cm} \\
Segue ent\~ao,
para a diferen\c ca
$ {\displaystyle
\;\!
\mbox{\boldmath $v$}(\cdot,t)
\;\!-\;\!
\tilde{\mbox{\boldmath $v$}}(\cdot,t)
} $,
sendo $\:\! t > \tilde{t}_0 $
qualquer,
que \\
\mbox{} \vspace{+0.050cm} \\
\mbox{} \hspace{+0.025cm}
$ {\displaystyle
D^{\alpha} \:\!
\tilde{\mbox{\boldmath $v$}}(\cdot,t)
\,-\,
D^{\alpha} \:\!
\mbox{\boldmath $v$}(\cdot,t)
\;\;\!=\;\;\;\!
D^{\alpha} \:\!
\bigl\{\,
e^{\;\!\nu \,\!\Delta \:\!(\;\!\mbox{\footnotesize $t$}
\;\!-\,
\mbox{\footnotesize $\tilde{t}_0$})} \;\!
[\,\mbox{\boldmath $u$}(\cdot,\tilde{t}_0) -
\:\!\mbox{\boldmath $u$}_{\mbox{}_{\scriptstyle \!\:\!\delta}}
\!\;\!(\cdot,\tilde{t}_0) \,]
\,\bigr\}
\;\,+
} $ \\
\mbox{} \vspace{-0.700cm} \\
\mbox{} \hfill (3.2) \\
\mbox{} \vspace{-0.550cm} \\
\mbox{} \hfill
$ {\displaystyle
\;-\;
D^{\alpha} \:\!\bigl\{\,
e^{\;\!\nu \,\!\Delta \:\!(\;\!\mbox{\footnotesize $t$} \;\!-\,
\mbox{\footnotesize $t_0$})} \;\!
[\,\mbox{\boldmath $u$}(\cdot,t_0) -
\:\!\mbox{\boldmath $u$}_{\mbox{}_{\scriptstyle \!\:\!\delta}}
\!\;\!(\cdot,t_0) \,]
\,\bigr\}
\:+\,
D^{\alpha} \!\!
\int_{\:\!\mbox{\footnotesize $t_0$}}
    ^{\;\!\mbox{\footnotesize $\tilde{t}_0$}}
\!\!
e^{\;\!\nu\,\!\Delta \:\!(\;\!\mbox{\footnotesize $t$} \;\!-\,
\mbox{\footnotesize $s$})} \;\!
\mbox{\boldmath $Q$}_{\mbox{}_{\scriptstyle \!\delta}}
\!\:\!(\cdot,s)
\; ds
} $. \\
\mbox{} \vspace{+0.150cm} \\
Portanto,
dado
$ {\displaystyle
\;\!\mathbb{K} \subset \mathbb{R}^{3}
\!\:\!
} $
compacto qualquer,
tem-se,
para cada
$ \;\! t > \tilde{t}_0 $,
$ \;\! \delta > 0 $: \\
\mbox{} \vspace{-0.150cm} \\
\mbox{} \hspace{-0.050cm}
$ {\displaystyle
\|\, D^{\alpha} \tilde{\mbox{\boldmath $v$}}(\cdot,t)
\,-\, D^{\alpha} \mbox{\boldmath $v$}(\cdot,t) \,
\|_{\mbox{}_{\scriptstyle L^{2}(\mathbb{K})}}
\,\leq\;
J_{\mbox{}_{\scriptstyle \!\alpha, \;\!\delta}}(t)
\;+\:\!
\int_{\:\!\mbox{\footnotesize $t_0$}}
    ^{\;\!\mbox{\footnotesize $\tilde{t}_0$}}
\!\!
\|\, D^{\alpha} \bigl\{\;\!
e^{\;\!\nu \Delta \:\!(\;\!\mbox{\footnotesize $t$} \;\!-\,
\mbox{\footnotesize $s$})} \;\!
\mbox{\boldmath $Q$}_{\mbox{}_{\scriptstyle \!\delta}}
\!\:\!(\cdot,s)
\;\!\bigr\}
\;\!\|_{\mbox{}_{\scriptstyle L^{2}(\mathbb{K})}}
\, ds
} $ \\
\mbox{} \vspace{+0.025cm} \\
\mbox{} \hfill
$ {\displaystyle
\leq\;
J_{\mbox{}_{\scriptstyle \!\alpha, \;\!\delta}}(t)
\;+\,
K\!\;\!(m)
\:
\nu^{\mbox{}^{\scriptstyle \!-\, \frac{m}{2} }}
\!\!\!\:\!
\int_{\:\!\mbox{\footnotesize $t_0$}}
    ^{\;\!\mbox{\footnotesize $\tilde{t}_0$}}
\!
\Bigl(\, \frac{t - s}{\mbox{\scriptsize $2$}}
\;\!\Bigr)^{\mbox{}^{\scriptstyle \!\!\!\!-\,\frac{m}{2}}}
\!\;\!
\|\: e^{\;\!\nu \Delta \:\!(\;\!\mbox{\footnotesize $t$} \;\!-\,
\mbox{\footnotesize $s$})/2} \,
\mbox{\boldmath $Q$}_{\mbox{}_{\scriptstyle \!\delta}}
\!\:\!(\cdot,s)
\;\!\bigr\}
\;\!\|_{\mbox{}_{\scriptstyle L^{2}(\mathbb{R}^{3})}}
\, ds
} $ \\
\mbox{} \vspace{+0.025cm} \\
\mbox{} \hfill
$ {\displaystyle
\leq\;
J_{\mbox{}_{\scriptstyle \!\alpha, \;\!\delta}}(t)
\;+\,
K\!\;\!(m)
\:
\nu^{\mbox{}^{\scriptstyle \!-\, \frac{m}{2} \,-\, \frac{3}{4} }}
\!\;\!
(\,\!t - \tilde{t}_0
)^{\mbox{}^{\scriptstyle \!\!-\,\frac{m}{2}}}
\!\!\!
\int_{\:\!\mbox{\footnotesize $t_0$}}
    ^{\;\!\mbox{\footnotesize $\tilde{t}_0$}}
\!
\Bigl(\, \frac{t - s}{\mbox{\footnotesize $2$}}
\;\!\Bigr)^{\mbox{}^{\scriptstyle \!\!\!-\,\frac{3}{4}}}
\|\,\mbox{\boldmath $u$}_{\mbox{}_{\scriptstyle \!\:\!\delta}}
\!\:\!(\cdot,s) \,
\|_{\mbox{}_{\scriptstyle L^{2}(\mathbb{R}^{3})}}
\:\!
\|\, D \mbox{\boldmath $u$}_{\mbox{}_{\scriptstyle \!\:\!\delta}}
\!\:\!(\cdot,s) \,
\|_{\mbox{}_{\scriptstyle L^{2}(\mathbb{R}^{3})}}
\;\! ds
} $ \\
\mbox{} \vspace{+0.025cm} \\
\mbox{} \hspace{+3.050cm}
$ {\displaystyle
\leq\;
J_{\mbox{}_{\scriptstyle \!\alpha, \;\!\delta}}(t)
\;+\;
K\!\;\!(m)
\:
\nu^{\mbox{}^{\scriptstyle \!-\, \frac{m}{2} \,-\, \frac{5}{4} }}
\!\;\!
\bigl(\;\!t - \tilde{t}_0
\bigr)^{\mbox{}^{\scriptstyle \!\!-\,\frac{m}{2} \,-\, \frac{3}{4}}}
\bigl(\;\! \tilde{t}_0 \!\;\!-\;\! t_0
\bigr)^{\mbox{}^{\scriptstyle \!\:\! \frac{1}{2} }}
\:\!
\|\, \mbox{\boldmath $u$}_0 \,
\|_{\mbox{}_{\scriptstyle L^{2}(\mathbb{R}^{3})}}
  ^{\mbox{}^{\scriptstyle \:\! 2}}
} $
\hspace{+0.100cm} (3.3) \\
\mbox{} \vspace{+0.075cm} \\
por (2.11), (2.13) e (2.14) acima
(e Cauchy-Schwarz),
onde \\
\mbox{} \vspace{+0.000cm} \\
\mbox{} \hspace{+0.750cm}
$ {\displaystyle
J_{\mbox{}_{\scriptstyle \!\alpha, \;\!\delta}}(t)
\;=\;
\|\, D^{\alpha} \bigl\{\;\!
e^{\;\!\nu \Delta \:\!(\;\!\mbox{\footnotesize $t$}
\;\!-\, \mbox{\footnotesize $\tilde{t}_0$})} \;\!
[\,\mbox{\boldmath $u$}(\cdot,\tilde{t}_0) -
\:\!\mbox{\boldmath $u$}_{\mbox{}_{\scriptstyle \!\:\!\delta}}
\!\;\!(\cdot,\tilde{t}_0) \,]
\;\!\bigr\}
\;\!\|_{\mbox{}_{\scriptstyle L^{2}(\mathbb{K})}}
\;+
} $ \\
\mbox{} \vspace{-0.100cm} \\
\mbox{} \hfill
$ {\displaystyle
+\;\;\!
\|\, D^{\alpha} \bigl\{\;\!
e^{\;\!\nu \Delta \:\!(\;\!\mbox{\footnotesize $t$}
\;\!-\, \mbox{\footnotesize $t_0$})} \;\!
[\,\mbox{\boldmath $u$}(\cdot,t_0) -
\:\!\mbox{\boldmath $u$}_{\mbox{}_{\scriptstyle \!\:\!\delta}}
\!\;\!(\cdot,t_0) \,]
\;\!\bigr\}
\;\!\|_{\mbox{}_{\scriptstyle L^{2}(\mathbb{K})}}
\!\:\!
} $. \\
\mbox{} \vspace{+0.050cm} \\
Tomando
$ \;\!\delta = \delta^{\prime} \!\rightarrow  0 \;\! $
conforme (2.2),
obt\'em-se
$ J_{\mbox{}_{\scriptstyle \!\alpha, \;\!\delta}}(t) \rightarrow 0 $,
pois,
pelo teorema da converg\^encia dominada
e (2.2),
tem-se,
para cada
$ \sigma \!\:\!, \, \tau > 0 \:\!$: \\
\mbox{} \vspace{-0.550cm} \\
\begin{equation}
\tag{3.4}
\|\: D^{\alpha} \bigl\{\;\!
e^{\;\! \nu \Delta \mbox{\footnotesize $\tau$}}
[\,\mbox{\boldmath $u$}(\cdot,\sigma) -
\:\!\mbox{\boldmath $u$}_{\mbox{}_{\scriptstyle \!\:\!\delta^{\prime}}}
\!\;\!(\cdot,\sigma) \,]
\;\!\bigr\}
\;\!\|_{\mbox{}_{\scriptstyle L^{2}(\mathbb{K})}}
\rightarrow\; 0
\qquad
\mbox{ao }
\;\; \delta^{\prime} \!\rightarrow 0.
\end{equation}
\mbox{} \vspace{-0.190cm} \\
\mbox{[}$\,$De fato,
sendo
$ {\displaystyle
\;\!
F_{\mbox{}_{\!\delta}}(\cdot,\tau) :=
D^{\alpha} \bigl\{\;\!
e^{\;\! \nu \Delta \mbox{\footnotesize $\tau$}}
[\,\mbox{\boldmath $u$}(\cdot,\sigma) -
\:\!\mbox{\boldmath $u$}_{\mbox{}_{\scriptstyle \!\:\!\delta}}
\!\;\!(\cdot,\sigma) \,]
\;\!\bigr\}
} $,
tem-se \\
\mbox{} \vspace{-0.750cm} \\
\begin{equation}
\notag
F_{\mbox{}_{\!\;\!\delta}}(\cdot,\tau)
\,=\,
H_{\mbox{}_{\!\;\!\alpha}}(\cdot,\tau) \ast
[\,\mbox{\boldmath $u$}(\cdot,\sigma) -
\:\!\mbox{\boldmath $u$}_{\mbox{}_{\scriptstyle \!\:\!\delta}}
\!\;\!(\cdot,\sigma) \,]
\end{equation}
\mbox{} \vspace{-0.300cm} \\
onde
$ {\displaystyle
H_{\mbox{}_{\!\alpha}}(\cdot,\tau)
\in
L^{1}(\mathbb{R}^{3}) \cap L^{\infty}(\mathbb{R}^{3})
\;\!
} $
\'e independente de $ \delta $.
Como
$ {\displaystyle
\;\!
\mbox{\boldmath $u$}(\cdot,\sigma) -
\mbox{\boldmath $u$}_{\mbox{}_{\scriptstyle \!\:\!\delta^{\prime}}}
\!\;\!(\cdot,\sigma)
\;\!\rightharpoonup \;\!0
\;\!
} $
em $ L^{2}(\mathbb{R}^{3}) $,
ver (2.2),
segue que
$ {\displaystyle
F_{\mbox{}_{\!\:\!\delta^{\prime}}}(x,\tau)
\rightarrow 0
} $
(ao $ \delta^{\prime} \rightarrow 0 $)
para cada $ x \in \mathbb{R}^{3} \!\;\!$;
por outro lado,
por (2.13) e Cauchy-Schwarz,
tem-se \\
\mbox{} \vspace{-0.550cm} \\
\begin{equation}
\notag
\begin{split}
|\, F_{\mbox{}_{\!\;\!\delta}}(x,\tau) \,|
\;&\leq\;
\|\, H_{\mbox{}_{\!\;\!\alpha}}(\cdot,\tau) \,
\|_{\mbox{}_{\scriptstyle L^{2}(\mathbb{R}^{3})}}
\;\!
\|\,\mbox{\boldmath $u$}(\cdot,\sigma) -
\:\!\mbox{\boldmath $u$}_{\mbox{}_{\scriptstyle \!\:\!\delta}}
\!\;\!(\cdot,\sigma) \,
\|_{\mbox{}_{\scriptstyle L^{2}(\mathbb{R}^{3})}} \\
&\leq\;
2 \:
\|\, H_{\mbox{}_{\!\;\!\alpha}}(\cdot,\tau) \,
\|_{\mbox{}_{\scriptstyle L^{2}(\mathbb{R}^{3})}}
\;\!
\|\, \mbox{\boldmath $u$}_0 \;\!
\|_{\mbox{}_{\scriptstyle L^{2}(\mathbb{R}^{3})}}
\end{split}
\end{equation}
\mbox{} \vspace{-0.050cm} \\
para todo $ x \in \mathbb{R}^{3} \!\;\!$.
Por
converg\^encia dominada,
segue ent\~ao
que
$ {\displaystyle
\;\!
\|\, F_{\mbox{}_{\!\:\!\delta^{\prime}}}(\cdot,\tau) \,
\|_{{\scriptstyle L^{2}(\mathbb{K})}}
\!\rightarrow\;\! 0
\;\!
} $
ao
$ \delta^{\:\!\prime} \!\;\!\rightarrow 0 $,
visto que $ \mathbb{K} $ \'e compacto,
o que conclui a
demonstra\c c\~ao da afirma\c c\~ao (3.4) acima.$\;\!$\mbox{]}
Assim,
fazendo
$ \;\! \delta = \delta^{\prime} \!\rightarrow 0 \;\!$
em (3.3),
resulta \\
\mbox{} \vspace{-0.650cm} \\
\begin{equation}
\notag
\|\, D^{\alpha} \:\! \tilde{\mbox{\boldmath $v$}}(\cdot,t)
\;\!-\;\! D^{\alpha} \:\!\mbox{\boldmath $v$}(\cdot,t) \,
\|_{\mbox{}_{\scriptstyle L^{2}(\mathbb{K})}}
\leq\:
K\!\;\!(m)
\:
\nu^{\mbox{}^{\scriptstyle \!-\, \frac{m}{2} \,-\, \frac{5}{4} }}
\;\!
\|\, \mbox{\boldmath $u$}_0 \,
\|_{\mbox{}_{\scriptstyle L^{2}(\mathbb{R}^{3})}}
  ^{\mbox{}^{\scriptstyle \:\! 2}}
\bigl(\;\! \tilde{t}_0 \!\;\!-\;\! t_0
\bigr)^{\mbox{}^{\scriptstyle \!\:\! \frac{1}{2} }}
\bigl(\;\!t - \tilde{t}_0
\bigr)^{\mbox{}^{\scriptstyle \!\!-\,\frac{m}{2} \,-\, \frac{3}{4}}}
\end{equation}
\mbox{} \vspace{-0.250cm} \\
para todo
$ \;\! t > \tilde{t}_0 $
(sendo
$ \;\!\mathbb{K} \subset \mathbb{R}^{3} \!\;\!$
compacto {\em arbitr\'ario\/}),
o que \'e equivalente a (3.1).
}
\mbox{} \hfill $\Box$ \\
%
%
\mbox{} \vspace{-0.675cm} \\

O Teorema 3.1 estabelece (1.11)
para todo $ s = m \geq 0 $ inteiro.
Utilizando (1.13), segue
a validade de (1.11)
para todo $ \;\!s \geq 0 $
($s$ real), como afirmado.

Para a obten\c c\~ao
das estimativas (1.12),
que ser\'a feito nas
Se\c c\~oes 4 e 5
a seguir,
precisaremos ter
anteriormente
estabelecido (1.12$a$)
nos casos particulares
$ m = 0, 1 $.
Estes dois resultados
j\'a foram obtidos
por outros autores;
as provas mais simples
s\~ao fornecidas em
\cite{SchutzZinganoZingano2014}
($ m = 0 $)
e
\cite{KreissHagstromLorenzZingano2002, %
KreissHagstromLorenzZingano2003}
($ m = 1 $),
repetidas abaixo por conveni\^encia.
O resultado mais f\'acil
\'e o segundo,
por seguir imediatamente
da desigualdade de energia
(1.2)
e do fato de
$ {\displaystyle
\;\!
\|\, D \mbox{\boldmath $u$}(\cdot,t) \,
\|_{\mbox{}_{\scriptstyle L^{2}(\mathbb{R}^{3})}}
\!\;\!
} $
ser eventualmente monot\^onica,
conforme a Proposi\c c\~ao 2.3
da Se\c c\~ao 2 acima: \\
\nl
%
%
%
%
%
{\bf Teorema 3.2.}
\textit{%
Sendo
$ \,\mbox{\boldmath $u$}(\cdot,t) $
solu\c c\~ao de Leray
para as equa\c c\~oes $\;\!(1.1)$,
tem-se
} \\
\mbox{} \vspace{-0.650cm} \\
\begin{equation}
\tag{3.5}
\lim_{\;t\,\rightarrow\,\infty}
\,
t^{\:\!1/2}
\,
\|\, D \mbox{\boldmath $u$}(\cdot,t) \,
\|_{\mbox{}_{\scriptstyle L^{2}(\mathbb{R}^{3})}}
=\; 0.
\end{equation}
%
%
%
\nl
\mbox{} \vspace{-0.475cm} \\
%
%
{\small
{\bf Prova:}
O argumento seguinte
\'e obtido de
\cite{KreissHagstromLorenzZingano2002, %
KreissHagstromLorenzZingano2003}.
Se a propriedade (3.5) fosse falsa,
existiria uma sequ\^encia
crescente
$ \;\! t_{\ell} \;\!\mbox{\scriptsize $\nearrow$}\;\! \infty $
(com
$ \:\! t_{\ell} \geq t_{\!\;\!\ast\ast} $
e
$ \;\! t_{\ell} \geq 2 \;\! t_{\ell - 1} $
para todo $\ell $,
digamos) \linebreak
e algum
$ \;\!\eta > 0 \;\!$
fixo
tais que \\
\mbox{} \vspace{-0.650cm} \\
\begin{equation}
\notag
\mbox{} \;\;\;
t_{\ell} \:
\|\, D \mbox{\boldmath $u$}(\cdot,t_{\ell}) \,
\|_{\mbox{}_{\scriptstyle L^{2}(\mathbb{R}^{3})}}^{\:\!2}
\:\!\geq\; \eta
\qquad
\forall \;\;\! \ell.
\end{equation}
\mbox{} \vspace{-0.240cm} \\
Em particular,
ter\'\i amos
de ter \\
\mbox{} \vspace{-0.600cm} \\
\begin{equation}
\notag
\int_{\mbox{}_{\mbox{\footnotesize $\!\:\!t_{\ell - 1}$}}}
    ^{\mbox{\footnotesize $\:\!t_{\ell}$}}
\!\!\!\!\!
\|\, D \mbox{\boldmath $u$}(\cdot,s) \,
\|_{\mbox{}_{\scriptstyle L^{2}(\mathbb{R}^{3})}}^{\:\!2}
\:\!dt
\:\,\!\geq\:
(\:\!t_{\ell} -\;\! t_{\ell - 1})
\,
\|\, D \mbox{\boldmath $u$}(\cdot,t_{\ell}) \,
\|_{\mbox{}_{\scriptstyle L^{2}(\mathbb{R}^{3})}}^{\:\!2}
\!\:\!\geq\,
\mbox{\footnotesize $ {\displaystyle \frac{1}{2} }$}
\: t_{\ell} \,
\|\, D \mbox{\boldmath $u$}(\cdot,t_{\ell}) \,
\|_{\mbox{}_{\scriptstyle L^{2}(\mathbb{R}^{3})}}^{\:\!2}
\!\:\!\geq\,
\mbox{\footnotesize $ {\displaystyle \frac{1}{2} }$}
\, \eta
\end{equation}
\mbox{} \vspace{-0.100cm} \\
para todo $\ell $,
em contradi\c c\~ao com (1.2), (2.3).
Portanto,
(3.5) tem de ser verdadeira.
}
\mbox{} \hfill $\Box$ \\
%
%
\nl
\mbox{} \vspace{-0.400cm} \\
%
%
%
%
%
{\bf Teorema 3.3}
({\small \sc Solu\c c\~ao do Problema Cl\'assico de Leray\/}).
\textit{%
Tem-se
} \\
\mbox{} \vspace{-0.675cm} \\
\begin{equation}
\tag{3.6}
\lim_{t\,\rightarrow\,\infty}
\,
\|\, \mbox{\boldmath $u$}(\cdot,t) \,
\|_{\mbox{}_{\scriptstyle L^{2}(\mathbb{R}^{3})}}
=\;0.
\end{equation}
%
\nl
\mbox{} \vspace{-0.500cm} \\
{\small
{\bf Prova:}
Seguindo o argumento em
\cite{SchutzZinganoZingano2014},
seja
$ t_{\ast\ast} \!\;\!$
definido na Proposi\c c\~ao 2.3
da Se\c c\~ao 2.
Dado $ \:\!\epsilon > 0 $,
tomemos
$ \;\!t_0 \!\;\!\geq t_{\ast\ast} $
suficientemente grande
tal que,
pelo Teorema 3.2,
tenha-se \\
\mbox{} \vspace{-0.600cm} \\
\begin{equation}
\tag{3.7}
\mbox{} \hspace{+0.500cm}
t^{\:\!1/2} \,
\|\, D \mbox{\boldmath $u$}(\cdot,t) \,
\|_{\mbox{}_{\scriptstyle L^{2}(\mathbb{R}^{3})}}
\leq\: \epsilon
\qquad
\forall \;\;\!
t \geq t_0.
\end{equation}
\mbox{} \vspace{-0.175cm} \\
Como
$ \:\!\mbox{\boldmath $u$}(\cdot,t) \:\!$
\'e suave em
$ [\,t_0, \infty) $,
podemos escrever
(pelo princ\'\i pio de Duhamel) \\
\mbox{} \vspace{-0.525cm} \\
\begin{equation}
\tag{3.8}
\mbox{\boldmath $u$}(\cdot,t)
\;=\;
e^{\;\!\nu \Delta (\mbox{\footnotesize $t$} \;\!-\, \mbox{\footnotesize $t_0$})}
\:\!
\mbox{\boldmath $u$}(\cdot,t_0)
\;+
\int_{\mbox{\footnotesize $\!\;\!t_0$}}
    ^{\mbox{\footnotesize $\:\!t$}}
\!\!\;\!
e^{\;\!\nu \Delta (\mbox{\footnotesize $t$} \;\!-\, \mbox{\footnotesize $s$})}
\:\!
\mbox{\boldmath $Q$}(\cdot,s)
\: ds,
\qquad
t \geq t_0,
\end{equation}
\mbox{} \vspace{-0.150cm} \\
onde
$ {\displaystyle
\;\!
\mbox{\boldmath $Q$}(\cdot,s)
\;\!
} $
\'e dada em (2.4), Se\c c\~ao 2..
Usando a representa\c c\~ao (3.8)
para $ \mbox{\boldmath $u$}(\cdot,t) $,
obt\'em-se \\
\mbox{} \vspace{-0.075cm} \\
\mbox{} \hspace{-0.200cm}
$ {\displaystyle
\|\, \mbox{\boldmath $u$}(\cdot,t) \,
\|_{\mbox{}_{\scriptstyle L^{2}(\mathbb{R}^{3})}}
\;\!\leq\;
\|\: e^{\;\!\nu \Delta (t \;\!-\, t_0)} \:\!
\mbox{\boldmath $u$}(\cdot,t_0) \,
\|_{\mbox{}_{\scriptstyle L^{2}(\mathbb{R}^{3})}}
\:\!+
\int_{\mbox{\scriptsize $t_0$}}
    ^{\mbox{\scriptsize $t$}}
\!
\|\: e^{\;\!\nu \Delta (t \;\!-\, s)} \:\!
\mbox{\boldmath $Q$}(\cdot,s) \,
\|_{\mbox{}_{\scriptstyle L^{2}(\mathbb{R}^{3})}}
\:\!
ds
} $ \\
\mbox{} \vspace{+0.050cm} \\
\mbox{} \hfill
$ {\displaystyle
\leq\:
\|\: e^{\;\!\nu \Delta (t \;\!-\, t_0)} \:\!
\mbox{\boldmath $u$}(\cdot,t_0) \,
\|_{\mbox{}_{\scriptstyle L^{2}(\mathbb{R}^{3})}}
\:\!+\:
K \;\!
\nu^{-\,3/4} \!\!
\int_{\mbox{\scriptsize $t_0$}}
    ^{\mbox{\scriptsize $t$}}
\!
(t - s)^{-\,3/4} \:
\|\, \mbox{\boldmath $u$}(\cdot,s) \,
\|_{\mbox{}_{\scriptstyle L^{2}(\mathbb{R}^{3})}}
\:\!
\|\, D \mbox{\boldmath $u$}(\cdot,s) \,
\|_{\mbox{}_{\scriptstyle L^{2}(\mathbb{R}^{3})}}
\:\!
ds
} $ \\
\mbox{} \vspace{+0.050cm} \\
\mbox{} \hfill
$ {\displaystyle
\leq\:
\|\: e^{\;\!\nu \Delta (t \;\!-\, t_0)} \:\!
\mbox{\boldmath $u$}(\cdot,t_0) \,
\|_{\mbox{}_{\scriptstyle L^{2}(\mathbb{R}^{3})}}
\:\!+\:
K \;\!
\nu^{-\,3/4} \,
\|\, \mbox{\boldmath $u$}_0 \,
\|_{\mbox{}_{\scriptstyle L^{2}(\mathbb{R}^{3})}}
\!\!\;\!
\int_{\mbox{\scriptsize $t_0$}}
    ^{\mbox{\scriptsize $t$}}
\!\;\!
(t - s)^{-\,3/4} \:
\|\, D \mbox{\boldmath $u$}(\cdot,s) \,
\|_{\mbox{}_{\scriptstyle L^{2}(\mathbb{R}^{3})}}
\:\!
ds
} $ \\
\mbox{} \vspace{+0.050cm} \\
\mbox{} \hspace{+0.100cm}
$ {\displaystyle
\leq\:
\|\: e^{\;\!\nu \Delta (\mbox{\scriptsize $t$} \;\!-\, t_0)} \:\!
\mbox{\boldmath $u$}(\cdot,t_0) \,
\|_{\mbox{}_{\scriptstyle L^{2}(\mathbb{R}^{3})}}
+\,
K \;\! \nu^{-\,3/4} \,
\|\, \mbox{\boldmath $u$}_0 \,
\|_{\mbox{}_{\scriptstyle L^{2}(\mathbb{R}^{3})}}
\;\!
\epsilon
\!\:\!
\int_{\mbox{\scriptsize $t_0$}}
    ^{\mbox{\scriptsize $t$}}
\!\;\!
(t - s)^{-\,3/4} \, s^{-\,1/2}
\,ds
} $
\mbox{} \hspace{+0.200cm} [$\,$por (3.7)$\;\!$] \\
\mbox{} \vspace{+0.150cm} \\
para todo $ \:\!t > t_0 $,
onde
$ \:\!K \!\:\!=\:\! (\:\!8 \:\!\pi )^{-\,3/4} \!\:\!$,
usando
(1.2) e (2.5).
Observando que  \\
\mbox{} \vspace{-0.575cm} \\
\begin{equation}
\notag
\int_{\mbox{\scriptsize $t_0$}}
    ^{\mbox{\scriptsize $t$}}
\!\;\!
(t - s)^{-\,3/4} \, s^{-\,1/2}
\,ds
\;\leq\;
6 \: \sqrt[4]{\;\!2\;}
\qquad
\forall
\;\,
t \;\!\geq\;\! t_0 + 1,
\end{equation}
\mbox{} \vspace{-0.075cm} \\
obt\'em-se, ent\~ao,
para todo $ \:\! t \geq\;\! t_0 +\;\! 1 $: \\
\mbox{} \vspace{-0.550cm} \\
\begin{equation}
\notag
\|\, \mbox{\boldmath $u$}(\cdot,t) \,
\|_{\mbox{}_{\scriptstyle L^{2}(\mathbb{R}^{3})}}
\;\!\leq\;
\|\: e^{\;\! \nu \Delta (t \;\!-\, t_0)} \:\!
\mbox{\boldmath $u$}(\cdot,t_0) \,
\|_{\mbox{}_{\scriptstyle L^{2}(\mathbb{R}^{3})}}
\:\!+\:
\nu^{-\,3/4} \,
\|\, \mbox{\boldmath $u$}_0 \,
\|_{\mbox{}_{\scriptstyle L^{2}(\mathbb{R}^{3})}}
\epsilon.
\end{equation}
\mbox{} \vspace{-0.150cm} \\
Como,
para o semigrupo do calor,
tem-se
que
$ {\displaystyle
\;\!
\|\: e^{\;\!\nu \Delta (t \;\!-\, t_0)} \:\!
\mbox{\boldmath $u$}(\cdot,t_0) \,
\|_{\mbox{}_{\scriptstyle L^{2}(\mathbb{R}^{3})}}
\!\rightarrow \;\!0
\;\!
} $
ao $ \;\!t \rightarrow \infty $, \linebreak
\mbox{} \vspace{-0.550cm} \\
segue que \\
\mbox{} \vspace{-0.900cm} \\
\begin{equation}
\notag
\|\, \mbox{\boldmath $u$}(\cdot,t) \,
\|_{\mbox{}_{\scriptstyle L^{2}(\mathbb{R}^{3})}}
\;\!\leq\;
(\;\! 1 \,+\, \nu^{-\,3/4} \,
\|\, \mbox{\boldmath $u$}_0 \,
\|_{\mbox{}_{\scriptstyle L^{2}(\mathbb{R}^{3})}})
\: \epsilon
\end{equation}
\mbox{} \vspace{-0.250cm} \\
para todo $ \;\! t \gg 1 $.
Dado que
$ \;\!\epsilon > 0 \;\!$
\'e arbitr\'ario,
isto mostra (3.6),
como afirmado.
}
\mbox{} \hfill $\Box$ \\
%
\nl
\mbox{} \vspace{-0.525cm} \\
%
%
%
%
%
{\bf Observa\c c\~ao 3.1.}
De modo similar,
poder\'\i amos tamb\'em
estabelecer (1.8),
ou seja, \\
\mbox{} \vspace{-0.550cm} \\
\begin{equation}
\tag{3.9}
\lim_{t\,\rightarrow\,\infty}
\,
t^{\;\! \frac{\scriptstyle 3}{\scriptstyle 4}}
\;\!
\|\, \mbox{\boldmath $u$}(\cdot,t) -\,
e^{\:\!\nu \:\! \Delta \:\!(t - t_0)}
\;\!\mbox{\boldmath $u$}(\cdot,t_0) \,
\|_{\mbox{}_{\scriptstyle L^{2}(\mathbb{R}^{3})}}
\:\!=\;
0,
\end{equation}
\mbox{} \vspace{-0.170cm} \\
mas este resultado
n\~ao ser\'a necess\'ario
na an\'alise a seguir.
(Para uma prova de (3.9), \linebreak
ver \cite{SchutzZinganoZingano2014}, Section 3.)
O mesmo vale para
as propriedades (1.9),
i.e., \\
\mbox{} \vspace{-0.650cm} \\
\begin{equation}
\tag{3.10$a$}
\lim_{t\,\rightarrow\,\infty} \,
t^{\;\! \frac{\scriptstyle 3}{\scriptstyle 4}}
\;\!
\|\, \mbox{\boldmath $u$}(\cdot,t) \,
\|_{\mbox{}_{\scriptstyle L^{\infty}(\mathbb{R}^{3})}}
\:\!=\;
0,
\end{equation}
\mbox{} \vspace{-0.750cm} \\
\begin{equation}
\tag{3.10$b$}
\lim_{t\,\rightarrow\,\infty}
\,
t
\:
\|\, \mbox{\boldmath $u$}(\cdot,t) -\,
e^{\:\!\nu \:\! \Delta \:\!(t - t_0)}
\;\!\mbox{\boldmath $u$}(\cdot,t_0) \,
\|_{\mbox{}_{\scriptstyle L^{\infty}(\mathbb{R}^{3})}}
\:\!=\;
0,
\end{equation}
\mbox{} \vspace{-0.100cm} \\
obtidas em
\cite{SchutzZinganoZingano2014}, Section 4.
Um ponto a destacar sobre as provas
de (3.6), (3.9), (3.10) \linebreak
apresentadas em
\cite{SchutzZinganoZingano2014}
\'e que
\textit{%
elas utilizam
apenas resultados
obtidos por Leray}
\cite{Leray1934}
{\em e} \linebreak
\textit{%
m\'etodos matem\'aticos conhecidos
at\'e aquela \'epoca}.
O mesmo {\em n\~ao\/} vale para
as estima\-tivas (bem mais dif\'\i ceis)
em (1.10) ou (1.12),
como se ver\'a
nas Se\c c\~oes 4 e 5 seguir.
%
%

%
%
\mbox{} \vspace{-1.250cm} \\

{\bf 4. Prova de (1.12{\em a\/})} \\
\setcounter{section}{4}

Nesta se\c c\~ao,
vamos demonstrar o seguinte resultado,
dada uma solu\c c\~ao de Leray \linebreak
(qualquer)
$ \mbox{\boldmath $u$}(\cdot,t) $
para as equa\c c\~oes de Navier-Stokes (1.1). \\
\mbox{} \vspace{-0.200cm} \\
%
%
%
%
%
%
\mbox{} \hspace{-0.800cm}
\fbox{%
\begin{minipage}[t]{16.000cm}
\mbox{} \vspace{-0.450cm} \\
\mbox{} \hspace{+0.300cm}
\begin{minipage}[t]{15.000cm}
\mbox{} \vspace{+0.100cm} \\
{\bf Teorema 4.1.}
\textit{%
Para todo
$\;\!m \geq 0 $,
tem-se
} \\
\mbox{} \vspace{-0.625cm} \\
\begin{equation}
\tag{4.1}
\mbox{} \;\;\,
\lim_{t\,\rightarrow\,\infty}
\;
t^{\mbox{}^{\scriptstyle
\frac{\scriptstyle m}{2} }}
\;\!
\|\, D^{m} \,\! \mbox{\boldmath $u$}(\cdot,t)\,
\|_{\mbox{}_{\scriptstyle L^{2}(\mathbb{R}^{3})}}
\;\!=\; 0.
\end{equation}
\mbox{} \vspace{-0.500cm} \\
\end{minipage}
\end{minipage}
}
%
%
%
%
\nl
\mbox{} \vspace{-0.100cm} \\
A prova do Teorema 4.1
ocupar\'a
toda a discuss\~ao
a seguir.
Novamente,
o ponto de partida
\'e dado pela representa\c c\~ao (3.8)
para $ \mbox{\boldmath $u$}(\cdot,t) $,
ou seja:
considerando
$ t_0 \geq \mbox{\small $T$}_{\!\:\!\ast\ast} $
($ \:\!\mbox{\small $T$}_{\!\:\!\ast\ast} $
dado em (1.3), Se\c c\~ao 1),
podemos escrever,
como
$ \mbox{\boldmath $u$}(\cdot,t) $
\'e suave em $ [\,t_0, \infty ) $: \\
\mbox{} \vspace{-0.525cm} \\
\begin{equation}
\tag{4.2}
\mbox{\boldmath $u$}(\cdot,t)
\;=\;
e^{\;\!\nu \Delta (\mbox{\footnotesize $t$} \;\!-\, \mbox{\footnotesize $t_0$})}
\:\!
\mbox{\boldmath $u$}(\cdot,t_0)
\;+
\int_{\mbox{\footnotesize $\!\;\!t_0$}}
    ^{\mbox{\footnotesize $\:\!t$}}
\!\!\;\!
e^{\;\!\nu \Delta (\mbox{\footnotesize $t$} \;\!-\, \mbox{\footnotesize $s$})}
\:\!
\mbox{\boldmath $Q$}(\cdot,s)
\: ds,
\qquad
\forall \;\,
t \geq t_0,
\end{equation}
\mbox{} \vspace{-0.050cm} \\
onde
$ {\displaystyle
\;\!
\mbox{\boldmath $Q$}(\cdot,s)
\;\!
} $
\'e dada em (2.4), Se\c c\~ao 2.
Em particular,
segue que,
para cada $\alpha$, \\
\mbox{} \vspace{-0.600cm} \\
\begin{equation}
\tag{4.3$a$}
\|\, D^{\alpha} \mbox{\boldmath $u$}(\cdot,t) \,
\|_{\mbox{}_{\scriptstyle L^{2}(\mathbb{R}^{3})}}
\,\leq\:
I_{\alpha}(\:\!t\:\!; t_0)
\;+
\int_{\mbox{\footnotesize $t_0$}}
    ^{\mbox{\footnotesize $\;\!t$}}
\!
\|\, D^{\alpha} \:\![\:
e^{\;\!\nu \Delta (\mbox{\footnotesize $t$} \;\!-\, \mbox{\footnotesize $s$})}
\:\!
\mbox{\boldmath $Q$}(\cdot,s)
\,] \:
\|_{\mbox{}_{\scriptstyle L^{2}(\mathbb{R}^{3})}}
\;\! ds
\end{equation}
\mbox{} \vspace{-0.150cm} \\
para todo $ \;\!t \geq t_0 $,
onde \\
\mbox{} \vspace{-0.625cm} \\
\begin{equation}
\tag{4.3$b$}
\mbox{} \;\;
I_{\alpha}(\:\!t\:\!; t_0)
\;=\;
\|\: D^{\alpha} \;\![\;\!\;\!
e^{\;\!\nu \Delta (\mbox{\footnotesize $t$} \;\!-\, \mbox{\footnotesize $t_0$})}
\:\!
\mbox{\boldmath $u$}(\cdot,t_0)
\,] \,
\|_{\mbox{}_{\scriptstyle L^{2}(\mathbb{R}^{3})}}
\!\:\!,
\quad \;\;
t \;\!\geq\;\! t_{0}.
\end{equation}
\mbox{} \vspace{-0.150cm} \\
Note-se que
(4.1) j\'a foi obtido
se $ m = 0 $ (Teorema 3.3)
e $ m = 1 $ (Teorema 3.2),
seguindo
\cite{KreissHagstromLorenzZingano2002, %
KreissHagstromLorenzZingano2003, SchutzZinganoZingano2014}.
Em particular,
tem-se \\
\mbox{} \vspace{-0.550cm} \\
\begin{equation}
\tag{4.4}
t^{\mbox{}^{\scriptstyle \frac{1}{2} }}
\:\!
\|\, D \:\! \mbox{\boldmath $u$}(\cdot,t) \,
\|_{\mbox{}_{\scriptstyle L^{2}(\mathbb{R}^{3})}}
\;\!\leq\:
M_{\mbox{}_{1}}
\!\:\!,
\qquad
\forall \;\,
t \,\geq\, \mbox{\small $T$}_{\!\:\!\ast\ast}
\end{equation}
\mbox{} \vspace{-0.225cm} \\
para certa constante $ M_{\mbox{}_{1}} \!\;\!> 0 $
dependendo da solu\c c\~ao
$ \mbox{\boldmath $u$}(\cdot,t) $
considerada.
Mais geralmente,
tem-se
o seguinte resultado. \\
\nl
%
%
%
%
%
{\bf Lemma 4.1.}
\textit{%
Sendo
$ {\displaystyle
\;\!
U_{m}(t) :=\,
t^{\mbox{}^{\scriptstyle \frac{\scriptstyle m}{2} }}
\!\;\!
\|\, D^{m} \,\!\mbox{\boldmath $u$}(\cdot,t) \,
\|_{\mbox{}_{\scriptstyle L^{2}(\mathbb{R}^{3})}}
\!\;\!
} $,
tem-se,
para cada $\;\!m \geq 0 \!\!\;\!:$
} \\
\mbox{} \vspace{-0.625cm} \\
\begin{equation}
\tag{4.5}
U_{m} \in
L^{\infty}(\:\![\,\mbox{\small $T$}_{\!\:\!\ast\ast}, \infty)\,\!).
\end{equation}.
%
%
%
\nl
{\small
{\bf Prova:}
Para $ m = 0, 1 $,
(4.5) segue de (1.2), (4.4);
assim,
resta provar (4.5)
para
$ m \geq 2 $.
Dados
$ \;\! t_0 \!\;\!\geq \mbox{\small $T$}_{\!\:\!\ast\ast} \:\!$
e $ \;\!\alpha \;\!$ multi-\'\i ndice com
$ \;\!|\;\!\alpha\;\!| = m $,
tem-se,
de (4.3) acima,
para
$ \;\! t \geq t_0 $: \\
\mbox{} \vspace{-0.500cm} \\
\begin{equation}
\tag{4.6}
V_{\alpha}(t) \;\equiv\;
t^{\mbox{}^{\scriptstyle \frac{\scriptstyle m}{2} }}
\!\;\!
\|\, D^{\alpha} \,\!\mbox{\boldmath $u$}(\cdot,t) \,
\|_{\mbox{}_{\scriptstyle L^{2}(\mathbb{R}^{3})}}
\;\!\leq\;
I_{1}(\alpha, \:\!t) \,+\,
I_{2}(\alpha, \:\!t) \,+\,
J_{\alpha}(t),
\end{equation}
\mbox{} \vspace{-0.350cm} \\
onde \\
\mbox{} \vspace{-0.950cm} \\
\begin{equation}
\tag{4.7$a$}
I_{1}(\alpha, \:\!t)
\;=\;\;\!
t^{\mbox{}^{\scriptstyle \frac{\scriptstyle m}{2} }}
\,\!
\|\: D^{\alpha} \;\![\;\!\;\!
e^{\;\!\nu \Delta (\mbox{\footnotesize $t$} \;\!-\, \mbox{\footnotesize $t_0$})}
\:\!
\mbox{\boldmath $u$}(\cdot,t_0)
\,] \,
\|_{\mbox{}_{\scriptstyle L^{2}(\mathbb{R}^{3})}}
\!\:\!,
\end{equation}
\mbox{} \vspace{-0.850cm} \\
\begin{equation}
\tag{4.7$b$}
I_{2}(\alpha, \:\!t)
\;=\;\;\!
t^{\mbox{}^{\scriptstyle \frac{\scriptstyle m}{2} }}
\!\!\!
\int_{\mbox{\footnotesize $\;\!t_0$}}
    ^{\mbox{\footnotesize $\;\!\mu(t)$}}
\!\!
\|\, D^{\alpha} \:\![\:
e^{\;\!\nu \Delta (\mbox{\footnotesize $t$} \;\!-\, \mbox{\footnotesize $s$})}
\:\!
\mbox{\boldmath $Q$}(\cdot,s)
\,] \:
\|_{\mbox{}_{\scriptstyle L^{2}(\mathbb{R}^{3})}}
\;\! ds,
\end{equation}
\mbox{} \vspace{-0.750cm} \\
\begin{equation}
\tag{4.7$c$}
J_{\alpha}(t)
\;=\;\;\!
t^{\mbox{}^{\scriptstyle \frac{\scriptstyle m}{2} }}
\!\!\!
\int_{\mbox{\footnotesize $\mu(t)$}}
    ^{\mbox{\footnotesize $\;\!t$}}
\!
\|\, D^{\alpha} \:\![\:
e^{\;\!\nu \Delta (\mbox{\footnotesize $t$} \;\!-\, \mbox{\footnotesize $s$})}
\:\!
\mbox{\boldmath $Q$}(\cdot,s)
\,] \:
\|_{\mbox{}_{\scriptstyle L^{2}(\mathbb{R}^{3})}}
\;\! ds,
\end{equation}
\mbox{} \vspace{+0.100cm} \\
sendo
$ {\displaystyle
\;\!
\mu(t) = (\:\!t_0 + \;\!t)/2
} $.
Estimar $ I_{1}(t) $, $ I_{2}(t) $ \'e simples:
obt\'em-se, de (2.14) e (4.7$a$), \\
\mbox{} \vspace{-0.500cm} \\
\begin{equation}
\tag{4.8}
|\: I_{1}(\alpha, \:\!t) \,|
\;\leq\:
K\!\;\!(m,\;\!\nu) \:
\|\, \mbox{\boldmath $u$}(\cdot,t_0) \,
\|_{\mbox{}_{\scriptstyle L^{2}(\mathbb{R}^{3})}}
\,\!
(\:\! t - t_0)^{\mbox{}^{\scriptstyle
\!\!-\, \frac{\scriptstyle m}{2} }}
\;\!
t^{\mbox{}^{\scriptstyle \frac{\scriptstyle m}{2} }}
\leq\;
K\!\;\!(m,\;\!\nu, \;\!t_0) \:
\|\, \mbox{\boldmath $u$}_0 \;\!
\|_{\mbox{}_{\scriptstyle L^{2}(\mathbb{R}^{3})}}
\end{equation}
\mbox{} \vspace{-0.150cm} \\
para todo
$ \;\! t \geq t_0 + 1 $,
enquanto,
por (2.15), (4.4) e (4.7$b$), \\
\mbox{} \vspace{-0.650cm} \\
\begin{equation}
\notag
\begin{split}
|\: I_{2}(\alpha, \:\!t) \,|
\;&\leq\:
K\!\;\!(m,\;\!\nu) \:
t^{\mbox{}^{\scriptstyle \frac{\scriptstyle m}{2} }}
\!\!\!
\int_{\mbox{\footnotesize $\;\!t_0$}}
    ^{\mbox{\footnotesize $\;\!\mu(t)$}}
\!\!
(\:\! t - s)^{\mbox{}^{\scriptstyle \!\!
-\, \frac{\scriptstyle m}{2} \,-\, \frac{3}{4} }}
\;\!
\|\, \mbox{\boldmath $u$}(\cdot,s) \,
\|_{\mbox{}_{\scriptstyle L^{2}(\mathbb{R}^{3})}}
\;\!
\|\, D \:\! \mbox{\boldmath $u$}(\cdot,s) \,
\|_{\mbox{}_{\scriptstyle L^{2}(\mathbb{R}^{3})}}
\;\! ds \\
&\leq\:
K\!\;\!(m,\;\!\nu) \:
\|\, \mbox{\boldmath $u$}_0 \;\!
\|_{\mbox{}_{\scriptstyle L^{2}(\mathbb{R}^{3})}}
\;\!
t^{\mbox{}^{\scriptstyle \frac{\scriptstyle m}{2} }}
\!\!\!
\int_{\mbox{\footnotesize $\;\!t_0$}}
    ^{\mbox{\footnotesize $\;\!\mu(t)$}}
\!\!
s^{\mbox{}^{\scriptstyle \!- \, \frac{1}{2} }}
\;\!
(\:\! t - s)^{\mbox{}^{\scriptstyle \!\!
-\, \frac{\scriptstyle m}{2} \,-\, \frac{3}{4} }}
\;\!
\bigl\{\, s^{\mbox{}^{\scriptstyle \frac{1}{2} }}
\;\!
\|\, D \:\! \mbox{\boldmath $u$}(\cdot,s) \,
\|_{\mbox{}_{\scriptstyle L^{2}(\mathbb{R}^{3})}}
\;\!\bigr\}
\, ds \\
&\leq\:
K\!\;\!(m,\;\!\nu) \:
\|\, \mbox{\boldmath $u$}_0 \;\!
\|_{\mbox{}_{\scriptstyle L^{2}(\mathbb{R}^{3})}}
\,
t^{\mbox{}^{\scriptstyle \frac{\scriptstyle m}{2} }}
\:\!
(\:\! t - t_0)^{\mbox{}^{\scriptstyle \!\!
-\, \frac{\scriptstyle m}{2} \,-\, \frac{3}{4} }}
\,
M_{\mbox{}_{1}}
\!\!\;\!
\int_{\mbox{\footnotesize $\;\!t_0$}}
    ^{\mbox{\footnotesize $\:\!\mu(t)$}}
\!\!
s^{\mbox{}^{\scriptstyle -\, \frac{1}{2} }}
\;\!
ds \\
&\leq\:
K\!\;\!(m,\;\!\nu, \;\!t_0)
\:
\|\, \mbox{\boldmath $u$}_0 \;\!
\|_{\mbox{}_{\scriptstyle L^{2}(\mathbb{R}^{3})}}
\,
M_{\mbox{}_{1}}
\;\!
(\:\! t - t_0)^{\mbox{}^{\scriptstyle \!\!-\, \frac{3}{4} }}
\:\!
(\:\! t + t_0)^{\mbox{}^{\scriptstyle \frac{1}{2} }} \\
\end{split}
\end{equation}
\mbox{} \vspace{-1.200cm} \\
\mbox{} \hfill (4.9$a$) \\
\mbox{} \vspace{+0.200cm} \\
para todo
$ \;\! t \geq t_0 + 1 $,
ou seja, \\
\mbox{} \vspace{-0.550cm} \\
\begin{equation}
\tag{4.9$b$}
|\: I_{2}(\alpha, \:\!t) \,|
\;\leq\:
K\!\;\!(m,\;\!\nu, \;\!t_0, \;\!M_{\mbox{}_{1}}\!\:\!, \;\!\mbox{\boldmath $u$}_0),
\qquad
\forall \;\;
t \,\geq\, t_0 \;\!+\;\! 1,
\end{equation}
\mbox{} \vspace{-0.175cm} \\
onde
$ {\displaystyle
\;\!
K\!\;\!(m,\;\!\nu, \;\!t_0, \:\!M_{\mbox{}_{1}}\!\:\!,
\:\!\mbox{\boldmath $u$}_0) > 0
\;\!
} $
depende dos dados
$ (m,\;\!\nu, \;\!t_0, \:\!M_{\mbox{}_{1}}\!\:\!,
\:\!\mbox{\boldmath $u$}_0) $.
Assim,
por (4.6), (4.8) e (4.9),
tem-se \\
\mbox{} \vspace{-0.600cm} \\
\begin{equation}
\tag{4.10}
V_{\alpha}(t)
\;\leq\:
K\!\:\!(m,\;\!\nu, \;\!t_0, \:\!M_{\mbox{}_{1}}\!\:\!,
\:\!\mbox{\boldmath $u$}_0)
\,+\;\!
J_{\alpha}(t),
\qquad
\forall \;\;
t \,\geq\, t_0 \;\!+\;\! 1,
\end{equation}
\mbox{} \vspace{-0.175cm} \\
onde
$ J_{\alpha}(t) $
\'e dada em (4.7$c$).
Para estimar $ J_{\alpha}(t) $,
que \'e o termo mais complicado,
podemos proceder
do modo seguinte.
Ilustraremos o m\'etodo
considerando inicialmente
o caso mais simples
$ |\,\alpha\,| = 2 $,
ou seja,
$ D^{\alpha} \!\;\!=\,\! D_{j} \:\!D_{\ell} $,
procedendo depois por indu\c c\~ao em
$ \;\!m = |\,\alpha\,| $.

Considerando
$ D^{\alpha} = D_{j} \:\!D_{\ell} $,
tem-se,
por (4.7$c$) e (2.14), \\
\mbox{} \vspace{-0.550cm} \\
\begin{equation}
\notag
\begin{split}
J_{j,\,\ell}(t)
\;&\equiv\;\;\!
t \!\;\!
\int_{\mbox{\footnotesize $\mu(t)$}}
    ^{\mbox{\footnotesize $\;\!t$}}
\!
\|\, D_{j} \:\!D_{\ell} \;\![\:
e^{\;\!\nu \Delta (\mbox{\footnotesize $t$} \;\!-\, \mbox{\footnotesize $s$})}
\:\!
\mbox{\boldmath $Q$}(\cdot,s)
\,] \:
\|_{\mbox{}_{\scriptstyle L^{2}(\mathbb{R}^{3})}}
\;\! ds, \\
&=\;\;\!
t \!\;\!
\int_{\mbox{\footnotesize $\mu(t)$}}
    ^{\mbox{\footnotesize $\;\!t$}}
\!
\|\, D_{j} \:\![\:
e^{\;\!\nu \Delta (\mbox{\footnotesize $t$} \;\!-\, \mbox{\footnotesize $s$})}
\:\!
D_{\ell} \;\! \mbox{\boldmath $Q$}(\cdot,s)
\,] \:
\|_{\mbox{}_{\scriptstyle L^{2}(\mathbb{R}^{3})}}
\;\! ds, \\
&\leq\,
K\!\;\!(\nu) \;\;\!
t \!\;\!
\int_{\mbox{\footnotesize $\mu(t)$}}
    ^{\mbox{\footnotesize $\;\!t$}}
\!
(\:\!t - s)^{\mbox{}^{\scriptstyle \!-\, \frac{7}{8} }}
\,\!
\|\, D_{\ell} \;\! \mbox{\boldmath $Q$}(\cdot,s) \,
\|_{\mbox{}_{\scriptstyle L^{4/3}(\mathbb{R}^{3})}}
\;\! ds.
\end{split}
\end{equation}
\mbox{} \vspace{+0.025cm} \\
Para prosseguir,
\'e preciso estimar
$ {\displaystyle
\;\!
\|\, D_{\ell} \;\! \mbox{\boldmath $Q$}(\cdot,s) \,
\|_{\scriptstyle L^{4/3}(\mathbb{R}^{3})}
\!\;\!
} $,
o que \'e feito usando
a teoria de Calderon-Zygmund
de operadores integrais singulares
(ver e.g.$\;\!$\cite{Stein1970}, Ch.$\;$2,
ou
\cite{DunfordSchwartz1963}, Ch.$\;$11).
Por conveni\^encia,
revisamos brevemente como isso \'e feito,
considerando o caso geral
em $ \mathbb{R}^{n} \!\;\!$. \linebreak
Por (2.4), Se\c c\~ao 2,
tem-se
$ {\displaystyle
\;\!
D_{\ell} \;\! \mbox{\boldmath $Q$}(\cdot,t)
\;\!=\;\!
-\, D_{\ell} \,[\, \mbox{\boldmath $u$} \cdot
\nabla \:\! \mbox{\boldmath $u$}\;\!(\cdot,t) \,]
\;\!-\;\! \nabla \:\! q_{\mbox{}_{\scriptstyle \ell}}(\cdot,t)
} $,
$ {\displaystyle
\;\!
q_{\mbox{}_{\scriptstyle \ell}}(\cdot,t)
=
D_{\ell} \,p(\cdot,t)
} $.
Tomando o divergente em $x$,
resulta que
$ \;\!q_{\mbox{}_{\scriptstyle \ell}}(\cdot,t) \;\!$
\'e a solu\c c\~ao (\'unica)
do proble\-ma de Poisson \\
\mbox{} \vspace{-0.875cm} \\
\begin{equation}
\notag
-\, \Delta \;\!
q_{\mbox{}_{\scriptstyle \ell}}(\cdot,t)
\;=\;
\mbox{div}\;
\bigl\{\, D_{\ell} \,[\, \mbox{\boldmath $u$} \cdot
\nabla \:\! \mbox{\boldmath $u$}\;\!(\cdot,t) \,]
\,\bigr\},
\quad \;\;
q_{\mbox{}_{\scriptstyle \ell}}(\cdot,t)
\in L^{2}(\mathbb{R}^{n}).
\end{equation}
\mbox{} \vspace{-0.250cm} \\
Aplicando a teoria de Calderon-Zygmund
(cf.$\;$\cite{Guterres2014}, Ch.$\;$5),
obt\'em-se,
para cada $ 1 < r < \infty $: \\
\mbox{} \vspace{-0.600cm} \\
\begin{equation}
\notag
\|\, \nabla \:\!
q_{\mbox{}_{\scriptstyle \ell}}(\cdot,t) \,
\|_{\mbox{}_{\scriptstyle L^{r}(\mathbb{R}^{n})}}
\,\leq\;
K\!\;\!(r,n) \:
\|\, D_{\ell} \,[\, \mbox{\boldmath $u$} \cdot
\nabla \:\! \mbox{\boldmath $u$}\;\!(\cdot,t) \,] \,
\|_{\mbox{}_{\scriptstyle L^{r}(\mathbb{R}^{n})}}
\end{equation}
\mbox{} \vspace{-0.240cm} \\
para certa constante
$ K\!\;\!(r,n) > 0 $
dependendo de $ \;\!r, \;\!n $
somente.
Isso implica,
por (2.4),
que \\
\mbox{} \vspace{-0.600cm} \\
\begin{equation}
\tag{4.11}
\|\, D_{\ell} \;\!
\mbox{\boldmath $Q$}(\cdot,t) \,
\|_{\mbox{}_{\scriptstyle L^{r}(\mathbb{R}^{n})}}
\,\leq\;
K\!\;\!(r,n) \:
\|\, D_{\ell} \,[\, \mbox{\boldmath $u$} \cdot
\nabla \:\! \mbox{\boldmath $u$}\;\!(\cdot,t) \,] \,
\|_{\mbox{}_{\scriptstyle L^{r}(\mathbb{R}^{n})}}
\end{equation}
\mbox{} \vspace{-0.225cm} \\
para cada $ 1 < r < \infty $,
onde, novamente,
$ K\!\;\!(r,n) $
depende de $ \;\!r, \;\!n $.
(Por exemplo, repetindo o argumento usado
na prova da Proposi\c c\~ao 2.1,
tem-se $ K\!\;\!(2,n) = 1 $, para todo $n$.)
Assim,
com $ \;\!r = 4/3 $,
$ n = 3 $,
obt\'em-se
$ {\displaystyle
\;\!
\|\, D_{\ell} \, \mbox{\boldmath $Q$}(\cdot,s) \,
\|_{\mbox{}_{\scriptstyle L^{4/3}(\mathbb{R}^{3})}}
\!\;\!\leq\;\!
K \:
\|\, D_{\ell} \,[\, \mbox{\boldmath $u$}(\cdot,s) \cdot
\nabla \:\! \mbox{\boldmath $u$}(\cdot,s) \,] \,
\|_{\mbox{}_{\scriptstyle L^{4/3}(\mathbb{R}^{3})}}
\!
} $,
de modo que \\
\mbox{} \vspace{-0.150cm} \\
\mbox{} \hspace{+1.500cm}
$ {\displaystyle
J_{j,\,\ell}(t)
\;\leq\;
K\!\;\!(\nu) \;
t \!
\int_{\mbox{\footnotesize $\mu(t)$}}
    ^{\mbox{\footnotesize $\;\!t$}}
\!
(\:\!t - s)^{\mbox{}^{\scriptstyle \!-\, \frac{7}{8} }}
\,\!
\|\, D_{\ell} \;\![\, \mbox{\boldmath $u$}(\cdot,s) \cdot
\nabla \:\! \mbox{\boldmath $u$}(\cdot,s) \,]\,
\|_{\mbox{}_{\scriptstyle L^{4/3}(\mathbb{R}^{3})}}
\;\! ds
} $ \\
\mbox{} \vspace{+0.150cm} \\
\mbox{} \hspace{+2.800cm}
$ {\displaystyle
\leq\;
K\!\;\!(\nu) \;
t \!
\int_{\mbox{\footnotesize $\mu(t)$}}
    ^{\mbox{\footnotesize $\;\!t$}}
\!
(\:\!t - s)^{\mbox{}^{\scriptstyle \!-\, \frac{7}{8} }}
\;\!
\Bigl\{\,
\|\, D_{\ell} \;\! \mbox{\boldmath $u$}(\cdot,s)
\;\!\cdot\;\!
\nabla \:\! \mbox{\boldmath $u$}(\cdot,s) \,
\|_{\mbox{}_{\scriptstyle L^{4/3}(\mathbb{R}^{3})}}
\:+
} $ \\
\mbox{} \vspace{-0.150cm} \\
\mbox{} \hfill
$ {\displaystyle
+\;\;\!
\|\, \mbox{\boldmath $u$}(\cdot,s)
\;\!\cdot\;\!
\nabla \:\! D_{\ell} \;\! \mbox{\boldmath $u$}(\cdot,s) \,
\|_{\mbox{}_{\scriptstyle L^{4/3}(\mathbb{R}^{3})}}
\;\! \Bigr\}
\;\! ds
} $ \\
\mbox{} \vspace{-0.050cm} \\
\mbox{} \hspace{+2.750cm}
$ {\displaystyle
\leq\;
K\!\;\!(\nu) \;
t \!
\int_{\mbox{\footnotesize $\mu(t)$}}
    ^{\mbox{\footnotesize $\;\!t$}}
\!
(\:\!t - s)^{\mbox{}^{\scriptstyle \!-\, \frac{7}{8} }}
\;\!
\Bigl\{\,
\|\, D \:\!\mbox{\boldmath $u$}(\cdot,s) \,
\|_{\mbox{}_{\scriptstyle L^{4}(\mathbb{R}^{3})}}
\:\!
\|\, D \:\!\mbox{\boldmath $u$}(\cdot,s) \,
\|_{\mbox{}_{\scriptstyle L^{2}(\mathbb{R}^{3})}}
\:+
} $ \\
\mbox{} \vspace{-0.150cm} \\
\mbox{} \hfill
$ {\displaystyle
+\;\;\!
\|\, \mbox{\boldmath $u$}(\cdot,s) \,
\|_{\mbox{}_{\scriptstyle L^{4}(\mathbb{R}^{3})}}
\:\!
\|\, D^{2} \,\! \mbox{\boldmath $u$}(\cdot,s) \,
\|_{\mbox{}_{\scriptstyle L^{2}(\mathbb{R}^{3})}}
\;\! \Bigr\}
\;\! ds
} $, \\
\mbox{} \vspace{+0.050cm} \\
pela
desigualdade de H\"older
(no \'ultimo passo acima).
Utilizando agora
a desigualdade
de Nirenberg-Gagliardo \\
\mbox{} \vspace{-0.750cm} \\
\begin{equation}
\tag{4.12}
\|\, \mbox{u} \,
\|_{\mbox{}_{\scriptstyle L^{4}(\mathbb{R}^{3})}}
\;\!\leq\;
K \,
\|\, \mbox{u} \,
\|_{\mbox{}_{\scriptstyle L^{2}(\mathbb{R}^{3})}}
  ^{1/4}
\:\!
\|\, D \:\! \mbox{u} \,
\|_{\mbox{}_{\scriptstyle L^{2}(\mathbb{R}^{3})}}
  ^{3/4}
\end{equation}
\mbox{} \vspace{-0.150cm} \\
para $ \;\!\mbox{u} \in H^{1}(\mathbb{R}^{3}) \;\!$
arbitr\'aria,
obt\'em-se,
ent\~ao, \\
\mbox{} \vspace{-0.025cm} \\
\mbox{} \hspace{+0.150cm}
$ {\displaystyle
J_{j,\,\ell}(t)
\;\leq\;
K\!\;\!(\nu) \;
t \!
\int_{\mbox{\footnotesize $\mu(t)$}}
    ^{\mbox{\footnotesize $\;\!t$}}
\!
(\:\!t - s)^{\mbox{}^{\scriptstyle \!-\, \frac{7}{8} }}
\;\!
\Bigl\{\,
\|\, D \:\!\mbox{\boldmath $u$}(\cdot,s) \,
\|_{\mbox{}_{\scriptstyle L^{4}(\mathbb{R}^{3})}}
  ^{\:\!5/4}
\:\!
\|\, D^{2} \,\!\mbox{\boldmath $u$}(\cdot,s) \,
\|_{\mbox{}_{\scriptstyle L^{2}(\mathbb{R}^{3})}}
  ^{\:\!3/4}
\:+
} $ \\
\mbox{} \vspace{-0.150cm} \\
\mbox{} \hfill
$ {\displaystyle
+\;\;\!
\|\, \mbox{\boldmath $u$}(\cdot,s) \,
\|_{\mbox{}_{\scriptstyle L^{4}(\mathbb{R}^{3})}}
  ^{\:\!1/4}
\:\!
\|\, D \:\!\mbox{\boldmath $u$}(\cdot,s) \,
\|_{\mbox{}_{\scriptstyle L^{2}(\mathbb{R}^{3})}}
  ^{\:\!3/4}
\:\!
\|\, D^{2} \,\! \mbox{\boldmath $u$}(\cdot,s) \,
\|_{\mbox{}_{\scriptstyle L^{2}(\mathbb{R}^{3})}}
\;\! \Bigr\}
\;\! ds
} $, \\
\mbox{} \vspace{+0.050cm} \\
\mbox{} \hspace{+1.500cm}
$ {\displaystyle
=\:
K\!\;\!(\nu) \;
t \!
\int_{\mbox{\footnotesize $\mu(t)$}}
    ^{\mbox{\footnotesize $\;\!t$}}
\!\!\!\!\;\!
s^{\mbox{}^{\scriptstyle \!-\, \frac{11}{8} }}
(\:\!t - s)^{\mbox{}^{\scriptstyle \!-\, \frac{7}{8} }}
\;\!
\Bigl\{\;\!
\bigl[\, s^{\mbox{}^{\scriptstyle \frac{1}{2} }}
\:\!
\|\, D \:\!\mbox{\boldmath $u$}(\cdot,s) \,
\|_{\mbox{}_{\scriptstyle L^{2}(\mathbb{R}^{3})}}
\:\!\bigr]^{\:\!5/4}
\;\!
\bigl[\, s \,
\|\, D^{2} \,\!\mbox{\boldmath $u$}(\cdot,s) \,
\|_{\mbox{}_{\scriptstyle L^{2}(\mathbb{R}^{3})}}
\bigr]^{\:\!3/4}
} $ \\
\mbox{} \vspace{-0.150cm} \\
\mbox{} \hfill
$ {\displaystyle
+\;\;\!
\|\, \mbox{\boldmath $u$}(\cdot,s) \,
\|_{\mbox{}_{\scriptstyle L^{4}(\mathbb{R}^{3})}}
  ^{\:\!1/4}
\:\!
\bigl[\, s^{\mbox{}^{\scriptstyle \frac{1}{2} }}
\:\!
\|\, D \:\!\mbox{\boldmath $u$}(\cdot,s) \,
\|_{\mbox{}_{\scriptstyle L^{2}(\mathbb{R}^{3})}}
\bigr]^{\:\!3/4}
\:\!
\bigl[\, s \,
\|\, D^{2} \,\! \mbox{\boldmath $u$}(\cdot,s) \,
\|_{\mbox{}_{\scriptstyle L^{2}(\mathbb{R}^{3})}}
\bigr]
\;\! \Bigr\}
\;\! ds
} $, \\
\mbox{} \vspace{+0.150cm} \\
\mbox{} \hspace{+1.500cm}
$ {\displaystyle
\leq \:
K\!\;\!(\nu) \,
M_{\mbox{}_{1}}^{\mbox{}^{\scriptstyle \frac{5}{4} }}
\;\!
t \;\!\;\!
(\:\!t + t_0)^{\mbox{}^{\scriptstyle \!\!-\, \frac{11}{8} }}
\!\!\!\!\;\!
\int_{\mbox{\footnotesize $\mu(t)$}}
    ^{\mbox{\footnotesize $\;\!t$}}
\!\!
(\:\!t - s)^{\mbox{}^{\scriptstyle \!-\, \frac{7}{8} }}
\;\!
U_{\mbox{}_{2}}(s)^{\mbox{}^{\scriptstyle \frac{3}{4} }}
\;\! ds
\;\,+
} $ \\
\mbox{} \vspace{-0.150cm} \\
\mbox{} \hfill
$ {\displaystyle
+\;\;\!
K\!\;\!(\nu) \,
M_{\mbox{}_{1}}^{\mbox{}^{\scriptstyle \frac{3}{4} }}
\;\!
\|\, \mbox{\boldmath $u$}_0 \;\!
\|_{\mbox{}_{\scriptstyle L^{2}(\mathbb{R}^{3})}}
  ^{\mbox{}^{\scriptstyle \frac{1}{4} }}
\:\!
t \;\!\;\!
(\:\!t + t_0)^{\mbox{}^{\scriptstyle \!\!-\, \frac{11}{8} }}
\!\!\!\!\;\!
\int_{\mbox{\footnotesize $\mu(t)$}}
    ^{\mbox{\footnotesize $\;\!t$}}
\!\!
(\:\!t - s)^{\mbox{}^{\scriptstyle \!-\, \frac{7}{8} }}
\;\!
U_{\mbox{}_{2}}(s)
\;\!\;\! ds
} $, \\
\mbox{} \vspace{+0.100cm} \\
onde
$ {\displaystyle
U_{\mbox{}_{2}}(t)
\;\!=\;\!
t \;
\|\,D^{2} \:\! \mbox{\boldmath $u$}(\cdot,t) \,
\|_{\mbox{}_{\scriptstyle L^{2}(\mathbb{R}^{3})}}
} $.
Portanto,
para cada
$ 1 \leq j, \;\!\ell \leq 3 $
(e cada $ t \geq t_0 + 1$),
obt\'em-se \\
\mbox{} \vspace{-0.200cm} \\
\mbox{} \hspace{+1.000cm}
$ {\displaystyle
J_{j,\,\ell}(t)
\;\leq \;\!\;\!
K\!\;\!(\nu) \,
M_{\mbox{}_{1}}^{\mbox{}^{\scriptstyle \frac{5}{4} }}
\,\!
(\:\!t + t_0)^{\mbox{}^{\scriptstyle \!\!-\, \frac{3}{8} }}
\!\!\!
\int_{\mbox{\footnotesize $\mu(t)$}}
    ^{\mbox{\footnotesize $\;\!t$}}
\!\!
(\:\!t - s)^{\mbox{}^{\scriptstyle \!-\, \frac{7}{8} }}
\;\!
U_{\mbox{}_{2}}(s)^{\mbox{}^{\scriptstyle \frac{3}{4} }}
\;\! ds
\;\,+
} $ \\
\mbox{} \vspace{-0.150cm} \\
\mbox{} \hfill
$ {\displaystyle
+\;\;\!
K\!\;\!(\nu) \,
M_{\mbox{}_{1}}^{\mbox{}^{\scriptstyle \frac{3}{4} }}
\;\!
\|\, \mbox{\boldmath $u$}_0 \;\!
\|_{\mbox{}_{\scriptstyle L^{2}(\mathbb{R}^{3})}}
  ^{\mbox{}^{\scriptstyle \frac{1}{4} }}
\:\!
(\:\!t + t_0)^{\mbox{}^{\scriptstyle \!\!-\, \frac{3}{8} }}
\!\!\!\!\;\!
\int_{\mbox{\footnotesize $\mu(t)$}}
    ^{\mbox{\footnotesize $\;\!t$}}
\!\!
(\:\!t - s)^{\mbox{}^{\scriptstyle \!-\, \frac{7}{8} }}
\;\!
U_{\mbox{}_{2}}(s)
\;\!\;\! ds
} $, \\
\mbox{} \vspace{+0.150cm} \\
\mbox{} \hspace{+2.300cm}
$ {\displaystyle
\leq \;\!\;\!
K\!\;\!(\nu) \,
M_{\mbox{}_{1}}^{\mbox{}^{\scriptstyle \frac{5}{4} }}
\,\!
(\:\!t + t_0)^{\mbox{}^{\scriptstyle \!\!-\, \frac{1}{4} }}
+\;
K\!\;\!(\nu, M_{\mbox{}_{1}}\!\:\!, \mbox{\boldmath $u$}_0)
\,
(\:\!t + t_0)^{\mbox{}^{\scriptstyle \!\!-\, \frac{3}{8} }}
\!\!\!
\int_{\mbox{\footnotesize $\mu(t)$}}
    ^{\mbox{\footnotesize $\;\!t$}}
\!\!
(\:\!t - s)^{\mbox{}^{\scriptstyle \!-\, \frac{7}{8} }}
\;\!
U_{\mbox{}_{2}}(s)
\: ds
} $, \\
\mbox{} \vspace{+0.150cm} \\
pelo fato de se ter
$ {\displaystyle
\;\!
U_{\mbox{}_{2}}(s)^{\mbox{}^{\scriptstyle \:\!\frac{3}{4} }}
\!\:\!\leq\,\! 1 + U_{\mbox{}_{2}}(s)
} $,
e onde
$ \;\!K\!\;\!(\nu, M_{\mbox{}_{1}}\!\:\!, \mbox{\boldmath $u$}_0) > 0 \;\!$
depende apenas de
$ \nu $, $ M_{\mbox{}_{1}}\!\:\! $ e
$ \;\!\|\,\mbox{\boldmath $u$}_0 \;\! \|_{L^{2}(\mathbb{R}^{3})} $.
Assim, por (4.6) e (4.10),
obt\'em-se,
para
$ {\displaystyle
\;\!
V_{j,\;\!\ell}(t)
\;\!=\;\!
t \:
\|\, D_{j} \:\! D_{\ell} \, \mbox{\boldmath $u$}(\cdot,t) \,
\|_{\mbox{}_{\scriptstyle L^{2}(\mathbb{R}^{3})}}
\!\:\!
} $: \\
\mbox{} \vspace{-0.600cm} \\
\begin{equation}
\notag
V_{j,\;\!\ell}(t)
\;\leq\;
K\!\;\!(m,\;\!\nu, \;\!t_0, \;\!M_{\mbox{}_{1}}\!\:\!, \;\!\mbox{\boldmath $u$}_0)
\;+\;
K\!\;\!(\nu, M_{\mbox{}_{1}}\!\:\!, \mbox{\boldmath $u$}_0)
\,
(\:\!t + t_0)^{\mbox{}^{\scriptstyle \!\!-\, \frac{3}{8} }}
\!\!\!
\int_{\mbox{\footnotesize $\mu(t)$}}
    ^{\mbox{\footnotesize $\;\!t$}}
\!\!
(\:\!t - s)^{\mbox{}^{\scriptstyle \!-\, \frac{7}{8} }}
\;\!
U_{\mbox{}_{2}}(s)
\: ds
\end{equation}
\mbox{} \vspace{-0.150cm} \\
para todos $ j, \;\!\ell $
(e todo $ \;\! t \;\!\geq\;\! t_0 + 1 $),
de modo que \\
\mbox{} \vspace{-0.600cm} \\
\begin{equation}
\tag{4.13}
U_{\mbox{}_{2}}(t)
\;\leq\;
K_{\mbox{}_{\scriptstyle \ast}}\!\;\!
(m,\;\!\nu, \;\!t_0, \;\!M_{\mbox{}_{1}}\!\:\!, \;\!\mbox{\boldmath $u$}_0)
\;+\;
K_{\mbox{}_{\scriptstyle \ast\ast}}\!\;\!
(\nu, M_{\mbox{}_{1}}\!\:\!, \mbox{\boldmath $u$}_0)
\,
(\:\!t + t_0)^{\mbox{}^{\scriptstyle \!\!-\, \frac{3}{8} }}
\!\!\!
\int_{\mbox{\footnotesize $\mu(t)$}}
    ^{\mbox{\footnotesize $\;\!t$}}
\!\!
(\:\!t - s)^{\mbox{}^{\scriptstyle \!-\, \frac{7}{8} }}
\;\!
U_{\mbox{}_{2}}(s)
\: ds
\end{equation}
\mbox{} \vspace{-0.150cm} \\
para todo $ \;\! t \geq t_0 + 1 $,
onde
$ {\displaystyle
\;\!
K_{\mbox{}_{\scriptstyle \ast\ast}}\!\;\!
(\nu, M_{\mbox{}_{1}}\!\:\!, \mbox{\boldmath $u$}_0) > 0
\;\!
} $
depende apenas de
$ \nu $, $ M_{\mbox{}_{1}}\!\:\! $ e
$ \;\!\|\,\mbox{\boldmath $u$}_0 \;\! \|_{L^{2}(\mathbb{R}^{3})} $.
Tomemos ent\~ao
$ \;\!t_{2} $,
$ \mathbb{M}_{2} $
dados por \\
\begin{equation}
\tag{4.14}
t_{2}\;\!:=\;
1 +\;\! t_0 +\, 2^{16} \;\!
K_{\mbox{}_{\scriptstyle \!\;\!\ast\ast}}^{\:\!4}
\!\;\!,
\qquad
\mathbb{M}_{2} \;\!\!:=\;
\sup \,\{\, U_{\mbox{}_{2}}(s) \!:\;
t_0 \leq s \leq t_2 \;\!\},
\end{equation}
\mbox{} \vspace{-0.200cm} \\
onde
$ {\displaystyle
\;\!
K_{\mbox{}_{\scriptstyle \!\;\!\ast\ast}}
\!\;\!=\;\!
K_{\mbox{}_{\scriptstyle \!\;\!\ast\ast}}\!\:\!
(\nu, M_{\mbox{}_{1}}\!\:\!, \mbox{\boldmath $u$}_0)
\;\!
} $
\'e a constante
definida em (4.13) acima.
Afirmamos que \\
\mbox{} \vspace{-0.500cm} \\
\begin{equation}
\tag{4.15}
U_{\mbox{}_{2}}(t) \;\leq\;
2 \, K_{\mbox{}_{\scriptstyle \ast}}\!\;\!
(m,\;\!\nu, \;\!t_0, \;\!M_{\mbox{}_{1}}\!\:\!, \;\!\mbox{\boldmath $u$}_0)
\;+\;
16 \: K_{\mbox{}_{\scriptstyle \ast\ast}}\!\;\!
(\nu, M_{\mbox{}_{1}}\!\:\!, \mbox{\boldmath $u$}_0) \,
\mathbb{M}_{\mbox{}_{2}},
\qquad
\forall \;\,
t \;\!\geq\;\! t_{2},
\end{equation}
\mbox{} \vspace{-0.150cm} \\
onde
$ \;\! K_{\mbox{}_{\scriptstyle \!\;\!\ast}} \!\;\! $,
$ K_{\mbox{}_{\scriptstyle \!\;\!\ast\ast}} $
s\~ao as constantes dadas em (4.13).
\mbox{[}$\,$De fato,
sendo $ t \geq t_{2} $,
definamos
$ \mathbb{U}_{\mbox{}_{2}}(t) $
pondo
$ {\displaystyle
\;\!
\mathbb{U}_{\mbox{}_{2}}(t) :=\,
\sup\,\{\;\! U_{\mbox{}_{2}}(s) \!:\; t_{2} \leq s \leq t \,\}
} $.
Se $ \;\!\mu(t) \geq t_{2} $,
ent\~ao,
por (4.13),
obt\'em-se \\
\mbox{} \vspace{-0.600cm} \\
\begin{equation}
\notag
\begin{split}
U_{\mbox{}_{2}}(t)
\;&\leq\;
K_{\mbox{}_{\scriptstyle \ast}}\!\;\!
(m,\;\!\nu, \;\!t_0, \;\!M_{\mbox{}_{1}}\!\:\!, \;\!\mbox{\boldmath $u$}_0)
\;+\;
K_{\mbox{}_{\scriptstyle \ast\ast}}\!\;\!
(\nu, M_{\mbox{}_{1}}\!\:\!, \mbox{\boldmath $u$}_0)
\,
(\:\!t + t_0)^{\mbox{}^{\scriptstyle \!\!-\, \frac{3}{8} }}
\!\!\!
\int_{\mbox{\footnotesize $\mu(t)$}}
    ^{\mbox{\footnotesize $\;\!t$}}
\!\!
(\:\!t - s)^{\mbox{}^{\scriptstyle \!-\, \frac{7}{8} }}
\;\!
ds
\;
\mathbb{U}_{\mbox{}_{2}}(t) \\
&\leq\;
K_{\mbox{}_{\scriptstyle \ast}}\!\;\!
(m,\;\!\nu, \;\!t_0, \;\!M_{\mbox{}_{1}}\!\:\!, \;\!\mbox{\boldmath $u$}_0)
\;+\;
8 \, K_{\mbox{}_{\scriptstyle \ast\ast}}\!\;\!
(\nu, M_{\mbox{}_{1}}\!\:\!, \mbox{\boldmath $u$}_0)
\,
(\:\!t + t_0)^{\mbox{}^{\scriptstyle \!\!-\, \frac{1}{4} }}
\;\!
\mathbb{U}_{\mbox{}_{2}}(t),
\end{split}
\end{equation}
\mbox{} \vspace{+0.100cm} \\
de modo que, pela defini\c c\~ao de $\;\!t_{2} $,
tem-se
$ {\displaystyle
\;\!
U_{\mbox{}_{2}}(t) \;\!\leq\;\!
K_{\mbox{}_{\scriptstyle \ast}}
\!\;\!+\;\! \mathbb{U}_{\mbox{}_{2}}(t)/2
} $.
Se $ \;\!\mu(t) < t_{2} $,
ent\~ao \\
\mbox{} \vspace{-0.575cm} \\
\begin{equation}
\notag
\begin{split}
U_{\mbox{}_{2}}(t)
\;&\leq\;
K_{\mbox{}_{\scriptstyle \ast}}\!\;\!
\;\!+\,
K_{\mbox{}_{\scriptstyle \ast\ast}}\!
\cdot
(\:\!t + t_0)^{\mbox{}^{\scriptstyle \!\!-\, \frac{3}{8} }}
\!\!\!
\int_{\mbox{\footnotesize $\mu(t)$}}
    ^{\mbox{\footnotesize $\;\!t$}}
\!\!
(\:\!t - s)^{\mbox{}^{\scriptstyle \!-\, \frac{7}{8} }}
\;\!
ds
\;
\bigl(\, \mathbb{M}_{\mbox{}_{2}} \;\!+\,
\mathbb{U}_{\mbox{}_{2}}(t) \;\!\bigr) \\
&\leq\;
K_{\mbox{}_{\scriptstyle \ast}}\!\;\!
\;\!+\;
8 \, K_{\mbox{}_{\scriptstyle \ast\ast}}\;\!
\mathbb{M}_{\mbox{}_{2}}
\,+\,
8 \, K_{\mbox{}_{\scriptstyle \ast\ast}}\;\!
(\:\!t + t_0)^{\mbox{}^{\scriptstyle \!\!-\, \frac{1}{4} }}
\;\!
\mathbb{U}_{\mbox{}_{2}}(t),
\end{split}
\end{equation}
\mbox{} \vspace{-0.200cm} \\
de modo que,
pela defini\c c\~ao de $\;\!t_{2} $,
obt\'em-se
neste caso
$ {\displaystyle
\;\!
U_{\mbox{}_{2}}(t)
\;\!\leq\;\!
K_{\mbox{}_{\scriptstyle \ast}}
\!\;\!+\;\!
8 \, K_{\mbox{}_{\scriptstyle \ast\ast}} \,
\mathbb{M}_{\mbox{}_{2}}
+\;\!
\mathbb{U}_{\mbox{}_{2}}(t)/2
} $.
Portanto,
tem-se sempre \\
\mbox{} \vspace{-0.650cm} \\
\begin{equation}
\notag
U_{\mbox{}_{2}}(t)
\;\leq\;
K_{\mbox{}_{\scriptstyle \ast}}
\!\;\!+\,
8 \, K_{\mbox{}_{\scriptstyle \ast\ast}}
\;\!
\mathbb{M}_{\mbox{}_{2}}
\;\!+\,
\frac{1}{2} \;
\mathbb{U}_{\mbox{}_{2}}(t),
\qquad
\forall \;\,
t \;\!\geq\;\! t_{2}.
\end{equation}
\mbox{} \vspace{-0.250cm} \\
Logo,
$ {\displaystyle
\;\!
\mathbb{U}_{\mbox{}_{2}}(t)
\;\!\leq\;\!
K_{\mbox{}_{\scriptstyle \ast}}
\!\,\!+\;\!
8 \, K_{\mbox{}_{\scriptstyle \ast\ast}}
\;\!
\mathbb{M}_{\mbox{}_{2}}
+\;\!
\mathbb{U}_{\mbox{}_{2}}(t)/2
\;\!
} $
para todo $ \;\!t \geq t_{2} $,
que \'e equivalente a (4.15).$\;\!$\mbox{]} \linebreak
De (4.15),
segue imediatamente
que
$ {\displaystyle
\;\!
U_{\mbox{}_{2}} \in L^{\infty}(\:\![\, t_2, \infty) )
} $.
Como,
por (1.3),
tem-se tamb\'em
$ {\displaystyle
\;\!
U_{\mbox{}_{2}} \in L^{\infty}(\:\![\, \mbox{\small $T$}_{\!\;\!\ast\ast}, \;\!t_2\;\!]\:\! )
} $,
resulta que
$ {\displaystyle
\;\!
U_{\mbox{}_{2}} \in L^{\infty}(\:\![\,\mbox{\small $T$}_{\!\;\!\ast\ast}, \infty) )
} $,
provando assim (4.5)
no caso $ m = 2 $. \\
\mbox{} \vspace{-0.750cm} \\
%
%

Mais geralmente,
podemos mostrar (4.5)
para $ m \geq 2 $ qualquer
usando indu\c c\~ao em $m$.
Assim,
dado $ m \geq 3 $,
suponha-se que (4.5)
j\'a tenha sido obtida para
os valores anteriores de $m$.
Tomando-se,
ent\~ao,
$ \alpha $ (multi-\'\i ndice)
com $ \;\!|\,\alpha\,| = m $,
e escrevendo-se
$ D^{\alpha} \!\;\!= D_{j} \;\! D^{\alpha^{\prime}} \!$
(para certo
multi-\'\i ndice $ \alpha^{\prime} $
com $ |\,\alpha^{\prime} \,| = m - 1 $),
tem-se,
lembrando (4.7$c$) acima, \\
\mbox{} \vspace{-0.525cm} \\
\begin{equation}
\notag
\begin{split}
J_{\alpha}(t)
\;&=\;\;\!
t^{\mbox{}^{\scriptstyle \frac{\scriptstyle m}{2} }}
\!\!\!
\int_{\mbox{\footnotesize $\mu(t)$}}
    ^{\mbox{\footnotesize $\;\!t$}}
\!
\|\, D^{\alpha} \;\![\:
e^{\;\!\nu \Delta (\mbox{\footnotesize $t$} \;\!-\, \mbox{\footnotesize $s$})}
\:\!
\mbox{\boldmath $Q$}(\cdot,s)
\,] \:
\|_{\mbox{}_{\scriptstyle L^{2}(\mathbb{R}^{3})}}
\;\! ds, \\
&=\;\;\!
t^{\mbox{}^{\scriptstyle \frac{\scriptstyle m}{2} }}
\!\!\!
\int_{\mbox{\footnotesize $\mu(t)$}}
    ^{\mbox{\footnotesize $\;\!t$}}
\!
\|\, D_{j} \,[\:
e^{\;\!\nu \Delta (\mbox{\footnotesize $t$} \;\!-\, \mbox{\footnotesize $s$})}
\:\!
D^{\alpha^{\prime}} \!\;\! \mbox{\boldmath $Q$}(\cdot,s)
\,] \:
\|_{\mbox{}_{\scriptstyle L^{2}(\mathbb{R}^{3})}}
\;\! ds, \\
&\leq\,
K\!\;\!(\nu) \;\;\!
t^{\mbox{}^{\scriptstyle \frac{\scriptstyle m}{2} }}
\!\!\!
\int_{\mbox{\footnotesize $\mu(t)$}}
    ^{\mbox{\footnotesize $\;\!t$}}
\!
(\:\!t - s)^{\mbox{}^{\scriptstyle \!-\, \frac{7}{8} }}
\,\!
\|\, D^{\alpha^{\prime}} \!\;\! \mbox{\boldmath $Q$}(\cdot,s) \,
\|_{\mbox{}_{\scriptstyle L^{4/3}(\mathbb{R}^{3})}}
\;\! ds \\
&\leq\,
K\!\;\!(\nu) \;\;\!
t^{\mbox{}^{\scriptstyle \frac{\scriptstyle m}{2} }}
\!\!\!
\int_{\mbox{\footnotesize $\mu(t)$}}
    ^{\mbox{\footnotesize $\;\!t$}}
\!
(\:\!t - s)^{\mbox{}^{\scriptstyle \!-\, \frac{7}{8} }}
\,\!
\|\, D^{\alpha^{\prime}} \!\;\![\,
\mbox{\boldmath $u$}(\cdot,s) \cdot
\nabla \:\! \mbox{\boldmath $u$}(\cdot,s) \,]\,
\|_{\mbox{}_{\scriptstyle L^{4/3}(\mathbb{R}^{3})}}
\;\! ds \\
\end{split}
\end{equation}
\mbox{} \vspace{+0.025cm} \\
usando (2.14)
e a estimativa
(4.16) abaixo
(que segue novamente por Calderon-Zygmund),

\mbox{} \vspace{-1.500cm} \\
\begin{equation}
\tag{4.16}
\|\, D^{\beta} \,\!\,\!\mbox{\boldmath $Q$}(\cdot,t) \,
\|_{\mbox{}_{\scriptstyle L^{r}(\mathbb{R}^{n})}}
\,\leq\;
K\!\;\!(r,n) \:
\|\, D^{\beta} \:\! [\,
\mbox{\boldmath $u$}(\cdot,t) \cdot
\nabla \:\! \mbox{\boldmath $u$}(\cdot,t) \,]\,
\|_{\mbox{}_{\scriptstyle L^{r}(\mathbb{R}^{n})}}
\!\;\!,
\end{equation}
\mbox{} \vspace{-0.200cm} \\
para cada $ 1 < r < \infty $,
e qualquer
multi-\'\i ndice $ \beta $,
que \'e mostrada
de modo an\'alogo a (4.11). \linebreak
Resulta,
ent\~ao,
por H\"older e (4.12),
como antes, \\
\mbox{} \vspace{-0.575cm} \\
\begin{equation}
\notag
\begin{split}
J_{\alpha}(t)
\;\;\!&\leq
\hspace{-1.250cm}
\sum_{\mbox{} \hspace{+1.250cm} |\,\beta\,| \,+\, |\,\gamma\,| \:=\:
m \,-\, {\scriptscriptstyle 1}}
\hspace{-1.400cm}
K\!\;\!(\nu) \;\;\!
t^{\mbox{}^{\scriptstyle \frac{\scriptstyle m}{2} }}
\!\!\!
\int_{\mbox{\footnotesize $\mu(t)$}}
    ^{\mbox{\footnotesize $\;\!t$}}
\!
(\:\!t - s)^{\mbox{}^{\scriptstyle \!-\, \frac{7}{8} }}
\,\!
\|\, D^{\beta} \,\! \mbox{\boldmath $u$}(\cdot,s) \cdot
\nabla \:\!  D^{\gamma} \,\! \mbox{\boldmath $u$}(\cdot,s) \,]\,
\|_{\mbox{}_{\scriptstyle L^{4/3}(\mathbb{R}^{3})}}
\;\! ds \\
&\leq\;\;
K\!\;\!(\nu) \!
\sum_{\ell \,=\,0}^{m\,-\,{\scriptscriptstyle 1}}
t^{\mbox{}^{\scriptstyle \frac{\scriptstyle m}{2} }}
\!\!\!
\int_{\mbox{\footnotesize $\mu(t)$}}
    ^{\mbox{\footnotesize $\;\!t$}}
\!
(\:\!t - s)^{\mbox{}^{\scriptstyle \!-\, \frac{7}{8} }}
\,\!
\|\, D^{\ell} \,\! \mbox{\boldmath $u$}(\cdot,s) \,
\|_{\mbox{}_{\scriptstyle L^{4}(\mathbb{R}^{3})}}
\:\!
\|\, D^{m - \ell} \,\! \mbox{\boldmath $u$}(\cdot,s) \,
\|_{\mbox{}_{\scriptstyle L^{2}(\mathbb{R}^{3})}}
\;\! ds \\
&\leq\;\;
K\!\;\!(\nu) \!
\sum_{\ell \,=\,0}^{m\,-\,{\scriptscriptstyle 1}}
t^{\mbox{}^{\scriptstyle \frac{\scriptstyle m}{2} }}
\!\!\!
\int_{\mbox{\footnotesize $\mu(t)$}}
    ^{\mbox{\footnotesize $\;\!t$}}
\!
(\:\!t - s)^{\mbox{}^{\scriptstyle \!-\, \frac{7}{8} }}
\,\!
\|\, D^{\ell} \,\! \mbox{\boldmath $u$}(\cdot,s) \,
\|_{\mbox{}_{\scriptstyle L^{2} }}
  ^{\mbox{}^{\scriptstyle \frac{1}{4} }}
\:\!
\|\, D^{\ell + 1} \,\! \mbox{\boldmath $u$}(\cdot,s) \,
\|_{\mbox{}_{\scriptstyle L^{2} }}
  ^{\mbox{}^{\scriptstyle \frac{3}{4} }}
\:\!
\|\, D^{m - \ell} \,\! \mbox{\boldmath $u$}(\cdot,s) \,
\|_{\mbox{}_{\scriptstyle L^{2} }}
\, ds, \\
\end{split}
\end{equation}
\mbox{} \vspace{+0.150cm} \\
ou seja,
$ {\displaystyle
\;\!
J_{\alpha}(t)
\;\!\leq\;\!
J_{\mbox{}_{1}}(t) \;\!+
J_{\mbox{}_{2}}(t) \;\!+
J_{\mbox{}_{3}}(t)
} $,
onde \\
\mbox{} \vspace{-0.500cm} \\
\begin{equation}
\tag{4.17$a$}
J_{\mbox{}_{1}}(t)
\;=\;
K\!\;\!(\nu) \;
t^{\mbox{}^{\scriptstyle \frac{\scriptstyle m}{2} }}
\!\!\!
\int_{\mbox{\footnotesize $\mu(t)$}}
    ^{\mbox{\footnotesize $\;\!t$}}
\!
(\:\!t - s)^{\mbox{}^{\scriptstyle \!-\, \frac{7}{8} }}
\,\!
\|\, \mbox{\boldmath $u$}(\cdot,s) \,
\|_{\mbox{}_{\scriptstyle L^{2} }}
  ^{\mbox{}^{\scriptstyle \frac{1}{4} }}
\:\!
\|\, D \:\! \mbox{\boldmath $u$}(\cdot,s) \,
\|_{\mbox{}_{\scriptstyle L^{2} }}
  ^{\mbox{}^{\scriptstyle \frac{3}{4} }}
\:\!
\|\, D^{m} \,\! \mbox{\boldmath $u$}(\cdot,s) \,
\|_{\mbox{}_{\scriptstyle L^{2} }}
\, ds,
\end{equation}
\mbox{} \vspace{-0.550cm} \\
\begin{equation}
\tag{4.17$b$}
J_{\mbox{}_{2}}(t)
\;=\;
K\!\;\!(\nu) \!\;\!
\sum_{\ell \,=\,1}^{m\,-\,{\scriptscriptstyle 2}}
t^{\mbox{}^{\scriptstyle \frac{\scriptstyle m}{2} }}
\!\!\!
\int_{\mbox{\footnotesize $\mu(t)$}}
    ^{\mbox{\footnotesize $\;\!t$}}
\!
(\:\!t - s)^{\mbox{}^{\scriptstyle \!-\, \frac{7}{8} }}
\,\!
\|\, D^{\ell} \,\! \mbox{\boldmath $u$}(\cdot,s) \,
\|_{\mbox{}_{\scriptstyle L^{2} }}
  ^{\mbox{}^{\scriptstyle \frac{1}{4} }}
\:\!
\|\, D^{\ell + 1} \,\! \mbox{\boldmath $u$}(\cdot,s) \,
\|_{\mbox{}_{\scriptstyle L^{2} }}
  ^{\mbox{}^{\scriptstyle \frac{3}{4} }}
\:\!
\|\, D^{m - \ell} \,\! \mbox{\boldmath $u$}(\cdot,s) \,
\|_{\mbox{}_{\scriptstyle L^{2} }}
\, ds,
\end{equation}
\mbox{} \vspace{-0.950cm} \\
\begin{equation}
\tag{4.17$c$}
J_{\mbox{}_{3}}(t)
\;=\;
K\!\;\!(\nu) \;
t^{\mbox{}^{\scriptstyle \frac{\scriptstyle m}{2} }}
\!\!\!
\int_{\mbox{\footnotesize $\mu(t)$}}
    ^{\mbox{\footnotesize $\;\!t$}}
\!
(\:\!t - s)^{\mbox{}^{\scriptstyle \!-\, \frac{7}{8} }}
\,\!
\|\, D^{m - 1} \,\! \mbox{\boldmath $u$}(\cdot,s) \,
\|_{\mbox{}_{\scriptstyle L^{2} }}
  ^{\mbox{}^{\scriptstyle \frac{1}{4} }}
\:\!
\|\, D \:\! \mbox{\boldmath $u$}(\cdot,s) \,
\|_{\mbox{}_{\scriptstyle L^{2} }}
\:\!
\|\, D^{m} \,\! \mbox{\boldmath $u$}(\cdot,s) \,
\|_{\mbox{}_{\scriptstyle L^{2} }}
  ^{\mbox{}^{\scriptstyle \frac{3}{4} }}
\, ds.
\end{equation}
\mbox{} \vspace{+0.050cm} \\
Escrevendo
$ {\displaystyle
\;\!
\|\, D \:\! \mbox{\boldmath $u$}(\cdot,s) \,
\|_{\mbox{}_{\scriptstyle L^{2} }}
\!=\;\!
s^{\mbox{}^{\scriptstyle -\, \frac{1}{2} }}
U_{\mbox{}_{1}}(s)
\;\!\leq\;\!
M_{\mbox{}_{1}} \,\!
s^{\mbox{}^{\scriptstyle -\, \frac{1}{2} }}
\!\!\;\!
%
%
} $,
$ {\displaystyle
\;\!
\|\, D^{m} \,\! \mbox{\boldmath $u$}(\cdot,s) \,
\|_{\mbox{}_{\scriptstyle L^{2} }}
\!\;\!=\;\!
s^{\mbox{}^{\scriptstyle \!-\, \frac{\scriptstyle m}{2} }}
U_{\mbox{}_{\scriptstyle m}}\!\;\!(s)
} $,
obt\'em-se \\
\mbox{} \vspace{-0.500cm} \\
\begin{equation}
\tag{4.18$a$}
\begin{split}
J_{\mbox{}_{1}}(t)
\;&\leq\;
K\!\;\!(\nu) \:
M_{\mbox{}_{1}}^{\mbox{}^{\scriptstyle \frac{3}{4} }}
\;\!
\|\, \mbox{\boldmath $u$}_0 \;\!
\|_{\mbox{}_{\scriptstyle L^{2}(\mathbb{R}^{3}) }}
  ^{\mbox{}^{\scriptstyle \frac{1}{4} }}
\,
t^{\mbox{}^{\scriptstyle \frac{\scriptstyle m}{2} }}
\!\!\!
\int_{\mbox{\footnotesize $\mu(t)$}}
    ^{\mbox{\footnotesize $\;\!t$}}
\!
(\:\!t - s)^{\mbox{}^{\scriptstyle \!-\, \frac{7}{8} }}
\;\!
s^{\mbox{}^{\scriptstyle \!-\, \frac{3}{8} \,-\, \frac{\scriptstyle m}{2} }}
\;\!
U_{\mbox{}_{\scriptstyle m}}(s)
\: ds \\
&\leq\;
K\!\;\!(\nu) \:
M_{\mbox{}_{1}}^{\mbox{}^{\scriptstyle \frac{3}{4} }}
\;\!
\|\, \mbox{\boldmath $u$}_0 \;\!
\|_{\mbox{}_{\scriptstyle L^{2}(\mathbb{R}^{3}) }}
  ^{\mbox{}^{\scriptstyle \frac{1}{4} }}
\,
(\:\!t + t_0)^{\mbox{}^{\scriptstyle \!-\, \frac{3}{8} }}
\!\!\!
\int_{\mbox{\footnotesize $\mu(t)$}}
    ^{\mbox{\footnotesize $\;\!t$}}
\!
(\:\!t - s)^{\mbox{}^{\scriptstyle \!-\, \frac{7}{8} }}
\;\!
U_{\mbox{}_{\scriptstyle m}}(s)
\: ds.
\end{split}
\end{equation}
\mbox{} \vspace{-0.100cm} \\
Analogamente,
tem-se \\
\mbox{} \vspace{-0.650cm} \\
\begin{equation}
\tag{4.18$b$}
\begin{split}
J_{\mbox{}_{2}}(t)
\;&\leq\;
K\!\;\!(\nu) \:\!
\sum_{\ell\,=\,{\scriptscriptstyle 1}}^{m - {\scriptscriptstyle 2}}
\;\!
M_{\mbox{}_{\scriptstyle \ell}}^{\mbox{}^{\scriptstyle \frac{1}{4} }}
\;\!
M_{\mbox{}_{{\scriptstyle \ell} + 1}}^{\mbox{}^{\scriptstyle \frac{3}{4} }}
\;\!
M_{\mbox{}_{{\scriptstyle m} - {\scriptstyle \ell}}}
\:
t^{\mbox{}^{\scriptstyle \frac{\scriptstyle m}{2} }}
\!\!\!
\int_{\mbox{\footnotesize $\mu(t)$}}
    ^{\mbox{\footnotesize $\;\!t$}}
\!
(\:\!t - s)^{\mbox{}^{\scriptstyle \!-\, \frac{7}{8} }}
\;\!
s^{\mbox{}^{\scriptstyle \!-\, \frac{3}{8} \,-\, \frac{\scriptstyle m}{2} }}
\, ds \\
&\leq\;
K\!\;\!(\nu)
\:\!
\sum_{\ell\,=\,{\scriptscriptstyle 1}}^{m - {\scriptscriptstyle 2}}
\;\!
M_{\mbox{}_{\scriptstyle \ell}}^{\mbox{}^{\scriptstyle \frac{1}{4} }}
\;\!
M_{\mbox{}_{{\scriptstyle \ell} + 1}}^{\mbox{}^{\scriptstyle \frac{3}{4} }}
\;\!
M_{\mbox{}_{{\scriptstyle m} - {\scriptstyle \ell}}}
\,
(\:\!t + t_0)^{\mbox{}^{\scriptstyle \!-\, \frac{1}{4} }}
\end{split}
\end{equation}
\mbox{} \vspace{-0.250cm} \\
e \\
\mbox{} \vspace{-0.900cm} \\
\begin{equation}
\tag{4.18$c$}
\begin{split}
J_{\mbox{}_{3}}(t)
\;&\leq\;
K\!\;\!(\nu) \:
M_{\mbox{}_{1}}
\,
M_{\mbox{}_{{\scriptstyle m} - 1}}^{\mbox{}^{\scriptstyle \frac{1}{4} }}
\,
t^{\mbox{}^{\scriptstyle \frac{\scriptstyle m}{2} }}
\!\!\!
\int_{\mbox{\footnotesize $\mu(t)$}}
    ^{\mbox{\footnotesize $\;\!t$}}
\!
(\:\!t - s)^{\mbox{}^{\scriptstyle \!-\, \frac{7}{8} }}
\;\!
s^{\mbox{}^{\scriptstyle \!-\, \frac{3}{8} \,-\, \frac{\scriptstyle m}{2} }}
\;\!
U_{\mbox{}_{\scriptstyle m}}(s)
 ^{\mbox{}^{\scriptstyle \frac{3}{4} }}
\: ds \\
&\leq\;
K\!\;\!(\nu) \:
M_{\mbox{}_{1}}
\,
M_{\mbox{}_{{\scriptstyle m} - 1}}^{\mbox{}^{\scriptstyle \frac{1}{4} }}
\,
(\:\!t + t_0)^{\mbox{}^{\scriptstyle \!-\, \frac{3}{8} }}
\!\!\!
\int_{\mbox{\footnotesize $\mu(t)$}}
    ^{\mbox{\footnotesize $\;\!t$}}
\!
(\:\!t - s)^{\mbox{}^{\scriptstyle \!-\, \frac{7}{8} }}
\;\!
U_{\mbox{}_{\scriptstyle m}}(s)
 ^{\mbox{}^{\scriptstyle \frac{3}{4} }}
\: ds.
\end{split}
\end{equation}

\mbox{} \vspace{-0.500cm} \\
Observando novamente que
$ {\displaystyle
\;\!
U_{\mbox{}_{\scriptstyle m}}(s)^{\mbox{}^{\scriptstyle \:\!\frac{3}{4} }}
\!\:\!\leq\,\! 1 + U_{\mbox{}_{\scriptstyle m}}(s)
} $,
resulta de (4.6), (4.10), (4.17) e (4.18)
que,
para todo $ \;\! t \geq t_0 + 1 $,
tem-se
($\:\!$sendo
$ {\displaystyle
\;\!
M_{\mbox{}_{0}} \!\;\!\equiv\;\!
\|\, \mbox{\boldmath $u$}_0 \;\!
\|_{\mbox{}_{\scriptstyle L^{2}(\mathbb{R}^{3})}}
\!\;\!
} $): \\
\mbox{} \vspace{-0.625cm} \\
\begin{equation}
\notag
V_{\alpha}(t)
\,\leq\,
K\!\;\!(m,\;\!\nu, \;\!t_0, M_{\mbox{}_{0}}\!\;\!,
M_{\mbox{}_{1}}\!\;\!, ...,
M_{\mbox{}_{{\scriptstyle m} - 1}}\!\;\!)
\,+\,
K\!\;\!(\nu, M_{\mbox{}_{0}}\!\;\!, M_{\mbox{}_{1}}\!\;\!)
\,
(t + t_0)^{\mbox{}^{\scriptstyle \!\!-\, \frac{3}{8} }}
\!\!\!
\int_{\mbox{\footnotesize $\mu(t)$}}
    ^{\mbox{\footnotesize $\;\!t$}}
\!\!
(t - s)^{\mbox{}^{\scriptstyle \!-\, \frac{7}{8} }}
\,\!
U_{\mbox{}_{\scriptstyle m}}(s)
\, ds
\end{equation}
\mbox{} \vspace{-0.150cm} \\
para todo $ \;\!\alpha \;\!$
com $ \;\!|\;\!\alpha\;\!|\;\!=\;\!m $,
$\;\!$onde
$ {\displaystyle
\;\!
V_{\alpha}(t)
\;\!=\,
t^{\mbox{}^{\scriptstyle \!\!-\, \frac{\scriptstyle m}{2} }}
\,\!
\|\, D^{\alpha} \,\!\mbox{\boldmath $u$}(\cdot,t) \,
\|_{\mbox{}_{\scriptstyle L^{2}(\mathbb{R}^{3})}}
\!\;\!
} $.
Portanto,
tem-se \\
\mbox{} \vspace{-0.575cm} \\
\begin{equation}
\tag{4.19}
U_{\mbox{}_{\scriptstyle m}}(t)
\,\leq\,
K_{\mbox{}_{\scriptstyle \ast}}\!\;\!
(m,\;\!\nu, \;\!t_0, M_{\mbox{}_{0}}\!\;\!, ...,
M_{\mbox{}_{{\scriptstyle m} - 1}}\!\;\!)
\;+\;
K_{\mbox{}_{\scriptstyle \ast\ast}}\!\;\!
(\nu, M_{\mbox{}_{0}}\!\;\!, M_{\mbox{}_{1}}\!\;\!)
\,
(\,\!t + t_0)^{\mbox{}^{\scriptstyle \!\!-\, \frac{3}{8} }}
\!\!\!
\int_{\mbox{\footnotesize $\mu(t)$}}
    ^{\mbox{\footnotesize $\;\!t$}}
\!\!
(\,\!t - s)^{\mbox{}^{\scriptstyle \!-\, \frac{7}{8} }}
\;\!
U_{\mbox{}_{\scriptstyle m}}(s)
\: ds
\end{equation}
\mbox{} \vspace{-0.450cm} \\
para todo $ \;\! t \geq t_0 + 1 $,
onde
$ {\displaystyle
\;\!
K_{\mbox{}_{\scriptstyle \ast\ast}}\!\;\!
(\nu, M_{\mbox{}_{0}}\!\:\!, M_{\mbox{}_{1}}\!\;\!) > 0
\;\!
} $
depende apenas de
$ \nu $, $ M_{\mbox{}_{1}}\!\:\! $ e
$ {\displaystyle
M_{\mbox{}_{0}} \!\;\!=\,\!
\|\,\mbox{\boldmath $u$}_0 \;\! \|_{L^{2}}
} $.
Como no caso $ m = 2 $,
tomemos agora
$ \;\!t_{m} $,
$ \mathbb{M}_{{\scriptstyle m}} \!\;\!$
dados por \\
\begin{equation}
\tag{4.20}
t_{m}\;\!:=\;
1 +\;\! t_0 +\, 2^{16} \;\!
K_{\mbox{}_{\scriptstyle \!\;\!\ast\ast}}^{\:\!4}
\!\;\!,
\qquad
\mathbb{M}_{{\scriptstyle m}} \;\!\!:=\;
\sup \,\{\, U_{\mbox{}_{\scriptstyle m}}(s) \!:\;
t_0 \leq s \leq t_{m} \;\!\},
\end{equation}
\mbox{} \vspace{-0.200cm} \\
onde
$ {\displaystyle
\;\!
K_{\mbox{}_{\scriptstyle \!\;\!\ast\ast}}
\!\,\!=\;\!
K_{\mbox{}_{\scriptstyle \!\;\!\ast\ast}}\!\:\!
(\nu, M_{\mbox{}_{0}}\!\:\!, M_{\mbox{}_{1}}\!\:\!)
\;\!
} $
\'e a constante
definida em (4.19) acima.
Obt\'em-se, ent\~ao, \\
\mbox{} \vspace{-0.500cm} \\
\begin{equation}
\tag{4.21}
U_{\mbox{}_{\scriptstyle m}}\!\;\!(t)
\:\leq\:
2 \, K_{\mbox{}_{\scriptstyle \ast}}\!\;\!
(m,\;\!\nu, \;\!t_0, M_{\mbox{}_{0}}\!\;\!,
M_{\mbox{}_{1}}\!\;\!, ...,
M_{\mbox{}_{{\scriptstyle m} - 1}}\!\,\!)
\;+\;
16 \: K_{\mbox{}_{\scriptstyle \ast\ast}}\!\;\!
(\nu, M_{\mbox{}_{0}}\!\;\!, M_{\mbox{}_{1}}\!\,\!) \,
\mathbb{M}_{\mbox{}_{\scriptstyle m}},
\quad \;\;
\forall \;\,
t \;\!\geq\;\! t_{m},
\end{equation}
\mbox{} \vspace{-0.200cm} \\
onde
$ \;\! K_{\mbox{}_{\scriptstyle \!\;\!\ast}} \!\;\! $,
$ K_{\mbox{}_{\scriptstyle \!\;\!\ast\ast}} $
s\~ao as constantes dadas
em (4.19).
\mbox{[}$\,$A prova de (4.21)
\'e exatamente an\'aloga \`a de (4.15).$\;\!$\mbox{]}
Isso mostra que
$ U_{\mbox{}_{\scriptstyle m}}\!\;\! $
\'e limitada em $\,[\,t_{m}, \infty) $;
como
$ U_{\mbox{}_{\scriptstyle m}}\!\;\! $
\'e limitada em $\,[\,t_{0}, \;\!t_{m}\;\!] $
(por (1.3)),
segue ent\~ao que
$ U_{\mbox{}_{\scriptstyle m}}\!\;\! \in L^{\infty}(\:\![\,t_0, \infty)\,\!) $,
concluindo a prova do Lema 4.1.
}
\mbox{} \hfill $\Box$ \\
%
\mbox{} \vspace{-0.650cm} \\

Usando o Lema 4.1,
obt\'em-se o Teorema 4.1
sem dificuldade.
De fato,
sendo
$ m \geq 2 $,
$ t_0 \geq \mbox{\small $T$}_{\!\;\!\ast\ast} \!\;\!$
quaisquer,
pode-se proceder
como segue.
Dado $ \;\!\epsilon > 0\;\!$
(arbitr\'ario), \linebreak
seja $ \;\! t_{\epsilon} \gg t_0 $
suficientemente grande
tal que,
por (1.6$c$) e (4.7$b$), (4.9$a$) acima,
tenha-se,
para todo $ \;\! t \geq t_{\epsilon} $: \\
\mbox{} \vspace{-0.625cm} \\
\begin{equation}
\tag{4.22$a$}
t^{\mbox{}^{\scriptstyle \:\!\frac{\scriptstyle m}{2} }}
\:\!
\|\, D^{m} \:\![\:
e^{\;\!\nu \Delta (\mbox{\footnotesize $t$} \;\!-\, \mbox{\footnotesize $t$}_0)}
\:\!
\mbox{\boldmath $u$}(\cdot,t_0)
\,] \:
\|_{\mbox{}_{\scriptstyle L^{2}(\mathbb{R}^{3})}}
\:\leq\;
\mbox{\small $ {\displaystyle \frac{1}{3} }$} \:
\epsilon
\end{equation}
\mbox{} \vspace{-0.800cm} \\
e \\
\mbox{} \vspace{-1.100cm} \\
\begin{equation}
\tag{4.22$b$}
t^{\mbox{}^{\scriptstyle \frac{\scriptstyle m}{2} }}
\!\!\!
\int_{\mbox{\footnotesize $\;\!t_0$}}
    ^{\mbox{\footnotesize $\;\!\mu(t)$}}
\!\!
\|\, D^{m} \;\![\:
e^{\;\!\nu \Delta (\mbox{\footnotesize $t$} \;\!-\, \mbox{\footnotesize $s$})}
\:\!
\mbox{\boldmath $Q$}(\cdot,s)
\,] \:
\|_{\mbox{}_{\scriptstyle L^{2}(\mathbb{R}^{3})}}
\;\! ds
\;\;\! \leq\;\;\!
\mbox{\small $ {\displaystyle \frac{1}{3} }$} \,
\epsilon,
\end{equation}
\mbox{} \vspace{-0.100cm} \\
onde,
como antes,
$ \;\!\mu(t) = (\:\!t + t_0)/2 $.
Por (4.5), Lema 4.1,
e lembrando
(4.6), (4.7$c$) e (4.18),
obtemos
(aumentando $ \:\!t_{\epsilon}$ se necess\'ario) \\
\mbox{} \vspace{-0.625cm} \\
\begin{equation}
\tag{4.22$c$}
t^{\mbox{}^{\scriptstyle \frac{\scriptstyle m}{2} }}
\!\!\!
\int_{\mbox{\footnotesize $\mu(t)$}}
    ^{\mbox{\footnotesize $\;\!t$}}
\!\!\;\!
\|\, D^{m} \;\![\:
e^{\;\!\nu \Delta (\mbox{\footnotesize $t$} \;\!-\, \mbox{\footnotesize $s$})}
\:\!
\mbox{\boldmath $Q$}(\cdot,s)
\,] \:
\|_{\mbox{}_{\scriptstyle L^{2}(\mathbb{R}^{3})}}
\;\! ds
\;\;\! \leq\;\;\!
\mbox{\small $ {\displaystyle \frac{1}{3} }$} \,
\epsilon
\end{equation}
\mbox{} \vspace{-0.100cm} \\
para todo $ \;\! t \geq t_{\epsilon} $.
$\!$Por (4.2),
v\^e-se que
(4.22) implica
termos
$ {\displaystyle
\,
t^{\:\!m/2} \;\!
\|\, D\:\! \mbox{\boldmath $u$}(\cdot,t) \,
\|_{\mbox{}_{\scriptstyle L^{2}(\mathbb{R}^{3})}}
\!\,\!\leq\,\! \epsilon
} $
para todo
$ \:\! t \geq t_{\epsilon} $,
estabelecendo (4.1),
como afirmado.
\hfill $\Box$

\nl
%
%
%
\mbox{} \vspace{-1.750cm} \\

{\bf 5. Prova de (1.12{\em b\/})} \\
\setcounter{section}{5}

Nesta se\c c\~ao,
provaremos a estimativa (1.12$b$),
para todo $ m \geq 0 $.
Este fato,
combinado com a propriedade
de interpola\c c\~ao (1.13)
\mbox{[}$\,$tomando-se, por exemplo,
$ s_{\mbox{}_{1}} = m $,
$ s_{\mbox{}_{2}} = m + 1\,$\mbox{]},
estabelece o resultado mais geral
dado em (1.10$b$),
v\'alido
para todo $ s \in \mathbb{R} $ n\~ao negativo.
Como sempre,
$ \mbox{\boldmath $u$}(\cdot,t) $
denota uma solu\c c\~ao de Leray
(qualquer) para as equa\c c\~oes
de Navier-Stokes (1.1);
tal solu\c c\~ao \'e suave
em
$ \mathbb{R}^{3} \!\times \!\;\!
[\,\mbox{\small $T$}_{\!\;\!\ast\ast}\!\;\!, \infty) $,
para certo
$ \mbox{\small $T$}_{\!\;\!\ast\ast} \!\;\! \gg 1 $,
e satisfaz:
$ {\displaystyle
\;\!
\mbox{\boldmath $u$}(\cdot,t)
\in L^{\infty}_{\tt loc}
(\:\![\,\mbox{\small $T$}_{\!\;\!\ast\ast} \!\;\!, \infty),
H^{m}(\mathbb{R}^{3})\:\!)
} $,
para todo $ m $
\mbox{[}$\,$ver (1.3)$\,$\mbox{]}. \linebreak
O ponto de partida para
obter (1.12$b$)
envolve v\'arios resultados
anteriores,
particularmente (2.5), (2.14), (2.15)
e os Teoremas 3.1, 3.2 e 4.1 acima. \\
\mbox{} \vspace{-0.200cm} \\
%
%
%
%
%
%
\mbox{} \hspace{-0.800cm}
\fbox{%
\begin{minipage}[t]{16.000cm}
\mbox{} \vspace{-0.450cm} \\
\mbox{} \hspace{+0.300cm}
\begin{minipage}[t]{15.000cm}
\mbox{} \vspace{+0.100cm} \\
{\bf Teorema 5.1.}
\textit{%
Para todo
$\;\!m \geq 0 $,
e todo $ t_0 \geq 0 $,
tem-se
} \\
\mbox{} \vspace{-0.625cm} \\
\begin{equation}
\tag{5.1}
\mbox{} \;\;\,
\lim_{t\,\rightarrow\,\infty}
\;
t^{\mbox{}^{\scriptstyle
\frac{\scriptstyle m}{2} \,+\, \frac{1}{4} }}
\:\!
\|\, D^{\alpha} \,\! \mbox{\boldmath $u$}(\cdot,t)
\;\!-\;\!
D^{\alpha} \:\![\:
e^{\:\! \nu \Delta (t - t_0)} \:\!
\mbox{\boldmath $u$}(\cdot,t_0) \,] \,
\|_{\mbox{}_{\scriptstyle L^{2}(\mathbb{R}^{3})}}
\;\!=\; 0
\end{equation}
\mbox{} \vspace{-0.200cm} \\
\textit{%
para cada
multi-\'\i ndice
$\;\!\alpha $
com
$ \;\!|\;\!\alpha\;\!| \:\!=\;\! m $.\\
}
\mbox{} \vspace{-0.500cm} \\
\end{minipage}
\end{minipage}
}
%
%
%
%
\nl
\mbox{} \vspace{-0.000cm} \\
%
%
%
{\small
{\bf Prova:}
Pelo Teorema 3.1,
\'e suficiente mostrar (5.1)
supondo-se
$ t_0 \geq \mbox{\small $T$}_{\!\;\!\ast\ast} \!\;\!$.
$\!$\mbox{[}$\,$Com efeito,
tendo-se j\'a estabelecido o resultado neste caso,
ent\~ao se poderia estend\^e-lo
para
$ t_0 < \mbox{\small $T$}_{\!\;\!\ast\ast} \!\;\!$
do seguinte modo:
tomando-se
$ t_0^{\;\!\prime} \geq \mbox{\small $T$}_{\!\;\!\ast\ast} \!\;\!$,
ter\'\i amos \\
\mbox{} \vspace{+0.100cm} \\
\mbox{} \hspace{+0.500cm}
$ {\displaystyle
\limsup_{t\,\rightarrow\,\infty}
\;
t^{\mbox{}^{\scriptstyle
\frac{\scriptstyle m}{2} \,+\, \frac{1}{4} }}
\:\!
\|\, D^{\alpha} \,\! \mbox{\boldmath $u$}(\cdot,t)
\;\!-\;\!
D^{\alpha} \:\![\:
e^{\:\! \nu \Delta (t - t_0)} \:\!
\mbox{\boldmath $u$}(\cdot,t_0) \,] \,
\|_{\mbox{}_{\scriptstyle L^{2}(\mathbb{R}^{3})}}
\;\!\leq
} $ \\
\mbox{} \vspace{+0.100cm} \\
\mbox{} \hspace{+2.500cm}
$ {\displaystyle
\leq\;
\limsup_{t\,\rightarrow\,\infty}
\;
t^{\mbox{}^{\scriptstyle
\frac{\scriptstyle m}{2} \,+\, \frac{1}{4} }}
\:\!
\|\, D^{\alpha} \,\! \mbox{\boldmath $u$}(\cdot,t)
\;\!-\;\!
D^{\alpha} \:\![\:
e^{\:\! \nu \Delta (t - t_0^{\;\!\prime})} \:\!
\mbox{\boldmath $u$}(\cdot,t_0^{\;\!\prime}) \,] \,
\|_{\mbox{}_{\scriptstyle L^{2}(\mathbb{R}^{3})}}
\;\!+
} $ \\
\mbox{} \vspace{-0.100cm} \\
\mbox{} \hspace{+3.000cm}
$ {\displaystyle
+\;
\limsup_{t\,\rightarrow\,\infty}
\;
t^{\mbox{}^{\scriptstyle
\frac{\scriptstyle m}{2} \,+\, \frac{1}{4} }}
\:\!
\|\,
D^{\alpha} \:\![\:
e^{\:\! \nu \Delta (t - t_0^{\;\!\prime})} \:\!
\mbox{\boldmath $u$}(\cdot,t_0^{\;\!\prime})
\;\!-\,
e^{\:\! \nu \Delta (t - t_0)} \:\!
\mbox{\boldmath $u$}(\cdot,t_0) \,] \,
\|_{\mbox{}_{\scriptstyle L^{2}(\mathbb{R}^{3})}}
} $ \\
\mbox{} \vspace{+0.100cm} \\
\mbox{} \hspace{+2.500cm}
$ {\displaystyle
=\;
\limsup_{t\,\rightarrow\,\infty}
\;
t^{\mbox{}^{\scriptstyle
\frac{\scriptstyle m}{2} \,+\, \frac{1}{4} }}
\:\!
\|\, D^{\alpha} \,\! \mbox{\boldmath $u$}(\cdot,t)
\;\!-\;\!
D^{\alpha} \:\![\:
e^{\:\! \nu \Delta (t - t_0^{\;\!\prime})} \:\!
\mbox{\boldmath $u$}(\cdot,t_0^{\;\!\prime}) \,] \,
\|_{\mbox{}_{\scriptstyle L^{2}(\mathbb{R}^{3})}}
} $
\mbox{} \hfill
\mbox{[}$\,$por (3.1)$\,$\mbox{]} \\
\mbox{} \vspace{-0.050cm} \\
\mbox{} \hspace{+2.500cm}
$ {\displaystyle
=\; 0
} $, \\
\mbox{} \vspace{-0.200cm} \\
onde no \'ultimo passo
se teria usado o resultado (5.1),
j\'a mostrado neste caso
($ \:\!t_0^{\;\!\prime} \geq
\mbox{\small $T$}_{\!\;\!\ast\ast}$).$\;\!$\mbox{]} \\
Supondo-se, ent\~ao,
$ \:\!t_0 \geq
\mbox{\small $T$}_{\!\;\!\ast\ast} \!\;\!$,
podemos escrever $ \mbox{\boldmath $u$}(\cdot,t) $
na forma (4.2), para $ t \geq t_0 $,
ou seja, \\
\mbox{} \vspace{-0.500cm} \\
\begin{equation}
\tag{5.2}
\mbox{\boldmath $u$}(\cdot,t)
\;=\;
e^{\;\!\nu \Delta (\mbox{\footnotesize $t$} \;\!-\, \mbox{\footnotesize $t_0$})}
\:\!
\mbox{\boldmath $u$}(\cdot,t_0)
\;+
\int_{\mbox{\footnotesize $\!\;\!t_0$}}
    ^{\mbox{\footnotesize $\:\!t$}}
\!\!\;\!
e^{\;\!\nu \Delta (\mbox{\footnotesize $t$} \;\!-\, \mbox{\footnotesize $s$})}
\:\!
\mbox{\boldmath $Q$}(\cdot,s)
\: ds,
\qquad
\forall \;\,
t \geq t_0,
\end{equation}
\mbox{} \vspace{-0.100cm} \\
onde
$ {\displaystyle
\;\!
\mbox{\boldmath $Q$}(\cdot,s)
\;\!
} $
\'e dada em (2.4), Se\c c\~ao 2.
%
%
%
Considerando
inicialmente
os casos $ m = 0 $
e $ m = 1 $,
podemos proceder
do seguinte modo.
Dado $ \;\!\epsilon > 0 $,
seja
$ \:\!t_{\!\;\!\epsilon} \!\:\!> t_{0} \!\;\!$
suficientemente grande
tal que,
pelo Teorema 3.2,
se tenha \\
\mbox{} \vspace{-0.525cm} \\
\begin{equation}
\tag{5.3}
\mbox{} \hspace{+0.500cm}
t^{\:\!1/2} \,
\|\, D \mbox{\boldmath $u$}(\cdot,t) \,
\|_{\mbox{}_{\scriptstyle L^{2}(\mathbb{R}^{3})}}
\leq\: \epsilon
\qquad
\forall \;\;\!
t \geq t_{\!\;\!\epsilon}.
\end{equation}
\mbox{} \vspace{-0.200cm} \\
Por (5.2)
e (1.2), (2.5),
tem-se \\
\mbox{} \vspace{-0.125cm} \\
\mbox{} \hspace{+0.500cm}
$ {\displaystyle
t^{\:\!1/4} \,
\|\, \mbox{\boldmath $u$}(\cdot,t) \,-\:
e^{\:\!\nu \Delta (t \;\!-\, t_0)} \:\!
\mbox{\boldmath $u$}(\cdot,t_0) \,
\|_{\mbox{}_{\scriptstyle L^{2}(\mathbb{R}^{3})}}
\,\leq\;\:\!
t^{\:\!1/4}
\!\!
\int_{\mbox{\scriptsize $\!\;\!t_0$}}
    ^{\mbox{\scriptsize $t$}}
\!\!\:\!
\|\: e^{\:\!\nu \Delta (t \;\!-\, s)} \:\!
\mbox{\boldmath $Q$}(\cdot,s) \,
\|_{\mbox{}_{\scriptstyle L^{2}(\mathbb{R}^{3})}}
\:\!
ds
} $ \\
\mbox{} \vspace{+0.050cm} \\
\mbox{} \hfill
$ {\displaystyle
\leq\;
I(\:\!t, \:\!t_{\!\;\!\epsilon})
\;+\:
K\!\;\!(\nu) \: t^{\:\!1/4}
\!\!
\int_{\mbox{\scriptsize $\!\;\!t_{\epsilon}$}}
    ^{\mbox{\scriptsize $t$}}
\!
(t - s)^{-\,3/4} \,
\|\, \mbox{\boldmath $u$}(\cdot,s) \,
\|_{\mbox{}_{\scriptstyle L^{2}(\mathbb{R}^{3})}}
\:\!
\|\, D \mbox{\boldmath $u$}(\cdot,s) \,
\|_{\mbox{}_{\scriptstyle L^{2}(\mathbb{R}^{3})}}
\:\!
ds
} $ \\
\mbox{} \vspace{-0.050cm} \\
\mbox{} \hspace{+2.900cm}
$ {\displaystyle
\leq\;
I(\:\!t, \:\!t_{\!\;\!\epsilon})
\;+\:
K\!\;\!(\nu) \:
\|\, \mbox{\boldmath $u$}_0 \,
\|_{\mbox{}_{\scriptstyle L^{2}(\mathbb{R}^{3})}}
\epsilon
\:
t^{\:\!1/4}
\!\!
\int_{\mbox{\scriptsize $\!\;\!t_{\epsilon}$}}
    ^{\mbox{\scriptsize $t$}}
\!
(t - s)^{-\,3/4} \, s^{-\,1/2}
\,ds
} $
\mbox{} \hfill [$\,$por (5.3)$\,$] \\
\mbox{} \vspace{+0.050cm} \\
\mbox{} \hspace{+2.900cm}
$ {\displaystyle
\leq\;
I(\:\!t, \:\!t_{\!\;\!\epsilon})
\;+\:
K\!\;\!(\nu) \:
\|\, \mbox{\boldmath $u$}_0 \,
\|_{\mbox{}_{\scriptstyle L^{2}(\mathbb{R}^{3})}}
\:\!
\epsilon
\;
t^{\:\!1/4}
\,
(\:\! t - t_{\!\;\!\epsilon})^{-\,1/4}
} $ \\
\mbox{} \vspace{+0.025cm} \\
para todo
$ \;\! t > t_{\!\;\!\epsilon}$,
onde
$ \:\!K\!\;\!(\nu) > 0 \;\! $
\'e independente de $ \epsilon $,
e \\
\mbox{} \vspace{+0.100cm} \\
\mbox{} \hspace{+2.130cm}
$ {\displaystyle
I(\:\!t, \:\!t_{\!\;\!\epsilon})
\;\!:=\:
K\!\;\!(\nu) \; t^{\:\!1/4}
\!\!
\int_{\mbox{\scriptsize $\!\;\!t_{0}$}}
    ^{\mbox{\scriptsize $t_{\epsilon}$}}
\!
(t - s)^{-\,3/4} \,
\|\, \mbox{\boldmath $u$}(\cdot,s) \,
\|_{\mbox{}_{\scriptstyle L^{2}(\mathbb{R}^{3})}}
\:\!
\|\, D \mbox{\boldmath $u$}(\cdot,s) \,
\|_{\mbox{}_{\scriptstyle L^{2}(\mathbb{R}^{3})}}
\:\!
ds
} $ \\
\mbox{} \vspace{-0.050cm} \\
\mbox{} \hspace{+3.450cm}
$ {\displaystyle
\leq\:
K\!\;\!(\nu) \;
t^{\:\!1/4}
\,
(\:\!t - t_{\!\;\!\epsilon})^{-\,3/4}
\!\!
\int_{\mbox{\scriptsize $\!\;\!t_{0}$}}
    ^{\mbox{\scriptsize $t_{\epsilon}$}}
\!
\|\, \mbox{\boldmath $u$}(\cdot,s) \,
\|_{\mbox{}_{\scriptstyle L^{2}(\mathbb{R}^{3})}}
\:\!
\|\, D \mbox{\boldmath $u$}(\cdot,s) \,
\|_{\mbox{}_{\scriptstyle L^{2}(\mathbb{R}^{3})}}
\:\!
ds
} $ \\
\mbox{} \vspace{-0.050cm} \\
\mbox{} \hspace{+3.450cm}
$ {\displaystyle
\leq\:
K\!\;\!(\nu) \:
\|\, \mbox{\boldmath $u$}_0 \;\!
\|_{\mbox{}_{\scriptstyle L^{2}(\mathbb{R}^{3})}}^{\:\!2}
\,
( t_{\epsilon} -\;\! t_0 )^{1/2}
\:
t^{\:\!1/4}
\,
(\:\!t - t_{\!\;\!\epsilon})^{-\,3/4}
} $ \\
\mbox{} \vspace{-0.120cm} \\
Portanto,
temos \\
\mbox{} \vspace{-0.550cm} \\
\begin{equation}
\notag
t^{\:\!1/4} \,
\|\, \mbox{\boldmath $u$}(\cdot,t) \,-\:
e^{\:\! \nu \,\!\Delta (t \;\!-\, t_0)} \:\!
\mbox{\boldmath $u$}(\cdot,t_0) \,
\|_{\mbox{}_{\scriptstyle L^{2}(\mathbb{R}^{3})}}
\:\!\leq\
K\!\;\!(\nu) \;
(\;\! 1 \,+\,
\|\, \mbox{\boldmath $u$}_0 \,
\|_{\mbox{}_{\scriptstyle L^{2}(\mathbb{R}^{3})}})
\: \epsilon
\end{equation}
\mbox{} \vspace{-0.300cm} \\
para todo $ \;\! t > t_{\epsilon} $
grande,
com
$ \:\!K\!\;\!(\nu) \:\!$
independente de $ \epsilon $.
Isto mostra (5.1)
no caso $ m = 0 $. \linebreak
%
%
\mbox{} \vspace{-0.750cm} \\

Considerando, agora,
$ m = 1 $,
seja novamente
$ \;\!t_{\epsilon} > t_0 $
como em (5.3),
para $ \;\!\epsilon > 0 \;\!$ dado.
Por (5.2),
tem-se,
para cada
$ \:\!1 \leq \ell \leq 3 \;\!$
e $ \;\! t > t_{\epsilon} $: \\
\mbox{} \vspace{-0.600cm} \\
\begin{equation}
\tag{5.4}
\;\!
{\cal V}_{\mbox{}_{\scriptstyle \ell}}(t)
\:\equiv\;
t^{\mbox{}^{\scriptstyle \frac{3}{4} }}
\!\;\!
\|\, D_{\ell} \;\!\mbox{\boldmath $u$}(\cdot,t)
\;\!-\;\!
D_{\ell} \;\![\;\!\;\!
e^{\;\!\nu \Delta (\mbox{\footnotesize $t$} \;\!-\, \mbox{\footnotesize $t_0$})}
\:\!
\mbox{\boldmath $u$}(\cdot,t_0) \,]
\,
\|_{\mbox{}_{\scriptstyle L^{2}(\mathbb{R}^{3})}}
\leq\;
{\cal I}_{\mbox{}_{\scriptstyle 1}}(t) \,+\,
{\cal I}_{\mbox{}_{\scriptstyle 2}}(t) \,+\,
{\cal J}_{\mbox{}_{\scriptstyle \ell}}(t),
\end{equation}
\mbox{} \vspace{-0.200cm} \\
onde \\
\mbox{} \vspace{-1.050cm} \\
\begin{equation}
\tag{5.5$a$}
{\cal I}_{\mbox{}_{\scriptstyle 1}}(t)
\;=\;\;\!
t^{\mbox{}^{\scriptstyle \frac{3}{4} }}
\!\!\!\;\!
\int_{\mbox{\footnotesize $\:\!t_0$}}
    ^{\mbox{\footnotesize $\;\!t_{\epsilon}$}}
\!\!\;\!
\|\, D_{\ell} \,[\:
e^{\:\!\nu \Delta (\mbox{\footnotesize $t$} \;\!-\, \mbox{\footnotesize $s$})}
\:\!
\mbox{\boldmath $Q$}(\cdot,s)
\,] \:
\|_{\mbox{}_{\scriptstyle L^{2}(\mathbb{R}^{3})}}
\;\! ds,
\end{equation}
\mbox{} \vspace{-0.850cm} \\
\begin{equation}
\tag{5.5$b$}
{\cal I}_{\mbox{}_{\scriptstyle 2}}(t)
\;=\;\;\!
t^{\mbox{}^{\scriptstyle \frac{3}{4} }}
\!\!\!\;\!
\int_{\mbox{\footnotesize $\:\!t_{\epsilon}$}}
    ^{\mbox{\footnotesize $\;\!\mu_{\epsilon}(t)$}}
\!\!
\|\, D_{\ell} \,[\:
e^{\:\!\nu \Delta (\mbox{\footnotesize $t$} \;\!-\, \mbox{\footnotesize $s$})}
\:\!
\mbox{\boldmath $Q$}(\cdot,s)
\,] \:
\|_{\mbox{}_{\scriptstyle L^{2}(\mathbb{R}^{3})}}
\;\! ds,
\end{equation}
\mbox{} \vspace{-0.750cm} \\
\begin{equation}
\tag{5.5$c$}
{\cal J}_{\mbox{}_{\scriptstyle \ell}}(t)
\;=\;\;\!
t^{\mbox{}^{\scriptstyle \frac{3}{4} }}
\!\!\!\;\!
\int_{\mbox{\footnotesize $\mu_{\epsilon}(t)$}}
    ^{\mbox{\footnotesize $\;\!t$}}
\!
\|\, D_{\ell} \,[\:
e^{\:\!\nu \Delta (\mbox{\footnotesize $t$} \;\!-\, \mbox{\footnotesize $s$})}
\:\!
\mbox{\boldmath $Q$}(\cdot,s)
\,] \:
\|_{\mbox{}_{\scriptstyle L^{2}(\mathbb{R}^{3})}}
\;\! ds,
\end{equation}
\mbox{} \vspace{+0.100cm} \\
sendo
$ {\displaystyle
\;\!
\mu_{\epsilon}(t) := (\:\!t_{\epsilon} + \;\!t)/2
} $.
%
%
Come\c cando
com $ {\cal I}_{\mbox{}_{\scriptstyle 1}}(t) $,
tem-se,
por (1.2), (2.3) e (2.15): \\
\mbox{} \vspace{-0.050cm} \\
\mbox{} \hspace{+1.500cm}
$ {\displaystyle
{\cal I}_{\mbox{}_{\scriptstyle 1}}(t)
\;\leq\:
K\!\;\!(\nu) \; t^{\:\!3/4}
\!\!\:\!
\int_{\mbox{\scriptsize $\!\;\!t_{0}$}}
    ^{\mbox{\scriptsize $t_{\epsilon}$}}
\!
(t - s)^{-\,5/4} \,
\|\, \mbox{\boldmath $u$}(\cdot,s) \,
\|_{\mbox{}_{\scriptstyle L^{2}(\mathbb{R}^{3})}}
\:\!
\|\, D \mbox{\boldmath $u$}(\cdot,s) \,
\|_{\mbox{}_{\scriptstyle L^{2}(\mathbb{R}^{3})}}
\:\!
ds
} $ \\
\mbox{} \vspace{-0.050cm} \\
\mbox{} \hspace{+2.600cm}
$ {\displaystyle
\leq\:
K\!\;\!(\nu) \;
t^{\:\!3/4}
\:
(\:\!t - t_{\!\;\!\epsilon})^{-\,5/4}
\!\!\:\!
\int_{\mbox{\scriptsize $\!\;\!t_{0}$}}
    ^{\mbox{\scriptsize $t_{\epsilon}$}}
\!
\|\, \mbox{\boldmath $u$}(\cdot,s) \,
\|_{\mbox{}_{\scriptstyle L^{2}(\mathbb{R}^{3})}}
\:\!
\|\, D \mbox{\boldmath $u$}(\cdot,s) \,
\|_{\mbox{}_{\scriptstyle L^{2}(\mathbb{R}^{3})}}
\:\!
ds
} $ \\
\mbox{} \vspace{+0.050cm} \\
\mbox{} \hspace{+2.625cm}
$ {\displaystyle
\leq\:
K\!\;\!(\nu) \:
\|\, \mbox{\boldmath $u$}_0 \;\!
\|_{\mbox{}_{\scriptstyle L^{2}(\mathbb{R}^{3})}}^{\:\!2}
\,\!
( t_{\epsilon} -\;\! t_0 )^{1/2}
\;
t^{\:\!3/4}
\,
(\:\!t - t_{\!\;\!\epsilon})^{-\,5/4}
} $
\mbox{} \hfill (5.6$a$) \\
\mbox{} \vspace{+0.050cm} \\
para todo $\;\!t > t_{\epsilon} $,
onde $ \:\!K\!\;\!(\nu) $ independe de $ \:\!\epsilon $.
%
%
%
Para
$ {\cal I}_{\mbox{}_{\scriptstyle 2}}(t) $,
tem-se,
por (1.2), (2.15) e (5.3): \\
\mbox{} \vspace{-0.025cm} \\
\mbox{} \hspace{+1.500cm}
$ {\displaystyle
{\cal I}_{\mbox{}_{\scriptstyle 2}}(t)
\;\leq\:
K\!\;\!(\nu) \; t^{\:\!3/4}
\!\!\:\!
\int_{\mbox{\scriptsize $ t_{\epsilon}$}}
    ^{\mbox{\scriptsize $\:\!\mu_{\epsilon}(t)$}}
\hspace{-0.350cm}
(t - s)^{-\,5/4} \,
\|\, \mbox{\boldmath $u$}(\cdot,s) \,
\|_{\mbox{}_{\scriptstyle L^{2}(\mathbb{R}^{3})}}
\:\!
\|\, D \mbox{\boldmath $u$}(\cdot,s) \,
\|_{\mbox{}_{\scriptstyle L^{2}(\mathbb{R}^{3})}}
\:\!
ds
} $ \\
\mbox{} \vspace{-0.050cm} \\
\mbox{} \hspace{+2.600cm}
$ {\displaystyle
\leq\:
K\!\;\!(\nu) \;
\|\, \mbox{\boldmath $u$}_0 \;\!
\|_{\mbox{}_{\scriptstyle L^{2}(\mathbb{R}^{3})}}
\:\!
\epsilon \;\;\!
t^{\:\!3/4}
\,
(\:\!t - t_{\!\;\!\epsilon})^{-\,3/4}
\!\!\:\!
\int_{\mbox{\scriptsize $\!\;\!t_{\epsilon}$}}
    ^{\mbox{\scriptsize $\mu_{\epsilon}(t)$}}
\hspace{-0.350cm}
( t - s )^{-\,1/2}
\:
s^{-\,1/2}
\:
ds
} $ \\
\mbox{} \vspace{+0.025cm} \\
\mbox{} \hspace{+2.625cm}
$ {\displaystyle
\leq\:
K\!\;\!(\nu) \;
\|\, \mbox{\boldmath $u$}_0 \;\!
\|_{\mbox{}_{\scriptstyle L^{2}(\mathbb{R}^{3})}}
\:\!
\epsilon \;\;\!
t^{\:\!3/4}
\,
(\:\!t - t_{\!\;\!\epsilon})^{-\,3/4}
} $
\mbox{} \hfill (5.6$b$) \\
\mbox{} \vspace{+0.050cm} \\
para todo $\;\!t > t_{\epsilon} $,
onde $ \:\!K\!\;\!(\nu) $ independe de $ \:\!\epsilon $.
%
%
%
%
Finalmente,
considerando
$ {\cal J}_{\mbox{}_{\scriptstyle \ell}}(t) $,
tem-se \\
\mbox{} \vspace{-0.025cm} \\
\mbox{} \hspace{+1.500cm}
$ {\displaystyle
{\cal J}_{\mbox{}_{\scriptstyle \ell}}(t)
\;=\;\;\!
t^{\:\!3/4}
\!\!\:\!
\int_{\mbox{\footnotesize $\mu_{\epsilon}(t)$}}
    ^{\mbox{\footnotesize $\;\!t$}}
\hspace{-0.325cm}
\|\, D_{\ell} \,[\:
e^{\;\!\nu \Delta (\mbox{\footnotesize $t$} \;\!-\, \mbox{\footnotesize $s$})}
\:\!
\mbox{\boldmath $Q$}(\cdot,s)
\,] \:
\|_{\mbox{}_{\scriptstyle L^{2}(\mathbb{R}^{3})}}
\;\! ds
} $ \\
\mbox{} \vspace{+0.050cm} \\
\mbox{} \hspace{+2.600cm}
$ {\displaystyle
\leq\:
K\!\;\!(\nu)
\;\:\!
t^{\:\!3/4}
\!\!\:\!
\int_{\mbox{\footnotesize $\mu_{\epsilon}(t)$}}
    ^{\mbox{\footnotesize $\;\!t$}}
\hspace{-0.300cm}
(\:\!t - s)^{-\, 7/8}
\:
\|\, \mbox{\boldmath $Q$}(\cdot,s) \,
\|_{\mbox{}_{\scriptstyle L^{4/3}(\mathbb{R}^{3})}}
\;\! ds
} $ \\
\mbox{} \vspace{+0.050cm} \\
\mbox{} \hspace{+2.600cm}
$ {\displaystyle
\leq\:
K\!\;\!(\nu)
\;
t^{\:\!3/4}
\!\!\:\!
\int_{\mbox{\footnotesize $\mu_{\epsilon}(t)$}}
    ^{\mbox{\footnotesize $\;\!t$}}
\hspace{-0.300cm}
(\:\!t - s)^{-\, 7/8}
\:
\|\, \mbox{\boldmath $u$}(\cdot,s) \cdot
\nabla \mbox{\boldmath $u$}(\cdot,s) \,
\|_{\mbox{}_{\scriptstyle L^{4/3}(\mathbb{R}^{3})}}
\;\! ds
} $ \\
\mbox{} \vspace{+0.050cm} \\
\mbox{} \hspace{+2.600cm}
$ {\displaystyle
\leq\:
K\!\;\!(\nu)
\;
t^{\:\!3/4}
\!\!\:\!
\int_{\mbox{\footnotesize $\mu_{\epsilon}(t)$}}
    ^{\mbox{\footnotesize $\;\!t$}}
\hspace{-0.300cm}
(\:\!t - s)^{-\, 7/8}
\:
\|\, \mbox{\boldmath $u$}(\cdot,s) \,
\|_{\mbox{}_{\scriptstyle L^{4}(\mathbb{R}^{3})}}
\;\!
\|\, D \:\!\mbox{\boldmath $u$}(\cdot,s) \,
\|_{\mbox{}_{\scriptstyle L^{2}(\mathbb{R}^{3})}}
\;\! ds
} $ \\
\mbox{} \vspace{+0.050cm} \\
\mbox{} \hspace{+2.600cm}
$ {\displaystyle
\leq\:
K\!\;\!(\nu)
\;
t^{\:\!3/4}
\!\!\:\!
\int_{\mbox{\footnotesize $\mu_{\epsilon}(t)$}}
    ^{\mbox{\footnotesize $\;\!t$}}
\hspace{-0.300cm}
(\:\!t - s)^{-\,7/8}
\:
\|\, \mbox{\boldmath $u$}(\cdot,s) \,
\|_{\mbox{}_{\scriptstyle L^{2}(\mathbb{R}^{3})}}
  ^{\:\!1/4}
\;\!
\|\, D \:\!\mbox{\boldmath $u$}(\cdot,s) \,
\|_{\mbox{}_{\scriptstyle L^{2}(\mathbb{R}^{3})}}
  ^{\:\!7/4}
\;\! ds
} $ \\
\mbox{} \vspace{+0.050cm} \\
\mbox{} \hspace{+2.600cm}
$ {\displaystyle
\leq\:
K\!\;\!(\nu)
\:
\|\, \mbox{\boldmath $u$}_0 \;\!
\|_{\mbox{}_{\scriptstyle L^{2}(\mathbb{R}^{3})}}
  ^{\:\!1/4}
\:\!
\epsilon^{\:\! 7/4}
\!\!\;\!
\int_{\mbox{\footnotesize $\mu_{\epsilon}(t)$}}
    ^{\mbox{\footnotesize $\;\!t$}}
\hspace{-0.300cm}
(\:\!t - s)^{-\,7/8}
\:
s^{-\,7/8}
\:
ds
} $ \\
\mbox{} \vspace{+0.050cm} \\
\mbox{} \hspace{+2.600cm}
$ {\displaystyle
\leq\:
K\!\;\!(\nu)
\:
\|\, \mbox{\boldmath $u$}_0 \;\!
\|_{\mbox{}_{\scriptstyle L^{2}(\mathbb{R}^{3})}}
  ^{\:\!1/4}
\:\!
\epsilon^{\:\! 7/4}
\;
t^{\:\! 3/4}
\,
(\:\! t + t_{\epsilon})^{-\,3/4}
\!\!\;\!
\int_{\mbox{\footnotesize $\mu_{\epsilon}(t)$}}
    ^{\mbox{\footnotesize $\;\!t$}}
\hspace{-0.300cm}
(\:\!t - s)^{-\,7/8}
\:
s^{-\,1/8}
\:
ds
} $ \\
\mbox{} \vspace{+0.050cm} \\
\mbox{} \hspace{+2.600cm}
$ {\displaystyle
\leq\:
K\!\;\!(\nu)
\;
\|\, \mbox{\boldmath $u$}_0 \;\!
\|_{\mbox{}_{\scriptstyle L^{2}(\mathbb{R}^{3})}}
  ^{\:\!1/4}
\;\!
\epsilon^{\:\! 7/4}
} $
\hfill (5.6$c$) \\
\mbox{} \vspace{+0.050cm} \\
para todo
$ \:\! t > t_{\epsilon} $,
usando (2.14), (4.16),
desigualdade de H\"older, (4.12)
e o fato de se ter \\
%
\mbox{} \vspace{-0.700cm} \\
\begin{equation}
\notag
\int_{t_0}^{\:\!t}
(t - s)^{-\,7/8}
\,
s^{-\,1/8}
\:
ds
\;\;\!\leq\;
K \,=\;
\frac{\pi}{\,\mbox{sen}\,(\pi/8)\,}
\;\!.
\end{equation}

\mbox{} \vspace{-0.750cm} \\
Claramente,
obtemos
de (5.4) e
(5.6$a$), (5.6$b$), (5.6$c$)
acima
que \\
\mbox{} \vspace{-0.550cm} \\
\begin{equation}
\notag
t^{\:\!3/4}
\,
\|\, D_{\ell} \,\mbox{\boldmath $u$}(\cdot,t)
\;\!-\;\!
D_{\ell} \,[\;\!\;\!
e^{\;\!\nu \Delta (\mbox{\footnotesize $t$} \;\!-\, \mbox{\footnotesize $t_0$})}
\:\!
\mbox{\boldmath $u$}(\cdot,t_0) \,]
\,
\|_{\mbox{}_{\scriptstyle L^{2}(\mathbb{R}^{3})}}
\leq\:
\bigl(\;\! 1 + \epsilon +
K\!\;\!(\nu)
\;
\|\, \mbox{\boldmath $u$}_0 \;\!
\|_{\mbox{}_{\scriptstyle L^{2}(\mathbb{R}^{3})}}
\;\!
\bigr)
\,
\epsilon
\end{equation}
\mbox{} \vspace{-0.100cm} \\
para todo $ t > t_{\epsilon} $
suficientemente grande,
para cada $ \;\! 1 \leq \ell \leq 3 $,
o que prova (5.1) se $ m = 1 $. \\
%
%
%
%
\mbox{} \vspace{-0.750cm} \\

Finalmente,
consideramos o caso geral $ m \geq 2 $,
procedendo de modo similar
ao caso anterior.
Dado $ \epsilon > 0 $ (arbitr\'ario),
tomamos $ t_{\epsilon} > t_0 $
suficientemente grande
tal que,
por (4.1) \linebreak
\mbox{[}$\,$Teorema 4.1$\,$\mbox{]},
tenhamos \\
\mbox{} \vspace{-0.600cm} \\
\begin{equation}
\tag{5.7}
t^{\mbox{}^{\scriptstyle
\frac{\scriptstyle k}{2} }}
\:\!
\|\, D^{k} \:\! \mbox{\boldmath $u$}(\cdot,t) \,
\|_{\mbox{}_{\scriptstyle L^{2}(\mathbb{R}^{3})}}
\,\leq\;
\epsilon,
\qquad
\forall \;\,
t \geq t_{\epsilon}
\end{equation}
\mbox{} \vspace{-0.220cm} \\
para todo $ \;\!0 \leq k \leq m $.
Dado $\;\!\alpha\;\!$
multi-\'\i ndice
com
$ \;\!|\;\!\alpha\;\!| = m $,
definimos
(para $ t > t_{\epsilon} $) \\
\mbox{} \vspace{-0.500cm} \\
\begin{equation}
\tag{5.8}
{\cal V}_{\alpha}(t)
\;\equiv\;\;\!
t^{\mbox{}^{\scriptstyle \frac{\scriptstyle m}{2} \,+\, \frac{1}{4} }}
\:\!
\|\, D^{\alpha} \:\!\mbox{\boldmath $u$}(\cdot,t)
\;\!-\;\!
D^{\alpha} \;\![\:
e^{\:\!\nu \Delta (\mbox{\footnotesize $t$} \;\!-\, \mbox{\footnotesize $t_0$})}
\:\!
\mbox{\boldmath $u$}(\cdot,t_0)
\,] \,
\|_{\mbox{}_{\scriptstyle L^{2}(\mathbb{R}^{3})}},
\end{equation}
\mbox{} \vspace{-0.200cm} \\
escrevendo
$ {\displaystyle
\;\!
{\cal V}_{\alpha}(t)
\;\!\leq\,
{\cal I}_{\mbox{}_{\scriptstyle 1}}(\alpha, \:\!t) \,+\,
{\cal I}_{\mbox{}_{\scriptstyle 2}}(\alpha, \:\!t) \,+\,
{\cal J}(\alpha, t)
} $,
onde \\
\mbox{} \vspace{-0.500cm} \\
\begin{equation}
\tag{5.9$a$}
{\cal I}_{\mbox{}_{\scriptstyle 1}}(\alpha, \:\!t)
\;=\;\;\!
t^{\mbox{}^{\scriptstyle \frac{\scriptstyle m}{2} \,+\, \frac{1}{4} }}
\!\!\!\;\!
\int_{\mbox{\footnotesize $\;\!t_0$}}
    ^{\mbox{\footnotesize $\;\!t_{\epsilon}$}}
\!\!
\|\, D^{\alpha} \:\![\:
e^{\;\!\nu \Delta (\mbox{\footnotesize $t$} \;\!-\, \mbox{\footnotesize $s$})}
\:\!
\mbox{\boldmath $Q$}(\cdot,s)
\,] \:
\|_{\mbox{}_{\scriptstyle L^{2}(\mathbb{R}^{3})}}
\;\! ds,
\end{equation}
\mbox{} \vspace{-0.650cm} \\
\begin{equation}
\tag{5.9$b$}
{\cal I}_{\mbox{}_{\scriptstyle 2}}(\alpha, \:\!t)
\;=\;\;\!
t^{\mbox{}^{\scriptstyle \frac{\scriptstyle m}{2} \,+\, \frac{1}{4} }}
\!\!\!\;\!
\int_{\mbox{\footnotesize $\;\!t_0$}}
    ^{\mbox{\footnotesize $\;\!\mu_{\epsilon}(t)$}}
\hspace{-0.570cm}
\|\, D^{\alpha} \:\![\:
e^{\;\!\nu \Delta (\mbox{\footnotesize $t$} \;\!-\, \mbox{\footnotesize $s$})}
\:\!
\mbox{\boldmath $Q$}(\cdot,s)
\,] \:
\|_{\mbox{}_{\scriptstyle L^{2}(\mathbb{R}^{3})}}
\;\! ds,
\end{equation}
\mbox{} \vspace{-0.650cm} \\
\begin{equation}
\tag{5.9$c$}
{\cal J}(\alpha, t)
\;=\;\;\!
t^{\mbox{}^{\scriptstyle \frac{\scriptstyle m}{2} \,+\, \frac{1}{4} }}
\!\!\!\;\!
\int_{\mbox{\footnotesize $\mu_{\epsilon}(t)$}}
    ^{\mbox{\footnotesize $\;\!t$}}
\hspace{-0.350cm}
\|\, D^{\alpha} \:\![\:
e^{\;\!\nu \Delta (\mbox{\footnotesize $t$} \;\!-\, \mbox{\footnotesize $s$})}
\:\!
\mbox{\boldmath $Q$}(\cdot,s)
\,] \:
\|_{\mbox{}_{\scriptstyle L^{2}(\mathbb{R}^{3})}}
\;\! ds,
\end{equation}
\mbox{} \vspace{+0.050cm} \\
sendo
$ {\displaystyle
\;\!
\mu_{\epsilon}(t) = (\:\!t_{\epsilon} + \;\!t)/2
} $.
%
%
%
%
Com rela\c c\~ao a
$ {\cal I}_{\mbox{}_{\scriptstyle 1}}(\alpha, t) $,
tem-se,
por (1.2), (2.3) e (2.15): \\
\mbox{} \vspace{-0.000cm} \\
\mbox{} \hspace{+1.000cm}
$ {\displaystyle
{\cal I}_{\mbox{}_{\scriptstyle 1}}(\alpha, t)
\;\leq\:
K\!\;\!(\nu) \;
t^{\mbox{}^{\scriptstyle \frac{\scriptstyle m}{2} \,+\, \frac{1}{4} }}
\!\!\!\;\!
\int_{\mbox{\scriptsize $\!\;\!t_{0}$}}
    ^{\mbox{\scriptsize $t_{\epsilon}$}}
\!
(t - s)^{\mbox{}^{\scriptstyle \!\! -\, \frac{\scriptstyle m}{2} \,-\, \frac{3}{4} }}
\,
\|\, \mbox{\boldmath $u$}(\cdot,s) \,
\|_{\mbox{}_{\scriptstyle L^{2}(\mathbb{R}^{3})}}
\:\!
\|\, D \mbox{\boldmath $u$}(\cdot,s) \,
\|_{\mbox{}_{\scriptstyle L^{2}(\mathbb{R}^{3})}}
\:\!
ds
} $ \\
\mbox{} \vspace{-0.050cm} \\
\mbox{} \hspace{+2.475cm}
$ {\displaystyle
\leq\:
K\!\;\!(\nu) \;\;\!
t^{\mbox{}^{\scriptstyle \frac{\scriptstyle m}{2} \,+\, \frac{1}{4} }}
\,\!
(t - t_{\epsilon})^{\mbox{}^{\scriptstyle \!\! -\,
\frac{\scriptstyle m}{2} \,-\, \frac{3}{4} }}
\!\!\!\;\!
\int_{\mbox{\scriptsize $\!\;\!t_{0}$}}
    ^{\mbox{\scriptsize $t_{\epsilon}$}}
\!
\|\, \mbox{\boldmath $u$}(\cdot,s) \,
\|_{\mbox{}_{\scriptstyle L^{2}(\mathbb{R}^{3})}}
\:\!
\|\, D \mbox{\boldmath $u$}(\cdot,s) \,
\|_{\mbox{}_{\scriptstyle L^{2}(\mathbb{R}^{3})}}
\:\!
ds
} $ \\
\mbox{} \vspace{+0.050cm} \\
\mbox{} \hspace{+2.475cm}
$ {\displaystyle
\leq\:
K\!\;\!(\nu) \:
\|\, \mbox{\boldmath $u$}_0 \;\!
\|_{\mbox{}_{\scriptstyle L^{2}(\mathbb{R}^{3})}}^{\:\!2}
\:\!
( t_{\epsilon} -\;\! t_0 )^{\mbox{}^{\scriptstyle \frac{1}{2} }}
\,
t^{\mbox{}^{\scriptstyle \frac{\scriptstyle m}{2} \,+\, \frac{1}{4} }}
\,\!
(t - t_{\epsilon})^{\mbox{}^{\scriptstyle \!\! -\,
\frac{\scriptstyle m}{2} \,-\, \frac{3}{4} }}
} $
\mbox{} \hfill (5.10$a$) \\
\mbox{} \vspace{+0.050cm} \\
para todo $\;\!t > t_{\epsilon} $,
onde $ \:\!K\!\;\!(\nu) $ independe de $ \:\!\epsilon $.
%
%
%
Para
$ {\cal I}_{\mbox{}_{\scriptstyle 2}}(\alpha, t) $,
tem-se,
de (1.2), (2.15) e (5.7): \\
\mbox{} \vspace{-0.025cm} \\
\mbox{} \hspace{+1.000cm}
$ {\displaystyle
{\cal I}_{\mbox{}_{\scriptstyle 2}}(\alpha, t)
\;\leq\:
K\!\;\!(\nu) \;\;\!
t^{\mbox{}^{\scriptstyle \frac{\scriptstyle m}{2} \,+\, \frac{1}{4} }}
\!\!\!\;\!
\int_{\mbox{\scriptsize $ t_{\epsilon}$}}
    ^{\mbox{\scriptsize $\:\!\mu_{\epsilon}(t)$}}
\hspace{-0.425cm}
(t - s)^{\mbox{}^{\scriptstyle \!\! -\,\frac{\scriptstyle m}{2} \,-\, \frac{3}{4} }}
\,
\|\, \mbox{\boldmath $u$}(\cdot,s) \,
\|_{\mbox{}_{\scriptstyle L^{2}(\mathbb{R}^{3})}}
\:\!
\|\, D \mbox{\boldmath $u$}(\cdot,s) \,
\|_{\mbox{}_{\scriptstyle L^{2}(\mathbb{R}^{3})}}
\:\!
ds
} $ \\
\mbox{} \vspace{-0.050cm} \\
\mbox{} \hspace{+2.475cm}
$ {\displaystyle
\leq\:
K\!\;\!(\nu) \;
\:\!
\epsilon^{\:\!2}
\;
t^{\mbox{}^{\scriptstyle \frac{\scriptstyle m}{2} \,+\, \frac{1}{4} }}
\;\!
(\:\!t - t_{\!\;\!\epsilon})^{\mbox{}^{\scriptstyle \!\!- \,
\frac{\scriptstyle m}{2} \,-\, \frac{1}{4} }}
\!\!\:\!
\int_{\mbox{\scriptsize $\!\;\!t_{\epsilon}$}}
    ^{\mbox{\scriptsize $\mu_{\epsilon}(t)$}}
\hspace{-0.350cm}
( t - s )^{-\,1/2}
\:
s^{-\,1/2}
\:
ds
} $ \\
\mbox{} \vspace{+0.025cm} \\
\mbox{} \hspace{+2.475cm}
$ {\displaystyle
\leq\:
K\!\;\!(\nu) \;
\:\!
\epsilon^{\:\!2}
\;
t^{\mbox{}^{\scriptstyle \frac{\scriptstyle m}{2} \,+\, \frac{1}{4} }}
\;\!
(\:\!t - t_{\!\;\!\epsilon})^{\mbox{}^{\scriptstyle \!\!- \,
\frac{\scriptstyle m}{2} \,-\, \frac{1}{4} }}
} $
\mbox{} \hfill (5.10$b$) \\
\mbox{} \vspace{+0.050cm} \\
para todo $\;\!t > t_{\epsilon} $,
onde $ \:\!K\!\;\!(\nu) $ independe de $ \:\!\epsilon $.
%
%
%
%
Finalmente,
para
$ {\displaystyle
\;\!
{\cal J}(\alpha,t)
} $,
tem-se,
escrevendo
$ D^{\alpha} \!\;\!= D_{j} \;\! D^{\alpha^{\prime}} \!\!$,
sendo
$ \alpha^{\prime} $ multi-\'\i ndice de ordem $ m - 1 $, \\
\mbox{} \vspace{-0.025cm} \\
\mbox{} \hspace{+0.085cm}
$ {\displaystyle
{\cal J}(\alpha, t)
\;=\;\:\!
t^{\mbox{}^{\scriptstyle \frac{\scriptstyle m}{2} \,+\, \frac{1}{4} }}
\!\!\!\;\!
\int_{\mbox{\footnotesize $\mu_{\epsilon}(t)$}}
    ^{\mbox{\footnotesize $\;\!t$}}
\hspace{-0.325cm}
\|\, D^{\alpha} \;\![\:
e^{\;\!\nu \Delta (\mbox{\footnotesize $t$} \;\!-\, \mbox{\footnotesize $s$})}
\:\!
\mbox{\boldmath $Q$}(\cdot,s)
\,] \:
\|_{\mbox{}_{\scriptstyle L^{2}(\mathbb{R}^{3})}}
\;\! ds
} $ \\
\mbox{} \vspace{+0.050cm} \\
\mbox{} \hspace{+1.500cm}
$ {\displaystyle
=\;\:\!
t^{\mbox{}^{\scriptstyle \frac{\scriptstyle m}{2} \,+\, \frac{1}{4} }}
\!\!\!\;\!
\int_{\mbox{\footnotesize $\mu_{\epsilon}(t)$}}
    ^{\mbox{\footnotesize $\;\!t$}}
\hspace{-0.325cm}
\|\, D_{j} \,[\:
e^{\;\!\nu \Delta (\mbox{\footnotesize $t$} \;\!-\, \mbox{\footnotesize $s$})}
\:\!
D^{\alpha^{\prime}} \!\:\!\mbox{\boldmath $Q$}(\cdot,s)
\,] \:
\|_{\mbox{}_{\scriptstyle L^{2}(\mathbb{R}^{3})}}
\;\! ds
} $ \\
\mbox{} \vspace{+0.050cm} \\
\mbox{} \hspace{+1.500cm}
$ {\displaystyle
\leq\:
K\!\;\!(\nu)
\;\:\!
t^{\mbox{}^{\scriptstyle \frac{\scriptstyle m}{2} \,+\, \frac{1}{4} }}
\!\!\!\;\!
\int_{\mbox{\footnotesize $\mu_{\epsilon}(t)$}}
    ^{\mbox{\footnotesize $\;\!t$}}
\hspace{-0.300cm}
(\:\!t - s)^{-\, 7/8}
\:
\|\, D^{\alpha^{\prime}} \!\:\! \mbox{\boldmath $Q$}(\cdot,s) \,
\|_{\mbox{}_{\scriptstyle L^{4/3}(\mathbb{R}^{3})}}
\;\! ds
} $ \\
\mbox{} \vspace{+0.050cm} \\
\mbox{} \hspace{+1.500cm}
$ {\displaystyle
\leq\:
K\!\;\!(\nu)
\;\;\!
t^{\mbox{}^{\scriptstyle \frac{\scriptstyle m}{2} \,+\, \frac{1}{4} }}
\!\!\!\;\!
\int_{\mbox{\footnotesize $\mu_{\epsilon}(t)$}}
    ^{\mbox{\footnotesize $\;\!t$}}
\hspace{-0.300cm}
(\:\!t - s)^{-\, 7/8}
\:
\|\, D^{\alpha^{\prime}} \!\;\![\,
\mbox{\boldmath $u$}(\cdot,s) \cdot
\nabla \mbox{\boldmath $u$}(\cdot,s) \,] \,
\|_{\mbox{}_{\scriptstyle L^{4/3}(\mathbb{R}^{3})}}
\;\! ds
} $ \\
\mbox{} \vspace{+0.050cm} \\
\mbox{} \hspace{+1.500cm}
$ {\displaystyle
\leq
\hspace{-1.350cm}
\sum_{\mbox{} \hspace{+1.250cm} |\,\beta\,| \,+\, |\,\gamma\,| \:=\:
m \,-\, {\scriptscriptstyle 1}}
\hspace{-1.400cm}
K\!\;\!(\nu) \;\;\!
t^{\mbox{}^{\scriptstyle \frac{\scriptstyle m}{2} \,+\, \frac{1}{4} }}
\!\!\!\;\!
\int_{\mbox{\footnotesize $\mu_{\epsilon}(t)$}}
    ^{\mbox{\footnotesize $\;\!t$}}
\hspace{-0.300cm}
(\:\!t - s)^{-\, 7/8}
\:
\|\, D^{\beta} \,\! \mbox{\boldmath $u$}(\cdot,s) \cdot
\nabla \:\!  D^{\gamma} \,\! \mbox{\boldmath $u$}(\cdot,s) \,]\,
\|_{\mbox{}_{\scriptstyle L^{4/3}(\mathbb{R}^{3})}}
\;\! ds
} $ \\
\mbox{} \vspace{+0.150cm} \\
\mbox{} \hspace{+1.500cm}
$ {\displaystyle
\leq\;
K\!\;\!(\nu) \!
\sum_{\ell \,=\,0}^{m\,-\,{\scriptscriptstyle 1}}
t^{\mbox{}^{\scriptstyle \frac{\scriptstyle m}{2} \,+\, \frac{1}{4}}}
\!\!\!\;\!
\int_{\mbox{\footnotesize $\mu_{\epsilon}(t)$}}
    ^{\mbox{\footnotesize $\;\!t$}}
\hspace{-0.300cm}
(\:\!t - s)^{-\, 7/8}
\:
\|\, D^{\ell} \,\! \mbox{\boldmath $u$}(\cdot,s) \,
\|_{\mbox{}_{\scriptstyle L^{4}(\mathbb{R}^{3})}}
\:\!
\|\, D^{m - \ell} \,\! \mbox{\boldmath $u$}(\cdot,s) \,
\|_{\mbox{}_{\scriptstyle L^{2}(\mathbb{R}^{3})}}
} $ \\
\mbox{} \vspace{+0.150cm} \\
\mbox{} \hspace{+1.500cm}
$ {\displaystyle
\leq\;
K\!\;\!(\nu) \!
\sum_{\ell \,=\,0}^{m\,-\,{\scriptscriptstyle 1}}
t^{\mbox{}^{\scriptstyle \frac{\scriptstyle m}{2} \,+\, \frac{1}{4} }}
\!\!\!
\int_{\mbox{\footnotesize $\mu_{\epsilon}(t)$}}
    ^{\mbox{\footnotesize $\;\!t$}}
\hspace{-0.300cm}
(\:\!t - s)^{- \,7/8}
\:
\|\, D^{\ell} \,\! \mbox{\boldmath $u$}(\cdot,s) \,
\|_{\mbox{}_{\scriptstyle L^{2} }}
  ^{\mbox{}^{\scriptstyle \frac{1}{4} }}
\:\!
\|\, D^{\ell + 1} \,\! \mbox{\boldmath $u$}(\cdot,s) \,
\|_{\mbox{}_{\scriptstyle L^{2} }}
  ^{\mbox{}^{\scriptstyle \frac{3}{4} }}
\:\!
\|\, D^{m - \ell} \,\! \mbox{\boldmath $u$}(\cdot,s) \,
\|_{\mbox{}_{\scriptstyle L^{2} }}
} $ \\
\mbox{} \vspace{+0.050cm} \\
\mbox{} \hspace{+1.500cm}
$ {\displaystyle
\leq\;
K\!\;\!(\nu) \!
\sum_{\ell \,=\,0}^{m\,-\,{\scriptscriptstyle 1}}
\!\;\!
\epsilon^{\:\!2}
\;\;\!
t^{\mbox{}^{\scriptstyle \frac{\scriptstyle m}{2} \,+\, \frac{1}{4} }}
\!\!\!\;\!
\int_{\mbox{\footnotesize $\mu_{\epsilon}(t)$}}
    ^{\mbox{\footnotesize $\;\!t$}}
\hspace{-0.300cm}
(\:\!t - s)^{-\,7/8}
\:
s^{-\,m/2 \,-\, 9 \,\ell/8 \,-\, 3/8}
\:
ds
} $ \\
\mbox{} \vspace{+0.050cm} \\
\mbox{} \hspace{+1.500cm}
$ {\displaystyle
\leq\:
K\!\;\!(m, \nu)
\;
\epsilon^{\:\! 2}
\;
t^{\mbox{}^{\scriptstyle \frac{\scriptstyle m}{2} \,+\, \frac{1}{4} }}
\,
(\:\! t + t_{\epsilon})^{\mbox{}^{\scriptstyle
\!\! -\, \frac{\scriptstyle m}{2} \,-\, \frac{1}{4} }}
\!\!\;\!
\int_{\mbox{\footnotesize $\mu_{\epsilon}(t)$}}
    ^{\mbox{\footnotesize $\;\!t$}}
\hspace{-0.300cm}
(\:\!t - s)^{-\,7/8}
\:
s^{-\,1/8}
\:
ds
} $ \\
\mbox{} \vspace{+0.050cm} \\
\mbox{} \hspace{+1.500cm}
$ {\displaystyle
\leq\:
K\!\;\!(m, \nu)
\;
\epsilon^{\:\! 2}
} $
\hfill (5.10$c$) \\
\mbox{} \vspace{+0.050cm} \\
para todo
$ \:\! t > t_{\epsilon} $,
usando-se (2.14), (4.16),
desigualdade de H\"older, (4.12)
e (5.7),
tendo-se que a constante
$ K\!\;\!(m, \nu) > 0 \;\!$
em (5.10$c$)
independe de $ \;\!\epsilon\:\! $ (e $\;\!t_{\epsilon}$).
\mbox{[}$\,$Foi suposto tamb\'em,
no pen\'ultimo passo acima,
que
$ \;\!t_{\epsilon} $ tenha sido tomado
em (5.7) de modo a satisfazer:
$ \:\!t_{\epsilon} \geq 1 $.$\;\!$\mbox{]}
\mbox{} \vspace{-0.525cm} \\
Por
(5.8) e
(5.10$a$), (5.10$b$), (5.10$c$),
segue
em particular
que
se tem \\
\mbox{} \vspace{-0.475cm} \\
\begin{equation}
\notag
t^{\mbox{}^{\scriptstyle \frac{\scriptstyle m}{2} \,+\, \frac{1}{4} }}
\:\!
\|\, D^{\alpha} \:\!\mbox{\boldmath $u$}(\cdot,t)
\;\!-\;\!
D^{\alpha} \;\![\:
e^{\:\!\nu \Delta (\mbox{\footnotesize $t$} \;\!-\, \mbox{\footnotesize $t_0$})}
\:\!
\mbox{\boldmath $u$}(\cdot,t_0)
\,] \,
\|_{\mbox{}_{\scriptstyle L^{2}(\mathbb{R}^{3})}}
\,\leq\;
(\;\! 1 + K\!\;\!(m,\nu) \,\epsilon \,)
\; \epsilon
\end{equation}
\mbox{} \vspace{-0.200cm} \\
para todo $ t > t_{\epsilon} $
suficientemente grande,
para cada $ \alpha  $
com $ |\;\!\alpha\;\!| = m $,
sendo $ \epsilon > 0 $
arbitr\'ario, \linebreak
e
$ m \geq 2 $
qualquer
(com
$ K\!\;\!(m,\nu) $
independente de $\:\!\epsilon $).
Somado aos casos
$m = 0$ e $ m = 1 $ \linebreak
considerados antes,
isso conclui a prova
de (5.1)
para todo $ m \geq 0 $,
como afirmado.
}
\mbox{} \hfill $\Box$ \\
%
%
%

Deixamos em aberto a obten\c c\~ao
de (1.10), (1.12)
para $ n \geq 4 $, $ m \geq 0 $
arbitr\'arios.

\newpage
%
%
\mbox{} \vspace{-1.500cm} \\
\begin{center}
{\large \sc Ap\^endice A} \\
\end{center}
\mbox{} \vspace{-0.700cm} \\

Neste ap\^endice,
vamos mostrar como obter
a estimativa (2.20)
para o valor
$ t_{\ast\ast} $
dado na Proposi\c c\~ao 2.3
da Se\c c\~ao 2.
O ponto de partida
\'e a seguinte estimativa,  \\
\mbox{} \vspace{-0.650cm} \\
\begin{equation}
\tag{A.1}
\int_{\mbox{}_{\scriptstyle \mathbb{R}^{3}}}
\!\!\;\!\Bigl\{\!\!\!\!\:\!
\sum_{\mbox{}\;\;\,i, \,j, \,\ell \,=\,1}^{3}
\!\!\!\!
|\, D_{\ell} \;\!u_{i} \;\!| \:
|\, D_{\ell} \;\!u_{j} \;\!| \:
|\, D_{j} \;\!u_{i} \;\!|
\,\Bigr\}
\;\!dx
\:\leq\,
K_{\mbox{}_{\!\:\!3}}^{3} \;\!
\|\, D \:\!\mbox{\boldmath $u$} \,
\|_{\mbox{}_{\scriptstyle L^{2}(\mathbb{R}^{3})}}
  ^{3/2}
\;\!
\|\, D^{2} \mbox{\boldmath $u$} \,
\|_{\mbox{}_{\scriptstyle L^{2}(\mathbb{R}^{3})}}
  ^{3/2}
\!\:\!,
\end{equation}
\mbox{} \vspace{-0.000cm} \\
onde
$ \:\!K_{\mbox{}_{\!\:\!3}} \!< 0.581\,862\,001\,307 $
(ver \cite{Agueh2008}, Theorem 2.1)
\'e a constante na
desigualdade de
Niren\-berg-Gagliardo
\cite{Agueh2008} \\
\mbox{} \vspace{-0.800cm} \\
\begin{equation}
\tag{A.2}
\|\: \mbox{u} \:
\|_{\mbox{}_{\scriptstyle L^{3}(\mathbb{R}^{3})}}
\!\leq
K_{\mbox{}_{\!\:\!3}}
\;\!
\|\: \mbox{u} \:
\|_{\mbox{}_{\scriptstyle L^{2}(\mathbb{R}^{3})}}^{1/2}
\|\, D \:\!\mbox{u} \,
\|_{\mbox{}_{\scriptstyle L^{2}(\mathbb{R}^{3})}}^{1/2}
\!\:\!.
\end{equation}
\mbox{} \vspace{-0.050cm} \\
%
%
%
\mbox{[}$\,$Prova de (A.1):
aplicando-se repetidamente
a desigualdade de Cauchy-Schwarz,
tem-se \\
\mbox{} \vspace{-0.150cm} \\
\mbox{} \hfill
$ {\displaystyle
\sum_{\mbox{}\;\;\,i, \,j, \,\ell \,=\,1}^{3}
\!\!\!\!\!\;\!
|\, D_{\ell} \;\!u_{i} \;\!| \:
|\, D_{\ell} \;\!u_{j} \;\!| \:
|\, D_{j} \;\!u_{i} \;\!|
\,\Bigr\}
\:\leq\!
\sum_{\mbox{}\;\, i, \,\ell \,=\,1}^{3}
\!\!\:\!
|\, D_{\ell} \;\!u_{i} \;\!|
\:
\Bigl\{\;\!
\sum_{j\,=\,1}^{3}
|\, D_{\ell} \;\!u_{j} \;\!|^{2}
\;\!\Bigr\}^{\!\!\:\!1/2}
\!\;\!
\Bigl\{\,
\sum_{j\,=\,1}^{3}
|\, D_{j} \:\!u_{i} \;\!|^{2}
\;\!\Bigr\}^{\!\!\:\!1/2}
} $ \\
\mbox{} \vspace{+0.050cm} \\
\mbox{} \hfill
$ {\displaystyle
\leq\:
\sum_{i\,=\,1}^{3}
\;
\Bigl\{\,
\sum_{j\,=\,1}^{3}
|\, D_{j} \:\!u_{i} \;\!|^{2}
\;\!\Bigr\}^{\!\!\:\!1/2}
\,\!
\Bigl\{\,
\sum_{\ell\,=\,1}^{3}
|\, D_{\ell} \;\!u_{i} \;\!|^{2}
\;\!\Bigr\}^{\!\!\:\!1/2}
\,\!
\Bigl\{\!
\sum_{\mbox{}\;\,j,\,\ell\,=\,1}^{3}
\!\!\!\;\!
|\, D_{\ell} \;\!u_{j} \;\!|^{2}
\;\!\Bigr\}^{\!\!\:\!1/2}
} $ \\
\mbox{} \vspace{+0.050cm} \\
\mbox{} \hfill
$ {\displaystyle
\leq\;
\Bigl\{\!
\sum_{\mbox{}\;j,\,\ell\,=\,1}^{3}
\!\!\!\;\!
|\, D_{\ell} \;\!u_{j} \;\!|^{2}
\;\!\Bigr\}^{\!\!\:\!1/2}
\:\!
\Bigl\{\!
\sum_{\mbox{}\;i,\,j\,=\,1}^{3}
\!\!\!\;\!
|\, D_{j} \:\!u_{i} \;\!|^{2}
\;\!\Bigr\}^{\!\!\:\!1/2}
\:\!
\Bigl\{\!
\sum_{\mbox{}\;i,\,\ell\,=\,1}^{3}
\!\!\!\;\!
|\, D_{\ell} \;\!u_{i} \;\!|^{2}
\;\!\Bigr\}^{\!\!\:\!1/2}
\!\!\!\!
} $, \mbox{}$\;\!$ \\
\mbox{} \vspace{-0.250cm} \\
de modo que \\
\mbox{} \vspace{-0.800cm} \\
\begin{equation}
\notag
\int_{\mbox{}_{\scriptstyle \mathbb{R}^{3}}}
\!\!\;\!\Bigl\{\!\!\!\!\:\!
\sum_{\mbox{}\;\;\,i, \,j, \,\ell \,=\,1}^{3}
\!\!\!\!
|\, D_{\ell} \;\!u_{i} \;\!| \:
|\, D_{\ell} \;\!u_{j} \;\!| \:
|\, D_{j} \;\!u_{i} \;\!|
\,\Bigr\}
\;\!dx
\:\leq\;
\|\: \mbox{w} \:
\|_{\mbox{}_{\scriptstyle L^{3}(\mathbb{R}^{3})}}^{3}
\!,
\quad \;\,
\mbox{w}(x) :=
\Bigl\{\!\!\;\!
\sum_{\mbox{} \;i, \,j, \,=\,1}^{3}
\!\!
|\, D_{j} \;\!u_{i} \;\!|^{\:\!2}
\,
\Bigr\}^{\!\!\;\!1/2}
\!\!\!\!.
\end{equation}
\mbox{} \vspace{+0.020cm} \\
Isso fornece \\
\mbox{} \vspace{-0.250cm} \\
\mbox{} \hspace{+0.500cm}
$ {\displaystyle
\notag
\int_{\mbox{}_{\scriptstyle \mathbb{R}^{3}}}
\!\!\;\!\Bigl\{\!\!\!\!\:\!
\sum_{\mbox{}\;\;\,i, \,j, \,\ell \,=\,1}^{3}
\!\!\!\!
|\, D_{\ell} \;\!u_{i} \;\!| \:
|\, D_{\ell} \;\!u_{j} \;\!| \:
|\, D_{j} \;\!u_{i} \;\!|
\,\Bigr\}
\;\!dx
\:\leq\:
K_{\mbox{}_{\!\:\!3}}^{3}
\;\!
\|\: \mbox{w} \:
\|_{\mbox{}_{\scriptstyle L^{2}(\mathbb{R}^{3})}}^{3/2}
\|\, D \:\!\mbox{w} \,
\|_{\mbox{}_{\scriptstyle L^{2}(\mathbb{R}^{3})}}^{3/2}
} $
\mbox{} \hfill
\mbox{[}$\;\!$por (A.2)$\;\!$\mbox{]} \\
\mbox{} \vspace{-0.150cm} \\
\mbox{} \hspace{+7.150cm}
$ {\displaystyle
\leq\;\!
K_{\mbox{}_{\!\:\!3}}^{3}
\;\!
\|\, D \:\!\mbox{\boldmath $u$} \,
\|_{\mbox{}_{\scriptstyle L^{2}(\mathbb{R}^{3})}}^{3/2}
\|\, D^{2} \mbox{\boldmath $u$} \,
\|_{\mbox{}_{\scriptstyle L^{2}(\mathbb{R}^{3})}}^{3/2}
} $ \\
\mbox{} \vspace{-0.010cm} \\
pois,
por (1.19),
tem-se
$ {\displaystyle
\;\!
\|\,\mbox{w}\,\|_{\mbox{}_{\scriptstyle L^{2}(\mathbb{R}^{3})}}
\!\!\;\!=
\|\,D\:\!\mbox{\boldmath $u$}\,\|_{\mbox{}_{\scriptstyle L^{2}(\mathbb{R}^{3})}}
} $
e
$ {\displaystyle
\,
\|\,D\:\!\mbox{w}\,\|_{\mbox{}_{\scriptstyle L^{2}(\mathbb{R}^{3})}}
\!\!\;\!\leq\;\!
\|\,D^{2}\mbox{\boldmath $u$}\,\|_{\mbox{}_{\scriptstyle L^{2}(\mathbb{R}^{3})}}
\!\;\!
} $.$\;\!$\mbox{]} \linebreak
%
%
%
%
\mbox{} \vspace{-0.600cm} \\

Agora,
considere
$ \;\!\hat{t} > 0 \;\!$
qualquer
satisfazendo \\
\mbox{} \vspace{-0.600cm} \\
\begin{equation}
\tag{A.3}
\hat{t} \;>\:
\frac{\;\!1\;\!}{2} \,
K_{\mbox{}_{\!\:\!3}}^{12}
\,
\nu^{-\,5}
\;\!
\|\, \mbox{\boldmath $u$}_{0} \;\!
\|_{\mbox{}_{\scriptstyle L^{2}(\mathbb{R}^{3})}}^{4}
\!\!:
\end{equation}
\mbox{} \vspace{-0.250cm} \\
Como
(por (1.2))
$ {\displaystyle
\!
\int_{0}^{\:\!\mbox{\footnotesize $\hat{t}$}}
\!\!\:\!
\|\, D \:\!\mbox{\boldmath $u$}(\cdot,t) \,
\|_{\mbox{}_{\scriptstyle L^{2}(\mathbb{R}^{3})}}^{2}
dt
\;\!\leq\;\!
\mbox{\small $ {\displaystyle \frac{1}{2\;\!\nu} }$}
\:
\|\, \mbox{\boldmath $u$}_0 \;\!
\|_{\mbox{}_{\scriptstyle L^{2}(\mathbb{R}^{3})}}^{2}
\!
} $,
$\;\!$tem de existir
$ \;\!t^{\prime} \!\in (\:\!0,\;\!\hat{t}\,] $ \linebreak
\mbox{} \vspace{-0.550cm} \\
tal que \\
\mbox{} \vspace{-0.900cm} \\
\begin{equation}
\tag{A.4}
\|\, D \:\!\mbox{\boldmath $u$}(\cdot,t^{\prime}) \,
\|_{\mbox{}_{\scriptstyle L^{2}(\mathbb{R}^{3})}}
\,\leq\;
\frac{1}{\sqrt{\:\!2\:\!\nu\,}\,}
\:
\|\, \mbox{\boldmath $u$}_{0} \;\!
\|_{\mbox{}_{\scriptstyle L^{2}(\mathbb{R}^{3})}}
\!\cdot\;\!
\frac{\mbox{} \;1\,}{\!\sqrt{\:\!\hat{t}^{\mbox{}}\:}\;}
.
\end{equation}
\mbox{} \vspace{-0.200cm} \\
Portanto,
por (A.3),
tem-se \\
\mbox{} \vspace{-0.650cm} \\
\begin{equation}
\tag{A.5}
K_{\mbox{}_{\!\:\!3}}^{3}
\,
\|\, \mbox{\boldmath $u$}(\cdot,s) \,
\|_{\mbox{}_{\scriptstyle L^{2}(\mathbb{R}^{3})}}^{1/2}
\:\!
\|\, D\:\!\mbox{\boldmath $u$}(\cdot,s) \,
\|_{\mbox{}_{\scriptstyle L^{2}(\mathbb{R}^{3})}}^{1/2}
<\, \nu
\end{equation}
\mbox{} \vspace{-0.220cm} \\
para todo $ \;\!s \geq t^{\prime} $
pr\'oximo
do ponto $ t^{\prime} \!\:\!$.
Isso fornece,
diferenciando (1.1$a$)
com respeito \`a vari\`avel $ x_{\ell} $,
tomando o produto escalar com
$ D_{\ell}\;\!\mbox{\boldmath $u$}(\cdot,t) $
e somando para $ \ell = 1, 2, 3$, \\
\mbox{} \vspace{-0.100cm} \\
\mbox{} \hspace{+3.500cm}
$ {\displaystyle
\|\, D \mbox{\boldmath $u$}(\cdot,t) \,
\|_{\mbox{}_{\scriptstyle L^{2}(\mathbb{R}^{3})}}^{\:\!2}
+\:
2 \,\nu \!\!\;\!
\int_{\mbox{\footnotesize $ t^{\prime} $}}
    ^{\mbox{\footnotesize $\:\!t$}}
\!
\|\, D^{2} \mbox{\boldmath $u$}(\cdot,s) \,
\|_{\mbox{}_{\scriptstyle L^{2}(\mathbb{R}^{3})}}^{\:\!2}
ds
} $ \\
\mbox{} \vspace{-0.625cm} \\
\mbox{} \hfill (A.6) \\
\mbox{} \vspace{-0.425cm} \\
\mbox{} \hfill
$ {\displaystyle
\leq\;\:\!
\|\, D \mbox{\boldmath $u$}(\cdot,t^{\prime}) \,
\|_{\mbox{}_{\scriptstyle L^{2}(\mathbb{R}^{3})}}^{\:\!2}
\!\;\!+\:
2
\sum_{i, \, j, \, \ell}
\int_{\mbox{\footnotesize $ t^{\prime} $}}
    ^{\mbox{\footnotesize $\:\!t$}}
\!
\int_{\mathbb{R}^{3}}
\!\!\!\;\!
|\, D_{\ell} \:\! u_{i}(x,s) \,|
\;
|\, D_{\ell} \:\! u_{j}(x,s) \,|
\;
|\, D_{\scriptstyle \!j} \:\!u_{i}(x,s) \,|
\;
dx \: ds
} $ \\
\mbox{} \vspace{+0.050cm} \\
\mbox{} \hspace{+0.130cm}
$ {\displaystyle
\leq\;
\|\, D \mbox{\boldmath $u$}(\cdot,t^{\prime}) \,
\|_{\mbox{}_{\scriptstyle L^{2}(\mathbb{R}^{3})}}^{\:\!2}
+\:
2 \!\!\;\!
\int_{\mbox{\footnotesize $ t^{\prime} $}}
    ^{\mbox{\footnotesize $\:\!t$}}
\!
K_{\mbox{}_{\!\:\!3}}^{3}
\:\!
\|\, D \mbox{\boldmath $u$}(\cdot,s) \,
\|_{\mbox{}_{\scriptstyle L^{2}(\mathbb{R}^{3})}}^{\:\!3/2}
\|\, D^{2} \mbox{\boldmath $u$}(\cdot,s) \,
\|_{\mbox{}_{\scriptstyle L^{2}(\mathbb{R}^{3})}}^{\:\!3/2}
\,\!
ds
} $
\mbox{} \hfill
\mbox{[}$\,$por (A.1)$\,$\mbox{]} \\
\mbox{} \vspace{+0.025cm} \\
\mbox{} \hfill
$ {\displaystyle
\leq\;
\|\, D \mbox{\boldmath $u$}(\cdot,t^{\prime}) \,
\|_{\mbox{}_{\scriptstyle L^{2}(\mathbb{R}^{3})}}^{\:\!2}
\!+\,
2 \!\!\;\!
\int_{\mbox{\footnotesize $ t^{\prime} $}}
    ^{\mbox{\footnotesize $\:\!t$}}
\!
\bigl[\,
K_{\mbox{}_{\!\:\!3}}^{3}
\,
\|\, \mbox{\boldmath $u$}(\cdot,s) \,
\|_{\mbox{}_{\scriptstyle L^{2}(\mathbb{R}^{3})}}^{\:\!1/2}
\|\, D \:\!\mbox{\boldmath $u$}(\cdot,s) \,
\|_{\mbox{}_{\scriptstyle L^{2}(\mathbb{R}^{3})}}^{\:\!1/2}
\:\!
\bigr]
\:
\|\, D^{2} \mbox{\boldmath $u$}(\cdot,s) \,
\|_{\mbox{}_{\scriptstyle L^{2}}}^{2}
ds
} $ \\
\mbox{} \vspace{-0.025cm} \\
\mbox{} \hspace{+2.950cm}
$ {\displaystyle
\leq\:
\|\, D \mbox{\boldmath $u$}(\cdot,t^{\prime}) \,
\|_{\mbox{}_{\scriptstyle L^{2}(\mathbb{R}^{3})}}^{\:\!2}
\!+\,
2 \, \nu \!\!\;\!
\int_{\mbox{\footnotesize $ t^{\prime} $}}
    ^{\mbox{\footnotesize $\:\!t$}}
\!
\|\, D^{2} \mbox{\boldmath $u$}(\cdot,s) \,
\|_{\mbox{}_{\scriptstyle L^{2}(\mathbb{R}^{3})}}^{2}
ds
} $
\mbox{} \hfill
\mbox{[}$\,$por (A.5)$\,$\mbox{]} \\
\mbox{} \vspace{+0.100cm} \\
para todo $ \;\!t \geq t^{\prime} \!\;\!$
pr\'oximo a $ t^{\prime} \!\:\!$,
onde
na quarta linha acima
usamos a desigualdade (2.16$b$)
da Se\c c\~ao 2.
Segue da\'{\i} que
$ {\displaystyle
\;\!
\|\, D \:\!\mbox{\boldmath $u$}(\cdot,t) \,
\|_{\scriptstyle L^{2}(\mathbb{R}^{3})}
\!
} $
\'e limitada por
$ {\displaystyle
\|\, D \:\!\mbox{\boldmath $u$}(\cdot,t^{\prime}) \,
\|_{\scriptstyle L^{2}(\mathbb{R}^{3})}
\!\;\!
} $,
e como
$ {\displaystyle
\|\, \mbox{\boldmath $u$}(\cdot,t) \,
\|_{\scriptstyle L^{2}(\mathbb{R}^{3})}
\!\;\!
} $
n\~ao pode crescer
(por (1.2)),
resulta que \\
\mbox{} \vspace{-0.600cm} \\
\begin{equation}
\tag{A.7}
K_{\mbox{}_{\!\:\!3}}^{3}
\,
\|\, \mbox{\boldmath $u$}(\cdot,t) \,
\|_{\mbox{}_{\scriptstyle L^{2}(\mathbb{R}^{3})}}^{1/2}
\:\!
\|\, D\:\!\mbox{\boldmath $u$}(\cdot,t) \,
\|_{\mbox{}_{\scriptstyle L^{2}(\mathbb{R}^{3})}}^{1/2}
\;\!<\: \nu,
\qquad
\forall \;\,
t \geq t^{\prime}
\!\:\!.
\end{equation}
\mbox{} \vspace{-0.200cm} \\
Em particular,
podemos repetir
a deriva\c c\~ao de
(A.6) acima
no intervalo
$ [\;\!t_0, t\;\!] $,
para
$ t_0 < t \in [\,t^{\prime}\!\:\!, \infty) $
arbitr\'ario,
obtendo \\
\mbox{} \vspace{-0.200cm} \\
\mbox{} \hspace{+3.500cm}
$ {\displaystyle
\|\, D \mbox{\boldmath $u$}(\cdot,t) \,
\|_{\mbox{}_{\scriptstyle L^{2}(\mathbb{R}^{3})}}^{\:\!2}
+\:
2 \, \nu \!\!\;\!
\int_{\mbox{\footnotesize $ t_0 $}}
    ^{\mbox{\footnotesize $\:\!t$}}
\!
\|\, D^{2} \mbox{\boldmath $u$}(\cdot,s) \,
\|_{\mbox{}_{\scriptstyle L^{2}(\mathbb{R}^{3})}}^{\:\!2}
ds
} $ \\
\mbox{} \vspace{-0.075cm} \\
\mbox{} \hfill
$ {\displaystyle
\leq\,
\|\, D \mbox{\boldmath $u$}(\cdot,t^{\prime}) \,
\|_{\mbox{}_{\scriptstyle L^{2}(\mathbb{R}^{3})}}^{\:\!2}
\!+\,
2 \!\!\;\!
\int_{\mbox{\footnotesize $ t_0 $}}
    ^{\mbox{\footnotesize $\:\!t$}}
\!\!\;\!
\bigl[\,
K_{\mbox{}_{\!\:\!3}}^{3}
\,
\|\, \mbox{\boldmath $u$}(\cdot,s) \,
\|_{\mbox{}_{\scriptstyle L^{2}(\mathbb{R}^{3})}}^{\:\!1/2}
\|\, D \:\!\mbox{\boldmath $u$}(\cdot,s) \,
\|_{\mbox{}_{\scriptstyle L^{2}(\mathbb{R}^{3})}}^{\:\!1/2}
\:\!
\bigr]
\:
\|\, D^{2} \mbox{\boldmath $u$}(\cdot,s) \,
\|_{\mbox{}_{\scriptstyle L^{2}}}^{2}
ds
} $ \\
\mbox{} \vspace{-0.100cm} \\
\mbox{} \hspace{+2.950cm}
$ {\displaystyle
\leq\:
\|\, D \mbox{\boldmath $u$}(\cdot,t_0) \,
\|_{\mbox{}_{\scriptstyle L^{2}(\mathbb{R}^{3})}}^{\:\!2}
\!+\,
2 \, \nu \!\!\;\!
\int_{\mbox{\footnotesize $ t_0 $}}
    ^{\mbox{\footnotesize $\:\!t$}}
\!\!\;\!
\|\, D^{2} \mbox{\boldmath $u$}(\cdot,s) \,
\|_{\mbox{}_{\scriptstyle L^{2}(\mathbb{R}^{3})}}^{2}
ds
} $.
\mbox{} \hfill
\mbox{[}$\,$por (A.7)$\,$\mbox{]} \\
\mbox{} \vspace{+0.100cm} \\
Portanto,
$ {\displaystyle
\|\, D\:\! \mbox{\boldmath $u$}(\cdot,t) \,
\|_{\mbox{}_{\scriptstyle L^{2}(\mathbb{R}^{3})}}
\!
} $
\'e monotonicamente decrescente
em
$ {\displaystyle
[\,t^{\prime}\!, \:\!\infty\:\!)
\supseteq
[\,\hat{t}, \:\!\infty\:\!)
} $,
de modo que,
pela teoria cl\'assica de Leray
\cite{Leray1934},
tem-se de ter
$ {\displaystyle
\;\!
\mbox{\boldmath $u$} \in
C^{\infty}(\mathbb{R}^{3} \!\times \!\;\!
(\;\!t^{\:\!\prime} \!\;\!, \:\! \infty))
} $.
Lembrando a desigualdade (A.3)
definindo $ \;\!t^{\:\!\prime} \!\;\!$,
isto completa a prova de (2.20),
visto que
tem-se
$ {\displaystyle
\;\!
1/2
\:
K_{\mbox{}_{\!\:\!3}}^{12}
\!\;\!<\:\!
0.000\,753\,026
} $.
\mbox{} \hfill $\Box$ \\
\mbox{} \vspace{-0.600cm} \\

Em resumo, pelo argumento acima,
tem-se demonstrado o seguinte resultado: \\
\mbox{} \vspace{-0.350cm} \\
%
%
%
%
%
%
\mbox{} \hspace{-0.800cm}
\fbox{%
\begin{minipage}[t]{16.000cm}
\mbox{} \vspace{-0.450cm} \\
\mbox{} \hspace{+0.300cm}
\begin{minipage}[t]{15.000cm}
\mbox{} \vspace{+0.100cm} \\
{\bf Teorema A.1.}
\textit{%
Dado um estado inicial
$\;\!\mbox{\boldmath $u$}_0 \!\in L^{2}_{\sigma}(\mathbb{R}^{3}) $
qualquer,
seja
$ \;\!\mbox{\boldmath $u$}(\cdot,t) $
uma so-\linebreak
lu\c c\~ao de Leray
para
as equa\c c\~oes de Navier-Stokes
$\:\!(1.1)$.
Ent\~ao,
existe
$\,0 \leq t_{\ast\ast} \!\:\!< $ \linebreak
$ 0.000\,753\,026 \:  \nu^{-\,5} \:\!
\|\, \mbox{\boldmath $u$}_0 \;\! \|_{L^{2}(\mathbb{R}^{3})}^{\:\!4} \!$
com
$ \;\!\mbox{\boldmath $u$} \in C^{\infty}(\:\!\mathbb{R}^{3} \!\!\;\!\times\!\:\!
[\,t_{\ast\ast} \!\;\!, \infty)) $
e
tal que
$ {\displaystyle
\;\!
\|\, D\:\!\mbox{\boldmath $u$}(\cdot,t) \,
\|_{\mbox{}_{\scriptstyle L^{2}(\mathbb{R}^{3})}}
\!
} $
\'e finita
e monotonicamente decrescente
em todo o intervalo
$ [\,t_{\ast\ast}\!\;\!, \infty) $.
} \\
\end{minipage}
\end{minipage}
}
\nl

\newpage

%
%

%
%
%
%

%
%
\mbox{} \vspace{-2.250cm} \\
%
%
%
%
%
\begin{center}
{\Large \sc Parte II} \\
\nl
\mbox{} \vspace{-0.000cm} \\
{\Large \bf Problema de Exist\^encia Global para} \\
\mbox{} \vspace{-0.150cm} \\
{\Large \bf Equa\c c\~oes de Advec\c c\~ao-Difus\~ao
Conservativas} \\
%
%
%
%
\mbox{} \vspace{-0.350cm} \\
\end{center}
\mbox{} \vspace{-0.200cm} \\
%
%
%
%
%
%

{\bf 1. Introdu\c c\~ao} \\

Na segunda parte deste trabalho,
estenderemos um procedimento de an\'alise
introduzido pelo autor
para a deriva\c c\~ao
de v\'arias estimativas b\'asicas
importantes
para as solu\c c\~oes
$ u(\cdot,t) $
de equa\c c\~oes
de advec\c c\~ao-difus\~ao
conservativas
em meios hetero\-g\^eneos
\cite{Zingano2010}.
O m\'etodo foi
inicialmente
aplicado a equa\c c\~oes
(ou sistemas de equa\c c\~oes)
em uma dimens\~ao espacial ($ n = 1 $),
no caso mais simples
de velocidades advectivas
limitadas
(i.e., $ \kappa = 0 $ em (1.1) a seguir),
ver \cite{BarrionuevoOliveiraZingano2014, %
BrazMeloZingano2015, Melo2011, Oliveira2013}.
Posteriormente,
o autor estendeu os resultados para
equa\c c\~oes mais gerais \\
\mbox{} \vspace{-0.525cm} \\
\begin{equation}
\notag
u_t \,+\;
\bigl(\;\! b(x,t,u) \, |\;\!u\;\!|^{\:\!\kappa} \:\! u \;\!
\bigr)_{\!x}
\;\!=\;
u_{xx},
\qquad
x \in \mathbb{R}, \;\, t > 0,
\end{equation}
\mbox{} \vspace{-0.175cm} \\
com $ \kappa > 0 $ constante,
$ b(x,t,u) $ limitada
\cite{Zingano2011},
tendo orientado trabalhos de doutorado \linebreak
na aplica\c c\~ao do m\'etodo
a equa\c c\~oes unidimensionais similares
no caso de difu\-s\~ao n\~ao linear
\cite{Chagas2015, Diehl2015, Fabris2013, Guidolin2015}.
No presente trabalho,
consideraremos finalmente
o desen\-volvimento
destas t\'ecnicas
%
%
em dimens\~ao $n$ arbitr\'aria,
adotando
(por simplicidade)
como prot\'otipo
o pro\-blema \\
\mbox{} \vspace{-0.550cm} \\
\begin{equation}
\tag{1.1$a$}
u_t \,+\;
\mbox{div}\,
\bigl(\;\! \mbox{\boldmath $b$}(x,t,u) \, |\;\!u\;\!|^{\:\!\kappa} \:\! u \;\!
\bigr)
\:+\:
\mbox{div} \,
\mbox{\boldmath $f$}(t,u)
\;=\;
\mbox{div}\,\bigl(\;\!
A(x,t,u) \;\! \nabla u \;\!\bigr),
%
%
\end{equation}
\mbox{} \vspace{-0.900cm} \\
\begin{equation}
\tag{1.1$b$}
u(\cdot,0) \,=\,
u_0 \in L^{1}(\mathbb{R}^{n})
\cap L^{\infty}(\mathbb{R}^{n}),
\end{equation}
\mbox{} \vspace{-0.175cm} \\
sendo
$ A(x,t,\mbox{u}) $
matriz suave
satisfazendo
$ A(x,t,\mbox{u}) \geq \mu(t) \:\!I $
para
$ \mu \in C^{0}([\,0, \infty)) $
positiva,
ou seja, \\
\mbox{} \vspace{-0.800cm} \\
\begin{equation}
\tag{1.2}
\left<\;\! A(x,t,\mbox{u}) \;\! \mbox{\bf v},
\, \mbox{\bf v} \;\!\right>
\;\geq\;
\mu(t) \, |\:\mbox{\bf v} \,|_{\mbox{}_{2}}^{\:\!2}
\quad \;\;
\forall \;\,
\mbox{\bf v} \in \mathbb{R}^{n}
\end{equation}
\mbox{} \vspace{-0.150cm} \\
para todo
$ x \in \mathbb{R}^{n} \!\;\!$,
$ t \geq 0 $,
$ \mbox{u} \in \mathbb{R} $,
e onde
$ {\displaystyle
\;\!
\mbox{\boldmath $b$} = (\;\! b_{\mbox{}_{1}} \!\;\!,
..., b_{\mbox{}_{\scriptstyle n}} )
} $,
$ {\displaystyle
\;\!
\mbox{\boldmath $f$} = (\:\! f_{\mbox{}_{1}} \!\;\!,
..., f_{\mbox{}_{\scriptstyle n}} )
} $
s\~ao fun\c c\~oes suaves,
com
$ {\displaystyle
\;\!
\mbox{\boldmath $b$}
\;\!
} $
satisfazendo \\
\mbox{} \vspace{-0.550cm} \\
\begin{equation}
\tag{1.3}
\mbox{} \hspace{+3.500cm}
\mbox{\boldmath $b$} \in
L^{\infty}(\mathbb{R}^{n} \!\times\!\;\!
[\,0,\mbox{\small $T$}\:\!] \!\;\!\times\:\! \mathbb{R})
\qquad
\mbox{\mbox{[}$\,$para cada $\;\!\mbox{\small $T$} > 0 \,]$.}
\end{equation}

\mbox{} \vspace{-1.250cm} \\

\'E conhecido
(ver e.g.$\;$\cite{LadyzhenskayaSolonnikovUralceva1968, Serre1999}
e Se\c c\~ao~2 abaixo)
que o problema (1.1)-(1.3)
possui solu\c c\~ao
(cl\'assica, limitada, \'unica)
$ {\displaystyle
\;\!
u(\cdot,t) \in
C^{0}([\,0, \mbox{\small $T$}\:\!],
L^{1}(\mathbb{R}^{n}))
\cap
L^{\infty}([\,0, \mbox{\small $T$}\:\!],
L^{\infty}(\mathbb{R}^{n}))
} $
para certo $ 0 < \mbox{\small $T$} < \infty $
(ou seja, exist\^encia {\em local\/}
est\'a bem estabelecida);
esta solu\c c\~ao pode
ser continuada (i.e., estendida)
a intervalos de exist\^encia
mais amplos enquanto \linebreak
permanecer
limitada.
Assim,
\'e importante
examinar o comportamento
das normas \linebreak
altas
(especialmente
$ {\displaystyle
\;\!
\|\, u(\cdot,t) \,
\|_{\scriptstyle L^{\infty}(\mathbb{R}^{n})}
} $)
no intervalo
de exist\^encia
da solu\c c\~ao.
Por\'em,
sob hip\'oteses t\~ao gerais
como (1.3) acima,
esta quest\~ao pode tornar-se
{\em muito\/} dif\'\i cil,
como explicamos
intuitivamente a seguir.
Considere-se,
por exemplo,
solu\c c\~oes $ v(\cdot,t) $
n\~ao negativas
da equa\c c\~ao \\
\mbox{} \vspace{-0.700cm} \\
\begin{equation}
\tag{1.4}
v_t \,+\; \mbox{div}\,(\;\!
\mbox{\boldmath $b$}(x) \, v^{\kappa \;\!+\;\! 1}
\:\!)
\,=\,
\Delta \;\!v,
\end{equation}
\mbox{} \vspace{-0.200cm} \\
que reescrevemos na forma \\
\mbox{} \vspace{-0.600cm} \\
\begin{equation}
\tag{$1.4^{\prime}$}
v_t \,+\: (\kappa + 1) \,
v^{\kappa} \,
\mbox{\boldmath $b$}(x) \cdot \nabla v
\;=\;
\Delta \;\!v
\,+\, \beta(x) \, v^{\kappa \;\!+\;\! 1}
\end{equation}
\mbox{} \vspace{-0.200cm} \\
onde
$ \beta(x) \!\:\!:= -\,\mbox{div}\:\mbox{\boldmath $b$}(x) $.
Supondo que
$ {\displaystyle
\Omega \equiv
\{\, x \in \mathbb{R}^{n} \!\!\;\!:\;\!\beta(x) > 0 \,\}
} $
seja n\~ao vazio,
v\^e-se
de (1.4$^\prime$)
que $ v(x,t) $
\'e estimulada a crescer
nos pontos $ x \in \Omega $,
particularmente
onde ocorrer $ \;\!\beta(x) \gg 1 $.
Como (1.4) conserva massa,
se $ v(\cdot,t) $ crescer
pronunciadamente
em alguma parte de $ \Omega $
ent\~ao o perfil de $ v(\cdot,t) $
ter\'a de afinar-se,
tornando-se assim mais suscet\'\i vel
aos efeitos dissipativos
do termo difusivo presente em (1.4).
O efeito final sobre a solu\c c\~ao
(i.e., ocorr\^encia de blow-up ou n\~ao,
supondo $ \kappa > 0 $)
resultante desta competi\c c\~ao
entre os termos do lado direito em (1.$4^{\prime}$)
\'e dif\'\i cil de ser previsto. \linebreak
A situa\c c\~ao pode
\`a primeira vista
parecer (equivocadamente)
similar \`a
das so\-lu\-\c c\~oes
n\~ao negativas da equa\c c\~ao \\
\mbox{} \vspace{-0.525cm} \\
\begin{equation}
\tag{1.5}
w_t \,=\; \Delta \:\!w
\,+\, w^{\kappa \;\!+\;\!1}
\!,
\qquad
x \in \mathbb{R}^{n} \!\;\!, \;\, t > 0,
\end{equation}
\mbox{} \vspace{-0.150cm} \\
examinada originalmente por Fujita \cite{Fujita1966}
e subsequentemente generalizada por outros
(ver e.g.$\;$\cite{BandleBrunner1998, DengLevine2000, %
Levine1990, Pinsky1997, QuittnerSouplet2007}),
onde {\em todas as solu\c c\~oes n\~ao negativas\/}
(exceto $ w(\cdot,t) \equiv 0 $)
{\em explodem em tempo finito se\/}
$ \:\!0 < \kappa \leq 2/n $
(e tamb\'em para $ \kappa > 2/n \:\!$
se $ w(\cdot,0) $
for apropriadamente grande)
\cite{Fujita1966, Hayakawa1973}.
Como ficar\'a mostrado nos resultados a seguir,
a \linebreak
situa\c c\~ao em (1.4) tem natureza oposta:
{\em todas as solu\c c\~oes de\/} (1.4)
{\em s\~ao globalmente de- \linebreak
finidas se\/} $ \;\!0 \leq \kappa < 1/n $
(e tamb\'em para $ \kappa \geq 1/n \;\!$
se $ v(\cdot,0) $
for apropriadamente pequena).
%
%
Esta diferen\c ca not\'avel
entre os dois sistemas
\'e devida ao fato de (1.4) conservar massa,
o que n\~ao acontece com (1.5).
Considera\c c\~oes an\'alogas
podem ser feitas no caso geral
do problema (1.1)-(1.3):
todas as solu\c c\~oes
v\~ao existir globalmente
se $ \;\!0 \leq \kappa < 1/n $
(pela raz\~ao de se ter
conserva\c c\~ao de massa
e, melhor ainda,
no caso de solu\c c\~oes
$ u(\cdot,t) $
sem sinal definido,
a propriedade
dada em (1.11)
abaixo).

\mbox{} \vspace{-1.500cm} \\

Como sugerido em (1.4$^{\prime}$),
a {\em magnitude\/} do coeficiente
$ \mbox{\boldmath $b$}(x,t,u) $
n\~ao deve desempenhar papel importante,
ao contr\'ario de suas derivadas ---
ou, mais propriamente,
a {\em varia\c c\~ao\/}
de $ \mbox{\boldmath $b$}(x,t,u) $
em $ \mathbb{R}^{n} \!\,\!$,
dada pela quantidade $ B(t) $
definida do seguinte modo.
Para cada $ 1 \leq j \leq n $,
introduzimos $ B_{j}(t) $
dada por \\
\mbox{} \vspace{-0.500cm} \\
\begin{equation}
\tag{1.6$a$}
B_{j}(t) \;\!:= \;
\mbox{\small $ {\displaystyle \frac{1}{2} } $}
\:
\Bigl[\,
\sup_{\;x \,\in\, \mathbb{R}^{n}}
\!b_{j}(x,t,u(x,t))
\;\;\!-\!\;\!
\inf_{\;x \,\in\, \mathbb{R}^{n}}
\!\!\:\!b_{j}(x,t,u(x,t))
\;\Bigr],
\qquad
0 \leq t < \mbox{\small $T$}_{\!\;\!\ast} \!\;\!,
\end{equation}
\mbox{} \vspace{-0.050cm} \\
e ent\~ao
definimos \\
\mbox{} \vspace{-0.550cm} \\
\begin{equation}
\tag{1.6$b$}
B(t) \:=\:
|\: (\:\!
B_{\mbox{}_{1}}\!\;\!(t), ..., B_{\mbox{}_{\scriptstyle n}}\!\;\!(t)
\:\!) \,|_{\mbox{}_{2}}
\;\!=\:
\Bigl\{\:\! B_{\mbox{}_{1}}\!\;\!(t)^{\:\!2} \!\;\!+ ... +
B_{\mbox{}_{\scriptstyle n}}\!\;\!(t)^{\:\!2} \:\!
\Bigr\}^{\!1/2}
\end{equation}
\mbox{} \vspace{-0.170cm} \\
para cada
$ \:\! 0 \leq t < \mbox{\small $T$}_{\!\;\!\ast} $
\mbox{[}$\,$acima,
e em todo o texto que segue,
$ [\,0, \:\!\mbox{\small $T$}_{\!\;\!\ast}\!\;\!) $
denota sempre o
intervalo m\'aximo de exist\^encia
da solu\c c\~ao
$ u(\cdot,t) $ considerada$\,$\mbox{]}.
Para a descri\c c\~ao
dos resultados principais
a serem obtidos neste trabalho,
precisamos ainda introduzir
as quantidades
$ \;\!\mathbb{B}_{\mu}\!\;\!(0\:\!; \:\!t) \;\!$
e
$ \;\!\mathbb{U}_{p}(0\:\!; \:\!t) $,
$ 1 \leq p \leq \infty $,
definidas
por \\
\mbox{} \vspace{-0.475cm} \\
\begin{equation}
\tag{1.7}
\mathbb{B}_{\mu}\!\;\!(0\:\!; \:\!t)
\,:=\;
\sup \;\Bigl\{\,
\frac{B(\tau)}{\mu(\tau)} \!\;\!:
\; 0 \,\leq\,\tau \,\leq\, t
\,\Bigr\},
\end{equation}
\mbox{} \vspace{-0.600cm} \\
\begin{equation}
\tag{1.8}
\mathbb{U}_{p}(0\:\!; \:\!t)
\,:=\;
\sup \;\Bigl\{\:
\|\, u(\cdot,\tau) \,
\|_{\mbox{}_{\scriptstyle L^{p}(\mathbb{R}^{n})}}
\!\!\;\!:
\; 0 \,\leq\,\tau \,\leq\, t
\,\Bigr\},
\end{equation}
\mbox{} \vspace{-0.050cm} \\
para
$ \;\!0 \leq t < \mbox{\small $T$}_{\!\;\!\ast} \!\;\!$,
$ \:\!1 \leq p \leq \infty $.
Na Se\c c\~ao~3,
ap\'os alguns preliminares
coletados \linebreak
(por conveni\^encia)
na Se\c c\~ao 2 anterior,
as seguintes
propriedades fundamentais
das solu\-\c c\~oes do
problema
(1.1)$\;\!$-$\;\!$(1.3)
ser\~ao mostradas.
A primeira delas esta\-belece
o impor\-tante fato de
suas solu\c c\~oes
serem todas globais
quando $ \kappa > 0 $
n\~ao for grande \linebreak
(sendo o valor cr\'\i tico,
no caso do problema (1.1),
dado por $ 1/n $). \\
\mbox{} \vspace{-0.250cm} \\
%
%
%
%
%
%
\mbox{} \hspace{-0.800cm}
\fbox{%
\begin{minipage}[t]{16.000cm}
\mbox{} \vspace{-0.450cm} \\
\mbox{} \hspace{+0.300cm}
\begin{minipage}[t]{15.000cm}
\mbox{} \vspace{+0.100cm} \\
{\bf Teorema A.}
\textit{%
$\;\!$Sendo
$ \;\! 0 \leq \kappa < 1/n $,
$\:\!$a solu\c c\~ao
do problema
$\,\!(1.1)\;\!$-$\;\!(1.3)\!\:\!$ acima
est\'a definida
globalmente
$($i.e.,
$ \mbox{\small $T$}_{\!\;\!\ast} \!\;\!= \infty )$,
para qualquer
$ \;\!u_0 \in L^{1}(\mathbb{R}^{n}) \cap L^{\infty}(\mathbb{R}^{n}) $,
tendo-se
} \linebreak
\mbox{} \vspace{-0.650cm} \\
\begin{equation}
\tag{1.9}
\|\, u(\cdot,t) \,
\|_{\mbox{}_{\scriptstyle L^{\infty}(\mathbb{R}^{n})}}
\;\!\leq\:
K\!\;\!(n,\:\!\kappa) \,\cdot\,
\max\;
\Bigl\{\:
\|\, u_0 \;\!
\|_{\mbox{}_{\scriptstyle L^{\infty}(\mathbb{R}^{n})}};
\;\:\!
\mathbb{B}_{\mu}\!\;\!(0\:\!;\:\!t)^{\mbox{}^{\scriptstyle \!
\frac{\scriptstyle n}{1 \,-\, {\scriptstyle n \:\!\kappa}} }}
\:\!
\|\, u_0 \;\!
\|_{\mbox{}_{\scriptstyle L^{1}(\mathbb{R}^{n})}}
^{\mbox{}^{\scriptstyle
\frac{1}{1 \,-\, {\scriptstyle n \:\!\kappa}} }}
\;\!\Bigr\}
\end{equation}
\mbox{} \vspace{-0.200cm} \\
\textit{%
para todo
$\;\! 0 \leq t < \infty $,
$\:\!$onde
$ {\displaystyle
\:\!
K\!\;\!(n,\:\!\kappa)
\:\!=\,
2^{\mbox{}^{\frac{\scriptstyle n}{1 \,-\, {\scriptstyle n \:\!\kappa}} }}
\!\!
} $. \\
}
\end{minipage}
\end{minipage}
}
%
%
\nl
\mbox{} \vspace{-0.200cm} \\
No caso $ \kappa \geq 1/n $,
uma solu\c c\~ao
ser\'a garantidamente global
quando conseguir ser mostrado
que alguma
(e ent\~ao todas)
de suas normas altas
$ {\displaystyle
\;\!
\|\, u(\cdot,t) \,
\|_{\mbox{}_{\scriptstyle L^{p}(\mathbb{R}^{n})}}
\!\;\!
} $,
$ \;\! p > n \:\!\kappa $,
n\~ao puder
explodir em tempo finito,
como consequ\^encia
do seguinte resultado: \\
\mbox{} \vspace{-0.500cm} \\
%
%
%
%
%
%
\mbox{} \hspace{-0.800cm}
\fbox{%
\begin{minipage}[t]{16.000cm}
\mbox{} \vspace{-0.450cm} \\
\mbox{} \hspace{+0.300cm}
\begin{minipage}[t]{15.000cm}
\mbox{} \vspace{+0.100cm} \\
{\bf Teorema B.}
\textit{%
$\;\!$Sejam
$ \;\! \kappa \geq 0 \;\!$
e
$\;\!u(\cdot,t) $,
$ 0 \leq t < \mbox{\small $T$}_{\!\;\!\ast}\!\;\!$,
solu\c c\~ao
do problema
$\,\!(1.1)\;\!$-$\;\!(1.3) $.
%
%
Para cada $\;\! p \geq 1 $
satisfazendo
$ \;\!p > n \:\!\kappa $,
tem-se \\
}
\mbox{} \vspace{-0.650cm} \\
\begin{equation}
\tag{1.10}
\mbox{} \!\!
\|\, u(\cdot,t) \,
\|_{\mbox{}_{\scriptstyle L^{\infty}(\mathbb{R}^{n})}}
\;\!\leq\,
K\!\;\!(n,\:\!\kappa, \:\!p) \,\cdot\,
\max\;
\Bigl\{\:
\|\, u_0 \;\!
\|_{\mbox{}_{\scriptstyle L^{\infty}(\mathbb{R}^{n})}};
\;\:\!
\mathbb{B}_{\mu}\!\;\!(0; t)^{\mbox{}^{\scriptstyle \!\;\!
\frac{\scriptstyle n}{{\scriptstyle p} \,-\, {\scriptstyle n \;\!\kappa}} }}
\;\!
\mathbb{U}_{p}\!\;\!(0; t)^{\mbox{}^{\scriptstyle \!\;\!
\frac{\scriptstyle p}{{\scriptstyle p} \,-\, {\scriptstyle n \;\!\kappa}} }}
\;\!\Bigr\}
\end{equation}
\mbox{} \vspace{-0.200cm} \\
\textit{%
para todo
$\;\! 0 \leq t < \mbox{\small $T$}_{\!\;\!\ast} \!\;\! $,
$\:\!$sendo
$ {\displaystyle
\:\!
K\!\;\!(n,\:\!\kappa, \:\!p)
\:\!=\,
\{\;\!2\;\!p \;\!\}^{\mbox{}^{\scriptstyle \!\:\!
\frac{\scriptstyle n}{{\scriptstyle p} \,-\, {\scriptstyle n \:\!\kappa}} }}
\!\!
} $. \\
}
\end{minipage}
\end{minipage}
}
%
%
\nl
\mbox{} \vspace{-0.050cm} \\
Note-se que o {\small \sc Teorema A}
\'e um corol\'ario do {\small \sc Teorema B} acima
(tomando-se $ p = 1 $),
em virtude da seguinte propriedade
b\'asica (conhecida)
das solu\c c\~oes da equa\c c\~ao (1.1): \\
\mbox{} \vspace{-0.600cm} \\
\begin{equation}
\tag{1.11}
\|\, u(\cdot,t) \,
\|_{\mbox{}_{\scriptstyle L^{1}(\mathbb{R}^{n})}}
\;\!\leq\;
\|\, u_0 \;\!
\|_{\mbox{}_{\scriptstyle L^{1}(\mathbb{R}^{n})}}
\!\;\!,
\qquad
\forall \;\,
t > 0,
\end{equation}
\mbox{} \vspace{-0.200cm} \\
ou seja,
$ {\displaystyle
\|\, u(\cdot,t) \,
\|_{\mbox{}_{\scriptstyle L^{1}(\mathbb{R}^{n})}}
\!\:\!
} $
decresce monotonicamente em $\:\!t$.
Esta propriedade
\'e tamb\'em satisfeita
pelas demais normas
$ {\displaystyle
\|\, u(\cdot,t) \,
\|_{\mbox{}_{\scriptstyle L^{p}(\mathbb{R}^{n})}}
\!\:\!
} $,
$ \,\!p > 1 $,
quando o termo $ \mbox{\boldmath $b$} $
em (1.1$a$)
n\~ao depender explicitamente de $\:\!x$,
ou,
mais geralmente,
se tivermos
\cite{SchutzZinganoZingano2014} \\
\mbox{} \vspace{-0.650cm} \\
\begin{equation}
\tag{1.12$a$}
\sum_{j\,=\,1}^{n}
\;\!
\frac{\partial \:\!b_{\:\!j}}
{\partial \:\!x_{\!\;\!j}}
\:\!
(x,t,\mbox{u})
\;\geq\;0,
\qquad
\forall \;\,
x \in \mathbb{R}^{n} \!, \;
t \geq 0, \;
\mbox{u} \in \mathbb{R}.
\end{equation}
\mbox{} \vspace{-0.075cm} \\
Neste caso,
n\~ao apenas
estar\~ao
as solu\c c\~oes de (1.1$a$), (1.1$b$)
definidas
para todo $t > 0 $, \linebreak
como tamb\'em
decair\~ao
ao $ \;\!t \rightarrow \infty $,
tendo-se
(\cite{BrazSchutzZingano2013}, Theorem 3.2) \\
\mbox{} \vspace{-0.600cm} \\
\begin{equation}
\tag{1.12$b$}
\|\, u(\cdot,t) \,
\|_{\mbox{}_{\scriptstyle L^{\infty}(\mathbb{R}^{n})}}
\;\!\leq\:
(\;\! 2 \:\!e \;\!)^{\mbox{}^{\scriptstyle \!-\,
\frac{\scriptstyle n}{2} }}
\:\!
\|\, u_0 \;\!
\|_{\mbox{}_{\scriptstyle L^{1}(\mathbb{R}^{n})}}
\;\!
t^{\mbox{}^{\scriptstyle \!-\,
\frac{\scriptstyle n}{2} }}
\qquad
\forall \;\,
t > 0.
\end{equation}
\mbox{} \vspace{-0.175cm} \\
Neste trabalho,
estamos justamente
interessados
na situa\c c\~ao
(muito mais dif\'\i cil)
em que
(1.12$a$) n\~ao \'e v\'alida,
quando (em geral)
n\~ao se tem decaimento,
podendo existir
solu\c c\~oes estacion\'arias,
etc.
Mesmo quando
$ {\displaystyle
\;\!
\|\, u(\cdot,t) \,
\|_{\mbox{}_{\scriptstyle L^{\infty}(\mathbb{R}^{n})}}
\!
\rightarrow 0
} $,
a taxa de decaimento
n\~ao \'e conhecida, em geral.
Experimentos num\'ericos
parecem indicar
o seguinte comportamento,
quando
$ \:\!\mbox{\boldmath $b$}, \,\!\mbox{\boldmath $f$} \,\!$
na equa\c c\~ao (1.1$a$)
independem do tempo $\:\!t\;\!$: \\
\mbox{} \vspace{-0.125cm} \\
%
%
{\it Conjectura A\/}:
as solu\c c\~oes estacion\'arias,
quando existem,
s\~ao est\'aveis
(atratoras); \linebreak
\mbox{} \vspace{-0.200cm} \\
%
%
{\it Conjectura B\/}:
na aus\^encia
de solu\c c\~oes estacion\'arias
(exceto $ u \equiv 0$)
tem-se
sempre
$ {\displaystyle
\;\!
\|\, u(\cdot,t) \,
\|_{\mbox{}_{\scriptstyle L^{\infty}(\mathbb{R}^{n})}}
\!
\rightarrow 0
\;\!
} $
ao
$ \:\! t \rightarrow \infty $. \\
\mbox{} \vspace{-0.050cm} \\
Como estas,
muitas quest\~oes
de interesse para (1.1)$\;\!$-$\;\!$(1.3)
permanecem em aberto. \linebreak

%
%
%
%
%
\mbox{} \vspace{-1.500cm} \\

{\bf 2. Preliminares} \\

Nesta se\c c\~ao,
revisaremos resumidamente
alguns resultados b\'asicos
para as solu\-\c c\~oes $ u(\cdot,t) $
do problema (1.1)$\;\!$-$\;\!$(1.3),
que ser\~ao usadas na an\'alise
a seguir
(Se\c c\~ao~3).
Estas propriedades
podem ser estabelecidas,
sem esfor\c co adicional,
para o problema
levemente mais geral \\
\mbox{} \vspace{-0.575cm} \\
\begin{equation}
\tag{2.1$a$}
u_t \,+\;
\mbox{div}\,
\bigl(\;\! \mbox{\boldmath $b$}(x,t,u) \, |\;\!u\;\!|^{\:\!\kappa} \:\! u \;\!
\bigr)
\:+\:
\mbox{div} \,
\mbox{\boldmath $f$}(t,u)
\;=\;
\mbox{div}\,\bigl(\;\!
A(x,t,u) \;\! \nabla u \;\!\bigr),
\end{equation}
\mbox{} \vspace{-0.900cm} \\
\begin{equation}
\tag{2.1$b$}
u(\cdot,0) \,=\,
u_0 \in L^{p_{\mbox{}_{0}}}(\mathbb{R}^{n})
\cap L^{\infty}(\mathbb{R}^{n}),
\end{equation}
\mbox{} \vspace{-0.175cm} \\
para
$ \:\! 1 \leq p_{\mbox{}_{0}} \!< \infty \;\!$
dado
(e n\~ao somente
$ \:\! p_{\mbox{}_{0}} \!= 1 $,
como em (1.1$b$)),
onde a matriz $A$ \linebreak
satisfaz a condi\c c\~ao
de elipticidade
(1.2) acima
(para certa $ \mu \in C^{0}([\,0,\infty)) $
positiva)
e
$ \mbox{\boldmath $b$} $
satisfaz (1.3).\footnote{%
%
%
Al\'em disso,
supo\~oe-se
$ \mbox{\boldmath $b$} $
suave
(mais precisamente:
$ \mbox{\boldmath $b$} $,
$ \mbox{\boldmath $b$}_{\mbox{\scriptsize $x_{\mbox{}_{\!\;\!1}}$}} \!\,\!,
\!\:\!..., \mbox{\boldmath $b$}_{\mbox{\scriptsize $x_{\mbox{}_{n}}$}} \!\!\!\,\!$
e
$ \mbox{\boldmath $b$}_{\mbox{\scriptsize u}} $
s\~ao supostas cont\'\i nuas).
Sobre o termo de fluxo
$ \mbox{\boldmath $f$} \!\;\!$,
por n\~ao depender de $x$,
s\'o ser\'a preciso supor
que
$ \mbox{\boldmath $f$} \!\:\!$,
$ \mbox{\boldmath $f$}_{\!\:\!\mbox{\scriptsize u}} $
sejam cont\'\i nuas.
}
%
%
%
$\!\!\;\!$Em (2.1),
a condi\c c\~ao
(2.1$b$)
\'e entendida
no sentido de
$ L^{1}_{\tt loc}(\mathbb{R}^{n}) $,
ou seja, \\
\mbox{} \vspace{-1.000cm} \\
\begin{equation}
\tag{2.2}
\|\, u(\cdot,t) \;\!-\;\! u_0 \;\!
\|_{\mbox{}_{\scriptstyle L^{1}(\mathbb{K})}}
\!\;\!\rightarrow\, 0
\quad \;\;
\mbox{ao }\,
t \rightarrow 0
\end{equation}
\mbox{} \vspace{-0.200cm} \\
para cada conjunto compacto
$ \:\!\mathbb{K} \subset \mathbb{R}^{n} \!\:\!$
considerado.
Por {\em solu\c c\~ao\/} de (2.1$a$), (2.1$b$) em \linebreak
um dado intervalo
$ \,\![\,0, \,\!\mbox{\small $T$}_{\!\,\!\ast}\!\;\!) $
entende-se uma fun\c c\~ao suave
$ {\displaystyle
u(\cdot,t)
\in
%
%
L^{\infty}_{\tt loc}(\:\!
[\;\!0, \:\!\mbox{\small $T$}_{\!\,\!\ast}\!\;\!),
\:\!
L^{\infty}(\mathbb{R}^{n})
\,\!)
} $
que satisfaz a equa\c c\~ao
(2.1$a$) classicamente
para
$ 0 < t < \mbox{\small $T$}_{\!\;\!\ast} \!\;\!$
e verifica (2.2)
ao $ t \rightarrow 0 $. \linebreak
A {\em exist\^encia\/} (local) de tais solu\c c\~oes
de\-corre da teoria geral de equa\c c\~oes parab\'olicas
(ver
e.g.$\;$\cite{LadyzhenskayaSolonnikovUralceva1968},
ou \cite{Serre1999}, Ch.$\;$7);
sabe-se tamb\'em que
as solu\c c\~oes s\~ao {\em \'unicas},
como pode ser mos\-trado usando princ\'\i pios de
compara\c c\~ao
(ver e.g.$\;$\cite{DiehlFabrisZingano2014}, Theorem~2.1). \\
%
%
\nl
%
%
%
%
{\bf Proposi\c c\~ao 2.1.}
\textit{%
Sendo
$ {\displaystyle
\;\!
u(\cdot,t)
\in
%
L^{\infty}_{\tt loc}(\:\!
[\,0, \:\!\mbox{\small $T$}_{\!\;\!\ast}),
\:\!
L^{\infty}(\mathbb{R}^{n})
\,\!)
} $
so\-lu\-\c c\~ao do
problema $\;\!(2.1)$,
onde
$ 0 < \mbox{\small $T$}_{\!\;\!\ast}\!\;\!\leq \infty $,
ent\~ao
tem-se
$ {\displaystyle
\;\!
u(\cdot,t) \in
C^{0}(\:\!
[\,0, \:\!\mbox{\small $T$}_{\!\;\!\ast}),
\:\!
L^{p}(\mathbb{R}^{n})
\,\!)
} $
para cada
$ \;\! p_{\mbox{}_{0}} \!\;\!\leq p < \infty $.
Al\'em disso,
tem-se,
para cada
$ \;\! p_{\mbox{}_{0}} \!\;\!\leq p < \infty \!\!\;\!:$
} \\
\mbox{} \vspace{-0.725cm} \\
\begin{equation}
\tag{2.3}
\|\, u(\cdot,t) \,
\|_{\mbox{}_{\scriptstyle L^{p}(\mathbb{R}^{n})}}
\;\!\leq\;
\|\, u_0 \;\!
\|_{\mbox{}_{\scriptstyle L^{p}(\mathbb{R}^{n})}}
\!\cdot\;\!\;\!
\exp\,\Bigl\{\:
\mbox{\small $ {\displaystyle \frac{1}{4} }$} \,
(p - 1) \:
\mathbb{B}_{\mu}\!\;\!(0\:\!;\:\!t)^{\:\!2} \:
\mathbb{U}_{p}\!\;\!(0\:\!;\:\!t)^{\:\!2\;\!\kappa}
\!\!\!\;\!
\int_{0}^{\;\!\mbox{\footnotesize $t$}}
\!\!\:\!
\mu(\tau) \, d\tau
\;\!\Bigr\}
\end{equation}
\mbox{} \vspace{-0.175cm} \\
\textit{%
para todo $\;\! 0 \leq t < \mbox{\small $T$}_{\!\;\!\ast} \!\;\!$,
com
$ \;\! \mathbb{B}_{\mu}\!\;\!(0\:\!;\:\!t) $,
$ \mathbb{U}_{p}\!\;\!(0\:\!;\:\!t) $
definidas em
$\;\!(1.7) $, $(1.8)$
acima. \\
}
%
%
\mbox{} \vspace{-0.750cm} \\

A prova da Proposi\c c\~ao 2.1
pode ser feita adaptando-se o argumento usado
em (\cite{BarrionuevoOliveiraZingano2014},
Theorem 1),
ou (\cite{BrazSchutzZingano2013}, Theorem 2.1).
Em particular,
com $ \;\!p_{\mbox{}_{0}} \!\;\!= 1 $,
$ p = 1 $ em (2.3),
obt\'em-se
a estimativa (1.11)
referida na Se\c c\~ao~1,
ou seja, \\
\mbox{} \vspace{-0.600cm} \\
\begin{equation}
\tag{2.4}
\mbox{} \hspace{+0.750cm}
\|\, u(\cdot,t) \,
\|_{\mbox{}_{\scriptstyle L^{1}(\mathbb{R}^{n})}}
\;\!\leq\;
\|\, u_0 \;\!
\|_{\mbox{}_{\scriptstyle L^{1}(\mathbb{R}^{n})}}
\!\;\!,
\qquad
\forall \;\,
t > 0.
\end{equation}

\mbox{} \vspace{-1.250cm} \\
\nl
%
%
%
%
{\bf Proposi\c c\~ao 2.2.}
\textit{%
Sendo
$ {\displaystyle
\;\!
u(\cdot,t)
\in
%
%
L^{\infty}_{\tt loc}(\:\!
[\,0, \:\!\mbox{\small $T$}_{\!\;\!\ast}),
\:\!
L^{\infty}(\mathbb{R}^{n})
\,\!)
} $
so\-lu\-\c c\~ao do
problema $\;\!(2.1) \!\;\!$
em um dado intervalo
$ \:\![\,0, \:\!\mbox{\small $T$}_{\!\;\!\ast}) $,
ent\~ao
tem-se,
para cada
$ \;\! q \;\!\geq\;\! p_{\mbox{}_{0}} \!\:\!+ 1 \!\!\:\!:$ \\
}
\mbox{} \vspace{-0.675cm} \\
\begin{equation}
\tag{2.5}
\int_{0}^{\;\!\mbox{\footnotesize $t$}}
\!\!
\int_{\mathbb{R}^{n}}
\!\!\;\!
|\, u(x,\tau)\,|^{\:\!q \;\!-\;\! 2}
\,
|\, \nabla u \,|^{\:\!2}
\:dx\, d\tau
\;<\; \infty,
\qquad
\forall \;\,
0 < t < \mbox{\small $T$}_{\!\;\!\ast}
\!\;\!.
\end{equation}
\mbox{} \vspace{-0.350cm} \\
%
%
%
\nl
{\small
{\bf Prova:}
Para $ q > 2 $,
considere
$ \;\!\Phi(\mbox{u}) \!\:\!:= L_{\delta}(\mbox{u})^{q} \!\;\!$,
onde
$ L_{\delta}(\cdot) $
\'e uma fun\c c\~ao sinal regularizada
(ver e.g.$\;$\cite{BrazSchutzZingano2013, %
DiehlFabrisZingano2014, KreissLorenz1989}),
sendo
$ \delta > 0 $ dado.
Seja tamb\'em,
para $ R > 0 $ grande,
$ \zeta_{\mbox{}_{R}} \in C^{\infty}(\mathbb{R}^{n}) $
uma fun\c c\~ao de corte
satisfazendo
$ \zeta_{\mbox{}_{R}}(x) = 1 $
se $ |\,x\,| \leq R - 1 $,
$ \zeta_{\mbox{}_{R}}(x) = 0 $
se $ |\,x\,| \geq R $,
$ 0 \leq \zeta_{\mbox{}_{R}}(x) \leq 1 $ para todo $x$,
$ |\,\nabla \zeta_{\mbox{}_{R}}(x) \,| \leq C $
para certa constante $C$ independente de $x$ e $R$. \linebreak
Multiplicando-se a equa\c c\~ao (2.1$a$)
por $ \;\!\Phi^{\prime}(u(x,t)) \, \zeta_{\mbox{}_{R}}(x) $
e integrando-se em $ [\,t_0, \:\!t\;\!] $
(dado $ 0 < t_0 < t $ arbitr\'ario),
obt\'em-se,
integrando-se por partes e
fazendo $ \delta \rightarrow 0 $: \\
\mbox{} \vspace{-0.050cm} \\
\mbox{} \hspace{-0.200cm}
$ {\displaystyle
\int_{\mbox{}_{\scriptstyle |\,x\,| \,<\,R}}
\hspace{-0.750cm}
|\, u(x,t) \,|^{q} \,
\zeta_{\mbox{}_{R}}(x)
\: dx
\:+\:
q \,(q - 1) \!
\int_{\mbox{}_{\!\;\!\mbox{\footnotesize $t$}_0}}
    ^{\,\mbox{\footnotesize $t$}}
\!\!\;\!
\int_{\mbox{}_{\scriptstyle |\,x\,| \,<\,R}}
\hspace{-0.750cm}
|\, u(x,\tau) \,|^{q - 2} \,
\bigl<\,A(x,\tau,u) \, \nabla u, \;\! \nabla u \;\!\bigr>
\: \zeta_{\mbox{}_{R}}(x)
\: dx \, d\tau
\;=
} $ \\
\mbox{} \vspace{-0.250cm} \\
\mbox{} \hfill (2.6) \\
\mbox{} \vspace{-0.350cm} \\
\mbox{} \hfill
$ {\displaystyle
=
\int_{\mbox{}_{\scriptstyle |\,x\,| \,<\,R}}
\hspace{-0.750cm}
|\, u(x,t_0) \,|^{q} \,
\zeta_{\mbox{}_{R}}(x)
\: dx
\:+\:
q \,(q - 1) \!
\int_{\mbox{}_{\!\;\!\mbox{\footnotesize $t$}_0}}
    ^{\,\mbox{\footnotesize $t$}}
\!\!\;\!
\int_{\mbox{}_{\scriptstyle |\,x\,| \,<\,R}}
\hspace{-0.750cm}
|\, u \,|^{q - 2 + \kappa} \, u \,
\bigl<\;\!\mbox{\boldmath $b$}(x,\tau,u) -
\mbox{\boldmath $\beta$}(\tau), \;\! \nabla u \;\!\bigr>
\: \zeta_{\mbox{}_{R}}(x)
\: dx \, d\tau
} $ \\
\mbox{} \vspace{+0.2000cm} \\
\mbox{} \hfill
$ {\displaystyle
-\;
q \,(q - 1) \!
\int_{\mbox{}_{\!\;\!\mbox{\footnotesize $t$}_0}}
    ^{\,\mbox{\footnotesize $t$}}
\!\!\;\!
\int_{\mbox{}_{\scriptstyle \!\:\!R - 1 \,<\, |\,x\,| \,<\,R}}
\hspace{-1.675cm}
G_{q}(u) \:
\bigl<\;\!\mbox{\boldmath $\beta$}(\tau),
\;\! \nabla \zeta_{\mbox{}_{R}}(x)
\,\bigr>
\: dx \, d\tau
\;\!+\,
q \!
\int_{\mbox{}_{\!\;\!\mbox{\footnotesize $t$}_0}}
    ^{\,\mbox{\footnotesize $t$}}
\!\!\;\!
\int_{\mbox{}_{\scriptstyle \!\:\!R - 1 \,<\,|\,x\,| \,<\,R}}
\hspace{-1.675cm}
|\,u\,|^{q + \kappa} \,
\bigl<\;\!\mbox{\boldmath $b$}(x,\tau,u),
\;\! \nabla \zeta_{\mbox{}_{R}}(x)
\,\bigr>
\: dx \, d\tau
} $ \\
\mbox{} \vspace{+0.200cm} \\
\mbox{} \hfill
$ {\displaystyle
-\;
q \,(q - 1) \!
\int_{\mbox{}_{\!\;\!\mbox{\footnotesize $t$}_0}}
    ^{\,\mbox{\footnotesize $t$}}
\!\!\;\!
\int_{\mbox{}_{\scriptstyle \!\:\!R - 1 \,<\, |\,x\,| \,<\,R}}
\hspace{-1.675cm}
\bigl<\;\!\mbox{\boldmath $F$}_{\!\!\;\!q}(\tau,u),
\;\! \nabla \zeta_{\mbox{}_{R}}(x)
\,\bigr>
\: dx \, d\tau
\;\!+\,
q \!
\int_{\mbox{}_{\!\;\!\mbox{\footnotesize $t$}_0}}
    ^{\,\mbox{\footnotesize $t$}}
\!\!\;\!
\int_{\mbox{}_{\scriptstyle \!\:\!R - 1 \,<\,|\,x\,| \,<\,R}}
\hspace{-1.675cm}
|\,u\,|^{q - 2} \;\!u \:
\bigl<\;\!\mbox{\boldmath $f$}(\tau,u),
\;\! \nabla \zeta_{\mbox{}_{R}}(x)
\,\bigr>
\: dx \, d\tau
} $ \\
\mbox{} \vspace{+0.200cm} \\
\mbox{} \hspace{+2.750cm}
$ {\displaystyle
-\;
q \!
\int_{\mbox{}_{\!\;\!\mbox{\footnotesize $t$}_0}}
    ^{\,\mbox{\footnotesize $t$}}
\!\!\;\!
\int_{\mbox{}_{\scriptstyle \!\:\!R - 1 \,<\, |\,x\,| \,<\,R}}
\hspace{-1.675cm}
|\,u\,|^{q - 2} \;\!u \:
\bigl<\,A(x,\tau,u) \,\nabla u,
\;\! \nabla \zeta_{\mbox{}_{R}}(x)
\,\bigr>
\: dx \, d\tau
} $, \\
\mbox{} \vspace{+0.050cm} \\
onde \\
\mbox{} \vspace{-0.900cm} \\
\begin{equation}
\tag{2.7$a$}
G_{q}(\mbox{u})
\,:=\:\!
\int_{0}^{\:\!\mbox{\scriptsize u}}
\!
|\,\mbox{v}\,|^{\:\!q - 2 + \kappa}
\;\! \mbox{v} \:
d\mbox{v},
\qquad
\mbox{\boldmath $F$}_{\!\!\;\!q}(t,\mbox{u})
\,:=\:\!
\int_{0}^{\:\!\mbox{\scriptsize u}}
\!
|\,\mbox{v}\,|^{q - 2} \;\!
\mbox{\boldmath $f$}(t,\mbox{v}) \:
d\mbox{v},
\end{equation}
\mbox{} \vspace{-0.150cm} \\
e
$ {\displaystyle
\,
\mbox{\boldmath $\beta$}(t)
\;\!=\;\!
(\;\!\beta_{\mbox{}_{1}}(t),\!\;\!..., \beta_{\mbox{}_{\scriptstyle n}}(t)
\:\!)
} $,
com
$ \beta_{\mbox{}_{\scriptstyle j}}(t) $,
$ 1 \leq j \leq n $,
dado por \\
\mbox{} \vspace{-0.600cm} \\
\begin{equation}
\tag{2.7$b$}
\beta_{\scriptstyle j}(t) \;\!:= \;
\mbox{\small $ {\displaystyle \frac{1}{2} } $}
\:
\Bigl[\,
\sup_{\;x \,\in\, \mathbb{R}^{n}}
\!b_{j}(x,t,u(x,t))
\;\;\!+\!\;\!
\inf_{\;x \,\in\, \mathbb{R}^{n}}
\!\!\:\!b_{j}(x,t,u(x,t))
\;\Bigr].
\end{equation}
\mbox{} \vspace{-0.100cm} \\
Para obter (2.6),
observamos que \\
\mbox{} \vspace{-0.100cm} \\
\mbox{} \hspace{+0.250cm}
$ {\displaystyle
\int_{\mbox{}_{\scriptstyle |\,x\,| \,<\,R}}
\hspace{-0.750cm}
|\, u \,|^{q - 2 + \kappa} \, u \,
\bigl<\,\mbox{\boldmath $b$}(x,\tau,u), \;\! \nabla u \;\!\bigr>
\: \zeta_{\mbox{}_{R}}(x)
\: dx
\;=
} $ \\
\mbox{} \vspace{+0.200cm} \\
\mbox{} \hfill
$ {\displaystyle
=\;\!
\int_{\mbox{}_{\scriptstyle |\,x\,| \,<\,R}}
\hspace{-0.750cm}
|\, u \,|^{q - 2 + \kappa} \, u \,
\bigl<\,\mbox{\boldmath $b$}(x,\tau,u) -
\mbox{\boldmath $\beta$}(\tau), \;\! \nabla u \;\!\bigr>
\: \zeta_{\mbox{}_{R}}(x)
\: dx
\;+\;\!
\int_{\mbox{}_{\scriptstyle |\,x\,| \,<\,R}}
\hspace{-0.750cm}
\bigl<\,\mbox{\boldmath $\beta$}(\tau),
\;\! \nabla G_{q}(u) \;\!\bigr>
\: \zeta_{\mbox{}_{R}}(x)
\: dx
} $ \\
\mbox{} \vspace{+0.200cm} \\
\mbox{} \hfill
$ {\displaystyle
=\;\!
\int_{\mbox{}_{\scriptstyle |\,x\,| \,<\,R}}
\hspace{-0.750cm}
|\, u \,|^{q - 2 + \kappa} \, u \,
\bigl<\,\mbox{\boldmath $b$}(x,\tau,u) -
\mbox{\boldmath $\beta$}(\tau), \;\! \nabla u \;\!\bigr>
\: \zeta_{\mbox{}_{R}}(x)
\: dx
\;-\;\!
\int_{\mbox{}_{\scriptstyle \!\:\!R - 1 \,<\, |\,x\,| \,<\,R}}
\hspace{-1.675cm}
G_{q}(u) \:
\bigl<\,\mbox{\boldmath $\beta$}(\tau),
\;\! \nabla \zeta_{\mbox{}_{R}}(x) \;\!\bigr>
\: dx
} $
\mbox{} \vspace{-1.000cm} \\
e
tamb\'em \\
\mbox{} \vspace{-0.100cm} \\
\mbox{} \hspace{+0.750cm}
$ {\displaystyle
\int_{\mbox{}_{\scriptstyle |\,x\,| \,<\,R}}
\hspace{-0.750cm}
|\, u \,|^{q - 2} \,
\bigl<\;\!\mbox{\boldmath $f$}(\tau,u), \;\! \nabla u \;\!\bigr>
\: \zeta_{\mbox{}_{R}}(x)
\: dx
\;=\;\!
\int_{\mbox{}_{\scriptstyle |\,x\,| \,<\,R}}
\hspace{-0.715cm}
[\; \mbox{div} \,
\mbox{\boldmath $F$}_{\!\!\;\!q}(\tau,u) \:]
\:
\zeta_{\mbox{}_{R}}(x)
\: dx
} $ \\
\mbox{} \vspace{+0.200cm} \\
\mbox{} \hspace{+6.550cm}
$ {\displaystyle
=\;
-\!\;\!
\int_{\mbox{}_{\scriptstyle \!\:\!R - 1 \,<\, |\,x\,| \,<\,R}}
\hspace{-1.650cm}
\bigl<\,\mbox{\boldmath $F$}_{\!\!\;\!q}(\tau,u),
\;\! \nabla \zeta_{\mbox{}_{R}}(x) \;\!\bigr>
\: dx
} $. \\
\mbox{} \vspace{+0.100cm} \\
Como,
por (1.6),
tem-se \\
\mbox{} \vspace{-0.000cm} \\
\mbox{} \hspace{+0.750cm}
$ {\displaystyle
\int_{\mbox{}_{\scriptstyle |\,x\,| \,<\,R}}
\hspace{-0.750cm}
|\, u \,|^{q - 2 + \kappa} \, u \,
\bigl<\,\mbox{\boldmath $b$}(x,\tau,u) -
\mbox{\boldmath $\beta$}(\tau), \;\! \nabla u \;\!\bigr>
\: \zeta_{\mbox{}_{R}}(x)
\: dx
} $ \\
\mbox{} \vspace{+0.100cm} \\
\mbox{} \hspace{+3.100cm}
$ {\displaystyle
\leq\;
B(\tau) \!
\int_{\mbox{}_{\scriptstyle |\,x\,| \,<\,R}}
\hspace{-0.750cm}
|\, u \,|^{q - 1 + \kappa} \,
|\, \nabla u \,|
\: \zeta_{\mbox{}_{R}}(x)
\: dx
} $ \\
\mbox{} \vspace{+0.200cm} \\
\mbox{} \hfill
$ {\displaystyle
\leq\;
\frac{1}{2} \;
\mu(\tau) \!
\int_{\mbox{}_{\scriptstyle |\,x\,| \,<\,R}}
\hspace{-0.750cm}
|\, u \,|^{q - 2} \,
|\, \nabla u \,|^{\:\!2}
\: \zeta_{\mbox{}_{R}}(x)
\: dx
\;+\;
\frac{1}{2} \:
\frac{\,B(\tau)^{2}}{\mu(\tau)}
\!
\int_{\mbox{}_{\scriptstyle |\,x\,| \,<\,R}}
\hspace{-0.750cm}
|\, u \,|^{q + 2\;\! \kappa}
\: \zeta_{\mbox{}_{R}}(x)
\: dx
} $, \\
\mbox{} \vspace{+0.300cm} \\
resulta
de (1.2), (2.6)
que \\
\mbox{} \vspace{-0.100cm} \\
\mbox{} \hspace{-0.200cm}
$ {\displaystyle
\int_{\mbox{}_{\scriptstyle |\,x\,| \,<\,R}}
\hspace{-0.750cm}
|\, u(x,t) \,|^{q} \,
\zeta_{\mbox{}_{R}}(x)
\: dx
\:+\:
\frac{1}{2} \;
q \,(q - 1) \!
\int_{\mbox{}_{\!\;\!\mbox{\footnotesize $t$}_0}}
    ^{\,\mbox{\footnotesize $t$}}
\!\!
\mu(\tau)
\!\!\;\!
\int_{\mbox{}_{\scriptstyle |\,x\,| \,<\,R}}
\hspace{-0.750cm}
|\, u(x,\tau) \,|^{q - 2} \:
|\, \nabla u \,|^{2}
\: \zeta_{\mbox{}_{R}}(x)
\: dx \, d\tau
} $ \\
\mbox{} \vspace{+0.300cm} \\
\mbox{} \hfill
$ {\displaystyle
\leq
\int_{\mbox{}_{\scriptstyle |\,x\,| \,<\,R}}
\hspace{-0.750cm}
|\, u(x,t_0) \,|^{q} \,
\zeta_{\mbox{}_{R}}(x)
\: dx
\:+\:
\frac{1}{2} \;
q \,(q - 1) \!
\int_{\mbox{}_{\!\;\!\mbox{\footnotesize $t$}_0}}
    ^{\,\mbox{\footnotesize $t$}}
\!\!\;\!
\frac{\,B(\tau)^{2}}{\mu(\tau)}
\!
\int_{\mbox{}_{\scriptstyle |\,x\,| \,<\,R}}
\hspace{-0.750cm}
|\, u \,|^{q + 2\;\! \kappa}
\: \zeta_{\mbox{}_{R}}(x)
\: dx \, d\tau
} $ \\
\mbox{} \vspace{+0.200cm} \\
\mbox{} \hfill
$ {\displaystyle
-\;
q \,(q - 1) \!
\int_{\mbox{}_{\!\;\!\mbox{\footnotesize $t$}_0}}
    ^{\,\mbox{\footnotesize $t$}}
\!\!\;\!
\int_{\mbox{}_{\scriptstyle \!\:\!R - 1 \,<\, |\,x\,| \,<\,R}}
\hspace{-1.675cm}
G_{q}(u) \:
\bigl<\;\!\mbox{\boldmath $\beta$}(\tau),
\;\! \nabla \zeta_{\mbox{}_{R}}(x)
\,\bigr>
\: dx \, d\tau
\;\!+\,
q \!
\int_{\mbox{}_{\!\;\!\mbox{\footnotesize $t$}_0}}
    ^{\,\mbox{\footnotesize $t$}}
\!\!\;\!
\int_{\mbox{}_{\scriptstyle \!\:\!R - 1 \,<\,|\,x\,| \,<\,R}}
\hspace{-1.675cm}
|\,u\,|^{q + \kappa} \,
\bigl<\;\!\mbox{\boldmath $b$}(x,\tau,u),
\;\! \nabla \zeta_{\mbox{}_{R}}(x)
\,\bigr>
\: dx \, d\tau
} $ \\
\mbox{} \vspace{+0.200cm} \\
\mbox{} \hfill
$ {\displaystyle
-\;
q \,(q - 1) \!
\int_{\mbox{}_{\!\;\!\mbox{\footnotesize $t$}_0}}
    ^{\,\mbox{\footnotesize $t$}}
\!\!\;\!
\int_{\mbox{}_{\scriptstyle \!\:\!R - 1 \,<\, |\,x\,| \,<\,R}}
\hspace{-1.675cm}
\bigl<\;\!\mbox{\boldmath $F$}_{\!\!\;\!q}(\tau,u),
\;\! \nabla \zeta_{\mbox{}_{R}}(x)
\,\bigr>
\: dx \, d\tau
\;\!+\,
q \!
\int_{\mbox{}_{\!\;\!\mbox{\footnotesize $t$}_0}}
    ^{\,\mbox{\footnotesize $t$}}
\!\!\;\!
\int_{\mbox{}_{\scriptstyle \!\:\!R - 1 \,<\,|\,x\,| \,<\,R}}
\hspace{-1.675cm}
|\,u\,|^{q - 2} \;\!u \:
\bigl<\;\!\mbox{\boldmath $f$}(\tau,u),
\;\! \nabla \zeta_{\mbox{}_{R}}(x)
\,\bigr>
\: dx \, d\tau
} $ \\
\mbox{} \vspace{+0.200cm} \\
\mbox{} \hspace{+2.750cm}
$ {\displaystyle
-\;
q \!
\int_{\mbox{}_{\!\;\!\mbox{\footnotesize $t$}_0}}
    ^{\,\mbox{\footnotesize $t$}}
\!\!\;\!
\int_{\mbox{}_{\scriptstyle \!\:\!R - 1 \,<\, |\,x\,| \,<\,R}}
\hspace{-1.675cm}
|\,u\,|^{q - 2} \;\!u \:
\bigl<\,A(x,\tau,u) \,\nabla u,
\;\! \nabla \zeta_{\mbox{}_{R}}(x)
\,\bigr>
\: dx \, d\tau
} $, \\
\mbox{} \vspace{+0.400cm} \\
para todo $ R > 1 $.
Fazendo
$ R \rightarrow \infty $,
obt\'em-se,
ent\~ao, \\
\mbox{} \vspace{-0.100cm} \\
\mbox{} \hspace{+2.100cm}
$ {\displaystyle
\|\, u(\cdot,t) \,
\|_{\mbox{}_{\scriptstyle L^{q}(\mathbb{R}^{n})}}^{\:\!q}
\,\!+\;
\frac{1}{2} \;
q \,(q - 1) \!
\int_{\:\!\mbox{\footnotesize $t$}_0}
    ^{\,\mbox{\footnotesize $t$}}
\!\!\;\!
\mu(\tau)
\!
\int_{\scriptstyle \;\!\mathbb{R}^{n}}
\!\!
|\, u(x,\tau) \,|^{q - 2} \:
|\, \nabla u \,|^{2}
\: dx \, d\tau
} $ \\
\mbox{} \vspace{-0.500cm} \\
\mbox{} \hfill (2.8) \\
\mbox{} \vspace{-0.350cm} \\
\mbox{} \hspace{+2.500cm}
$ {\displaystyle
\leq\;
\|\, u(\cdot,t_0) \,
\|_{\mbox{}_{\scriptstyle L^{q}(\mathbb{R}^{n})}}^{\:\!q}
\,\!+\;
\frac{1}{2} \;
q \,(q - 1) \!
\int_{\:\!\mbox{\footnotesize $t$}_0}
    ^{\,\mbox{\footnotesize $t$}}
\!
\frac{\,B(\tau)^{2}}{\mu(\tau)}
\!\;\!
\int_{\scriptstyle \;\!\mathbb{R}^{n}}
\!\!
|\, u \,|^{q + 2\;\! \kappa}
\: dx \, d\tau
} $, \\
\mbox{} \vspace{+0.400cm} \\
de onde
(2.5) pode ser
derivado
sem dificuldade.
De fato,
(2.8)
produz \\
\mbox{} \vspace{-0.600cm} \\
\begin{equation}
\notag
\int_{\:\!\mbox{\footnotesize $t$}_0}
    ^{\,\mbox{\footnotesize $t$}}
\!
\mu(\tau)
\!
\int_{\scriptstyle \;\!\mathbb{R}^{n}}
\!\!\;\!
|\, u(x,\tau) \,|^{q - 2} \,
|\, \nabla u \,|^{2}
\: dx \, d\tau
\;\leq\;\!
\int_{\:\!\mbox{\footnotesize $t$}_0}
    ^{\,\mbox{\footnotesize $t$}}
\!
\frac{\,B(\tau)^{2}}{\mu(\tau)}
\!\;\!
\int_{\scriptstyle \;\!\mathbb{R}^{n}}
\!\!\;\!
|\, u \,|^{q + 2\;\! \kappa}
\: dx \, d\tau,
\end{equation}
de modo que,
fazendo
$\:\!t_0 \!\;\!\rightarrow\;\!0 $,
tem-se \\
\mbox{} \vspace{-0.100cm} \\
\mbox{} \hspace{+0.750cm}
$ {\displaystyle
\int_{\:\!0}
    ^{\,\mbox{\footnotesize $t$}}
\!
\mu(\tau)
\!
\int_{\scriptstyle \;\!\mathbb{R}^{n}}
\!\!\;\!
|\, u \,|^{q - 2} \,
|\, \nabla u \,|^{2}
\: dx \, d\tau
\;\leq\;\!
\int_{\:\!0}
    ^{\,\mbox{\footnotesize $t$}}
\!\!\;\!
\frac{\,B(\tau)^{2}}{\mu(\tau)}
\!\;\!
\int_{\scriptstyle \;\!\mathbb{R}^{n}}
\!\!\;\!
|\, u(x,\tau) \,|^{q + 2\;\! \kappa}
\: dx \, d\tau
} $ \\
\mbox{} \vspace{+0.200cm} \\
\mbox{} \hspace{+6.150cm}
$ {\displaystyle
\leq\;\:\!
\mathbb{U}_{\infty}(0\:\!;\:\!t)^{\:\!2 \;\!\kappa}
\!\!
\int_{\:\!0}
    ^{\,\mbox{\footnotesize $t$}}
\!\!\;\!
\frac{\,B(\tau)^{2}}{\mu(\tau)}
\!\;\!
\int_{\scriptstyle \;\!\mathbb{R}^{n}}
\!\!
|\, u(x,\tau) \,|^{q}
\: dx \, d\tau
} $ \\
\mbox{} \vspace{+0.200cm} \\
\mbox{} \hspace{+6.150cm}
$ {\displaystyle
\leq\;\:\!
\mathbb{U}_{q}(0\:\!;\:\!t)^{\:\!q} \;
\mathbb{U}_{\infty}(0\:\!;\:\!t)^{\:\!2 \;\!\kappa}
\!\!
\int_{\:\!0}
    ^{\,\mbox{\footnotesize $t$}}
\!\!\;\!
\frac{\,B(\tau)^{2}}{\mu(\tau)}
\; d\tau
} $, \\
\mbox{} \vspace{+0.200cm} \\
para cada
$ \;\!0 < t < T_{\!\;\!\ast} \!\;\!$,
o que
conclui a prova de (2.5)
se $ \;\!q > 2 $.
Considerando,
agora,
$ q = 2 $,
pode-se
obter (2.5)
de modo inteiramente an\'alogo,
com a diferen\c ca de se
multiplicar desta vez
a equa\c c\~ao
(2.1$a$)
por
$ \;\! u(x,t) \: \zeta_{\mbox{}_{R}}(x) $,
em vez de
$ \;\! \Phi_{\delta}^{\prime}(u(x,t)) \: \zeta_{\mbox{}_{R}}(x) $
como feito antes.
Integrando-se, ent\~ao,
o resultado em
$\;\![\,t_0, \;\!t\;\!] $,
para $ 0 < t_0 < t $ dado,
obt\'em-se,
seguindo marcha
inteiramente similar
a (2.6)$\;\!$-$\;\!$(2.8)
acima, \\
\mbox{} \vspace{-0.600cm} \\
\begin{equation}
\tag{2.$8^{\prime}$}
\|\, u(\cdot,t) \,
\|_{\mbox{}_{\scriptstyle L^{2}(\mathbb{R}^{n})}}^{\:\!2}
\!\;\!+
\int_{\:\!\mbox{\footnotesize $t$}_0}
    ^{\,\mbox{\footnotesize $t$}}
\!\!\;\!
\mu(\tau)
\!
\int_{\scriptstyle \;\!\mathbb{R}^{n}}
\!\!
|\, \nabla u \,|^{2}
\: dx \, d\tau
\;\leq\;
\|\, u(\cdot,t_0) \,
\|_{\mbox{}_{\scriptstyle L^{2}(\mathbb{R}^{n})}}^{\:\!2}
\!\;\!+
\int_{\:\!\mbox{\footnotesize $t$}_0}
    ^{\,\mbox{\footnotesize $t$}}
\!
\frac{\,B(\tau)^{2}}{\mu(\tau)}
\!\;\!
\int_{\scriptstyle \;\!\mathbb{R}^{n}}
\!\!
|\, u \,|^{2 + 2\;\! \kappa}
\: dx \, d\tau
\end{equation}
\mbox{} \vspace{-0.525cm} \\
(sendo $ q = 2 $, $ p_{\mbox{}_{0}} \!\:\!= 1 $),
de onde resulta,
como antes, \\
\mbox{} \vspace{-0.025cm} \\
\mbox{} \hspace{+2.025cm}
$ {\displaystyle
\int_{\:\!0}
    ^{\,\mbox{\footnotesize $t$}}
\!
\mu(\tau)
\!
\int_{\scriptstyle \;\!\mathbb{R}^{n}}
\!\!\;\!
|\, \nabla u \,|^{2}
\: dx \, d\tau
\;\;\!\leq\;\!
\int_{\:\!0}
    ^{\,\mbox{\footnotesize $t$}}
\!\!\;\!
\frac{\,B(\tau)^{2}}{\mu(\tau)}
\!\;\!
\int_{\scriptstyle \;\!\mathbb{R}^{n}}
\!\!\;\!
|\, u(x,\tau) \,|^{2 + 2\;\! \kappa}
\: dx \, d\tau
} $ \\
\mbox{} \vspace{+0.200cm} \\
\mbox{} \hspace{+6.300cm}
$ {\displaystyle
\leq\;\:\!
\mathbb{U}_{\scriptscriptstyle 2}(0\:\!;\:\!t)^{\:\!2} \;
\mathbb{U}_{\infty}(0\:\!;\:\!t)^{\:\!2 \;\!\kappa}
\!\!
\int_{\:\!0}
    ^{\,\mbox{\footnotesize $t$}}
\!\!\;\!
\frac{\,B(\tau)^{2}}{\mu(\tau)}
\; d\tau
} $, \\
\mbox{} \vspace{+0.100cm} \\
o que estabelece (2.5)
no caso $ \:\!q = 2 $.
Isso conclui
a demonstra\c c\~ao
da Proposi\c c\~ao 2.2.
}
\mbox{} \hfill $\Box$ \\
%
\nl
Observe-se que,
da Proposi\c c\~ao 2.2,
segue imediatamente
que,
para todo
$ \;\!0 < t < \mbox{\small $T$}_{\!\;\!\ast} \!\;\!$: \\
\mbox{} \vspace{-0.525cm} \\
\begin{equation}
\tag{2.9$a$}
\int_{0}^{\;\!\mbox{\footnotesize $t$}}
\!\!\;\!
\int_{\mathbb{R}^{n}}
\!\!\;\!
\bigl<\, A(x,\tau,u) \, \nabla u, \;\! \nabla u \,\bigr>
\:dx\, d\tau
\;<\; \infty
\end{equation}
\mbox{} \vspace{-0.025cm} \\
se $ \;\!p_{\mbox{}_{0}} \!\:\!= 1 $
($ q = 2 $),
e \\
\mbox{} \vspace{-0.550cm} \\
\begin{equation}
\tag{2.9$b$}
\int_{0}^{\;\!\mbox{\footnotesize $t$}}
\!\!\;\!
\int_{\mathbb{R}^{n}}
\!\!\;\!
|\, u(x,\tau)\,|^{\:\!q \;\!-\;\! 2}
\,
\bigl<\, A(x,\tau,u) \, \nabla u, \;\! \nabla u \,\bigr>
\:dx\, d\tau
\;<\; \infty
\end{equation}
\mbox{} \vspace{+0.025cm} \\
sendo
$ \:\! q \;\!\geq\;\! p_{\mbox{}_{0}} \!\:\!+ 1 $,
$ \:\!q > 2 $.
Outra consequ\^encia importante
de (2.5) acima
\'e dada na Proposi\c c\~ao 2.3
a seguir,
que representa o ponto
de partida
para a an\'alise
na Se\c c\~ao~3
estabelecendo
o resultado principal
({\sc Teorema B})
anunciado na Se\c c\~ao~1.

\mbox{} \vspace{-1.000cm} \\
%
%
%
%
%
{\bf Proposi\c c\~ao 2.3.}
\textit{%
Sendo
$ {\displaystyle
\;\!
u(\cdot,t)
\in
%
%
L^{\infty}_{\tt loc}(\:\!
[\,0, \:\!\mbox{\small $T$}_{\!\;\!\ast}),
\:\!
L^{\infty}(\mathbb{R}^{n})
\,\!)
} $
so\-lu\-\c c\~ao do
problema $\;\!(2.1) \!\;\!$
em um dado intervalo
$ \:\![\,0, \:\!\mbox{\small $T$}_{\!\;\!\ast}) $,
tem-se,
para cada
$ \;\! q \;\!\geq\;\! p_{\mbox{}_{0}}\!\:\!+ 1 \!\!\:\!:$ \\
}
\mbox{} \vspace{-0.000cm} \\
\mbox{} \hspace{+1.000cm}
$ {\displaystyle
\frac{d}{d \:\!t} \;
\|\, u(\cdot,t) \,\|_{\mbox{}_{\scriptstyle L^{q}(\mathbb{R}^{n})}}^{\:\!q}
\:\!+\;
q \, (q - 1)
\!\!\;\!
\int_{\mathbb{R}^{n}}
\!\!\;\!
|\, u(x,t) \,|^{\:\!q - 2}
\:
\bigl<\, A(x,t,u) \, \nabla u, \;\! \nabla u \,\bigr>
\:dx
} $ \\
\mbox{} \vspace{-0.450cm} \\
\mbox{} \hfill (2.10) \\
\mbox{} \vspace{-0.600cm} \\
\mbox{} \hspace{+1.600cm}
$ {\displaystyle
=\;
q \, (q - 1)
\!\!\;\!
\int_{\mathbb{R}^{n}}
\!\!\;\!
|\, u(x,t) \,|^{\:\!q - 2 + \kappa}
\;\!
u(x,t) \:
\bigl<\, \mbox{\boldmath $b$}(x,t,u) \;\!-\;\!
\mbox{\boldmath $\beta$}(t), \;\! \nabla u \,\bigr>
\:dx
} $ \\
\mbox{} \vspace{+0.100cm} \\
\textit{%
para todo
$ \;\! t \in (\,0, \;\!\mbox{\small $T$}_{\!\;\!\ast})
\,\mbox{\footnotesize $\setminus$}\, E_{q} $,
sendo
$ E_{q} \!\;\!\subset\!\;\! (\;\!0, \infty) $
de medida zero
e $ \mbox{\boldmath $\beta$}(t) $
dada em $(2.7b)$.
}
\mbox{} \vspace{-0.500cm} \\
%
%
%
\nl
{\small
{\bf Prova:}
Na nota\c c\~ao da prova anterior,
multiplicando-se
(2.1$a$)
por
$ u(x,t) \, \zeta_{\mbox{}_{R}}(x) $
se $ q = 2 $,
e por
$ \;\! \Phi_{\delta}^{\prime}(u(x,t)) \:
\zeta_{\mbox{}_{R}}(x) \;\!$
se $ q > 2 $,
e integrando-se o resultado
em $ \:\![\,t_0, \;\!t\;\!] $,
obt\'em-se,
fazendo
$ \delta \rightarrow 0 $,
$ t_0 \!\;\!\rightarrow 0 $
e $ R \rightarrow \infty $,
por (2.5), (2.6) e (2.7), \\
\mbox{} \vspace{-0.550cm} \\
\begin{equation}
\notag
\|\, u(\cdot,t) \,
\|_{\mbox{}_{\scriptstyle L^{q}(\mathbb{R}^{n})}}^{\:\!q}
\:\!+\;
q \, (q - 1)
\!\!\;\!
\int_{\:\!0}^{\;\!\mbox{\footnotesize $t$}}
\!\!\;\!
\int_{\mathbb{R}^{n}}
\!\!\;\!
|\, u \,|^{\:\!q - 2}
\:
\bigl<\, A(x,\tau,u) \, \nabla u, \;\! \nabla u \,\bigr>
\;dx
\, d\tau
\;=
\end{equation}
\mbox{} \vspace{-0.775cm} \\
\begin{equation}
\notag
=\;
\|\, u_0 \;\!
\|_{\mbox{}_{\scriptstyle L^{q}(\mathbb{R}^{n})}}^{\:\!q}
\:\!+\;
q \, (q - 1)
\!\!\;\!
\int_{\:\!0}^{\;\!\mbox{\footnotesize $t$}}
\!\!\;\!
\int_{\mathbb{R}^{n}}
\!\!\;\!
|\, u(x,\tau) \,|^{\:\!q - 2 + \kappa}
\;\!
u(x,\tau) \:
\bigl<\, \mbox{\boldmath $b$}(x,\tau,u) \;\!-\;\!
\mbox{\boldmath $\beta$}(\tau), \;\! \nabla u \,\bigr>
\;dx
\,d\tau
\end{equation}
\mbox{} \vspace{-0.050cm} \\
se
$ \;\!q > 2 $,
para todo
$ \;\! 0 < t < T_{\!\;\!\ast} \!\;\!$.
Sendo $\;\!q = 2 $,
o resultado
correspondente
obtido
\'e \\
\mbox{} \vspace{-0.550cm} \\
\begin{equation}
\notag
\|\, u(\cdot,t) \,
\|_{\mbox{}_{\scriptstyle L^{2}(\mathbb{R}^{n})}}^{\:\!2}
+\;
2
\!\!\;\!
\int_{\:\!0}^{\;\!\mbox{\footnotesize $t$}}
\!\!\;\!
\int_{\mathbb{R}^{n}}
\!
\bigl<\, A(x,\tau,u) \, \nabla u, \;\! \nabla u \,\bigr>
\;dx
\, d\tau
\;=
\end{equation}
\mbox{} \vspace{-0.800cm} \\
\begin{equation}
\notag
=\;
\|\, u_0 \;\!
\|_{\mbox{}_{\scriptstyle L^{2}(\mathbb{R}^{n})}}^{\:\!2}
\:\!+\;
2
\!\!\;\!
\int_{\:\!0}^{\;\!\mbox{\footnotesize $t$}}
\!\!\;\!
\int_{\mathbb{R}^{n}}
\!\!\;\!
|\, u(x,\tau) \,|^{\:\!\kappa}
\;\!
u(x,\tau) \:
\bigl<\, \mbox{\boldmath $b$}(x,\tau,u) \;\!-\;\!
\mbox{\boldmath $\beta$}(\tau), \;\! \nabla u \,\bigr>
\;dx
\,d\tau,
\end{equation}
\mbox{} \vspace{-0.050cm} \\
para todo
$ \;\! 0 < t < T_{\!\;\!\ast} \!\;\!$.
Nas express\~oes acima,
todas as integrais
s\~ao bem definidas e finitas \linebreak
\mbox{[}$\,$por (2.5), Proposi\c c\~ao 2.2, e (2.9)$\,$\mbox{]},
envolvendo fun\c c\~oes integr\'aveis
(no sentido de Lebesgue)
nas regi\~oes indicadas.
$\!$($\,\!$Para o termo
$ {\displaystyle
\:\!
z(x,t)
\:\!=\:\!
|\, u(x,t) \,|^{\:\!q - 2 + \kappa}
\,\!
u(x,t) \,
\bigl<\, \mbox{\boldmath $b$}(x,t,u) \;\!-\;\!
\mbox{\boldmath $\beta$}(t), \,\! \nabla u \,\bigr>
} $,
por exemplo,
tem-se,
por (1.6) e (2.7$b$): \\
\mbox{} \vspace{-0.150cm} \\
\mbox{} \hspace{+1.200cm}
$ {\displaystyle
|\, z(x,t) \,|
\;\leq\;
|\, u(x,t) \,|^{\:\!q - 1 + \kappa}
\;\!
|\, \mbox{\boldmath $b$}(x,t,u) \;\!-\;\!
\mbox{\boldmath $\beta$}(t) \,|_{\mbox{}_{2}}
\;\!
|\, \nabla u \,|_{\mbox{}_{2}}
} $ \\
\mbox{} \vspace{-0.200cm} \\
\mbox{} \hspace{+2.800cm}
$ {\displaystyle
\leq\;
B(t) \:
|\, u(x,t) \,|^{\:\!q - 1 + \kappa}
\;\!
|\, \nabla u \,|_{\mbox{}_{2}}
} $ \\
\mbox{} \vspace{-0.150cm} \\
\mbox{} \hspace{+2.800cm}
$ {\displaystyle
\leq\;
\frac{\,B(t)^{2}}{\mu(t)} \;
|\, u(x,t) \,|^{\:\!q + 2 \;\!\kappa}
\,+\;
\mu(t) \:
|\, u(x,t) \,|^{\:\!q - 2}
\;\!
|\, \nabla u \,|_{\mbox{}_{2}}^{\:\!2}
} $ \\
\mbox{} \vspace{-0.150cm} \\
\mbox{} \hspace{+2.800cm}
$ {\displaystyle
\leq\;\:\!
\mu(t) \:
\mathbb{B}_{\mu}\!\;\!(0\:\!;\:\!t)^{2}
\:
\mathbb{U}_{\infty}\!\;\!(0\:\!;\:\!t)^{2\;\!\kappa}
\,
|\, u(x,t) \,|^{\:\!q}
\,+\:
\mu(t) \:
|\, u(x,t) \,|^{\:\!q - 2}
\,
|\, \nabla u \,|_{\mbox{}_{2}}^{\:\!2}
\!\:\!
} $, \\
\mbox{} \vspace{-0.100cm} \\
de modo que
$ {\displaystyle
z \in L^{1}(\mathbb{R}^{n}
\!\!\;\!\times\!\;\! [\,0, \;\!t\;\!])
} $
para cada
$ \;\! 0 < t < T_{\!\;\!\ast} \!\;\!$,
como afirmado.)
O resultado (2.10) \linebreak
segue, ent\~ao,
aplicando o teorema de diferencia\c c\~ao de Lebesgue,
para cada
$ \;\!q \:\!\geq\;\! p_{\mbox{}_{0}}\!\:\!+ 1 $.
}
\hfill $\Box$ \\
%
\nl
%
%
%
%
{\bf Observa\c c\~ao 3.1.}
Na verdade,
por (2.1),
tem-se (2.10)
v\'alida
para todo
$ 0 < t < \mbox{\small $T$}_{\!\;\!\ast} \!\;\!$.

%
%
%
\mbox{} \vspace{-1.500cm} \\

{\bf 3. Prova de (1.10)} \\

Nesta se\c c\~ao,
vamos obter o resultado principal
({\small \sc Teorema B}),
reescrito na forma mais geral
a seguir
para as solu\c c\~oes
do problema (2.1),
sob as hip\'oteses de trabalho
(1.2) e (1.3)
descritas na Se\c c\~ao 1. \\
\mbox{} \vspace{-0.250cm} \\
%
%
%
%
%
%
\mbox{} \hspace{-0.800cm}
\fbox{%
\begin{minipage}[t]{16.000cm}
\mbox{} \vspace{-0.450cm} \\
\mbox{} \hspace{+0.300cm}
\begin{minipage}[t]{15.000cm}
\mbox{} \vspace{+0.100cm} \\
{\bf Teorema 3.1.}
\textit{%
$\!\;\!$Sejam
$ \;\! \kappa \geq 0 $,
$ \;\!p_{\mbox{}_{0}} \!\;\!\geq 1 $
$\!\;\!$e
$\;\!u(\cdot,t) $,
$ 0 \leq t < \mbox{\small $T$}_{\!\;\!\ast}\!\;\!$,
solu\c c\~ao
do problema
$\,\!(2.1) $. \linebreak
%
%
Para cada $\;\! p \geq p_{\mbox{}_{0}} \!\,\!$
satisfazendo
$ \;\!p > n \:\!\kappa $,
tem-se \\
}
\mbox{} \vspace{-0.650cm} \\
\begin{equation}
\tag{3.1}
\mbox{} \!\!
\|\, u(\cdot,t) \,
\|_{\mbox{}_{\scriptstyle L^{\infty}(\mathbb{R}^{n})}}
\;\!\leq\,
K\!\;\!(n,\:\!\kappa, \:\!p) \,\cdot\,
\max\;
\Bigl\{\:
\|\, u_0 \;\!
\|_{\mbox{}_{\scriptstyle L^{\infty}(\mathbb{R}^{n})}};
\;\:\!
\mathbb{B}_{\mu}\!\;\!(0; t)^{\mbox{}^{\scriptstyle \!\;\!
\frac{\scriptstyle n}{{\scriptstyle p} \,-\, {\scriptstyle n \;\!\kappa}} }}
\;\!
\mathbb{U}_{p}\!\;\!(0; t)^{\mbox{}^{\scriptstyle \!\;\!
\frac{\scriptstyle p}{{\scriptstyle p} \,-\, {\scriptstyle n \;\!\kappa}} }}
\;\!\Bigr\}
\end{equation}
\mbox{} \vspace{-0.200cm} \\
\textit{%
para todo
$\:\! 0 \leq t < \mbox{\small $T$}_{\!\;\!\ast} \!\;\! $,
$\:\!$sendo
$ {\displaystyle
\:\!
K\!\;\!(n,\:\!\kappa, \:\!p)
\:\!=\,
\{\;\!2\;\!p \;\!\}^{\mbox{}^{\scriptstyle \!\:\!
\frac{\scriptstyle n}{{\scriptstyle p} \,-\, {\scriptstyle n \:\!\kappa}} }}
\!\!
} $
e
$ {\displaystyle
\,
\mathbb{B}_{\mu}\!\;\!(0\:\!;\:\!t),
\,
\mathbb{U}_{p}\!\;\!(0\:\!;\:\!t)
} $
definidas em
$ \:\!(1.7) $, $(1.8)$
acima. \\
}
\end{minipage}
\end{minipage}
}
%
%
\nl
\mbox{} \vspace{-0.050cm} \\
A prova do Teorema 3.1 ser\'a feita
a partir da Proposi\c c\~ao 2.3,
com o aux\'\i lio
de v\'arios resultados auxiliares
apresentados nos lemas abaixo.
Tamb\'em ser\~ao necess\'arias
diversas desigualdades de Nirenberg-Gagliardo,
incluindo a desigualdade de Nash
\cite{Nash1958} \\
\mbox{} \vspace{-0.550cm} \\
\begin{equation}
\label{Nash}
\tag{3.2}
\|\:\mbox{v}\:
\|_{\mbox{}_{\scriptstyle L^{2}(\mathbb{R}^{n})}}
\;\!\leq\;\!\;\!
K\!\;\!(n) \:
\|\:\mbox{v}\:
\|_{\mbox{}_{\scriptstyle L^{1}(\mathbb{R}^{n})}}
^{\mbox{}^{\scriptstyle \frac{2}{{\scriptstyle n} \;\!+\;\! 2} }}
\;\!
\|\, \nabla \mbox{v}\:
\|_{\mbox{}_{\scriptstyle L^{2}(\mathbb{R}^{n})}}
^{\mbox{}^{\scriptstyle \frac{\scriptstyle n}{{\scriptstyle n} \;\!+\;\! 2} }}
\qquad
\forall \;\,
\mbox{v} \in
L^{1}(\mathbb{R}^{n})
\cap
\dot{H}^{1}(\mathbb{R}^{n}),
\end{equation}
\mbox{} \vspace{-0.200cm} \\
para certa constante
$ \:\!K\!\;\!(n) < 1 $
(ver \cite{CarlenLoss1993}
para a determina\c c\~ao de seu valor optimal). \linebreak
Esta desigualdade ser\'a importante
mais adiante para se estimar
$ {\displaystyle
\;\!
\|\, u(\cdot,t) \,
\|_{\scriptstyle L^{q}(\mathbb{R}^{n})}
} $
em termos de
$ {\displaystyle
\;\!
\|\, u(\cdot,t) \,
\|_{\scriptstyle L^{q/2}(\mathbb{R}^{n})}
\!\;\!
} $,
dado $ \:\!q \:\!\geq\:\! p_{\mbox{}_{0}}/2 $
(Lema 3.2),
na vers\~ao
desenvolvida
pelo \linebreak
autor
para o cl\'assico m\'etodo
$ L^{p}$-$\,\!L^{q}$
em ordem a se aplicar
\`a investiga\c c\~ao de
pro\-blemas da forma (2.1) acima.
Em particular,
dado
$ \;\!q \;\!\geq\;\! 2 \;\!p_{\mbox{}_{0}}\!\;\! $,
resulta desde j\'a conve- \linebreak
niente
introduzir a fun\c c\~ao auxiliar
$ {\displaystyle
\;\!
v^{[\,q\,]}(\cdot,t)
\in
L^{1}(\mathbb{R}^{n})
\cap
L^{\infty}(\mathbb{R}^{n})
\:\!
} $
definida por \\
\mbox{} \vspace{-0.450cm} \\
\begin{equation}
\label{vq}
\tag{3.3}
v^{[\,q\,]}(x,t)
\,:=\:
\left\{
\begin{array}{ll}
\!\;\,u(x,t) \; &
\mbox{se }\, q = 2, \\
\mbox{} \vspace{-0.200cm} \\
\! |\, u(x,t) \,|^{\mbox{}^{\scriptstyle q/2}} \; &
\mbox{se }\, q > 2.
\end{array}
\right.
\end{equation}
\mbox{} \vspace{+0.100cm} \\
Em termos de
$ v^{[\,q\,]}(\cdot,t) $,
tem-se \\
\mbox{} \vspace{-0.550cm} \\
\begin{equation}
\tag{3.4$a$}
\|\, u(\cdot,t) \,
\|_{\mbox{}_{\scriptstyle L^{q}(\mathbb{R}^{n})}}^{\;\!q}
\:\!=\;
\|\, v^{[\,q\,]}(\cdot,t) \,
\|_{\mbox{}_{\scriptstyle L^{2}(\mathbb{R}^{n})}}^{\;\!2}
\!\;\!,
\end{equation}
\mbox{} \vspace{-0.750cm} \\
\begin{equation}
\tag{3.4$b$}
\|\, u(\cdot,t) \,
\|_{\mbox{}_{\scriptstyle L^{q/2}(\mathbb{R}^{n})}}^{\;\!q/2}
=\;
\|\, v^{[\,q\,]}(\cdot,t) \,
\|_{\mbox{}_{\scriptstyle L^{1}(\mathbb{R}^{n})}}
\!\;\!,
\end{equation}
\mbox{} \vspace{-0.550cm} \\
e tamb\'em \\
\mbox{} \vspace{-0.750cm} \\
\begin{equation}
\tag{3.5}
\int_{\mathbb{R}^{n}}
\!\!\;\!
|\, u(x,t) \,|^{\:\!q - 2}
\:
|\, \nabla u\;\!(x,t) \,|_{\mbox{}_{2}}^{\:\!2}
\: dx
\;=\;
\mbox{\small $ {\displaystyle \frac{4}{\,q^{2}} }$}
\:
\|\, \nabla v^{[\,q\,]}(\cdot,t) \,
\|_{\mbox{}_{\scriptstyle L^{2}(\mathbb{R}^{n})}}^{\:\!2}
\!.
\end{equation}
\mbox{} \vspace{-0.400cm} \\
\nl
%
%
%
%
{\bf Lema 3.1.}
\textit{%
Seja
$\;\! q \:\!\geq\:\!2 \;\! p_{\mbox{}_{0}} \!\;\!$.
Sendo
$ {\displaystyle
\;\!
v^{[\,q\,]}(\cdot,t)
\in
L^{1}(\mathbb{R}^{n})
\cap
L^{\infty}(\mathbb{R}^{n})
} $
dada em
$\;\!(3.3) $ acima,
tem-se \\
}
\mbox{} \vspace{-0.350cm} \\
\mbox{} \hspace{+2.200cm}
$ {\displaystyle
\frac{d}{d \:\!t} \:
\|\, v^{[\,q\,]}(\cdot,t) \,
\|_{\mbox{}_{\scriptstyle L^{2}(\mathbb{R}^{n})}}^{\:\!2}
\!\:\!+\;
\mbox{\small $ {\displaystyle
4 \;\! \Bigl(\:\! 1 - \frac{1}{q} \,\Bigr)
\, \mu(t)
} $}
\:
\|\, \nabla v^{[\,q\,]}(\cdot,t) \,
\|_{\mbox{}_{\scriptstyle L^{2}(\mathbb{R}^{n})}}^{\:\!2}
} $ \\
\mbox{} \vspace{-0.525cm} \\
\mbox{} \hfill (3.6) \\
\mbox{} \vspace{-0.550cm} \\
\mbox{} \hspace{+2.000cm}
$ {\displaystyle
\leq\;
\mbox{\small $ {\displaystyle
2 \, q \;\!
\Bigl(\;\! 1 - \frac{1}{q} \,\Bigr) \,
B(t)
} $}
\:
\|\, v^{[\,q\,]}(\cdot,t) \,
\|_{\scriptstyle L^{2 \;\!+\;\!\frac{4\:\!\kappa}{q}}\!(\mathbb{R}^{n})}
  ^{\mbox{}^{\scriptstyle 1 \;\!+\;\! \frac{2\:\!\kappa}{q} }}
\;\!
\|\, \nabla v^{[\,q\,]}(\cdot,t) \,
\|_{\mbox{}_{\scriptstyle L^{2}(\mathbb{R}^{n})}}
} $ \\
\mbox{} \vspace{+0.100cm} \\
\textit{%
para todo
$ \;\! t \in (\,0, \;\!\mbox{\small $T$}_{\!\;\!\ast})
\,\mbox{\footnotesize $\setminus$}\, E_{q} $,
sendo
$ E_{q} \!\;\!\subset\!\;\! (\;\!0, \infty) $,
$ |\,E_{q} \;\!|_{\mbox{}_{1}} \!\;\!=\;\! 0 $,
dado
na Proposi\c c\~ao~$2.3$. \linebreak
}
%
%
\nl
{\small
{\bf Prova:}
De (2.10),
tem-se,
por (1.2) e (1.3), (1.6), (2.7$b$),
para
$ t \in (\,0, \;\!\mbox{\small $T$}_{\!\;\!\ast})
\,\mbox{\footnotesize $\setminus$}\, E_{q} $: \\
\mbox{} \vspace{+0.000cm} \\
\mbox{} \hspace{+0.700cm}
$ {\displaystyle
\frac{d}{d \:\!t} \;
\|\, u(\cdot,t) \,\|_{\mbox{}_{\scriptstyle L^{q}(\mathbb{R}^{n})}}^{\:\!q}
\:\!+\;
q \, (q - 1)
\:
\mu(t)
\!\!\;\!
\int_{\mathbb{R}^{n}}
\!\!\;\!
|\, u(x,t) \,|^{\:\!q - 2}
\:
|\, \nabla u \,|_{\mbox{}_{2}}^{\:\!2}
\:dx
} $ \\
\mbox{} \vspace{+0.050cm} \\
\mbox{} \hspace{+2.350cm}
$ {\displaystyle
\leq\;
q \,(q - 1)\,
B(t) \!
\int_{\mathbb{R}^{n}}
\!\!\;\!
|\, u(x,t) \,|^{\:\!q - 1 + \kappa}
\:
|\, \nabla u \,|_{\mbox{}_{2}}
\:dx
} $
\mbox{} \hfill
\mbox{[}$\,$por (1.2), (1.6), (2.10)$\,$\mbox{]} \\
\mbox{} \vspace{-0.050cm} \\
\mbox{} \hfill
$ {\displaystyle
\leq\;
q \,(q - 1)\,
B(t) \;
\biggl\{\:\!
\int_{\mathbb{R}^{n}}
\!\!\;\!
|\, u(x,t) \,|^{\:\!q \;\!+\;\! 2\;\! \kappa}
\:dx
\;\!\biggr\}^{\!\!\;\!1/2}
\;\!
\biggl\{\:\!
\int_{\mathbb{R}^{n}}
\!\!\;\!
|\, u(x,t) \,|^{\:\!q \;\!-\;\! 2}
\,
|\, \nabla u \,|_{\mbox{}_{2}}^{\:\!2}
\:dx
\;\!\biggr\}^{\!\!\;\!1/2}
\!\!\!\!,
} $ \\
\mbox{} \vspace{+0.100cm} \\
que,
em termos da
fun\c c\~ao
$ \:\!v^{[\,q\,]}(\cdot,t) \:\!$
definida em (3.3),
equivale a (3.6),
como afirmado.
}
\hfill $\Box$ \\
%
%
\nl
De (3.6),
resulta a importante estimativa (3.8) abaixo,
para $ \:\!q \:\!\geq\:\! 2 \;\!p_{\mbox{}_{0}} \!\;\!$
com $ \:\! q > 2 \;\!n\;\!\kappa $,
usando-se as desigualdades
de Nirenberg-Gagliardo \\
%
%
\mbox{} \vspace{-0.625cm} \\
\begin{equation}
\tag{3.7}
\|\:\mbox{v}\:
\|_{\mbox{}_{\scriptstyle L^{r}(\mathbb{R}^{n})}}
\:\!\leq\;
\|\:\mbox{v}\:
\|_{\mbox{}_{\scriptstyle L^{1}(\mathbb{R}^{n})}}
  ^{1 - \theta}
\|\, \nabla \mbox{v}\:
\|_{\mbox{}_{\scriptstyle L^{2}(\mathbb{R}^{n})}}
  ^{\:\!\theta}
\!\;\!,
\qquad
\theta \;=\;
\frac{\mbox{\small 1} \:\mbox{\small $-$}\: \mbox{\small 1/}r}
{\mbox{\small 1/2} \;\mbox{\small +}\:\mbox{\small 1/}n}
\end{equation}
\mbox{} \vspace{-0.125cm} \\
(ver e.g.$\;$\cite{Friedman1969}, p.$\;$24),
nos casos
$ \;\!\mbox{\small $2$} \leq r < \mbox{\small $2$} + \mbox{\small $2$}/n $.
Aplicando-se~(3.7)
para
estimar
$ {\displaystyle
\|\, v^{[\,q\,]}(\cdot,t) \,
\|_{\scriptstyle L^{2 \;\!+\;\!\frac{4\:\!\kappa}{q}}\!(\mathbb{R}^{n})}
\!
} $
no lado direito
da desigualdade (3.6), Lema~3.1,
%
%
obt\'em-se \\
\mbox{} \vspace{+0.150cm} \\
\mbox{} \hspace{+2.200cm}
$ {\displaystyle
\frac{d}{d \:\!t} \:
\|\, v^{[\,q\,]}(\cdot,t) \,
\|_{\mbox{}_{\scriptstyle L^{2}(\mathbb{R}^{n})}}^{\:\!2}
\!\:\!+\;
\mbox{\small $ {\displaystyle
4 \;\! \Bigl(\:\! 1 - \frac{1}{q} \,\Bigr)
\, \mu(t)
} $}
\:
\|\, \nabla v^{[\,q\,]}(\cdot,t) \,
\|_{\mbox{}_{\scriptstyle L^{2}(\mathbb{R}^{n})}}^{\:\!2}
} $ \\
\mbox{} \vspace{-0.550cm} \\
\mbox{} \hfill (3.8) \\
\mbox{} \vspace{-0.425cm} \\
\mbox{} \hspace{+0.900cm}
$ {\displaystyle
\leq\;
\mbox{\small $ {\displaystyle
2 \, q \;\!
\Bigl(\:\! 1 - \frac{1}{q} \,\Bigr)
\;\!
B(t)
} $}
\:
\|\, v^{[\,q\,]}(\cdot,t) \,
\|_{\mbox{}_{\scriptstyle L^{1}(\mathbb{R}^{n})}}
  ^{\mbox{}^{\scriptstyle \:\!
\frac{2}{{\scriptstyle n}\;\!+\;\!2} \:
\frac{{\scriptstyle q} \;\!-\;\!
({\scriptstyle n} - 2) \;\!
{\scriptstyle \kappa}}
{{\scriptstyle q}} }}
\|\, \nabla v^{[\,q\,]}(\cdot,t) \,
\|_{\mbox{}_{\scriptstyle L^{2}(\mathbb{R}^{n})}}
  ^{\mbox{}^{\scriptstyle \:\!
\frac{2}{{\scriptstyle n}\;\!+\;\!2} \:
\frac{({\scriptstyle n} + 1)\,{\scriptstyle q}
\,+\, 2 \;\!{\scriptstyle n \:\!\kappa}}
{{\scriptstyle q}} }}
\hspace{-1.625cm},
} $ \\
\mbox{} \vspace{+0.150cm} \\
para
$ \:\! q \:\!\geq\:\! 2 \;\!p_{\mbox{}_{0}} \!\:\!$
satisfazendo $ \;\! q > 2 \;\!n \;\!\kappa $.
(Esta \'ultima condi\c c\~ao
foi feita de modo a
(3.8) poder ser \'util:
ela torna
o expoente
do termo
$ {\displaystyle
\;\!
\|\, \nabla v^{[\,q\,]}(\cdot,t) \,
\|_{\mbox{}_{\scriptstyle L^{2}(\mathbb{R}^{n})}}
\!
} $
no lado direito da
express\~ao (3.8)
menor que 2
\mbox{[}$\,$que \'e o expoente
%
%
do mesmo termo no lado esquerdo$\,$\mbox{]}.
Posto de outra forma
(equivalente):
vamos deste ponto em diante
supor sempre
que se tenha $ q $
verificando
$ \;\!q \geq 2 \;\!p $,
com
$ p \:\!\geq\:\! p_{\mbox{}_{0}} \!\;\!$
dado
(fixo)
satisfazendo$\;\!$\footnote{%
%
%
A condi\c c\~ao
$ \;\! p > \:\! n \:\! \kappa \;\! $
imposta em (3.9) acima
n\~ao \'e resultado de limita\c c\~ao
do m\'etodo de an\'alise \linebreak
apresentado,
mas uma condi\c c\~ao {\em natural},
prevista por argumentos de escala
aplicados a (2.1).
%
%
}
%
%
%
\nl
\mbox{} \vspace{-0.750cm} \\
\begin{equation}
\tag{3.9}
p \;\!\geq\;\! p_{\mbox{}_{0}}
\!\;\!,
\quad
\;\,
p \;\!>\;\!
n \;\! \kappa.
\end{equation}
\mbox{} \vspace{-0.300cm} \\
Esta condi\c c\~ao sobre $q$
permite que se
prossiga a an\'alise
al\'em da estimativa (3.8),
como mostram os seguintes resultados. \\
\nl
%
%
%
%
{\bf Lema 3.2.}
\textit{%
Seja
$\;\! q \:\!\geq\:\!2 \;\! p$,
com
$ \;\!p $ dado em $\;\!(3.9)$.
Sendo
$ {\displaystyle
\;\!
v^{[\,q\,]}(\cdot,t)
\in
L^{1}(\mathbb{R}^{n})
\cap
L^{\infty}(\mathbb{R}^{n})
} $
definida em
$\;\!(3.3) $,
tem-se \\
}
\mbox{} \vspace{-0.000cm} \\
\mbox{} \hspace{+0.800cm}
$ {\displaystyle
\frac{d}{d \:\!t} \;
\|\, v^{[\,q\,]}(\cdot,t) \,
\|_{\mbox{}_{\scriptstyle L^{2}(\mathbb{R}^{n})}}^{\:\!2}
\!\:\!+\;
\mbox{\small $ {\displaystyle
\frac{4}{n + 2} \,
\Bigl(\:\! 1 - \frac{1}{q} \,\Bigr) \;\!
\Bigl(\:\! 1 - \frac{2 \;\!n \:\! \kappa}{q} \;\!\Bigr)
} $}
\,
\mbox{\small $\mu(t)$}
\,
\|\, \nabla v^{[\,q\,]}(\cdot,t) \,
\|_{\mbox{}_{\scriptstyle L^{2}(\mathbb{R}^{n})}}^{\:\!2}
} $ \\
\mbox{} \vspace{-0.400cm} \\
\mbox{} \hfill (3.10) \\
\mbox{} \vspace{-0.425cm} \\
\mbox{} \hspace{+0.150cm}
$ {\displaystyle
\leq\;
\mbox{\small $ {\displaystyle
\frac{4}{n + 2} \,
\Bigl(\:\! 1 - \frac{1}{q} \,\Bigr) \;\!
\Bigl(\:\! 1 - \frac{2 \;\!n \:\! \kappa}{q} \;\!\Bigr)
} $}
\,
\mbox{\small $ {\displaystyle
\Bigl(\;\! \frac{\;\!q\;\!}{2} \;\! \Bigr)
} $}^{\mbox{}^{\scriptstyle \!\!
\frac{({\scriptstyle n} + 2) \,{\scriptstyle q}}
{{\scriptstyle q} \;\!-\;\! 2 \;\!{\scriptstyle n \;\!\kappa}} }}
\hspace{-0.700cm}
\cdot
\hspace{+0.300cm}
\mbox{\small $\mu(t)$}
\;\!\;\!
\mbox{\small $ {\displaystyle
\Bigl(\, \frac{B(t)}{\mu(t)} \,\Bigr)
} $}^{\mbox{}^{\scriptstyle \!\!\!
\frac{({\scriptstyle n} + 2) \;\!{\scriptstyle q}}
{{\scriptstyle q} \;\!-\;\! 2 \;\!{\scriptstyle n \:\!\kappa}} }}
\hspace{-0.650cm}
\cdot
\hspace{+0.300cm}
\|\, v^{[\,q\,]}(\cdot,t) \,
\|_{\mbox{}_{\scriptstyle L^{1}(\mathbb{R}^{n})}}
  ^{\mbox{}^{\scriptstyle \:\!  2 \,\cdot\,
\left[\,
\frac{{\scriptstyle q} \;\!-\;\!
({\scriptstyle n} - 2) \;\!{\scriptstyle \kappa}}
{{\scriptstyle q} \;\!-\;\! 2 \;\!{\scriptstyle n \:\!\kappa}}
\,\right]}}
\hspace{-1.200cm},
} $ \\
\mbox{} \vspace{+0.200cm} \\
\textit{%
para todo
$ \;\! t \in (\,0, \;\!\mbox{\small $T$}_{\!\;\!\ast})
\,\mbox{\footnotesize $\setminus$}\, E_{q} $,
sendo
$ E_{q} \!\;\!\subset\!\;\! (\;\!0, \infty) $,
$ |\,E_{q} \;\!|_{\mbox{}_{1}} \!\;\!=\;\! 0 $,
dado
na Proposi\c c\~ao~$2.3$. \linebreak
}
%
%
\nl
{\small
{\bf Prova:}
Considerando o termo
no lado direito de (3.8),
tem-se \\
\mbox{} \vspace{-0.000cm} \\
\mbox{} \hspace{+0.250cm}
$ {\displaystyle
\mbox{\small $ {\displaystyle
2 \, q \;\!
\Bigl(\:\! 1 - \frac{1}{q} \,\Bigr)
\;\!
B(t)
} $}
\:
\|\, v^{[\,q\,]}(\cdot,t) \,
\|_{\mbox{}_{\scriptstyle L^{1}(\mathbb{R}^{n})}}
  ^{\mbox{}^{\scriptstyle \:\!
\frac{2}{{\scriptstyle n}\;\!+\;\!2} \:
\frac{{\scriptstyle q} \,-\,
({\scriptstyle n} - 2) \;\!
{\scriptstyle \kappa}}
{{\scriptstyle q}} }}
\|\, \nabla v^{[\,q\,]}(\cdot,t) \,
\|_{\mbox{}_{\scriptstyle L^{2}(\mathbb{R}^{n})}}
  ^{\mbox{}^{\scriptstyle \:\!
\frac{2}{{\scriptstyle n}\;\!+\;\!2} \:
\frac{({\scriptstyle n} + 1)\,{\scriptstyle q}
\,+\, 2 \;\!{\scriptstyle n \:\!\kappa}}
{{\scriptstyle q}} }}
\hspace{-1.625cm}
} $ \\
\mbox{} \vspace{+0.150cm} \\
\mbox{} \hspace{+2.500cm}
$ {\displaystyle
=\;
\mbox{\small $ {\displaystyle
4 \,
\Bigl(\:\! 1 - \frac{1}{q} \,\Bigr)
} $}
\;
\Bigl[\;\,
\mbox{\small $ {\displaystyle
\frac{q}{2} \:
B(t) \;
\mu(t)}$}^{\scriptstyle \!\!-\,
\frac{({\scriptstyle n} + 1)\;\!{\scriptstyle q}
\,+\, 2 \;\!{\scriptstyle n \:\!\kappa}}
{({\scriptstyle n} + 2)\;\!{\scriptstyle q}} }
\,
\|\, v^{[\,q\,]}(\cdot,t) \,
\|_{\mbox{}_{\scriptstyle L^{1}(\mathbb{R}^{n})}}
  ^{\mbox{}^{\scriptstyle \:\!
\frac{2}{{\scriptstyle n}\;\!+\;\!2} \:
\frac{{\scriptstyle q} \,-\,
({\scriptstyle n} - 2) \;\!
{\scriptstyle \kappa}}
{{\scriptstyle q}} }}
\;\Bigr]
\;\times
} $ \\
\mbox{} \vspace{-0.000cm} \\
\mbox{} \hspace{+5.700cm}
$ {\displaystyle
\times \;
\Bigl[\;\;\!
\mbox{\small $\mu(t)$}^{\mbox{}^{\scriptstyle \!\!
\frac{({\scriptstyle n} + 1)\;\!{\scriptstyle q}
\,+\, 2 \;\!{\scriptstyle n \:\!\kappa}}
{({\scriptstyle n} + 2)\;\!{\scriptstyle q}} }}
\,
\|\, \nabla v^{[\,q\,]}(\cdot,t) \,
\|_{\mbox{}_{\scriptstyle L^{2}(\mathbb{R}^{n})}}
  ^{\mbox{}^{\scriptstyle \:\!
\frac{2}{{\scriptstyle n}\;\!+\;\!2} \:
\frac{({\scriptstyle n} + 1)\,{\scriptstyle q}
\,+\, 2 \;\!{\scriptstyle n \:\!\kappa}}
{{\scriptstyle q}} }}
\;\Bigr]
} $ \\
\mbox{} \vspace{+0.050cm} \\
\mbox{} \hspace{+0.800cm}
$ {\displaystyle
\leq \;
\mbox{\small $ {\displaystyle
\frac{4}{n + 2} \,
\Bigl(\:\! 1 - \frac{1}{q} \,\Bigr) \;\!
\Bigl(\:\! 1 - \frac{2 \;\!n \:\! \kappa}{q} \;\!\Bigr)
} $}
\,
\mbox{\small $ {\displaystyle
\Bigl(\;\! \frac{\;\!q\;\!}{2} \;\! \Bigr)
} $}^{\mbox{}^{\scriptstyle \!\!
\frac{({\scriptstyle n} + 2) \,{\scriptstyle q}}
{{\scriptstyle q} \;\!-\;\! 2 \;\!{\scriptstyle n \;\!\kappa}} }}
\hspace{-0.700cm}
\cdot
\hspace{+0.300cm}
\mbox{\small $\mu(t)$}
\;\!\;\!
\mbox{\small $ {\displaystyle
\Bigl(\, \frac{B(t)}{\mu(t)} \,\Bigr)
} $}^{\mbox{}^{\scriptstyle \!\!\!
\frac{({\scriptstyle n} + 2) \;\!{\scriptstyle q}}
{{\scriptstyle q} \;\!-\;\! 2 \;\!{\scriptstyle n \:\!\kappa}} }}
\hspace{-0.650cm}
\cdot
\hspace{+0.300cm}
\|\, v^{[\,q\,]}(\cdot,t) \,
\|_{\mbox{}_{\scriptstyle L^{1}(\mathbb{R}^{n})}}
  ^{\mbox{}^{\scriptstyle \:\!  2 \,\cdot\,
\left[\,
\frac{{\scriptstyle q} \;\!-\;\!
({\scriptstyle n} - 2) \;\!{\scriptstyle \kappa}}
{{\scriptstyle q} \;\!-\;\! 2 \;\!{\scriptstyle n \:\!\kappa}}
\,\right]}}
} $ \\
\mbox{} \vspace{+0.050cm} \\
%
\mbox{} \hspace{+5.450cm}
$ {\displaystyle
+ \;\,
\mbox{\small $ {\displaystyle
4 \,
\Bigl(\:\! 1 - \frac{1}{q} \,\Bigr)
\;
\frac{(n + 1)\;\!q \,+\, 2 \;\!n\;\!\kappa}
{(n + 2)\;\!q}
\;
\mu(t)
} $}
\:
\|\, \nabla v^{[\,q\,]}(\cdot,t) \,
\|_{\mbox{}_{\scriptstyle L^{2}(\mathbb{R}^{n})}}^{\:\!2}
} $ \\
\mbox{} \vspace{+0.100cm} \\
aplicando-se (no \'ultimo passo)
a desigualdade elementar de Young
(ver e.g.$\;$\cite{Evans2002}, p.$\;$622).
Estimando-se em (3.8)
seu termo direito
como realizado nesta prova, obt\'em-se (3.10).
}
\hfill $\Box$
%
\nl
\mbox{} \vspace{-1.000cm} \\
%
%
%
%
%
{\bf Lema 3.3.}
\textit{%
Seja
$\;\! q \:\!\geq\:\!2 \;\! p$,
com
$ \;\!p $ dado em $\;\!(3.9)$,
e seja
$ {\displaystyle
\;\!
v^{[\,q\,]}(\cdot,t)
\in
L^{1}(\mathbb{R}^{n})
\cap
L^{\infty}(\mathbb{R}^{n})
} $
dada em
$\;\!(3.3) $.
Se
$ \;\! \hat{t} \in (\,0, \;\!\mbox{\small $T$}_{\!\;\!\ast})
\,\mbox{\footnotesize $\setminus$}\, E_{q} $
for tal que \\
}
\mbox{} \vspace{-0.700cm} \\
\begin{equation}
\tag{3.11$a$}
\frac{d}{d \:\!t} \;
\|\, v^{[\,q\,]}(\cdot,t) \,
\|_{\mbox{}_{\scriptstyle L^{2}(\mathbb{R}^{n})}}^{\:\!2}
{\mbox{}_{\bigr|}}_{\mbox{}_{\mbox{\footnotesize $t = \:\!\hat{t}$}}}
\hspace{-0.700cm}
\geq\: 0,
\end{equation}
\mbox{} \vspace{-0.550cm} \\
\textit{%
ent\~ao
} \\
\mbox{} \vspace{-0.950cm} \\
\begin{equation}
\tag{3.11$b$}
\|\, v^{[\,q\,]}(\cdot,\hat{t}) \,
\|_{\mbox{}_{\scriptstyle L^{2}(\mathbb{R}^{n})}}
\;\!\leq\;
\mbox{\small $ {\displaystyle
K\!\;\!(n) } $}
\:
\mbox{\small $ {\displaystyle
\Bigl(\;\! \frac{\;\!q\;\!}{2} \;\! \Bigr)
} $}^{\mbox{}^{\scriptstyle \!\!\!\;\!
\frac{\scriptstyle n}{2} \,
\frac{\scriptstyle q}
{{\scriptstyle q} \;\!-\;\! 2 \;\!{\scriptstyle n \;\!\kappa}} }}
\hspace{-0.810cm}
\cdot
\hspace{+0.445cm}
\mbox{\small $ {\displaystyle
\Bigl(\, \frac{B(\hat{t})}{\mu(\hat{t})} \,\Bigr)
} $}^{\scriptstyle \!\!\!\;\!
\frac{\scriptstyle n}{2} \,
\frac{\scriptstyle q}
{{\scriptstyle q} \;\!-\;\! 2 \;\!{\scriptstyle n \;\!\kappa}} }
\hspace{-0.775cm}
\cdot
\hspace{+0.500cm}
\|\, v^{[\,q\,]}(\cdot,\hat{t}) \,
\|_{\mbox{}_{\scriptstyle L^{1}(\mathbb{R}^{n})}}
  ^{\mbox{}^{\scriptstyle \:\!
\frac{{\scriptstyle q} \;\!-\;\!{\scriptstyle n \:\! \kappa}}
{{\scriptstyle q} \;\!-\;\! 2 \;\!{\scriptstyle n \:\!\kappa}} }}
\hspace{-0.100cm},
\end{equation}
\mbox{} \vspace{-0.150cm} \\
\textit{%
onde $\,\!K\!\;\!(n) > 0 $
\'e a constante de Nash dada em $\;\!(3.2) $. \\
}
%
%
\nl
{\small
{\bf Prova:}
De (3.10), Lema 3.2,
obt\'em-se,
usando a
hip\'otese (3.11$a$)
acima,
a estimativa \\
\mbox{} \vspace{-0.650cm} \\
\begin{equation}
\tag{3.12}
\|\, \nabla v^{[\,q\,]}(\cdot,\hat{t}) \,
\|_{\mbox{}_{\scriptstyle L^{2}(\mathbb{R}^{n})}}
\;\!\leq\;
\mbox{\small $ {\displaystyle
\Bigl(\;\! \frac{\;\!q\;\!}{2} \;\! \Bigr)
} $}^{\scriptstyle \!
\frac{{\scriptstyle n} + 2}{2} \,
\frac{\scriptstyle q}
{{\scriptstyle q} \;\!-\;\! 2 \;\!{\scriptstyle n \;\!\kappa}} }
\;\!
\mbox{\small $ {\displaystyle
\Bigl(\, \frac{B(\hat{t})}{\mu(\hat{t})} \,\Bigr)
} $}^{\scriptstyle \!\!\!\!
\frac{{\scriptstyle n} + 2}{2} \,
\frac{\scriptstyle q}
{{\scriptstyle q} \;\!-\;\! 2 \;\!{\scriptstyle n \;\!\kappa}} }
\;\!
\|\, v^{[\,q\,]}(\cdot,\hat{t}) \,
\|_{\mbox{}_{\scriptstyle L^{1}(\mathbb{R}^{n})}}
  ^{\mbox{}^{\scriptstyle \:\!
\frac{{\scriptstyle q} \;\!-\;\!({\scriptstyle n} - 2) \;\! {\scriptstyle \kappa} }
{{\scriptstyle q} \;\!-\;\! 2 \;\!{\scriptstyle n \:\!\kappa}} }}
\hspace{-0.500cm},
\end{equation}
\mbox{} \vspace{-0.150cm} \\
de onde
segue
o resultado (3.11$b$)
aplicando-se a desigualdade (3.2)
para
$ {\displaystyle
\;\!
\mbox{v} = v^{[\,q\,]}(\cdot,\hat{t})
} $.
}
\hfill $\Box$ \\
%
\nl
Em termos da solu\c c\~ao $ u(\cdot,t) $
do problema (2.1),
o Lema 3.3 \'e escrito
como segue. \\
\nl
%
%
%
%
{\bf Lema 3.\mbox{\boldmath $3^{\prime}$}\!.}
\textit{%
Seja
$\;\! q \:\!\geq\:\!2 \;\! p$,
com
$ \;\!p $ dado em $\;\!(3.9)$.
Se
$ \;\! \hat{t} \in (\,0, \;\!\mbox{\small $T$}_{\!\;\!\ast})
\,\mbox{\footnotesize $\setminus$}\, E_{q} $
for tal que \\
}
\mbox{} \vspace{-0.700cm} \\
\begin{equation}
\tag{3.13$a$}
\frac{d}{d \:\!t} \;
\|\, u(\cdot,t) \,
\|_{\mbox{}_{\scriptstyle L^{q}(\mathbb{R}^{n})}}^{\:\!q}
{\mbox{}_{\bigr|}}_{\mbox{}_{\mbox{\footnotesize $t = \:\!\hat{t}$}}}
\hspace{-0.700cm}
\geq\: 0,
\end{equation}
\mbox{} \vspace{-0.550cm} \\
\textit{%
ent\~ao
} \\
\mbox{} \vspace{-0.950cm} \\
\begin{equation}
\tag{3.13$b$}
\|\, u(\cdot,\hat{t}) \,
\|_{\mbox{}_{\scriptstyle L^{q}(\mathbb{R}^{n})}}
\;\!\leq\:
\mbox{\small $ {\displaystyle
K\!\;\!(n) } $}^{\mbox{}^{\scriptstyle \!
\frac{2}{\scriptstyle q} }}
\!\;\!
\mbox{\small $ {\displaystyle
\Bigl(\;\! \frac{\;\!q\;\!}{2} \;\! \Bigr)
} $}^{\mbox{}^{\scriptstyle \!\!\!\;\!
\frac{\scriptstyle n}
{{\scriptstyle q} \;\!-\;\! 2 \;\!{\scriptstyle n \;\!\kappa}} }}
\hspace{-0.700cm}
\cdot
\hspace{+0.300cm}
\mbox{\small $ {\displaystyle
\Bigl(\, \frac{B(\hat{t})}{\mu(\hat{t})} \,\Bigr)
} $}^{\scriptstyle \!\!\!\;\!
\frac{\scriptstyle n}
{{\scriptstyle q} \;\!-\;\! 2 \;\!{\scriptstyle n \;\!\kappa}} }
\hspace{-0.630cm}
\cdot
\hspace{+0.350cm}
\|\, u(\cdot,\hat{t}) \,
\|_{\mbox{}_{\scriptstyle L^{q/2}(\mathbb{R}^{n})}}
  ^{\mbox{}^{\scriptstyle \:\!
\frac{{\scriptstyle q} \;\!-\;\!{\scriptstyle n \:\! \kappa}}
{{\scriptstyle q} \;\!-\;\! 2 \;\!{\scriptstyle n \:\!\kappa}} }}
\hspace{-0.050cm},
\end{equation}
\mbox{} \vspace{-0.150cm} \\
\textit{%
onde $\,\!K\!\;\!(n) > 0 $
\'e a constante de Nash dada em $\;\!(3.2) $. \\
}
%
%
\mbox{} \vspace{-0.650cm} \\

Os lemas acima indicam
intuitivamente
um caminho  b\'asico
para a obten\c c\~ao
de resultados como (3.1)
usando argumentos
tipo $ L^{p}$-$\;\!L^{q}$:
para cada
$ \;\!q \geq 2 \,p_{\mbox{}_{0}} \!\;\!$,
$ \:\!q > 2 \;\! n \:\!\kappa $,
examina-se o comportamento (local)
em $t$
de
$ {\displaystyle
\;\!
\|\, u(\cdot,t) \,
\|_{L^{q}(\mathbb{R}^{n})}
\!\;\!
} $.
Caso esteja {\em crescendo},
ent\~ao
$ {\displaystyle
\;\!
\|\, u(\cdot,t) \,
\|_{L^{q}(\mathbb{R}^{n})}
} $
pode ser estimada
(localmente) em termos
de
$ {\displaystyle
\;\!
\|\, u(\cdot,t) \,
\|_{L^{q/2}(\mathbb{R}^{n})}
\!\;\!
} $,
por meio da desigualdade de energia (3.10),
de acordo com o
Lema 3.3$^{\prime}$;
se estiver {\em decrescendo},
ent\~ao (3.10) torna-se
neste caso in\'util,
mas possivelmente
esta situa\c c\~ao
possa ser
compensada
pelo fato de se saber
que
$ {\displaystyle
\;\!
\|\, u(\cdot,t) \,
\|_{L^{q/2}(\mathbb{R}^{n})}
\!\;\!
} $
n\~ao esteja crescendo
(momentaneamente, pelo menos).
Em qualquer dos casos,
sempre se possui
alguma informa\c c\~ao
aparentemente importante
sobre
$ {\displaystyle
\;\!
\|\, u(\cdot,t) \,
\|_{L^{q}(\mathbb{R}^{n})}
\!\;\!
} $,
embora varie
(de modo com\-plicado, provavelmente)
com $q $ e $t$.
Assim,
n\~ao \'e evidente
uma estrat\'egia simples
que indique como
utilizar a informa\c c\~ao dispon\'\i vel
de modo eficaz.
O pr\'oximo lema
mostra precisamente
como isso pode
ser feito.

\mbox{} \vspace{-1.000cm} \\
%
%
%
%
{\bf Lema 3.4.}
\textit{%
$\!$Seja
$\;\! q \:\!\geq 2 \;\! p$,
com
$ \:\!p $ dado em $\:\!(3.9)$.
$\!$Sendo
$\:\!u(\cdot,t) $,
$ 0 \leq t < \mbox{\small $T$}_{\!\;\!\ast}\!\;\!$,
$\!\;\!$solu\c c\~ao
do
problema
$\,\!(2.1) $,
tem-se
} \\
\mbox{} \vspace{-0.720cm} \\
\begin{equation}
\tag{3.14}
\mathbb{U}_{q}(0; t)
\;\!\leq\,
\max\,\biggl\{\,
\|\, u_0 \;\!\|_{\mbox{}_{\scriptstyle L^{q}(\mathbb{R}^{n})}}
\!\,\!;
\;\!\,\!
K\!\;\!(n)^{\mbox{}^{\scriptstyle \!\frac{2}{\scriptstyle q} }}
\!\;\!
\Bigl(\;\! \mbox{\small $ {\displaystyle \frac{q}{2} }$}
\;\!\Bigr)^{\!\!\;\! \frac{\scriptstyle n}
{{\scriptstyle q} \;\!-\;\! 2\;\!{\scriptstyle n \:\!\kappa}} }
\,
\mathbb{B}_{\mu}\!\;\!(0; t)^{\mbox{}^{\scriptstyle \!
\frac{{\scriptstyle n}}{{\scriptstyle q} \;\!-\;\! 2 \;\!
{\scriptstyle n \:\!\kappa}} }}
\;\!
\mathbb{U}_{\mbox{}_{\!\frac{\scriptstyle q}{2}}}\!
(0; t)^{\mbox{}^{\scriptstyle \!
\frac{{\scriptstyle q} \;\!-\;\! {\scriptstyle n\:\!\kappa} }
{{\scriptstyle q} \;\!-\;\! 2 \;\!
{\scriptstyle n \:\!\kappa}} }}
\:\!
\biggr\}
\end{equation}
\mbox{} \vspace{-0.170cm} \\
\textit{%
para todo
$ \:\!0 \leq t < \mbox{\small $T$}_{\!\;\!\ast}\!\;\!$,
$\!\;\!$onde
$ \;\!\mathbb{B}_{\mu}\!\;\!(0\:\!; \:\!t) $,
$ \!\;\!\mathbb{U}_{q}\!\;\!(0\:\!; \:\!t) $,
$ \!\;\!\mathbb{U}_{\frac{\scriptstyle q}{2}}\!\;\!(0\:\!; \:\!t) $
s\~ao definidas em
$\:\!(1.7)$ e $(\:\!1.8)$, \linebreak
\mbox{} \vspace{-0.570cm} \\
e $K\!\;\!(n) $
\'e dada na desigualdade $\:\!(3.2)$. \\
}
%
%
\nl
\mbox{} \vspace{-0.475cm} \\
%
%
%
{\small
{\bf Prova:}
Dado
$ \;\!0 \leq t < \mbox{\small $T$}_{\!\;\!\ast}$,
fixo no que segue,
seja
(por conveni\^encia)
$ \,\!\gamma_{q} \in \mathbb{R}^{\mbox{}^{+}}\!\!\;\!$
definido por \\
\mbox{} \vspace{-0.600cm} \\
\begin{equation}
\notag
\gamma_{q}
\,=\:
K\!\;\!(n)^{\mbox{}^{\scriptstyle \!\frac{2}{\scriptstyle q} }}
\!\;\!
\Bigl(\;\! \mbox{\small $ {\displaystyle \frac{q}{2} }$}
\;\!\Bigr)^{\!\!\;\! \frac{\scriptstyle n}
{{\scriptstyle q} \;\!-\;\! 2\;\!{\scriptstyle n \:\!\kappa}} }
\,
\mathbb{B}_{\mu}\!\;\!(0; t)^{\mbox{}^{\scriptstyle \!
\frac{{\scriptstyle n}}{{\scriptstyle q} \;\!-\;\! 2 \;\!
{\scriptstyle n \:\!\kappa}} }}
\;\!
\mathbb{U}_{\mbox{}_{\!\frac{\scriptstyle q}{2}}}\!
(0; t)^{\mbox{}^{\scriptstyle \!
\frac{{\scriptstyle q} \;\!-\;\! {\scriptstyle n\:\!\kappa} }
{{\scriptstyle q} \;\!-\;\! 2 \;\!
{\scriptstyle n \:\!\kappa}} }}
\end{equation}
\mbox{} \vspace{-0.270cm} \\
e consideremos
os casos poss\'\i veis
para os valores de
$ {\displaystyle
\;\!
\|\, u(\cdot,\tau) \,
\|_{\mbox{}_{\scriptstyle L^{q}(\mathbb{R}^{n})}}
\!\:\!
} $
no intervalo
$ 0 \leq \tau \leq t $: \\
\mbox{} \vspace{-0.100cm} \\
%
%
\mbox{\tt Caso} \mbox{\tt I:}
$ {\displaystyle
\;\!
\|\, u(\cdot,\tau) \,
\|_{\mbox{}_{\scriptstyle L^{q}(\mathbb{R}^{n})}}
\!\:\!
>\;\!
\gamma_{q}
} $
para todo
$ \;\!0 \leq \tau < t $. \\
\mbox{} \vspace{-0.150cm} \\
Neste caso,
segue do Lema 3.$3^{\prime}$
que temos de ter
$ {\displaystyle
\,
d/d\tau \,
\|\, u(\cdot,\tau) \,
\|_{\mbox{}_{\scriptstyle L^{q}(\mathbb{R}^{n})}}^{\:\!q}
\!\!\;\!<\:\!0
\;\!
} $
para quase todo \linebreak
\mbox{} \vspace{-0.570cm} \\
$ \tau \in I \equiv [\,0, \:\!t\,] $,
de modo que
$ {\displaystyle
\,
\|\, u(\cdot,\tau) \,
\|_{\mbox{}_{\scriptstyle L^{q}(\mathbb{R}^{n})}}
\!\!\;\!
} $
\'e (estritamente) decrescente
neste intervalo.
Em particular,
segue que
$ {\displaystyle
\,
\|\, u(\cdot,\tau) \,
\|_{\mbox{}_{\scriptstyle L^{q}(\mathbb{R}^{n})}}
\!\;\!\leq\;\!
\|\, u_0 \;\!
\|_{\mbox{}_{\scriptstyle L^{q}(\mathbb{R}^{n})}}
\!\;\!
} $
para todo
$ \tau $ em $I$,
ou seja,
tem-se \linebreak
\mbox{} \vspace{-0.570cm} \\
neste caso \\
\mbox{} \vspace{-1.150cm} \\
\begin{equation}
\notag
\mathbb{U}_{q}(0; t)
\,=\,
\|\, u_0 \;\!
\|_{\mbox{}_{\scriptstyle L^{q}(\mathbb{R}^{n})}}
\!\;\!.
\end{equation}
\mbox{} \vspace{-0.200cm} \\
%
%
\mbox{\tt Caso} \mbox{\tt II:}
$ {\displaystyle
\;\!
\|\, u_0 \;\!
\|_{\mbox{}_{\scriptstyle L^{q}(\mathbb{R}^{n})}}
\!\:\!
>\;\!
\gamma_{q}
} $,
$\;\!$tendo-se
$ {\displaystyle
\;\!
\|\, u(\cdot,t_{\mbox{}_{\scriptstyle \ast}}\!\;\!) \,
\|_{\mbox{}_{\scriptstyle L^{q}(\mathbb{R}^{n})}}
\!\:\!
\leq\;\!
\gamma_{q}
} $
para algum
$ \;\!0 < t_{\mbox{}_{\scriptstyle \ast}} < t $. \\
\mbox{} \vspace{-0.150cm} \\
Neste caso,
existe
$ \;\!0 < t_{\mbox{}_{1}} \!\:\! \leq t_{\mbox{}_{\scriptstyle \ast}} $
tal que
$ {\displaystyle
\;\!
\|\, u(\cdot,\tau) \,
\|_{\mbox{}_{\scriptstyle L^{q}(\mathbb{R}^{n})}}
\!\:\!
>
\gamma_{q}
} $
em
$ [\,0, \:\!t_{\mbox{}_{1}} \!\;\!) $,
$ {\displaystyle
\;\!
\|\, u(\cdot, t_{\mbox{}_{1}} \!\;\!)\,
\|_{\mbox{}_{\scriptstyle L^{q}(\mathbb{R}^{n})}}
\!\:\!
=
\gamma_{q}
} $.
Pelo Lema 3.$3^{\prime}\!\;\!$,
segue repetindo o argumento acima
que
$ {\displaystyle
\;\!
\|\, u(\cdot,\tau) \,
\|_{\mbox{}_{\scriptstyle L^{q}(\mathbb{R}^{n})}}
\!
} $
tem de ser decrescente
em $[\,0, \:\!t_{\mbox{}_{1}}] $.
Por outro lado,
no intervalo
$ J \equiv [\, t_{\mbox{}_{1}} \!\;\!, \:\! t\;\!] \;\!$
temos de ter
$ {\displaystyle
\;\!
\|\, u(\cdot,\tau) \,
\|_{\mbox{}_{\scriptstyle L^{q}(\mathbb{R}^{n})}}
\!\leq
\gamma_{q}
} $
para todo
$ \tau \in J $.
\mbox{[}$\,$De fato,
se assim n\~ao fosse,
teriam de existir
$ \;\!t_{\mbox{}_{2}} \!\:\!< t_{\mbox{}_{3}} \in
[\;\!\:\!t_{\mbox{}_{1}}\!\;\!, \:\!t\,] $
tais que
$ {\displaystyle
\|\, u(\cdot, \tau) \,
\|_{\mbox{}_{\scriptstyle L^{q}(\mathbb{R}^{n})}}
\!>
\gamma_{q}
} $
para todo
$ t_{\mbox{}_{2}} \!< \tau \leq t_{\mbox{}_{3}} \!\;\!$,
tendo-se
$ {\displaystyle
\;\!
\|\, u(\cdot,t_{\mbox{}_{2}}\!\;\!) \,
\|_{\mbox{}_{\scriptstyle L^{q}(\mathbb{R}^{n})}}
\!=
\gamma_{q}
} $.
Assim, teria de existir
$ t_{\mbox{}_{\ast\ast}} \!\!\:\!\in
(\,t_{\mbox{}_{2}}\!\;\!, \;\! t_{\mbox{}_{3}})
\setminus E_{q} $
com
$ {\displaystyle
\,
d/d\tau \,
\|\, u(\cdot, \tau) \,
\|_{\mbox{}_{\scriptstyle L^{q}(\mathbb{R}^{n})}}
  ^{\:\!q}
\!
} $
positiva em $ \:\!\tau = t_{\mbox{}_{\ast\ast}} \!\:\!$,
de modo que,
pelo \linebreak
\mbox{} \vspace{-0.560cm} \\
Lema 3.$3^{\prime}\!\;\!$,
valeria
$ {\displaystyle
\;\!
\|\, u(\cdot, t_{\mbox{}_{\ast\ast}}\!) \,
\|_{\mbox{}_{\scriptstyle L^{q}(\mathbb{R}^{n})}}
\!\leq \gamma_{q}
\!\;\!
} $,
contradizendo o fato
de ter-se
$ {\displaystyle
\;\!
\|\, u(\cdot, \tau) \,
\|_{\mbox{}_{\scriptstyle L^{q}(\mathbb{R}^{n})}}
\!> \gamma_{q}
\!\;\!
} $
em todo o intervalo
$ (\;\!t_{\mbox{}_{2}}\!\;\!, \:\!t_{\mbox{}_{3}}) $.$\,$\mbox{]}
$\;\!$Portanto,
tem-se tamb\'em aqui
$ {\displaystyle
\;\!
\|\, u(\cdot,\tau) \,
\|_{\mbox{}_{\scriptstyle L^{q}(\mathbb{R}^{n})}}
\!\;\!\leq\;\!
\|\, u_0 \;\!
\|_{\mbox{}_{\scriptstyle L^{q}(\mathbb{R}^{n})}}
\!\;\!
} $
para todo
$ \tau \in [\,0, \:\!t\;\!] $,
ou seja,
obt\'em-se novamente
$ {\displaystyle
\;\!
\mathbb{U}_{q}(0; t)
\,=\,
\|\, u_0 \;\!
\|_{\mbox{}_{\scriptstyle L^{q}(\mathbb{R}^{n})}}
\!\;\!
} $. \\
\mbox{} \vspace{-0.050cm} \\
%
%
\mbox{\tt Caso} \mbox{\tt III:}
$ {\displaystyle
\;\!
\|\, u_0 \;\!
\|_{\mbox{}_{\scriptstyle L^{q}(\mathbb{R}^{n})}}
\!\:\!
\leq\;\!
\gamma_{q}
} $. \\
\mbox{} \vspace{-0.150cm} \\
Neste caso,
repetindo-se o argumento
aplicado no {\tt Caso} {\tt II} acima
para o intervalo $J\!\;\!$,
resulta que se tem
$ {\displaystyle
\;\!
\|\, u(\cdot,\tau) \,
\|_{\mbox{}_{\scriptstyle L^{q}(\mathbb{R}^{n})}}
\!\,\!\leq
\gamma_{q}
\:\!
} $
em todo o intervalo
$ \;\![\;\!0, \;\!t\;\!\,\!] $,
de modo que,
neste caso,
tem-se \linebreak

\mbox{} \vspace{-1.150cm} \\
\begin{equation}
\notag
\mathbb{U}_{q}\!\;\!(0\:\!; \:\!t)
\:\leq\:
K\!\;\!(n)^{\mbox{}^{\scriptstyle \!\frac{2}{\scriptstyle q} }}
\!\;\!
\Bigl(\;\! \mbox{\small $ {\displaystyle \frac{q}{2} }$}
\;\!\Bigr)^{\!\!\;\! \frac{\scriptstyle n}
{{\scriptstyle q} \;\!-\;\! 2\;\!{\scriptstyle n \:\!\kappa}} }
\,
\mathbb{B}_{\mu}\!\;\!(0; t)^{\mbox{}^{\scriptstyle \!
\frac{{\scriptstyle n}}{{\scriptstyle q} \;\!-\;\! 2 \;\!
{\scriptstyle n \:\!\kappa}} }}
\;\!
\mathbb{U}_{\mbox{}_{\!\frac{\scriptstyle q}{2}}}\!
(0; t)^{\mbox{}^{\scriptstyle \!
\frac{{\scriptstyle q} \;\!-\;\! {\scriptstyle n\:\!\kappa} }
{{\scriptstyle q} \;\!-\;\! 2 \;\!
{\scriptstyle n \:\!\kappa}} }}
\hspace{-0.960cm}.
\end{equation}
\mbox{} \vspace{-0.300cm} \\
Em todos os casos,
tem-se sempre (3.14) acima,
o que conclui a prova do Lema 3.4.
}
\mbox{} \hfill $\Box$ \\
%

%
\mbox{} \vspace{-1.500cm} \\

Para os resultados seguintes,
ser\'a conveniente
introduzir,
para cada
$ 1 \leq j \leq m $,
$ p \geq p_{\mbox{}_{0}} \!\;\!$,
$ p > n \:\!\kappa $,
a constante
$ \:\!C(j,m) \,\equiv\, C(j,m\:\!; n,p,\kappa) > 0 \;\!$
definida por \\
\mbox{} \vspace{-0.750cm} \\
\begin{equation}
\tag{3.15$a$}
C(j,m)
\,:=\;
\prod_{\ell\,=\,j}^{m}
\,
\lambda(\,\!2^{\ell}p\,\!)^{\mbox{}^{\scriptstyle \!\!\;\!
\frac{\,{\scriptstyle p} \:-\:
2^{\mbox{}^{-\;\!m}} \!\:\!{\scriptstyle n \:\!\kappa}}
{{\scriptstyle p} \;-\;
2^{\mbox{}^{-\;\!\ell}} \!\;\!{\scriptstyle n \:\!\kappa}}
}}
\end{equation}
\mbox{} \vspace{-0.550cm} \\
sendo \\
\mbox{} \vspace{-1.350cm} \\
\begin{equation}
\tag{3.15$b$}
\lambda(q)
\,\equiv\,
\lambda(n, \:\!\kappa, \:\!q)
\;\!:=\;
K\!\;\!(n)^{\mbox{}^{\scriptstyle \!\frac{2}{\scriptstyle q} }}
\Bigl(\;\! \mbox{\small $ {\displaystyle \frac{\mbox{\normalsize $q$}}{2} }$}
\;\!\Bigr)^{\mbox{}^{\scriptstyle \!\!
\frac{\scriptstyle n}
{{\scriptstyle q} \;\!-\;\! 2\;\!{\scriptstyle n \:\!\kappa}}
}}
\end{equation}
\mbox{} \vspace{-0.175cm} \\
para
todo
$ q \geq 2 \,\!\,\!p $,
onde
$ K\!\;\!(n) > 0 $
denota a constante de Nash
na desigualdade (3.2). \\
\nl
%
%
%
%
{\bf Lema 3.5.}
\textit{%
$\!$Seja
$\;\! p \:\!\geq\;\! p_{\mbox{}_{0}} \!\;\!$,
$ \:\!p > n \:\! \kappa $.
$\!$Sendo
$\:\!u(\cdot,t) $,
$ 0 \leq t < \mbox{\small $T$}_{\!\;\!\ast}\!\;\!$,
$\!\;\!$solu\c c\~ao
de
$\,\!(2.1) $,
tem-se
} \\
\mbox{} \vspace{-0.720cm} \\
\begin{equation}
\tag{3.16$a$}
\mathbb{U}_{2\:\!p}(0; t)
\;\!\leq\,
\max\,\biggl\{\,
\|\, u_0 \;\!\|_{\mbox{}_{\scriptstyle L^{2\:\!p}(\mathbb{R}^{n})}}
\!\,\!;
\;\!\,\!
\lambda(2\:\!p)
\:
\mathbb{B}_{\mu}\!\;\!(0; t)^{\mbox{}^{\scriptstyle \!
\frac{{\scriptstyle n}}{2\;\!{\scriptstyle p} \;\!-\;\! 2 \;\!
{\scriptstyle n \:\!\kappa}} }}
\;\!
\mathbb{U}_{\mbox{}_{\scriptstyle \!p}}
(0; t)^{\mbox{}^{\scriptstyle \!
\frac{2 \;\!{\scriptstyle p} \;\!-\;\! {\scriptstyle n\:\!\kappa} }
{2\;\!{\scriptstyle p} \;\!-\;\! 2 \;\!
{\scriptstyle n \:\!\kappa}} }}
\:\!
\biggr\}
\end{equation}
\mbox{} \vspace{-0.225cm} \\
\textit{%
e, mais geralmente,
para todo $ \;\! m \geq 2 \!:$
} \\
\mbox{} \vspace{-0.250cm} \\
\mbox{} \hfill
$ {\displaystyle
\mathbb{U}_{2^{\mbox{}^{m}}\!p}(0; t)
\;\!\leq\,
\max\;\biggl\{\;\!\;\!
\|\, u_0 \;\!
\|_{\mbox{}_{\scriptstyle L^{2^{m}p}(\mathbb{R}^{n})}}
\!\:\!;
\;\;\!
C(j,m)
\;\:\!
\mathbb{B}_{\mu}\!\;\!(0; t)^{\mbox{}^{\scriptstyle
\!\!\!
\frac{\;\!{\scriptstyle p} \:-\;
{\scriptstyle n \:\!\kappa}\:\!/\;\!2^{m}}
{{\scriptstyle p}}
\;\!
\left[\;
\frac{2\;\!{\scriptstyle n}}{\;\!2^{j}{\scriptstyle p}
\,-\, 2{\scriptstyle n \:\!\kappa}}
\:-\;
\frac{\scriptstyle n}{2^{m}{\scriptstyle p}
\,-\, {\scriptstyle n \:\!\kappa}} \,
\right] }}
\,
\times
} $ \\
\mbox{} \vspace{-0.250cm} \\
\mbox{} \hspace{+6.650cm}
$ {\displaystyle
\times
\;\;\!
\|\, u_0 \;\!
\|_{{\scriptstyle L^{2^{j}p/2}(\mathbb{R}^{n})}}
  ^{\mbox{}^{\scriptstyle
\frac{\;\!{\scriptstyle p} \:-\;
{\scriptstyle n \:\!\kappa}\:\!/\;\!2^{m}}
{\;\!{\scriptstyle p} \:-\; 2 \;\!
{\scriptstyle n \:\!\kappa}\:\!/\;\!2^{\:\!j}}
}}
\!\!\:\!,
\hspace{+0.400cm}
2 \;\!\leq\;\! j \;\!\leq\;\! m
\;\!;
} $
\mbox{} \hfill (3.16$b$) \\
\mbox{} \vspace{-0.050cm} \\
\mbox{} \hspace{+6.700cm}
$ {\displaystyle
C(1,m)
\;\:\!
\mathbb{B}_{\mu}\!\;\!(0; t)^{\mbox{}^{\scriptstyle
\!\!\!\;\!
\frac{\;\!{\scriptstyle n}
\,(\;\!1 \:-\; 2^{-\;\!m})}
{\;\!{\scriptstyle p} \,-\: {\scriptstyle n \:\!\kappa}}
}}
\:\!
\mathbb{U}_{p}(0; t)^{\mbox{}^{\scriptstyle
\!\!\!\;\!
\frac{\;\!{\scriptstyle p} \:-\;
{\scriptstyle n \:\!\kappa}\:\!/\;\!2^{m}}
{{\scriptstyle p} \:-\; {\scriptstyle n \:\!\kappa}}
}}
\;
\biggr\}
\:\!
} $, \\
\mbox{} \vspace{-0.050cm} \\
\textit{%
onde
$\:\! C(j,m) $, $ 1 \leq j \leq m $,
s\~ao as constantes dadas em $\;\!(3.15)$,
com
$ \;\!\mathbb{B}_{\mu}\!\;\!(0; t) $,
$ \mathbb{U}_{q}\!\;\!(0; t) $
definidas em
$\;\!(1.7) $, $(1.8)$. \\
}
%
%
\nl
{\small
{\bf Prova:}
A express\~ao (3.16$a$)
corresponde a (3.14) do
Lema 3.4 acima,
tomando-se simples\-mente
$ q = 2 p $;
a prova de (3.16$b$)
\'e realizada por indu\c c\~ao em $m$,
como indicado a seguir.
Para $ m = 2 $,
(3.16$b$)
\'e obtida
combinando-se (3.16$a$)
e \mbox{[}$\,$(3.14), com $ q = 4 \:\!p \,$\mbox{]}.
Dado
$ m \geq 3 $ arbitr\'ario,
supondo-se que
(3.16$b$) seja v\'alida
para inteiros menores que $m$,
obt\'em-se,
pelo Lema 3.4,
tomando-se
$ q = 2^{m} p $
em (3.14), \\
\mbox{} \vspace{-0.720cm} \\
\begin{equation}
\notag
\mathbb{U}_{2^{m}p}(0; t)
\;\!\leq\,
\max\,\biggl\{\,
\|\, u_0 \;\!\|_{\mbox{}_{\scriptstyle L^{2^{m}p}(\mathbb{R}^{n})}}
\!\,\!;
\;\:\!
\lambda(2^{m}p)
\;
\mathbb{B}_{\mu}\!\;\!(0\:\!; \:\!t)^{\mbox{}^{\scriptstyle \!
\frac{{\scriptstyle n}}{\;\!2^{m}p \;\!-\;\! 2 \;\!
{\scriptstyle n \:\!\kappa}} }}
\;\!
\mathbb{U}_{\mbox{}_{\scriptstyle \!2^{m-1}p}}\!
(0\:\!; \:\!t)^{\mbox{}^{\scriptstyle \!
\frac{\;\!{\scriptstyle p} \:-\;
{\scriptstyle n \:\!\kappa}\:\!/\;\!2^{m}}
{\;\!{\scriptstyle p} \:-\; 2 \;\!
{\scriptstyle n \:\!\kappa}\:\!/\;\!2^{m}}
}}
\:\!
\biggr\},
\end{equation}
\mbox{} \vspace{-0.190cm} \\
de onde a express\~ao
(3.16$b$)
segue
aplicando-se a hip\'otese de indu\c c\~ao
para
$ {\displaystyle
\;\!
\mathbb{U}_{\mbox{}_{\scriptstyle \!2^{m-1}p}}\!
(0\:\!; \:\!t)
} $.
}
\hfill $\Box$ \\
%
%

%
\mbox{} \vspace{-1.050cm} \\

Antes de prosseguir,
ser\'a \'util estimar
as constantes
$ C(j,m) $ definidas em
(3.15).
Observando
as somas elementares
abaixo, \\
\mbox{} \vspace{-0.650cm} \\
\begin{equation}
\tag{3.17$a$}
\sum_{\ell \,=\, 1}^{m}
\;\!
\mbox{\small $ {\displaystyle
\frac{2^{\ell} \:\!p}
{\;\!(\:\!2^{\ell}\:\!p \:\!-\;\! 2 \;\! n \:\! \kappa\:\!)
\, (\:\!2^{\ell}\:\!p \:\!-\;\! n \:\! \kappa\:\!)\;\!}
\;\;\!=\;\;\!
\frac{1}{\;\!p \:\!-\;\! n \:\!\kappa \;\!}
\;-\;
\frac{1}{\;\!2^{m}\:\! p \:\!-\;\! n \:\!\kappa \;\!}
} $}
\;\!,
\mbox{} \hspace{+0.200cm}
\end{equation}

\mbox{} \vspace{-1.500cm} \\
\begin{equation}
\tag{3.17$b$}
\sum_{\ell \,=\, 1}^{m}
\;\!
\mbox{\small $ {\displaystyle
\frac{\ell \;\!\cdot\, 2^{\ell} \:\!p}
{\;\!(\:\!2^{\ell}\:\!p \:\!-\;\! 2 \;\! n \:\! \kappa\:\!)
\, (\:\!2^{\ell}\:\!p \:\!-\;\! n \:\! \kappa\:\!)\;\!}
\;\;\!=\;\;\!
\frac{1}{\;\!p \:\!-\;\! n \:\!\kappa \;\!}
\;-\;
\frac{m \;\!+\;\!1}{\;\!2^{m}\:\! p \:\!-\;\! n \:\!\kappa \;\!}
} $}
\;+\,
\sum_{\ell \,=\, 1}^{m}
\;\!
\mbox{\small $ {\displaystyle
\frac{1}{\;\!2^{j}\:\! p \:\!-\;\! n \:\!\kappa \;\!}
} $}
\;\!,
\mbox{} \hspace{+0.200cm}
\end{equation}
\mbox{} \vspace{+0.050cm} \\
resulta
a seguinte estimativa
b\'asica
para os coeficientes
$ C(j,m) $
em (3.15), (3.16). \linebreak
\mbox{} \vspace{+0.050cm} \\
%
%
%
%
%
{\bf Lema 3.6.}
\textit{%
$\!$Sejam
$ \;\!p > n \:\! \kappa $,
$ m \geq 2 $.
Ent\~ao,
para
$\:\!C(j,m) $
dado em $\;\!(3.15)$,
tem-se \\
}
\mbox{} \vspace{-0.700cm} \\
\begin{equation}
\tag{3.18}
C(j,m)
\;\leq\;
\bigl(\;\! 2 \;\!p \;\!
\bigr)^{\mbox{}^{\scriptstyle
\!
\frac{{\scriptstyle n}}
{\;\!{\scriptstyle p} \;\!-\;\! {\scriptstyle n \:\!\kappa}}
}}
\equiv\:
K\!\;\!(n, \,\!\kappa, \,\!p),
\qquad \,
\forall \;\,
1 \;\!\leq\;\!j \;\!\leq\;\! m.
\end{equation}
%
%
%
\mbox{} \vspace{+0.050cm} \\
%
%
%
{\small
{\bf Prova:}
Como
$ K\!\;\!(n) < 1 $
para todo $n$
$\;\!$(cf.$\;$\cite{CarlenLoss1993}, p.$\;$213)
obt\'em-se,
de (3.15): \\
\mbox{} \vspace{-0.600cm} \\
\begin{equation}
\notag
C(j,m)
\;\;\!\leq\;\;\!
p^{\mbox{}_{\scriptstyle
\!
\frac{\,\!(\:\!2^{m}\:\!{\scriptstyle p}
\,-\, {\scriptstyle n \:\!\kappa}\:\!)\,\!}
{2^{\mbox{}^{m}} \!\;\!{\scriptstyle p}}
\;\!
\mbox{\scriptsize
$ {\displaystyle
\sum_{\ell\,=\,j}^{m}
\;\!
\frac{2^{\ell} \:\!{\scriptstyle p \, n}}
{\,\!(\:\!2^{\ell} {\scriptstyle p} -
2 \;\! {\scriptstyle n \:\! \kappa}\:\!)
\, (\:\!2^{\ell} {\scriptstyle p} -
{\scriptstyle n \:\! \kappa}\:\!)\,\!}
} $}
}}
\;\!\times\;
2^{\mbox{}_{\scriptstyle
\!\!
\frac{\,\!(\:\!2^{m}\:\!{\scriptstyle p}
\,-\, {\scriptstyle n \:\!\kappa}\:\!)\,\!}
{2^{\mbox{}^{m}} \!\;\!{\scriptstyle p}}
\;\!
\mbox{\scriptsize
$ {\displaystyle
\sum_{\ell\,=\,j}^{m}
\;\!
\frac{{\scriptstyle \ell} \cdot
2^{\ell} \:\!{\scriptstyle p \, n}}
{\,\!(\:\!2^{\ell} {\scriptstyle p} -
2 \;\! {\scriptstyle n \:\! \kappa}\:\!)
\, (\:\!2^{\ell} {\scriptstyle p} -
{\scriptstyle n \:\! \kappa}\:\!)\,\!}
} $}
}}
,
\end{equation}
\mbox{} \vspace{-0.300cm} \\
de onde segue a estimativa (3.18),
visto que \\
\mbox{} \vspace{-0.025cm} \\
\mbox{} \hspace{+0.000cm}
$ {\displaystyle
\frac{\,\!2^{m} \,\!p - n \:\!\kappa \,\!}
{2^{m} p}
\,
\sum_{\ell \,=\, j}^{m}
\;\!
\mbox{\small $ {\displaystyle
\frac{2^{\ell} \:\!p}
{\;\!(\:\!2^{\ell}\:\!p \:\!-\;\! 2 \;\! n \:\! \kappa\:\!)
\, (\:\!2^{\ell}\:\!p \:\!-\;\! n \:\! \kappa\:\!)\;\!}
\;\,\!\leq\,
\sum_{\ell \,=\, 1}^{m}
\;\!
\mbox{\small $ {\displaystyle
\frac{2^{\ell} \:\!p}
{\;\!(\:\!2^{\ell}\:\!p \:\!-\;\! 2 \;\! n \:\! \kappa\:\!)
\, (\:\!2^{\ell}\:\!p \:\!-\;\! n \:\! \kappa\:\!)\;\!}
} $}
\;\,\!\leq\;\,\!
\frac{1}{\,\!p - n \:\!\kappa\,\!}
} $}
} $ \\
\mbox{} \vspace{-0.325cm} \\
e \\
\mbox{} \vspace{-0.275cm} \\
\mbox{} \hspace{+0.000cm}
$ {\displaystyle
\frac{\,\!2^{m} \,\!p - n \:\!\kappa \,\!}
{2^{m} p}
\,
\sum_{\ell \,=\, j}^{m}
\;\!
\mbox{\small $ {\displaystyle
\frac{\ell \;\!\cdot\, 2^{\ell} \:\!p}
{\;\!(\:\!2^{\ell}\:\!p \:\!-\;\! 2 \;\! n \:\! \kappa\:\!)
\, (\:\!2^{\ell}\:\!p \:\!-\;\! n \:\! \kappa\:\!)\;\!}
\;\,\!\leq\,
\sum_{\ell \,=\, 1}^{m}
\;\!
\mbox{\small $ {\displaystyle
\frac{\ell \;\!\cdot\, 2^{\ell} \:\!p}
{\;\!(\:\!2^{\ell}\:\!p \:\!-\;\! 2 \;\! n \:\! \kappa\:\!)
\, (\:\!2^{\ell}\:\!p \:\!-\;\! n \:\! \kappa\:\!)\;\!}
} $}
\;\,\!\leq\;\,\!
\frac{1}{\,\!p - n \:\!\kappa\,\!}
} $}
} $ \\
\mbox{} \vspace{+0.175cm} \\
(devido a (3.17) acima).
\mbox{[}$\,$Note-se que
estimativas mais finas
para $ C(j,m) $
tamb\'em podem ser
obtidas,
de modo an\'alogo,
mas este ponto n\~ao \'e essencial
no argumento a seguir.$\,$\mbox{]}
}
\mbox{} \hfill $\Box$ \\
%
\mbox{} \vspace{-0.500cm} \\

Usando-se as express\~oes (3.15), (3.16$b$) e (3.18) acima,
obt\'em-se
uma estimativa mais simples
para
$ {\displaystyle
\;\!
\mathbb{U}_{2^{\mbox{}^{m}}\!p}(0\:\!; \:\!t)
} $,
descrita em (3.19) abaixo.
Esta estimativa
representa o passo final
para estabelecermos o Teorema 3.1. \\
\mbox{} \vspace{-0.150cm} \\
%
%
%
%
%
%
\mbox{} \hspace{-0.800cm}
\fbox{%
\begin{minipage}[t]{16.000cm}
\mbox{} \vspace{-0.450cm} \\
\mbox{} \hspace{+0.300cm}
\begin{minipage}[t]{15.000cm}
\mbox{} \vspace{+0.100cm} \\
{\bf Lema 3.7.}
\textit{%
$\!$Seja
$\;\! p \,\!\geq\,\! p_{\mbox{}_{0}} \!\;\!$,
$\!\;\!$com
$ \;\!p > n \:\! \kappa $.
$\!$Sendo
$\:\!u(\cdot,t) $,
$ 0 \leq t < \mbox{\small $T$}_{\!\;\!\ast}\!\;\!$,
$\!\;\!$solu\c c\~ao
do pro\-blema
dado em
$\;\!(2.1) $
acima,
tem-se
} \\
\mbox{} \vspace{-0.050cm} \\
\mbox{} \hfill
$ {\displaystyle
\mathbb{U}_{2^{\mbox{}^{m}}\!p}(0; t)
\,\leq\,
K\!\;\!(n,\kappa,p)
\;\!\cdot\;\!
\max\;\biggl\{\;\!\;\!
\|\, u_0 \;\!
\|_{\mbox{}_{\scriptstyle L^{2^{m}p}(\mathbb{R}^{n})}}
\!\:\!;
\;\;\!
\mathbb{B}_{\mu}\!\;\!(0; t)^{\mbox{}^{\scriptstyle
\!\!\!\;\!
\frac{\;\!{\scriptstyle n}
\,(\;\!1 \:-\; 2^{-\;\!m})}
{\;\!{\scriptstyle p} \,-\: {\scriptstyle n \:\!\kappa}}
}}
\:\!
\mathbb{U}_{p}(0; t)^{\mbox{}^{\scriptstyle
\!\!\!\;\!
\frac{\;\!{\scriptstyle p} \:-\;
{\scriptstyle n \:\!\kappa}\:\!/\;\!2^{m}}
{{\scriptstyle p} \:-\; {\scriptstyle n \:\!\kappa}}
}}
\,
\biggr\}
} $ \\
\mbox{} \vspace{-0.300cm} \\
\mbox{} \hfill (3.19) \\
\mbox{} \vspace{-0.450cm} \\
\textit{%
para todo
$ \;\!0 \leq t < \mbox{\small $T$}_{\!\;\!\ast}\!\;\!$,
e todo $ \;\! m \geq 1 $,
onde
$ \;\!\mathbb{B}_{\mu}\!\;\!(0; t) $,
$ \mathbb{U}_{q}\!\;\!(0; t) $
s\~ao dadas em
$\;\!(1.7) $, $(1.8)$, \linebreak
e onde
$ K\!\;\!(n,\kappa,p) $
\'e definida em $\;\!(3.18)$. \\
%
%
%
%
%
}
\end{minipage}
\end{minipage}
}
%
%
\nl
\mbox{} \vspace{-0.500cm} \\
%
%
%
{\small
{\bf Prova:}
Se $ m = 1 $,
o resultado segue
imediatamente de
(3.15$b$), (3.16$a$),
j\'a que
$ \lambda(2\:\!p) \leq K\!\;\!(n,\kappa,p) $.
Consideremos,
assim,
$ m \geq 2 $.
Dado $ 2 \leq j \leq m $,
obt\'em-se,
ent\~ao,
estimando-se \linebreak
$ {\displaystyle
\|\, u_0 \;\!
\|_{L^{2^{j}p/2}(\mathbb{R}^{n})}
\!\;\!
} $
por interpola\c c\~ao
com respeito \`as
duas normas
$ {\displaystyle
\;\!
\|\, u_0 \;\!
\|_{\mbox{}_{\scriptstyle L^{p}(\mathbb{R}^{n})}}
\!\;\!
} $
e
$ {\displaystyle
\;\!
\|\, u_0 \;\!
\|_{\mbox{}_{\scriptstyle L^{2^{m}p}(\mathbb{R}^{n})}}
\!\;\!
} $
(e usando (3.18) acima): \\
\mbox{} \vspace{+0.100cm} \\
\mbox{} \hspace{-0.200cm}
$ {\displaystyle
C(j,m)
\;\:\!
\mathbb{B}_{\mu}\!\;\!(0; t)^{\mbox{}^{\scriptstyle
\!\!\!
\frac{\;\!{\scriptstyle p} \:-\;
{\scriptstyle n \:\!\kappa}\:\!/\;\!2^{m}}
{{\scriptstyle p}}
\;\!
\left[\;
\frac{2\;\!{\scriptstyle n}}{\;\!2^{j}{\scriptstyle p}
\,-\, 2{\scriptstyle n \:\!\kappa}}
\:-\;
\frac{\scriptstyle n}{2^{m}{\scriptstyle p}
\,-\, {\scriptstyle n \:\!\kappa}} \,
\right] }}
\,
\|\, u_0 \;\!
\|_{{\scriptstyle L^{2^{j}p/2}(\mathbb{R}^{n})}}
  ^{\mbox{}^{\scriptstyle
\frac{\;\!{\scriptstyle p} \:-\;
{\scriptstyle n \:\!\kappa}\:\!/\;\!2^{m}}
{\;\!{\scriptstyle p} \:-\; 2 \;\!
{\scriptstyle n \:\!\kappa}\:\!/\;\!2^{\:\!j}}
}}
} $ \\
\mbox{} \vspace{+0.010cm} \\
\mbox{} \hspace{+1.750cm}
$ {\displaystyle
\leq\;
K\!\;\!(n,\kappa,p)
\;
\biggl[\;\,
\|\, u_0 \;\!
\|_{{\scriptstyle L^{2^{m}p}(\mathbb{R}^{n})}}
  ^{\mbox{}^{\scriptstyle
\frac{\;\!{\scriptstyle p} \:-\;
{\scriptstyle n \:\!\kappa}\:\!/\;\!2^{m}}
{\;\!{\scriptstyle p} \:-\; 2 \;\!
{\scriptstyle n \:\!\kappa}\:\!/\;\!2^{\:\!j}}
\,
\frac{\,\!1 - 2^{-j + 1}}
{1 \,-\, 2^{-m}}
}}
\;
\biggr]
\;\times
} $ \\
\mbox{} \vspace{+0.050cm} \\
\mbox{} \hspace{+2.600cm}
$ {\displaystyle
\times \;
\biggl[\;\,
\mathbb{B}_{\mu}\!\;\!(0; t)^{\mbox{}^{\scriptstyle
\!\!\!
\frac{\;\!{\scriptstyle p} \:-\;
{\scriptstyle n \:\!\kappa}\:\!/\;\!2^{m}}
{{\scriptstyle p}}
\;\!
\left[\;
\frac{2\;\!{\scriptstyle n}}{\;\!2^{j}{\scriptstyle p}
\,-\, 2{\scriptstyle n \:\!\kappa}}
\:-\;
\frac{\scriptstyle n}{2^{m}{\scriptstyle p}
\,-\, {\scriptstyle n \:\!\kappa}} \,
\right]
}}
\;\!
\|\, u_0 \;\!
\|_{{\scriptstyle L^{p}(\mathbb{R}^{n})}}
  ^{\mbox{}^{\scriptstyle
\frac{\;\!{\scriptstyle p} \:-\;
{\scriptstyle n \:\!\kappa}\:\!/\;\!2^{m}}
{\;\!{\scriptstyle p} \:-\; 2 \;\!
{\scriptstyle n \:\!\kappa}\:\!/\;\!2^{\:\!j}}
\,
\frac{\,\!2^{-j + 1} -\, 2^{-m}}
{1 \,-\, 2^{-m}}
}}
\;
\biggr]
} $ \\
\mbox{} \vspace{+0.100cm} \\
\mbox{} \hspace{+0.250cm}
$ {\displaystyle
\leq\;
\theta \,\!\cdot\;\!
K\!\;\!(n,\kappa,p) \;
\|\, u_0 \;\!
\|_{\mbox{}_{\scriptstyle
L^{2^{m}p}(\mathbb{R}^{n})}}
\!\!\;\!+\,
(1 - \theta)
\,\!\cdot\;\!
K\!\;\!(n,\kappa,p) \;
\mathbb{B}_{\mu}\!\;\!(0; t)^{\mbox{}^{\scriptstyle
\!\!\!\;\!
\frac{\;\!{\scriptstyle n}
\,(\;\!1 \:-\; 2^{-\;\!m})}
{\;\!{\scriptstyle p} \,-\: {\scriptstyle n \:\!\kappa}}
}}
\:\!
\|\, u_0 \;\!
\|_{\mbox{}_{\scriptstyle
L^{p}(\mathbb{R}^{n})}}
^{\mbox{}^{\scriptstyle
\!\!\!\;\!
\frac{\;\!{\scriptstyle p} \:-\;
{\scriptstyle n \:\!\kappa}\:\!/\;\!2^{m}}
{{\scriptstyle p} \:-\; {\scriptstyle n \:\!\kappa}}
}}
} $ \\
\mbox{} \vspace{+0.125cm} \\
pela desigualdade de Young
(\cite{Evans2002}, p.$\;$622),
onde
$ \;\! \theta \in (\;\!0, 1\,\!) \;\!$
\'e dado por \\
\mbox{} \vspace{-0.500cm} \\
\begin{equation}
\notag
\theta
\;=\;
\frac{\,\! 1 - 2^{-\;\!j + 1}}
{1 - 2^{-\;\!m}}
\:
\frac{\,\!p \;\!-\;\! n \:\!\kappa/2^{m}}
{\,\!p \;\!-\;\! 2 \;\!n \:\!\kappa/2^{j}}
\:\!.
\end{equation}
\mbox{} \vspace{-0.150cm} \\
Assim,
de (3.16$b$) e (3.18),
segue que,
denotando
$ K \!\;\!\equiv K\!\;\!(n,\kappa,p) $: \\
\mbox{} \vspace{-0.100cm} \\
\mbox{} \hspace{+0.500cm}
$ {\displaystyle
\mathbb{U}_{2^{m}p}(0\:\!; \:\!t)
\;\leq\;
\max \;
\biggl\{\;\:\!
\|\, u_0 \;\!
\|_{\mbox{}_{\scriptstyle L^{2^{m}p}(\mathbb{R}^{n})}}
\:\!;
} $ \\
\mbox{} \vspace{-0.300cm} \\
\mbox{} \hspace{+3.250cm}
$ {\displaystyle
\theta \,\!\cdot\;\!
K \;
\|\, u_0 \;\!
\|_{\mbox{}_{\scriptstyle
L^{2^{m}p}(\mathbb{R}^{n})}}
\!\!\;\!+\,
(1 - \theta)
\,\!\cdot\;\!
K \;
\mathbb{B}_{\mu}\!\;\!(0; t)^{\mbox{}^{\scriptstyle
\!\!\!\;\!
\frac{\;\!{\scriptstyle n}
\,(\;\!1 \:-\; 2^{-\;\!m})}
{\;\!{\scriptstyle p} \,-\: {\scriptstyle n \:\!\kappa}}
}}
\:\!
\|\, u_0 \;\!
\|_{\mbox{}_{\scriptstyle
L^{p}(\mathbb{R}^{n})}}
^{\mbox{}^{\scriptstyle
\!\!\!\;\!
\frac{\;\!{\scriptstyle p} \:-\;
{\scriptstyle n \:\!\kappa}\:\!/\;\!2^{m}}
{{\scriptstyle p} \:-\; {\scriptstyle n \:\!\kappa}}
}}
\:\!;
} $ \\
\mbox{} \vspace{-0.000cm} \\
\mbox{} \hspace{+7.900cm}
$ {\displaystyle
K \;
\mathbb{B}_{\mu}\!\;\!(0; t)^{\mbox{}^{\scriptstyle
\!\!\!\;\!
\frac{\;\!{\scriptstyle n}
\,(\;\!1 \:-\; 2^{-\;\!m})}
{\;\!{\scriptstyle p} \,-\: {\scriptstyle n \:\!\kappa}}
}}
\:\!
\mathbb{U}_{p}(0\:\!;\:\!t)
^{\mbox{}^{\scriptstyle
\!\!\!\;\!
\frac{\;\!{\scriptstyle p} \:-\;
{\scriptstyle n \:\!\kappa}\:\!/\;\!2^{m}}
{{\scriptstyle p} \:-\; {\scriptstyle n \:\!\kappa}}
}}
\;
\biggr\}
} $ \\
\mbox{} \vspace{+0.150cm} \\
\mbox{} \hspace{+2.450cm}
$ {\displaystyle
\leq\;
K
\cdot\;
\max \;
\biggl\{\;\:\!
\|\, u_0 \;\!
\|_{\mbox{}_{\scriptstyle L^{2^{m}p}(\mathbb{R}^{n})}}
\:\!;
\;\;\!
\mathbb{B}_{\mu}\!\;\!(0; t)^{\mbox{}^{\scriptstyle
\!\!\!\;\!
\frac{\;\!{\scriptstyle n}
\,(\;\!1 \:-\; 2^{-\;\!m})}
{\;\!{\scriptstyle p} \,-\: {\scriptstyle n \:\!\kappa}}
}}
\:\!
\mathbb{U}_{p}(0\:\!;\:\!t)
^{\mbox{}^{\scriptstyle
\!\!\!\;\!
\frac{\;\!{\scriptstyle p} \:-\;
{\scriptstyle n \:\!\kappa}\:\!/\;\!2^{m}}
{{\scriptstyle p} \:-\; {\scriptstyle n \:\!\kappa}}
}}
\;
\biggr\}
} $ \\
\mbox{} \vspace{+0.125cm} \\
para todo
$ \;\! 0 \leq t < T_{\!\;\!\ast} \!\;\! $,
como afirmado.
Isso conclui a prova do Lema 3.7.
}
\mbox{} \hfill $\Box$ \\
%
\mbox{} \vspace{-0.500cm} \\

Do Lema 3.7,
pode-se finalmente obter
a estimativa (3.1),
sendo apenas necess\'ario
que se tome
$ \;\!m \rightarrow \infty \;\!$
em (3.19),
pelo fato de se ter \\
\mbox{} \vspace{-0.550cm} \\
\begin{equation}
\tag{3.20}
\lim_{q\,\rightarrow\,\infty}
\mathbb{U}_{\!\;\!q}(0\:\!;\:\!t)
\;=\;
\mathbb{U}_{\infty}\!\;\!(0\:\!;\:\!t).
\end{equation}
\mbox{} \vspace{-0.150cm} \\
Isso conclui a prova do Teorema 3.1,
que ocupou toda a discuss\~ao
da presente se\c c\~ao.

%
%
%
%
%
\mbox{} \vspace{-1.500cm} \\

{\bf 4. Condi\c c\~oes de exist\^encia global} \\

Nesta se\c c\~ao,
vamos aplicar a an\'alise acima
de modo a obter condi\c c\~oes
garantindo exist\^encia global
(i.e., $ \mbox{\small $T$}_{\!\;\!\ast} \!\;\! = \infty $)
das solu\c c\~oes $ u(\cdot,t) $
do problema (1.1)$\;\!$-$\;\!$(1.3),
ou seja, \linebreak
\mbox{} \vspace{-0.550cm} \\
\begin{equation}
\tag{4.1$a$}
u_t \,+\;
\mbox{div}\,
\bigl(\;\! \mbox{\boldmath $b$}(x,t,u) \, |\;\!u\;\!|^{\:\!\kappa} \:\! u \;\!
\bigr)
\:+\:
\mbox{div} \,
\mbox{\boldmath $f$}(t,u)
\;=\;
\mbox{div}\,\bigl(\;\!
A(x,t,u) \;\! \nabla u \;\!\bigr),
%
%
\end{equation}
\mbox{} \vspace{-0.900cm} \\
\begin{equation}
\tag{4.1$b$}
u(\cdot,0) \,=\,
u_0 \in L^{1}(\mathbb{R}^{n})
\cap L^{\infty}(\mathbb{R}^{n}),
\end{equation}
\mbox{} \vspace{-0.175cm} \\
sendo
$ A $ (tensor difusivo) e
$ \;\!\mbox{\boldmath $b$}, \mbox{\boldmath $f$} \!\;\!$
(campos vetoriais)
suaves
satisfazendo (1.2) e (1.3).
Um exemplo de condi\c c\~oes
de exist\^encia global
\'e dado no Teorema 4.1 a seguir,
onde
(ver (1.7), Se\c c\~ao 1) \\
\mbox{} \vspace{-0.950cm} \\
\begin{equation}
\tag{4.2}
\mathbb{B}_{\!\;\!\mu}\!\;\!(0\:\!; \:\!\infty)
\;=\;\;\!
\sup_{t \,\geq\,0} \;
\mbox{\small $ {\displaystyle
\frac{B(t)}{\mu(t)}
} $}
\:\leq\;\! \infty.
\end{equation}
\mbox{} \vspace{-0.100cm} \\
%
%
%
%
%
%
\mbox{} \hspace{-0.800cm}
\fbox{%
\begin{minipage}[t]{16.000cm}
\mbox{} \vspace{-0.400cm} \\
\mbox{} \hspace{+0.300cm}
\begin{minipage}[t]{15.000cm}
\mbox{} \vspace{+0.100cm} \\
{\bf Teorema 4.1.}
\textit{%
Na nota\c c\~ao acima,
tem-se,
a respeito das solu\c c\~oes}
{\em do problema}$\;\!$: \\
\mbox{} \vspace{-0.250cm} \\
({\em i\/})
\textit{%
se
$ \;\! 0 \leq \kappa < 1/n $,
ent\~ao
$ \:\!u(\cdot,t) $
est\'a definida para todo
$ \,t > 0 \!\!\;\!$
}
\mbox{({\em para todo dado} $\:\!u_0$)}; \\
\mbox{} \vspace{-0.275cm} \\
({\em ii\/})
\textit{%
se $ \;\!\kappa = 1/n $,
as solu\c c\~oes s\~ao globais
sempre
que}
$ {\displaystyle
\;\!
\|\, u_0 \;\!
\|_{\mbox{}_{\scriptstyle L^{1}(\mathbb{R}^{n})}}
\!\:\!\leq\;\!
\mathbb{B}_{\!\;\!\mu}\!\;\!(0\:\!; \infty)
^{\mbox{}^{-\,\mbox{\scriptsize $n$}}}
\!\;\!
} $; \\
\mbox{} \vspace{-0.275cm} \\
({\em iii\/})
\textit{%
se $ \;\!\kappa > 1/n $,
as solu\c c\~oes s\~ao globais
sempre que
o dado inicial
satisfizer
} \\
\mbox{} \vspace{-0.750cm} \\
\begin{equation}
\tag{4.3}
\|\, u_0 \;\!
\|_{\mbox{}_{\scriptstyle L^{1}(\mathbb{R}^{n})}}
\:\!
\|\, u_0 \;\!
\|_{\mbox{}_{\scriptstyle L^{\infty}(\mathbb{R}^{n})}}
  ^{\:\! n \:\!\kappa \,-\,1}
\:\!\leq\,
\bigl\{\;\! n \:\! \kappa \;
\mathbb{B}_{\!\;\!\mu}\!\;\!(0\:\!; \infty)
\,\bigr\}^{\!\:\!-\,\mbox{\scriptsize $n$}}
\end{equation}
\mbox{} \vspace{-0.575cm} \\
\end{minipage}
\end{minipage}
}
%
%
\nl
\nl
%
%
{\small
{\bf Prova:}
O caso ({\em i\/}),
j\'a considerado no {\sc Teorema A}
da Se\c c\~ao~1,
\'e consequ\^encia imediata
da propriedade (2.4) e
do Teorema 3.1
(tomando-se $ \;\!p = 1 $ em (3.1)).
Nos casos
({\em ii\/}) e ({\em iii\/}), \linebreak
podemos proceder
do seguinte modo.
Da desigualdade (3.8)
\mbox{[}$\,$reescrita em termos de $ u(\cdot,t) $,
usando (3.4), (3.5)$\,$\mbox{]},
obt\'em-se,
considerando
$ \;\! q = 2 \;\! n \:\! \kappa \geq 2 $, \\
\mbox{} \vspace{+0.040cm} \\
\mbox{} \hspace{+0.100cm}
$ {\displaystyle
\frac{d}{d\:\!t} \;
\|\, u(\cdot,t) \,
\|_{\mbox{}_{\scriptstyle L^{2\;\!n\:\!\kappa}(\mathbb{R}^{n})}}
  ^{\:\! 2 \;\! n \:\!\kappa}
\!+\,
 2\;\!n \:\!\kappa \,
(2\;\!n \:\!\kappa - 1) \;
\mu(t)
\!
\int_{\mathbb{R}^{n}}
\!\!\:\!
|\, u(x,t) \,|^{\:\!2\;\!n \:\!\kappa \;\!-\;\!2}
\,
|\, \nabla u \,|^{\:\!2}
\: dx
} $
\mbox{} \hfill (4.4) \\
\mbox{} \vspace{+0.100cm} \\
\mbox{} \hspace{+0.900cm}
$ {\displaystyle
\leq\;
n \:\! \kappa \;
\frac{B(t)}{\mu(t)} \;
\|\, u(\cdot,t) \,
\|_{\mbox{}_{\scriptstyle L^{n\:\!\kappa}(\mathbb{R}^{n})}}
  ^{\:\!\kappa}
\;\!
\Bigl\{\,
 2\;\!n \:\!\kappa \,
(2\;\!n \:\!\kappa - 1) \;
\mu(t)
\!
\int_{\mathbb{R}^{n}}
\!\!\:\!
|\, u(x,t) \,|^{\:\!2\;\!n \:\!\kappa \;\!-\;\!2}
\,
|\, \nabla u \,|^{\:\!2}
\: dx
\;\!\Bigr\}
} $ \\
\mbox{} \vspace{+0.100cm} \\
\mbox{} \hfill
$ {\displaystyle
\leq\;
n \:\! \kappa \;
\mathbb{B}_{\!\;\!\mu}\!\;\!(0\:\!; \infty)
\:
\|\, u(\cdot,t) \,
\|_{\mbox{}_{\scriptstyle L^{n\:\!\kappa}(\mathbb{R}^{n})}}
  ^{\:\!\kappa}
\;\!
\Bigl\{\,
 2\;\!n \:\!\kappa \,
(2\;\!n \:\!\kappa - 1) \;
\mu(t)
\!
\int_{\mathbb{R}^{n}}
\!\!\:\!
|\, u(x,t) \,|^{\:\!2\;\!n \:\!\kappa \;\!-\;\!2}
\,
|\, \nabla u \,|^{\:\!2}
\: dx
\;\!\Bigr\}
} $ \\
\mbox{} \vspace{+0.100cm} \\
para todo
$ \;\!t \in (\:\!0, \:\!T_{\!\;\!\ast}\!\;\!) \setminus
E_{2\;\!n \:\!\kappa} $,
de modo que temos
$ {\displaystyle
\,
\|\, u(\cdot,t) \,
\|_{\mbox{}_{\scriptstyle L^{2 \;\!n\:\!\kappa}(\mathbb{R}^{n})}}
\!
} $
decrescente em
$\;\![\,0, \:\!T\;\!] $

\mbox{} \vspace{-0.750cm} \\
sempre que tivermos \\
\mbox{} \vspace{-0.650cm} \\
\begin{equation}
\tag{4.5}
n \:\! \kappa \:\,\!
\mathbb{B}_{\!\;\!\mu}\!\;\!(0\:\!; \infty)
\:
\|\, u(\cdot,t) \,
\|_{\mbox{}_{\scriptstyle L^{n\:\!\kappa}(\mathbb{R}^{n})}}
  ^{\:\!\kappa}
\;\!\leq\:
1,
\qquad
\forall \;\,
t \in [\,0, \;\!T\;\!].
\end{equation}
\mbox{} \vspace{-0.250cm} \\
No caso ({\em ii\/}),
esta condi\c c\~ao \'e simplesmente
$ {\displaystyle
\;\!
\mathbb{B}_{\mu}\!\;\!(0\:\!; \infty)^{\:\!n}
\:\!
\|\, u(\cdot,t) \,
\|_{\mbox{}_{\scriptstyle L^{1}(\mathbb{R}^{n})}}
\!\;\!\leq\:\!
1
} $,
que
\'e automati\-ca\-mente
satisfeita
em qualquer intervalo
$ \,\![\,0, \:\!T\:\!]\,\! $
se for satisfeita em $ \:\!t = 0 $,
devido a (2.4). \linebreak
Isso mostra que
$ {\displaystyle
\;\!
\|\, u(\cdot,t) \,
\|_{\mbox{}^{\scriptstyle L^{2}(\mathbb{R}^{n})}}
\!\;\!
} $
\'e monotonicamente
decrescente em
$ \,\![\,0, \:\!T_{\!\;\!\ast}\!\;\!) \,\! $
caso se tenha
$ {\displaystyle
\mathbb{B}_{\mu}\!\;\!(0\:\!; \infty)^{\:\!n}
\:\!
\|\, u_0 \;\!
\|_{\mbox{}_{\scriptstyle L^{1}(\mathbb{R}^{n})}}
\!\;\!\leq\:\!
1
} $.
Pelo Teorema 3.1,
$ {\displaystyle
\|\, u(\cdot,t) \,
\|_{\mbox{}_{\scriptstyle L^{\infty}(\mathbb{R}^{n})}}
\!\!\;\!
} $
\'e controlada
(como
$ \:\!n \:\!\kappa = 1 $) \linebreak
\mbox{} \vspace{-0.550cm} \\
por
$ {\displaystyle
\;\!
\|\, u(\cdot,t) \,
\|_{\mbox{}^{\scriptstyle L^{2}(\mathbb{R}^{n})}}
\!\;\!
} $,
e,
assim sendo,
$ {\displaystyle
\|\, u(\cdot,t) \,
\|_{\mbox{}_{\scriptstyle L^{\infty}(\mathbb{R}^{n})}}
\!\!\;\!
} $
tem de permanecer limitada
em qualquer intervalo limitado.
Logo, n\~ao se pode ter
$ \;\!T_{\!\;\!\ast} < \infty \:\!$
neste caso,
como afirmado em~({\em ii\/}). \linebreak
\mbox{} \vspace{-0.700cm} \\

Finalmente,
consideremos
o caso ({\em iii\/}).
Observando que
(por interpola\c c\~ao)
tem-se \\
\mbox{} \vspace{-0.100cm} \\
\mbox{} \hspace{+0.600cm}
$ {\displaystyle
n \:\! \kappa \:\,\!
\mathbb{B}_{\!\;\!\mu}\!\;\!(0\:\!; \infty)
\:
\|\, u(\cdot,t) \,
\|_{\mbox{}_{\scriptstyle L^{n\:\!\kappa}(\mathbb{R}^{n})}}
  ^{\:\!\kappa}
\:\!\leq\,\:\!
n \:\! \kappa \:\,\!
\mathbb{B}_{\!\;\!\mu}\!\;\!(0\:\!; \infty)
\:
\|\, u(\cdot,t) \,
\|_{\mbox{}_{\scriptstyle L^{1}(\mathbb{R}^{n})}}
  ^{\mbox{}^{\scriptstyle
  \frac{1}{\,\!2 \;\!{\scriptstyle n} \;\!-\;\!1/{\scriptstyle \kappa}}}}
\|\, u(\cdot,t) \,
\|_{\mbox{}_{\scriptstyle L^{2\:\!n\:\!\kappa}(\mathbb{R}^{n})}}
  ^{\mbox{}^{\scriptstyle
  \frac{2 \;\!{\scriptstyle n \:\!\kappa} \;\!-\;\!2}
       {\,\!2 \;\!{\scriptstyle n} \;\!-\;\!1/{\scriptstyle \kappa}}}}
} $
\mbox{} \hfill (4.6) \\
\mbox{} \vspace{-0.000cm} \\
e tamb\'em
(novamente, por interpola\c c\~ao) \\
\mbox{} \vspace{-0.100cm} \\
\mbox{} \hspace{-0.200cm}
$ {\displaystyle
n \:\! \kappa \:\,\!
\mathbb{B}_{\!\;\!\mu}\!\;\!(0\:\!; \infty)
\:
\|\, u(\cdot,t) \,
\|_{\mbox{}_{\scriptstyle L^{1}(\mathbb{R}^{n})}}
  ^{\mbox{}^{\scriptstyle
  \frac{1}{\,\!2 \;\!{\scriptstyle n} \;\!-\;\!1/{\scriptstyle \kappa}}}}
\|\, u(\cdot,t) \,
\|_{\mbox{}_{\scriptstyle L^{2\:\!n\:\!\kappa}(\mathbb{R}^{n})}}
  ^{\mbox{}^{\scriptstyle
  \frac{2 \;\!{\scriptstyle n \:\!\kappa} \;\!-\;\!2}
       {\,\!2 \;\!{\scriptstyle n} \;\!-\;\!1/{\scriptstyle \kappa}}}}
\!\;\!<\:
n \:\! \kappa \:\,\!
\mathbb{B}_{\!\;\!\mu}\!\;\!(0\:\!; \infty)
\:
\|\, u(\cdot,t) \,
\|_{\mbox{}_{\scriptstyle L^{1}(\mathbb{R}^{n})}}
  ^{\mbox{}^{\scriptstyle \frac{1}{\scriptstyle n} }}
\:\!
\|\, u(\cdot,t) \,
\|_{\mbox{}_{\scriptstyle L^{\infty}}}
  ^{\mbox{}^{\scriptstyle \kappa \;\!-\;\! \frac{1}{\scriptstyle n} }}
} $ \\
\mbox{} \vspace{-0.425cm} \\
\mbox{} \hfill (4.7) \\
\mbox{} \vspace{-0.55cm} \\
para todo
$\;\!t \in [\,0, T_{\!\;\!\ast}\!\;\!) $,
obtemos,
por (4.3) e (4.7), \\
\mbox{} \vspace{-0.600cm} \\
\begin{equation}
\tag{4.8}
n \:\! \kappa \:\,\!
\mathbb{B}_{\!\;\!\mu}\!\;\!(0\:\!; \infty)
\:
\|\, u(\cdot,t) \,
\|_{\mbox{}_{\scriptstyle L^{1}(\mathbb{R}^{n})}}
  ^{\mbox{}^{\scriptstyle
  \frac{1}{\,\!2 \;\!{\scriptstyle n} \;\!-\;\!1/{\scriptstyle \kappa}}}}
\|\, u(\cdot,t) \,
\|_{\mbox{}_{\scriptstyle L^{2\:\!n\:\!\kappa}(\mathbb{R}^{n})}}
  ^{\mbox{}^{\scriptstyle
  \frac{2 \;\!{\scriptstyle n \:\!\kappa} \;\!-\;\!2}
       {\,\!2 \;\!{\scriptstyle n} \;\!-\;\!1/{\scriptstyle \kappa}}}}
\!\;\!<\:
1
\end{equation}
\mbox{} \vspace{-0.200cm} \\
para todo
$\;\!t \in [\,0, T_{\!\;\!\ast}\!\;\!) $
suficientemente pr\'oximo de zero.
Afirmamos que (4.8) acima
tem de ser verdadeira
para todo
$\;\!t \in [\,0, T_{\!\;\!\ast}\!\;\!) $.
\mbox{[}$\,$De fato,
se n\~ao fosse,
existiria
$\;\!T_{\mbox{}_{1}} \!\in (\:\!0, T_{\!\;\!\ast}\!\;\!) $
tal que
se teria \\
\mbox{} \vspace{-0.400cm} \\
\mbox{} \hspace{+1.400cm}
$ {\displaystyle
n \:\! \kappa \:\,\!
\mathbb{B}_{\!\;\!\mu}\!\;\!(0\:\!; \infty)
\:
\|\, u(\cdot,t) \,
\|_{\mbox{}_{\scriptstyle L^{1}(\mathbb{R}^{n})}}
  ^{\mbox{}^{\scriptstyle
  \frac{1}{\,\!2 \;\!{\scriptstyle n} \;\!-\;\!1/{\scriptstyle \kappa}}}}
\|\, u(\cdot,t) \,
\|_{\mbox{}_{\scriptstyle L^{2\:\!n\:\!\kappa}(\mathbb{R}^{n})}}
  ^{\mbox{}^{\scriptstyle
  \frac{2 \;\!{\scriptstyle n \:\!\kappa} \;\!-\;\!2}
       {\,\!2 \;\!{\scriptstyle n} \;\!-\;\!1/{\scriptstyle \kappa}}}}
\!\;\!<\:
1,
\quad \;\;\,
\forall \;\,
0 \;\!\leq\;\! t \;\!<\;\! T_{\mbox{}_{1}}\!\;\!
} $,
\mbox{} \hfill (4.9) \\
\mbox{} \vspace{+0.100cm} \\
enquanto
$ {\displaystyle
\,
n \:\! \kappa \:\,\!
\mathbb{B}_{\!\;\!\mu}\!\;\!(0\:\!; \infty)
\:
\|\, u(\cdot,T_{\mbox{}_{1}}) \,
\|_{\scriptstyle L^{1}(\mathbb{R}^{n})}
  ^{\mbox{}^{\scriptstyle
  \frac{1}{\,\!2 \;\!{\scriptstyle n} \;\!-\;\!1/{\scriptstyle \kappa}}}}
\|\, u(\cdot,T_{\mbox{}_{1}}) \,
\|_{\scriptstyle L^{2\:\!n\:\!\kappa}(\mathbb{R}^{n})}
  ^{\mbox{}^{\scriptstyle
  \frac{2 \;\!{\scriptstyle n \:\!\kappa} \;\!-\;\!2}
       {\,\!2 \;\!{\scriptstyle n} \;\!-\;\!1/{\scriptstyle \kappa}}}}
\!\;\!=\;\!
1
} $.
Por (4.6),
ter\'\i amos ent\~ao
(4.5) satisfeita para $ T = T_{\mbox{}_{1}} $,
de modo que, por (4.4),
$ {\displaystyle
\;\!
\|\, u(\cdot,t) \,
\|_{L^{2\;\!n\;\!\kappa}(\mathbb{R}^{n})}
\!\;\!
} $
seria decrescente
no intervalo
$ \;\![\,0, \;\!T_{\mbox{}_{1}}\:\!] $.
Assim,
ter\'\i amos \\
\mbox{} \vspace{-0.100cm} \\
\mbox{} \hspace{+2.875cm}
$ {\displaystyle
1
\;\;\!=\;\;\!
n \:\! \kappa \:\,\!
\mathbb{B}_{\!\;\!\mu}\!\;\!(0\:\!; \infty)
\:
\|\, u(\cdot,T_{\mbox{}_{1}}) \,
\|_{\scriptstyle L^{1}(\mathbb{R}^{n})}
  ^{\mbox{}^{\scriptstyle
  \frac{1}{\,\!2 \;\!{\scriptstyle n} \;\!-\;\!1/{\scriptstyle \kappa}}}}
\|\, u(\cdot,T_{\mbox{}_{1}}) \,
\|_{\scriptstyle L^{2\:\!n\:\!\kappa}(\mathbb{R}^{n})}
  ^{\mbox{}^{\scriptstyle
  \frac{2 \;\!{\scriptstyle n \:\!\kappa} \;\!-\;\!2}
       {\,\!2 \;\!{\scriptstyle n} \;\!-\;\!1/{\scriptstyle \kappa}}}}
} $ \\
\mbox{} \vspace{-0.050cm} \\
\mbox{} \hspace{+3.350cm}
$ {\displaystyle
\leq\;
n \:\! \kappa \:\,\!
\mathbb{B}_{\!\;\!\mu}\!\;\!(0\:\!; \infty)
\:
\|\, u_0 \;\!
\|_{\scriptstyle L^{1}(\mathbb{R}^{n})}
  ^{\mbox{}^{\scriptstyle
  \frac{1}{\,\!2 \;\!{\scriptstyle n} \;\!-\;\!1/{\scriptstyle \kappa}}}}
\|\, u_0 \;\!
\|_{\scriptstyle L^{2\:\!n\:\!\kappa}(\mathbb{R}^{n})}
  ^{\mbox{}^{\scriptstyle
  \frac{2 \;\!{\scriptstyle n \:\!\kappa} \;\!-\;\!2}
       {\,\!2 \;\!{\scriptstyle n} \;\!-\;\!1/{\scriptstyle \kappa}}}}
} $
\mbox{} \hfill
\mbox{[}$\,$por (2.4)$\,$\mbox{]} \\
\mbox{} \vspace{-0.050cm} \\
\mbox{} \hspace{+3.350cm}
$ {\displaystyle
\leq\;
n \:\! \kappa \:\,\!
\mathbb{B}_{\!\;\!\mu}\!\;\!(0\:\!; \infty)
\:
\|\, u_0 \;\!
\|_{\scriptstyle L^{1}(\mathbb{R}^{n})}
  ^{\mbox{}^{\scriptstyle
    \frac{1}{\scriptstyle n}}}
\:\!
\|\, u_0 \;\!
\|_{\scriptstyle L^{\infty}(\mathbb{R}^{n})}
  ^{\mbox{}^{\scriptstyle
    {\scriptstyle \kappa} \;\!-\;\! \frac{1}{\scriptstyle n}}}
\;<\, 1
} $.
\mbox{} \hfill
\mbox{[}$\,$por (4.7), (4.3)$\,$\mbox{]} \\
\mbox{} \vspace{+0.050cm} \\
Esta contradi\c c\~ao mostra que
(4.8) tem ser v\'alida
para todo $ \;\!0 \leq t < T_{\!\;\!\ast} \!\;\!$,
como afirmado.$\,$\mbox{]}
Sendo (4.8)
verdadeira para todo
$ \;\!0 \leq t < T_{\!\;\!\ast} \!\;\!$,
resulta ent\~ao,
por (4.6),
que

\mbox{} \vspace{-0.750cm} \\
\begin{equation}
\tag{4.10}
n \:\! \kappa \:\,\!
\mathbb{B}_{\!\;\!\mu}\!\;\!(0\:\!; \infty)
\:
\|\, u(\cdot,t) \,
\|_{\mbox{}_{\scriptstyle L^{n\:\!\kappa}(\mathbb{R}^{n})}}
  ^{\:\!\kappa}
\,<\; 1,
\qquad
\forall \;\,
t \in [\,0, \;\!T_{\!\;\!\ast}).
\end{equation}
\mbox{} \vspace{-0.150cm} \\
Isso mostra,
por (4.4),
que
$ {\displaystyle
\;\!
\|\, u(\cdot,t) \,
\|_{\mbox{}_{\scriptstyle L^{2\:\!n\;\!\kappa}(\mathbb{R}^{n})}}
\!\,\!
} $
\'e decrescente
no intervalo de exist\^encia
$\;\![\,0, \;\!T_{\!\;\!\ast}) $. \linebreak
Portanto,
pelo Teorema 3.1
\mbox{[}$\,$aplicado a $ \;\!p = 2\:\!n\:\!\kappa \,$\mbox{]},
temos de ter
$ {\displaystyle
\;\!
\|\, u(\cdot,t) \,
\|_{\mbox{}_{\scriptstyle L^{\infty}(\mathbb{R}^{n})}}
\!\:\!
} $
limitada
em todo o intervalo
$\;\![\,0, \;\!T_{\!\;\!\ast}) $,
de modo que,
como no caso ({\em ii\/}),
$ T_{\!\;\!\ast} $
n\~ao pode ser finito.
}
\mbox{} \hfill $\Box$ \\
%
%

%

\newpage

%
%
%
%
%
\mbox{} \vspace{-0.750cm} \\
%
%
%
%
%
%

%
%

%
%


\begin{thebibliography}{999}

%
%

\bibitem{Agueh2008}
\textsc{M. Agueh},
{\it Gagliardo-Nirenberg inequalities involving
the gradient $L^{2}\!$-norm},
C. R. Acad. Sci. Paris, Ser. I {\bf 346} (2008),
757-762.

\bibitem{BandleBrunner1998}
{\sc C. Bandle and H. Brunner},
{\it Blow-up in diffusion equations: a survey},
J. Comp. Appl. Math. {\bf 97} (1998), 3-22.

\bibitem{BarrionuevoOliveiraZingano2014}
\textsc{J.$\;$A.$\;$Barrionuevo, L.$\;$S.$\;$Oliveira
and P.$\;$R.$\;$Zingano},
{\it General asymptotic supnorm estimates for solutions
of one-dimensional advection-diffusion equations
in heterogeneous media},
Intern. J. Partial Diff. Equations (2014), 1-8
(freely available at:
http://www.hindawi.com/journals/ijpde/2014/450417).

\bibitem{Benameur2010}
\textsc{J. Benameur},
{\it On the blow-up criterion of
\mbox{\it 3D} Navier-Stokes equations},
J. Math. Anal. Appl. {\bf 371} (2010),
719-727.

\bibitem{BenameurSelmi2012}
{\sc J. Benameur and R. Selmi},
{\em Long time decay to the Leray solution of the
two-dimensional Navier-Stokes equations},
Bull. London Math. Soc. {\bf 44} (2012),
1001-1019.

\bibitem{BrazMeloZingano2015}
\textsc{P.$\;$Braz e Silva, W.$\;$G.$\;$Melo and
P.$\;$R.$\;$Zingano},
{\it An asymptotic supnorm estimate for solutions of
\mbox{\em 1-D} systems of convection-diffusion equations},
J. Diff. Eqs. {\bf 258} (2015),
2806-2822.

\bibitem{BrazLorenzMeloZingano2014}
{\sc P. Braz e Silva, J. Lorenz, W. G. Melo and P. R. Zingano},
{\it On the large time approximation of the Navier-Stokes equations
in $\:\!\mathbb{R}^{n}\!$ by Stokes flows\/}
(in preparation).

\bibitem{BrazSchutzZingano2013}
{\sc P. Braz e Silva, L. Sch\"utz and P. R. Zingano},
{\it On some energy inequalities and supnorm estimates
for advection-diffusion equations in $\mathbb{R}^{n}\!\;\!$},
Nonlin. Anal. {\bf 93} (2013),
90-96.

\bibitem{CarlenLoss1993}
{\sc E. A. Carlen and M. Loss},
{\it Sharp constant in Nash's inequality},
Intern. Math. Res. Notices, 1993,
213-215.

\bibitem{Chagas2015}
{\sc J.$\;$Q.$\;$Chagas},
\textit{Some results for doubly nonlinear
equations with advection},
Universidade Federal do Rio Grande do Sul,
Brazil
(work in progress).

\bibitem{Constantin2001}
{\sc P. Constantin},
{\it Some open problems and research directions
in the mathematical study of fluid dynamics},
in: B. Engquist and W. Schmid (Eds.),
{\sf Mathematics Unlimited --- 2001 and Beyond},
Springer, New York, 2001,
pp.$\;$353-360.

\bibitem{DengLevine2000}
{\sc K. Deng and H. A. Levine},
{\it The role of critical exponents in blow-up theorems:
the sequel},
J. Math. Anal. Appl. {\bf 243} (2000),
85-126.

\bibitem{DiehlFabrisZingano2014}
{\sc N. L. Diehl, L. Fabris and P. R. Zingano},
{\it Comparison results for smooth solutions of
quasilinear parabolic equations},
Adv. Diff. Eqs. Control Proc. {\bf 14} (2014),
11-22.

\bibitem{Diehl2015}
{\sc N.$\;$L.$\;$Diehl},
\textit{Some results for unsigned
porous medium equations with advec\-tion},
Universidade Federal do Rio Grande do Sul,
Brazil
(work in progress).

\bibitem{DunfordSchwartz1963}
{\sc N.$\;$Dunford and J.$\;$T.$\;$Schwartz},
{\sf Linear Opeartors}, vol.$\;$2,
Interscience, New York, 1963.

\bibitem{Evans2002}
{\sc L. C. Evans},
{\sf Partial Differential Equations},
American Mathematical Society, Providence, 2002.

\bibitem{Fabris2013}
{\sc L.$\;$Fabris},
\textit{On global existence and supnorm results
for nonnegative solutions of the porous medium equation
with arbitrary advection terms\/} (Portuguese),
PhD Thesis,
Universidade Federal do Rio Grande do Sul,
Porto Alegre, RS, Brazil, October 2013
(available at:
http://hdl.handle.net/10183/88277).

\bibitem{Fefferman2006}
{\sc C. L. Fefferman},
{\it Existence and smoothness of the Navier--Stokes equations}, \linebreak
in: A.$\;$M.$\;$Jaffe and A.$\;$J.$\;$Wiles (Eds.),
{\sf The Millenium Prize Problems},
American \linebreak
Mathematical
Society, Providence, 2006,
pp.$\;$57-70
($\,\!$freely available at http:// \linebreak
www.claymath.org/millenium/Navier-Stokes\_Equations/NavierStokes.pdf.)

\bibitem{Friedman1969}
{\sc A. Friedman},
{\sf Partial Differential Equations},
Holt, Rinehart and Winston,
New York, 1969.

\bibitem{Fujita1966}
{\sc H.$\;$Fujita},
{\it On the blowing up of solutions of the
Cauchy problem for $ \,\!u_t \!= \Delta u + u^{1 + \alpha}\!\;\!$},
J. Fac. Sci. Univ. Tokyo,
Sect. IA Math., {\bf 13} (1966), 109-124.

\bibitem{Galdi2000}
{\sc G. P. Galdi},
{\it An introduction to the Navier--Stokes
initial--boundary problem},
in:
G. P. Galdi, J. G. Heywood and R. Rannacher (Eds.),
{\sf Fundamental Directions in Mathematical Fluid Dynamics},
Birkhauser, Basel, 2000, pp.$\;$1-70.

\bibitem{Guidolin2015}
{\sc P.$\;$L.$\;$Guidolin},
\textit{Some results for p-Laplacian
evolution equations with advection},
Universidade Federal do Rio Grande do Sul,
Brazil
(work in progress).

\bibitem{Guterres2014}
{\sc R. H. Guterres},
{\it On singular integral operators,
with applications to partial differential equations\/} (Portuguese),
M. Sc. Dissertation,
Universidade Federal do Rio Grande do Sul,
Porto Alegre, RS, Brazil,
September 2014
(available at http://www.lume.ufrgs.br/handle/10183/115175).

\bibitem{Hayakawa1973}
{\sc K. Hayakawa},
{\it On the non-existence of global solutions
of some semilinear parabolic differential equations},
Proced. Japan Acad. Sci., Ser. A,
{\bf 49} (1973), 503-505.

\bibitem{KajikiyaMiyakawa1986}
{\sc R. Kajikiya and T. Miyakawa},
{\it On the $L^2$ decay of weak solutions of the
Navier--Stokes equations in $\mathbb{R}^{n}\!\;\!$},
Math. Z. {\bf 192} (1986), 135-148.

\bibitem{Kato1984}
{\sc T. Kato},
{\it Strong $L^p$-solutions of the Navier--Stokes equations
in $\mathbb{R}^{m}\!$, with applications to weak solutions},
Math. Z. {\bf 187} (1984), 471-480.

\bibitem{KreissHagstromLorenzZingano2002}
{\sc H.-O. Kreiss, T. Hagstrom, J. Lorenz and P. R. Zingano},
{\it Decay in time of the solutions of the Navier--Stokes equations
for incompressible flows},
unpublished note,
University of New Mexico, Albuquerque, NM, 2002.

\bibitem{KreissHagstromLorenzZingano2003}
{\sc H.-O. Kreiss, T. Hagstrom, J. Lorenz and P. R. Zingano},
{\it Decay in time of incompressible flows},
J. Math. Fluid Mech. {\bf 5} (2003), 231-244.

\bibitem{KreissLorenz1989}
{\sc H.-O. Kreiss and J. Lorenz},
{\sf Initial--boundary value problems and the Navier--Stokes equations},
Academic Press,
New York, 1989.
(Reprinted in the series SIAM Classics in Applied Mathematics,
 Vol. 47, 2004.)

\bibitem{LadyzhenskayaSolonnikovUralceva1968}
\textsc{O.$\;$A.$\;$Ladyzhenskaya, V.$\;$A.$\;$Solonnikov
and N.$\;$N.$\;$Uralceva},
{\sf Linear and Quasilinear Equations of Parabolic Type},
American Mathematical Society, Providence, 1968.

\bibitem{Levine1990}
{\sc H. A. Levine},
{\it The role of critical exponents in blow-up theorems},
SIAM Rev. {\bf 32} (1990),
262-288.

\bibitem{LorenzZingano2012}
{\sc J.$\;$Lorenz and P.$\;$R.$\;$Zingano},
{\it The Navier-Stokes equations for incompressible flows:
solution properties at potential blow-up times},
Universidade Federal do Rio Grande do Sul,
Porto Alegre, RS, July 2012
(available at
http://www.arXiv.org/abs/1503.01767).

\bibitem{Leray1934}
{\sc J. Leray},
{\it Essai sur le mouvement d'un fluide visqueux emplissant l'espace},
Acta Math. {\bf 63} (1934), 193-248.

\bibitem{Masuda1984}
{\sc K.$\;$Masuda},
{\it Weak solutions of the Navier-Stokes equations},
T\^ohoku Math. Journal {\bf 36} (1984), 623-646.

\bibitem{Melo2011}
\textsc{W.$\;$G.$\;$Melo},
{\it A-priori estimates for systems of
advection-diffusion equations\/} (Portuguese),
PhD Thesis,
Universidade Federal de Pernambuco,
Recife, PE, Brazil, 2011
(available at:
http://hdl.handle.net/123456789/7355).

\bibitem{Nash1958}
{\sc J.$\;$Nash},
{\it Continuity of solutions of parabolic and elliptic equations},
Amer. J. Math. {\bf 80} (1958),
931-954.

\bibitem{Oliveira2013}
\textsc{L.$\;$S.$\;$Oliveira},
{\it Two results in classical analysis\/} (Portuguese),
PhD Thesis,
Universidade Federal do Rio Grande do Sul,
Porto Alegre, RS, Brazil, 2013
(available at:
http://hdl.handle.net/10183/70212).

\bibitem{OliverTiti2000}
\textsc{M. Oliver and E. S. Titi},
{\it Remark on the rate of decay of higher order derivatives
for solutions to the Navier-Stokes equations in $ \mathbb{R}^{n} \!$},
J. Funct. Anal. {\bf 172} (2000), 1-18.

\bibitem{Perusato2014}
{\sc C. F. Perusato},
{\it On the Leray's problem for the Navier-Stokes equations
and some generalizations\/} (Portuguese),
M. Sc. Dissertation,
Universidade Federal do Rio Grande do Sul,
Porto Alegre, RS, Brazil,
September 2014
(available at http://www.lume.ufrgs.br/handle/10183/115208).

\bibitem{Pinsky1997}
{\sc R. G. Pinsky},
{\it Existence and non-existence of global solutions for
$u_t = \Delta u + a(x) \;\!u^p \!\;\!$ in $ \mathbb{R}^{d}\!$},
J. Diff. Eqs. {\bf 133} (1997), 152-177.

\bibitem{QuittnerSouplet2007}
{\sc P.$\;$Quittner and Ph.$\;$Souplet},
{\sf Superlinear Parabolic Problems:
blow-up, global existence and steady states},
Birkh\"auser, Basel, 2007.

\bibitem{RobinsonSadowskiSilva2012}
{\sc J. C. Robinson, W. Sadowski and R. P. Silva},
{\it Lower bounds on blow up solutions of the
three-dimensional Navier-Stokes equations
in homogeneous Sobolev spaces},
J. Math. Phys.
{\bf 53} (2012), no.$\;$11, 115618, 15 pp.

\bibitem{SchonbekWiegner1996}
{\sc M. E. Schonbek and M. Wiegner},
{\it On the decay of higher-order norms of the
solutions of Navier-Stokes equations},
Proc. Roy. Soc. Edinburgh {\bf 126A} (1996), 677-685.

\bibitem{Schutz2008}
{\sc L. Sch\"utz},
{\it Some results on advection-diffusion equations,
with applications to the Navier-Stokes equations\/}
(in Portuguese), Doctorate Thesis,
Graduate Program in Mathematics
(http://www.mat.ufrgs.br/${\sim}$ppgmat),
Universidade Federal do Rio Grande do Sul,
Porto Alegre, RS, June 2008
(available at http://hdl.handle.net/10183/13714).

\bibitem{SchutzZinganoZingano2014}
{\sc L. Sch\"utz, J. P. Zingano and P. R. Zingano},
{\it On the supnorm form of Leray's problem
for the incompressible Navier-Stokes equations\/}
(submitted).

\bibitem{Seregin2012}
{\sc G. Seregin},
{\it A certain necessary condition of potential blow--up
for Navier-Stokes equations},
Comm. Math. Phys. {\bf 312} (2012), 833-845.

\bibitem{Serre1999}
\textsc{D. Serre},
{\sf Systems of Conservation Laws}, vol.$\;$1,
Cambridge University Press, Cambridge, 1999.

\bibitem{Serrin1963}
{\sc J. Serrin},
{\it The initial value problem for the Navier--Stokes equations},
in:
R. Langer (Ed.), {\sf Nonlinear Problems},
University of Wisconsin Press, Madison, 1963,
pp.$\;$69-98.

\bibitem{Stein1970}
{\sc E. Stein},
{\sf Singular Integrals and
Differentiability Properties of Functions},
Princeton Univ. Press, Princeton, 1970.

\bibitem{Taylor2011}
{\sc M. E. Taylor},
{\sf Partial Differential Equations} (2nd ed.),
vol.$\;$III,
Springer, New York, 2011.

\bibitem{Yuan2008}
{\sc J.$\;$Yuan},
{\it Existence theorem and blow-up criterion
of the strong solutions to the magneto-micropolar
fluid equations},
Math. Meth. Appl. Sci. {\bf 31} (2008),
1113-1130.

\bibitem{Wiegner1987}
{\sc M. Wiegner},
{\it Decay results for weak solutions of the Navier-Stokes
equations on $ \mathbb{R}^{n}\!$},
J. London Math. Soc. {\bf 35} (1987), 303-313.

\bibitem{Zingano2010}
{\sc P.$\;$R.$\;$Zingano},
{\it $\,\!$New $\,\!L^{p}\,\!$-$\;\!L^{q}\!\,\!$
procedures for advection-diffusion equations, I\/}
($n = 1, \;\!\kappa = 0$),
unpublished notes,
Porto Alegre, RS, 2010.

%
%

\bibitem{Zingano2011}
{\sc P.$\;$R.$\;$Zingano},
{\it $\,\!$New $\,\!L^{p}\:\!$-$\;\!L^{q}\!\,\!$
procedures for advection-diffusion equations, II\/}
($n = 1, \;\!\kappa > 0$),
unpublished notes,
Porto Alegre, RS, 2011.

\end{thebibliography}
\end{document}